\documentclass{article}
\usepackage[T2A]{fontenc}
\usepackage[cp1251]{inputenc} %for Windows
\usepackage[english,russian]{babel}
\usepackage[tbtags]{amsmath}
\usepackage{amsfonts,amssymb,mathrsfs,amscd}

\newtheorem{theorem}{\sc \;\;\; ╥хюЁхьр}[section]
\newtheorem{lemma}{\sc \;\;\; ╦хььр}[section]

\newtheorem{remark}{\sc \;\;\; ╟рьхўрэшх}[section]

\newcommand{\supp}{\mathop{\mathrm{supp}}\nolimits}

\makeatletter
\renewcommand{\@begintheorem}[2]%
{\begin{trivlist}\itemindent=\parindent
\item[\hspace{\labelsep}{\bf #1\ #2}]\it} %
\renewcommand{\@endtheorem}
{\end{trivlist}}
\renewcommand{\@opargbegintheorem}[3]%
{\begin{trivlist}\item[\hspace{\labelsep}{\bf #1\ #2 \
\protect{#3}.}]\rm}
\newcommand{\subsect}{%
\@startsection{subsection}{2}{\parindent}{1\baselineskip plus 3ex minus .2ex}%
{2.0ex plus .2ex}{\Large\bf} }
\newcommand{\subsubsect}{%
\@startsection{subsubsection}{3}{\parindent}{2ex plus 1ex minus .2ex}%
{-2.0ex plus .2ex}{\Large\sl} }

\@addtoreset{equation}{section}

\renewcommand{\abstractname}{}

\makeatother
\begin{document}
%%%%%%%%%%%%%%%%%%%%%%%%%%%%%%%%%%%%%%%%%%%%%%%%%%%%%%%%%%%%%%%%%%%%%%%%%%%%
\renewcommand{\abstractname}{}

%\title{\bf └схыхт√ ш Єрє\-сх\-Ёют√ ЄхюЁхь√ фы  шэЄхуЁрыют}
%\author{└.\,╘.~├Ёш°шэ, ╚.\,┬.~╧юхфшэЎхтр}

\begin{center}
{\bf \large A.\,F.~Grishin, I.\,V.~Poedintseva  \\ \vskip 0.5cm
Abelian and Tauberian theorems for integrals}
\end{center}

%\date{17.11.2012}
%\maketitle

%%%%%%%%%%%%%%%%%%%%%%%%%%%%%%%%%%%%%%%%%%%%%%%%%%%%%%%%%%%%%%%%%%%%%%%%%%%%

\begin{abstract}
\\
A new method of obtaining Abelian and Tauberian theorems  for the
integral of the form $\int\limits_0^\infty K\left(\frac{t}{r}\right)
d\mu(t)$ is proposed. It is based on the use of limit sets of the
measures. A version of Azarin's sets is constructed for Radon's
measures on the ray $(0,\infty)$. Abelian theorems of a new type are
proved in which asymptotic behavior of the integral is described in
terms of these limit sets. Using these theorems together with an
improved version of the well-known Carleman's theorem on analytic
continuation, a substantial improvement of the second Wiener
Tauberian theorem is obtained.
\\
Reference: 25 units.
\\
{\bf Comments:} 73 pages, in Russian. \vskip0.5cm ╧Ёхф\-ырур\-хЄё 
эют√щ ьхЄюф яюыєўх\-эш  рсхых\-т√ї ш Єрє\-сх\-Ёю\-т√ї Єхю\-Ёхь фы 
шэЄху\-Ёр\-ыют тшфр $\int\limits_0^\infty K\left(\frac{t}{r}\right)
d\mu(t)$. ╬э срчш\-Ёє\-хЄё  эр шё\-яюы№чю\-тр\-эшш ётющёЄт
яЁх\-фхы№\-э√ї ьэю\-цхёЄт ьхЁ. ─ы  ¤Єюую ёЄЁюшЄ\-ё  трЁш\-рэЄ
Єхю\-Ёшш яЁх\-фхы№\-э√ї ьэю\-цх\-ёЄт └чр\-Ёш\-эр фы  Ёр\-фю\-эю\-т√ї
ьхЁ эр яюыє\-юёш $(0,\infty)$. ─юърч√\-тр■Єё  рсхых\-т√ Єхю\-Ёх\-ь√
эютюую Єшяр, т ъюЄю\-Ё√ї рёшья\-Єю\-Єш\-ўхёъюх яю\-тхфх\-эшх
т√°х\-эр\-чтрэ\-э√ї шэЄху\-Ёр\-ыют юяшё√\-тр\-хЄё  т ЄхЁьш\-эрї
яЁх\-фхы№\-э√ї ьэю\-цх\-ёЄт ьхЁ $\mu$. ╚ёяюы№чє  ¤Єш Єхю\-Ёх\-ь√, р
Єръ\-цх фю\-ърчрэ\-э√щ т ёЄрЄ№х єёшыхэ\-э√щ трЁш\-рэЄ шч\-тхёЄ\-эющ
ыхьь√ ╩рЁых\-ьрэр юс рэрыш\-Єш\-ўхёъюь яЁю\-фюы\-цх\-эшш,
фю\-ърч√\-тр\-хЄё  чэрўш\-Єхы№\-эюх єёшых\-эшх тЄюЁющ
Єрє\-сх\-Ёю\-тющ Єхю\-Ёх\-ь√ ┬шэхЁр.
\\
┴шсышюуЁрЇш : 25 эрчтрэш .
\\{\bf ╩ы■ўхт√х ёыютр:} єЄюўэ╕э\-э√щ яю\-Ё фюъ ┬рыш\-Ёюэр, Ёр\-фю\-эютр ьхЁр,
яЁх\-фхы№\-эюх ьэю\-цхёЄ\-тю └чрЁш\-эр ьхЁ√, Ёхуєы Ё\-эр  ьхЁр
└чрЁш\-эр, ЄрєсхЁю\-тр ЄхюЁх\-ьр ┬шэхЁр.
\end{abstract}

\tableofcontents   % ╤╬─┼╨╞└═╚┼

 %%%%%%%%%%%%%%%%%%%%%%%%%%%%%%%%%%%%%%%%%%%%%%%%%%%%%%%%%%%%%%%%%%%%%%%%%%%%%%%%%%%%%%%%%
\section{┬ёЄєяыхэшх}
%\addcontentsline{toc}{section}{\bf\quad\:  ┬ёЄєяыхэшх}
%%%%%%%%%%%%%%%%%%%%%%%%%%%%%%%%%%%%%%%%%%%%%%%%%%%%%%%%%%%%%%%%%%%%%%%%%%%%%%%%%%%%%%%%%
\qquad ╧ю\-ёъюы№ъє єЄюў\-э╕э\-э√щ яю\-Ё фюъ тёЄЁхўр\-хЄё  т
ЇюЁ\-ьєыш\-Ёют\-ърї сюы№\-°шэ\-ёЄтр єЄ\-тхЁ\-ц\-фх\-эшщ
яЁхф\-ырурх\-ьющ ЁрсюЄ√, Єю ь√ эрўэ╕ь ё эх\-ъюЄюЁ√ї юсю\-чэрўх\-эшщ
ш Ёх\-чєы№\-Єр\-Єют, ёт чрэ\-э√ї ё єЄюў\-э╕э\-э√ь яю\-Ё фъюь.

\par ╙Єюў\-э╕э\-э√щ яю\-Ё фюъ $\rho(r)$ шуЁрхЄ трцэє■ Ёюы№ т ЄхюЁшш
рсхых\-т√ї ш Єрє\-сх\-Ёю\-т√ї Єхю\-Ёхь, т ЄхюЁшш ЁюёЄр
ёєс\-урЁьюэш\-ўхёъшї Їєэъ\-Ўшщ, т ЄхюЁшш тхЁю Є\-эюё\-Єхщ. ╚ч Єюую,
ўЄю $\rho(r)$ -- єЄюў\-э╕э\-э√щ яю\-Ё фюъ, ёых\-фє\-хЄ, ўЄю
Їєэъ\-Ўш  $V(r)=r^{\rho(r)}$  ты \-хЄё  Ёхуєы Ё\-эю ьхэ ■\-∙хщ\-ё 
Їєэъ\-Ўшхщ т ёь√ёых ╩рЁр\-ьрЄ√. ╚чыюцх\-эшх эхъюЄю\-Ё√ї ётющёЄт
єЄюў\-э╕э\-эю\-ую яю\-Ё ф\-ър ьюц\-эю эрщЄш т  \cite{Levin}.
╧ю тшт\-°р ё  яючцх ъэшур \cite{Bingham}  ты \-хЄё 
¤эЎш\-ъыюяхфш\-ўхёъшь ЄЁєфюь яю Ёхує\-ы Ё\-эю ьхэ ■\-∙шь\-ё 
Їєэъ\-Ўш ь ш шї яЁш\-ыюцх\-эш\- ь.

\par ─ю чэръюь\-ёЄтр ё ЇюЁьры№\-э√ь юяЁх\-фхых\-эш\-хь єЄюў\-э╕э\-эю\-ую
яю\-Ё ф\-ър ўшЄр\-Єхы№ ью\-цхЄ ёўш\-ЄрЄ№, ўЄю $V(r)=r^\rho a(r)$,
$$
a(r)=\left\{\begin{array}{lcl}
\ln^\alpha(er),\quad r\geq 1,\\
\ln^\alpha \frac{e}{r},\qquad r\in (0,1),
\end{array}\right.
$$
уфх $\rho$ ш $\alpha$ -- яЁюшч\-тюы№эю т√сЁрэ\-э√х тх∙хёЄ\-тхэ\-э√х
ўшёыр.

\par ╧єёЄ№ $\rho=\rho(\infty)=\lim\limits_{r\to\infty}\rho(r)$. ┼ёыш
юсю\-чэрўшЄ№

\begin{equation}\label{def_gamma}
\gamma(t)=\sup\limits_{r>0}\frac{V(rt)}{t^\rho V(r)},
\end{equation}
Єю юўхтшфэю, ўЄю яЁш $r>0$ ш $t>0$ т√\-яюыэ \-хЄё  эх\-Ёртхэ\-ёЄтю
$V(rt)\leq t^\rho \gamma(t)V(r)$.

\par ┬ Єхю\-Ёхьх  2.5 (эєьхЁр\-Ўш  Єхю\-Ёхь тю ттхфх\-эшш
ёют\-ярфр\-хЄ ё шї эєьхЁр\-Ўшхщ т ЄхъёЄх ЁрсюЄ√) фю\-ърч√\-тр\-хЄё ,
ўЄю $\gamma(t)$  -- эх\-яЁх\-Ё√тэр  Їєэъ\-Ўш  ш ўЄю т√\-яюыэ \-■Єё 
Ёртхэ\-ёЄ\-тр
$$
\lim\limits_{t\to\infty} \frac{\ln\gamma(t)}{\ln t}=0, \qquad
\lim\limits_{t\to\infty} \frac{\ln\gamma
\left(\frac{1}{t}\right)}{\ln t}=0.
$$
\par ▌Єр Єхю\-Ёхьр эхьэю\-ую єёшыш\-трхЄ Ёх\-чєы№\-ЄрЄ ╧юЄЄхЁр  \cite{Potter}
ш чэрўш\-Єхы№\-эю єяЁю\-∙р\-хЄ хую ЇюЁ\-ьє\-ыш\-Ёютъє. ╬эр  --
трцэ√щ шэёЄЁє\-ьхэЄ яЁш Ёрсю\-Єх ё єЄюў\-э╕э\-э√ь яю\-Ё ф\-ъюь. ┬
фры№\-эхщ\-°хь ЄхъёЄх ЁрсюЄ√ ёшьтюы $\gamma(t)$ эх сєфхЄ
шё\-яюы№чю\-трЄ№\-ё  фы  юсю\-чэрўх\-эш  фЁєушї Їєэъ\-Ўшщ.

\par ╤хщўрё ь√ юяЁх\-фхышь фтр ъырёёр Ёр\-фю\-эю\-т√ї ьхЁ эр яюыє\-юёш  $(0,\infty)$.
┬ Ёрсю\-Єх, т юёэют\-эюь, шчєўр\-■Єё  ьхЁ√ шч ¤Єшї ъырёёют.

\par ╫хЁхч $\frak{M}_\infty$ $(\rho(r))$ ь√ юсю\-чэрўшь ьэю\-цх\-ёЄ\-тю
Ёр\-фю\-эю\-т√ї ьхЁ эр яюыє\-юёш  $(0,\infty)$, фы  ъюЄю\-Ё√ї
т√\-яюыэ \-хЄё  эх\-Ёртхэ\-ёЄтю
$$
\sup\limits_{r\geq 1}\frac{|\mu|([r,er])}{V(r)}<\infty.
$$
\par ╫хЁхч $\frak{M}(\rho(r))$ ь√ юсю\-чэрўшь ьэю\-цх\-ёЄ\-тю Ёр\-фю\-эю\-т√ї
ьхЁ эр яюыє\-юёш  $(0,\infty)$, фы  ъюЄю\-Ё√ї т√\-яюыэ \-хЄё 
эх\-Ёртхэ\-ёЄтю
$$
\sup\limits_{r>0}\frac{|\mu|([r,er])}{V(r)}<\infty.
$$
\par ┼ёыш т√\-яюыэ \-хЄё  ёююЄ\-эю°х\-эшх
$$
\mathop{\overline\lim}\limits_{r\to\infty}
\frac{|\mu|([r,er])}{V(r)}\in(0,\infty),
$$
Єю єЄюў\-э╕э\-э√щ яю\-Ё фюъ $\rho(r)$ сєфхь эр\-ч√трЄ№ {\it
єЄюў\-э╕э\-э√ь яю\-Ё ф\-ъюь} ьхЁ√ $\mu$, р ўшёыю $\rho=\rho(\infty)$
сєфхь эр\-ч√\-трЄ№ яю\-Ё ф\-ъюь Ёр\-фю\-эю\-тющ ьхЁ√ $\mu$. ╥ръюх
юяЁх\-фхых\-эшх яю\-Ё ф\-ър яюч\-тюы хЄ Ёрё\-ёьрЄЁш\-трЄ№
Ёр\-фю\-эют√ ьхЁ√ ы■сю\-ую тх∙хёЄ\-тхэ\-эю\-ую яю\-Ё ф\-ър.

\par {\it ╧Ёх\-фхы№\-эюх ьэю\-цх\-ёЄтю $Fr[\mu]=Fr[\rho(r),\mu]$}
Ёр\-фю\-эю\-тющ ьхЁ√ $\mu$ ь√ юяЁх\-фхы \-хь тёыхф чр └чрЁш\-э√ь ъръ
ьэю\-цх\-ёЄ\-тю ьхЁ $\nu$ тшфр $\nu=\lim\limits_{n\to\infty}
\mu_{t_n}$, уфх $t_n\to\infty$, р ьхЁр $\mu_t$ юяЁх\-фхы \-хЄ\-ё 
Ёртхэ\-ёЄ\-тюь
$$
\mu_t(E)=\frac{\mu(tE)}{V(t)}.
$$
╨ртхэёЄтю $\nu=\lim\limits_{n\to\infty}\mu_{t_n}$ ючэрўр\-хЄ, ўЄю
яю\-ёыхфю\-тр\-Єхы№\-эюёЄ№ ьхЁ $\mu_{t_n}$ {\it °ш\-Ёю\-ъю
ёїю\-фшЄ\-ё } ъ ьхЁх  $\nu$. ╬яЁх\-фхых\-эшх °ш\-Ёю\-ъющ
ёїюфш\-ьюёЄш фрэю т Ёрч\-фхых  3.

\par ─ю\-ърч√тр\-хЄё , ўЄю хёыш $\mu\in$ $\frak{M}_\infty$ $(\rho(r))$,
Єю ьэю\-цх\-ёЄ\-тю $Fr[\mu]$  юсырфр\-хЄ ётющёЄ\-трьш,
яхЁх\-ўшёыхэ\-э√ьш т ёююЄ\-тхЄ\-ёЄ\-тє■\-∙хщ Єхю\-Ёхьх └чрЁш\-эр
\cite{Azarin_79}. ▌Єш ётющ\-ёЄтр Єръ\-цх яхЁх\-ўшёых\-э√ т
Єхю\-Ёх\-ьх 3.10 эр°хщ ЁрсюЄ√. └чрЁшэ Ёрё\-ёьрЄ\-Ёш\-тры Єюы№ъю
яюыю\-цш\-Єхы№\-э√х ьхЁ√ т $\mathbb{R}^n$ ш юяЁх\-фхы ы
єЄюў\-э╕э\-э√щ яю\-Ё фюъ ьхЁ√ эх\-ёъюы№ъю яю-шэюьє, ўЄю тюч\-ьюц\-эю
Єюы№ъю фы  ёыєўр  $\rho>0$. ┬  ¤Єюь ёыєўрх эр°х юяЁх\-фхых\-эшх
єЄюў\-э╕э\-эю\-ую яю\-Ё ф\-ър ьхЁ√, ¤ътш\-трыхэЄ\-эю Єюьє, ъюЄю\-Ёюх
шё\-яюы№\-чю\-тры └чрЁшэ. ╬ЄьхЄшь х∙╕, ўЄю └чрЁшэ Ёрё\-ёьрЄ\-Ёш\-тры
ёїюфш\-ьюёЄ№ ьхЁ юЄышў\-эє■ юЄ Єющ, ъюЄю\-Ёр 
Ёрё\-ёьрЄ\-Ёш\-тр\-хЄё  т ёЄрЄ№х.

\par ╠хЁр  $\mu$ эр\-ч√тр\-хЄё  {\it Ёхує\-ы Ё\-эющ} ьхЁющ шыш {\it ьхЁющ
Ёхує\-ы Ё\-эющ т ёь√ёых └чрЁш\-эр}, хёыш яЁх\-фхы№\-эюх
ьэю\-цх\-ёЄтю $Fr[\mu]$ ёю\-ёЄюшЄ шч хфшэ\-ёЄтхэ\-эющ ьхЁ√ $\nu$. ┬
¤Єюь ёыєўрх ьхЁр  $\nu$ ё эхюс\-їюфш\-ьюё\-Є№■ шьххЄ тшф
$d\nu(t)=ct^{\rho-1}dt$ (ЄхюЁх\-ьр \ref{azarin_th_8}).

┬ ёЄрЄ№х сєфхЄ шё\-яюы№чю\-трЄ№ё  ъръ ёЄрэ\-фрЁЄ\-э√х Єръшх
юсю\-чэрўх\-эш .
\\{\bf 1.} ╘єэъЎш  $\Psi(r)$,
\begin{equation}\label{introduction_Psi(r)}
\Psi(r)=\int\limits_0^\infty K\left(\frac{t}{r}\right)d\mu(t).
\end{equation}
┬ эх\-ъюЄю\-Ё√ї ёыєўр ї тьхёЄю $\Psi(r)$ ь√ сєфхь шё\-яюы№чю\-трЄ№
сюыхх шэ\-ЇюЁьр\-Єшт\-эюх юсю\-чэрўх\-эшх $\Psi(K,r)$.
\\{\bf 2.} ╘єэъЎш  $J(r)=\frac{1}{V(r)}\Psi(r)$.
\\{\bf 3.}  ╠хЁр $s$, $ds(t)=\Psi(t)dt$.

\par ╫хЁхч $\mu(t)$ ь√ сєфхь юсю\-чэрўрЄ№  Їєэъ\-Ўш■ Ёрё\-яЁх\-фхых\-эш 
ьхЁ√ $\mu$, Єръ ўЄю $\mu((a,b])$ $=\mu(b)-\mu(a)$.

\par {\it ╧Ёхфхы№\-эюх ьэю\-цх\-ёЄтю} Їєэъ\-Ўшш $f(r)$ яю эр\-яЁрт\-ых\-эш■ $r\to\infty$
(ьэю\-цх\-ёЄтю яЁх\-фхыют тшфр $\lim\limits_{n\to\infty} f(r_n)$,
уфх $r_n\to\infty$) ь√ сєфхь юсю\-чэрўрЄ№ $L(f,\infty)$.

\par ╘єэъЎш■ $f(r)$ ь√ сєфхь эр\-ч√трЄ№ {\it ЇшэшЄ\-эющ} Їєэъ\-Ўшхщ
эр яюыє\-юёш $(0,\infty)$, хёыш $\supp f\subset
[a,b]\subset(0,\infty)$.

\par ╚ёяюы№чю\-тр\-эшх ётющёЄт ьэю\-цх\-ёЄтр $Fr[\mu]$ яю\-чтюыш\-ыю
яю\-ыєўшЄ№ эх\-ёъюы№\-ъю эю\-т√ї Єхю\-Ёхь рсхыхтр ш Єрє\-сх\-Ёю\-тр
Єшяют фы  шэ\-Єху\-Ёр\-ыют тшфр (\ref{introduction_Psi(r)}).

\par ╥хю\-Ёх\-ь√ {\it рсхыхтр} Єшяр  -- ¤Єю Єхю\-Ёх\-ь√, т ъюЄю\-Ё√ї шё\-їюф  шч
ётющёЄт ьхЁ√ $\mu$ эр\-їюф Є\-ё  ётющ\-ёЄтр Їєэъ\-Ўшш  $\Psi$.

\par ╥хю\-Ёх\-ь√ {\it Єрє\-сх\-Ёютр} Єшяр  -- ¤Єю Єхю\-Ёх\-ь√, т ъюЄю\-Ё√ї шёїюф  шч
ётющёЄт Їєэъ\-Ўшш  $\Psi$ эр\-їюф Є\-ё  ётющ\-ёЄтр ьхЁ√ $\mu$.

\par ┬ю ьэюушї шчтхёЄ\-э√ї Єхю\-Ёх\-ьрї рсхыхтр Єшяр єЄ\-тхЁ\-ц\-фр\-хЄё , ўЄю
шч ЇюЁ\-ьєы√ $\mu'(t)$ $\sim$ $A\frac{V(t)}{t}$ $(t\to\infty)$
т√\-Єхър\-хЄ рёшья\-Єю\-Єш\-ўхё\-ър  ЇюЁ\-ьє\-ыр $\Psi(r)\sim
BV(r)$. ╥хю\-Ёх\-ь√ Єръюую Єшяр ьюц\-эю эрщЄш, т ўрёЄ\-эюё\-Єш, т
ъэшурї \cite{Widder}, \cite{Hardy}, \cite{Bingham}, \cite{Korevaar}.

\par ╤ЇюЁ\-ьєыш\-Ёєхь яЁюёЄхщ\-°є■ Єхю\-Ёхьє рсхыхтр Єшяр шч эр°хщ ЁрсюЄ√.
% % % % % % % % % % % % % % % % % % % % % % % % % % % % % % % % % % % % %
\par{\sc  ╥хюЁхьр  4.1.} \:{\it ╧єёЄ№ $\mu\in
\frak{M}_\infty$ $(\rho(r))$, $K$  -- эх\-яЁх\-Ё√т\-эюх ЇшэшЄ\-эюх
 фЁю эр яюыє\-юёш $(0,\infty)$. ╥юуфр т√\-яюыэ \-хЄ\-ё 
Ёртхэ\-ёЄ\-тю
$$
L(J,\infty)=\left\{\int\limits_0^\infty K(t)d\nu(t):\quad \nu\in
Fr[\mu]\right\}.
$$}
% % % % % % % % % % % % % % % % % % % % % % % % % % % % % % % % % % % % %
\par ─Ёєушх Єхю\-Ёх\-ь√ рсхыхтр Єшяр т яЁхф\-ёЄрт\-ыхэ\-эющ Ёрсю\-Єх --
¤Єю рэрыюуш Єхю\-Ёх\-ь√  4.1, ъюЄюЁ√х яюыєўр\-■Є\-ё  яЁш
Ёрч\-ышў\-э√ї юуЁрэш\-ўх\-эш\- ї эр  фЁю  $K$ ш ьхЁє  $\mu$. ┬ юфэшї
ёыєўр ї ь√ юЄърч√\-тр\-хь\-ё  юЄ ЄЁхсютр\-эш  ЇшэшЄ\-эюё\-Єш  фЁр
$K$. ┬ фЁєушї сюыхх ёыюц\-э√ї ёыєўр ї ь√ Єръ\-цх юЄърч√\-трхь\-ё  юЄ
ЄЁхсю\-тр\-эш  эх\-яЁх\-Ё√т\-эюё\-Єш  фЁр  $K$. ┬ юс∙хь ёыєўрх  фЁю
$K$ -- ¤Єю сюЁх\-ыхт\-ёър  Їєэъ\-Ўш  эр яюыє\-юёш $(0,\infty)$.

\par ╬ЄьхЄшь трцэюх Ёрч\-ышўшх ьхцфє єяюь \-эє\-Є√ьш т√°х
яЁхф\-°хёЄ\-тє■\-∙ш\-ьш Ёх\-чєы№\-Єр\-Єрьш ш Ёх\-чєы№\-Єр\-Єрьш
яЁхф\-ёЄрт\-ыхэ\-э√ьш т ¤Єющ Ёрсю\-Єх. ┬ яЁхф\-°хёЄ\-тє■\-∙шї
Єхю\-Ёх\-ьрї шё\-ёых\-фє\-■Є\-ё  ёыєўрш, ъюуфр  $\mu$  ты \-хЄ\-ё 
Ёхує\-ы Ё\-эющ ьхЁющ. ┬ эр°шї Єхю\-Ёх\-ьрї шё\-ёых\-фє\-хЄё 
чэрўш\-Єхы№\-эю сюыхх °ш\-Ёю\-ъшщ ъырёё ьхЁ. ═ряЁш\-ьхЁ, т
Єхю\-Ёх\-ьх  4.1 ¤Єю ъырёё $\frak{M}_\infty$ $(\rho(r))$.

\par ╧Ёш юЄърчх юЄ эх\-яЁх\-Ё√т\-эюё\-Єш  фЁр  $K$ т ёыєўрх юс∙шї Ёр\-фю\-эю\-т√ї
ьхЁ яЁх\-фхы№\-эюх ьэю\-цх\-ёЄ\-тю $Fr[\mu]$ єцх эх юяЁх\-фхы \-хЄ
ьэю\-цх\-ёЄ\-тю $L(J,\infty)$, ъръ ¤Єю тшф\-эю шч Єхю\-Ёх\-ь√  4.2.

\par ┬рцэ√ь эют√ь Ёх\-чєы№\-ЄрЄюь  ты \-хЄё  єЄ\-тхЁ\-ц\-фх\-эшх, ўЄю т
ёыєўрх Ёрч\-Ё√т\-э√ї  фхЁ  $K$ яЁх\-фхы№\-эюх ьэю\-цх\-ёЄ\-тю
$Fr[\mu]$ ьхЁ√ $\mu$ юфэю\-чэрў\-эю юяЁх\-фхы \-хЄ ьэю\-цх\-ёЄ\-тю
$Fr[s]$. ╬фэшь шч фю\-ёЄшцх\-эшщ ЁрсюЄ√  ты \-хЄё  Єю, ўЄю
ттюфшЄ\-ё  т Ёрё\-ёьюЄЁх\-эшх яЁх\-фхы№\-эюх ьэю\-цх\-ёЄ\-тю $Fr[s]$
ьхЁ√ $s$. ╤ыхфє\-■\-∙р  Єхю\-Ёх\-ьр -- ¤Єю юфшэ шч юёэют\-э√ї
Ёх\-чєы№\-Єр\-Єют ЁрсюЄ√. ═ряюь\-эшь, ўЄю Їєэъ\-Ўш  $\gamma(t)$
юяЁх\-фхы \-хЄё  Ёртхэ\-ёЄ\-тюь (\ref{def_gamma}).
%%%%%%%%%%%%%%%%%%%%%%%%%%%%%%%%%%%%%%%%%%%%%%%%%%%%%%%%%%%%%%%%%%%%%%%%%%
\par{\sc  ╥хюЁхьр  4.8.} \: {\it ╧єёЄ№ $\rho(r)$ -- яЁюшч\-тюы№\-э√щ
єЄюў\-э╕э\-э√щ яю\-Ё фюъ, $\mu\in$ $\frak{M}(\rho(r))$. ╧єёЄ№ $K$ --
сюЁхыхт\-ёър  Їєэъ\-Ўш  эр яюыє\-юёш $(0,\infty)$ Єр\-ър , ўЄю
$t^{\rho-1}\gamma(t)K(t)\in L_1(0,\infty)$, $\Psi$  -- Їєэъ\-Ўш 
юяЁх\-фхы \-хьр  Ёртхэ\-ёЄ\-тюь (\ref{introduction_Psi(r)}). ╥юуфр
ьхЁр  $s$, $ds(t)=\Psi(t)dt$, яЁш\-эрф\-ыхцшЄ ъырёёє
$\frak{M}(\rho(r)+1)$ ш х╕ яЁх\-фхы№\-эюх ьэю\-цх\-ёЄ\-тю
$Fr[\rho(r)+1,s]$ ёюёЄю\-шЄ шч рсёюы■Є\-эю эх\-яЁх\-Ё√т\-э√ї ьхЁ,
ьэю\-цх\-ёЄ\-тю яыюЄ\-эюёЄхщ ъюЄю\-Ё√ї шьххЄ тшф
$$
\left\{\int\limits_0^\infty K\left(\frac{t}{u}\right) d\nu(t):\:
\nu\in Fr[\mu]\right\}.
$$}
%%%%%%%%%%%%%%%%%%%%%%%%%%%%%%%%%%%%%%%%%%%%%%%%%%%%%%%%%%%%%%%%%%%%%%%%%%

\par ┬ Єхї ёыєўр ї, ъюуфр шч эр°шї Єхю\-Ёхь яюыєўр\-хЄ\-ё  Ёртхэ\-ёЄтю
$L(J,\infty)=\{0\}$, тюяЁюё ю яю\-Ё ф\-ъх ЁюёЄр Їєэъ\-Ўшш $\Psi(r)$
эр схё\-ъю\-эхў\-эюё\-Єш юёЄр╕Є\-ё  юЄъЁ√\-Є√ь.

\par ╨рё\-ёьюЄЁшь ёыєўрщ, ъюуфр $K$  -- схё\-ъю\-эхўэю фшЇ\-Їх\-Ёхэ\-Ўш\-Ёєх\-ьюх
ЇшэшЄ\-эюх  фЁю эр яюыє\-юёш $(0,\infty)$. ┬ ¤Єюь ёыєўрх эрЁ фє ё
ЇюЁ\-ьє\-ыющ  (\ref{introduction_Psi(r)}) фы  Їєэъ\-Ўшш $\Psi(r)$
ёяЁртхф\-ышт√ ЇюЁ\-ьє\-ы√
\begin{equation}\label{formula_for_Psi}
(-1)^{n+1}r^{n+1}\Psi(r)=\int\limits_0^\infty
K^{(n+1)}\left(\frac{t}{r}\right)F_n(t)dt,\quad n=0,1,...,
\end{equation}
уфх $F_0(t)=\mu(t)$, $F_{n+1}'(t)=F_n(t)$.

\par ╤ЄртшЄё  тюяЁюё: яючтюы ■Є ыш Єхю\-Ёх\-ьр  4.1 ш
Ёртхэ\-ёЄ\-тр (\ref{formula_for_Psi}) юяЁх\-фх\-ышЄ№ яю\-Ё фюъ ЁюёЄр
Їєэъ\-Ўшш $\Psi(r)$ эр схё\-ъю\-эхў\-эюё\-Єш? ▌ЄюЄ тюяЁюё
шё\-ёых\-фє\-хЄё  т ъюэЎх Ёрч\-фхыр  4. ╬ЄтхЄ ёых\-фє\-■\-∙шщ. ┬
Ёрё\-ёьрЄЁш\-трх\-ьюь ёыєўрх Єхю\-Ёх\-ьр  4.1 ш Ёртхэ\-ёЄ\-тр
(\ref{formula_for_Psi}) ўрёЄю фр■Є юЄтхЄ ю яю\-Ё ф\-ъх ЁюёЄр
Їєэъ\-Ўшш $\Psi(r)$ эр схё\-ъю\-эхў\-эюё\-Єш. ╬фэръю, шьх■Є\-ё 
Ёрч\-ышў\-э√х шё\-ъы■ўш\-Єхы№\-э√х ёыєўрш. ┬ ¤Єшї ёыєўр\- ї тюяЁюё ю
яю\-Ё ф\-ъх ЁюёЄр Їєэъ\-Ўшш $\Psi(r)$ эр схё\-ъю\-эхў\-эюё\-Єш
т√\-їюфшЄ чр Ёрьъш эр°хщ ЁрсюЄ√.

\par ╬сЁрЄшьё  ЄхяхЁ№ ъ Єхю\-Ёх\-ьрь Єрє\-сх\-Ёютр Єшяр. ┬рцэ√ьш
Ёх\-чєы№\-Єр\-Єр\-ьш Єрє\-сх\-Ёю\-тющ ЄхюЁшш  ты \-■Єё 
Єрє\-сх\-Ёю\-т√ Єхю\-Ёх\-ь√ ┬шэхЁр. ╤ЇюЁ\-ьє\-ыш\-Ёєхь ¤Єш
Єхю\-Ёх\-ь√.
%%%%%%%%%%%%%%%%%%%%%%%%%%%%%%%%%%%%%%%%%%%%%%%%%%%%%%%%%%%%%%%%%%%%%%%%%%
\par {\sc  ╥хюЁхьр 1.1.}\: {\it ╧єёЄ№ $F(x)\in L_1(-\infty,\infty)$,
$\int\limits_{-\infty}^\infty F(x)e^{-i\lambda x}dx\neq 0$ яЁш
$\lambda\in(-\infty,\infty)$, $g(x)$  -- шч\-ьхЁш\-ьр 
юуЁрэш\-ўхэ\-эр  Їєэъ\-Ўш  эр тх∙хёЄ\-тхэ\-эющ юёш. ╧єёЄ№
$$
\lim\limits_{x\to +\infty} \int\limits_{-\infty}^\infty
F(x-y)g(y)dy= A \int\limits_{-\infty}^\infty F(y)dy.
$$
╥юуфр фы  ы■сющ Їєэъ\-Ўшш $H\in L_1(-\infty,\infty)$ т√\-яюыэ \-хЄё 
Ёртхэ\-ёЄтю
$$
\lim\limits_{x\to +\infty}  \int\limits_{-\infty}^\infty
H(x-y)g(y)dy=A \int\limits_{-\infty}^\infty H(y)dy.
$$}
%%%%%%%%%%%%%%%%%%%%%%%%%%%%%%%%%%%%%%%%%%%%%%%%%%%%%%%%%%%%%%%%%%%%%%%%%%
\par ╬сю\-чэрўшь ўхЁхч $M$ яЁюёЄ\-Ёрэ\-ёЄтю эх\-яЁх\-Ё√т\-э√ї
Їєэъ\-Ўшщ  $F$ эр юёш $(-\infty,\infty)$ ё ъю\-эхў\-эющ
рьры№\-урь-эюЁьющ
$$
\|F\|_M=\sum\limits_{n=-\infty}^\infty \max\{|F(x)|:\;
x\in[n,n+1]\}.
$$
\par ╤ыхфє\-■∙є■ Єхю\-Ёхьє ўрёЄю эр\-ч√тр\-■Є тЄюЁющ Єрє\-сх\-Ёю\-тющ
Єхю\-Ёх\-ьющ ┬шэхЁр.
%%%%%%%%%%%%%%%%%%%%%%%%%%%%%%%%%%%%%%%%%%%%%%%%%%%%%%%%%%%%%%%%%%%%%%%%%%
\par {\sc  ╥хюЁхьр 1.2.}\: {\it ╧єёЄ№ $F\in M$,
$\int\limits_{-\infty}^\infty F(x)e^{-i\lambda x}dx\neq 0$ яЁш
$\lambda\in(-\infty,\infty)$, ш яєёЄ№ $\nu$  -- Ёр\-фю\-эютр ьхЁр эр
тх∙хёЄ\-тхэ\-эющ юёш, Єр\-ър , ўЄю фы  эх\-ъюЄю\-Ёюую $B>0$ ш ы■сюую
Ўхыю\-ую  $n$ т√\-яюыэ \-хЄё  эх\-Ёртхэ\-ёЄтю $|\nu|([n,n+1])\leq
B$. ╧єёЄ№
$$
\lim\limits_{x\to +\infty}  \int\limits_{-\infty}^\infty
F(x-y)d\nu(y)=A \int\limits_{-\infty}^\infty F(y)dy.
$$
╥юуфр фы  ы■сющ Їєэъ\-Ўшш $H\in M$ т√\-яюыэ \-хЄё  Ёртхэ\-ёЄ\-тю
$$
\lim\limits_{x\to +\infty}  \int\limits_{-\infty}^\infty
H(x-y)d\nu(y)=A \int\limits_{-\infty}^\infty H(y)dy.
$$}
%%%%%%%%%%%%%%%%%%%%%%%%%%%%%%%%%%%%%%%%%%%%%%%%%%%%%%%%%%%%%%%%%%%%%%%%%%
\par ═рь сєфхЄ єфюсэю шьхЄ№ фхыю ё ьєы№Єш\-яышър\-Єшт\-э√ь трЁш\-рэ\-Єюь
Єхю\-Ёх\-ь√  1.2, ъюЄю\-Ё√щ яю\-ыєўр\-хЄё  шч Єхю\-Ёх\-ь√  1.2
чр\-ьхэющ яхЁх\-ьхэ\-эющ $y=\ln t$ т шэ\-ЄхуЁр\-ырї ш
шё\-яюы№чю\-тр\-эш\-хь юсю\-чэрўх\-эшщ $x=\ln r$,
$K(t)=\frac{1}{t}F(-\ln t)$, $d\mu(t)=td\nu(\ln t)$.

\par ╬сючэрўшь  ўхЁхч $M_1$ яЁюёЄ\-Ёрэ\-ёЄтю эх\-яЁх\-Ё√т\-э√ї Їєэъ\-Ўшщ  $K$
эр яюыє\-юёш $(0,\infty)$, фы  ъюЄю\-Ё√ї ёїю\-фшЄ\-ё  Ё ф
$\sum\limits_{n=-\infty}^\infty K_n e^n$, уфх
$K_n=\max\{|K(x)|:\:x\in[e^n,e^{n+1}]\}$.
% % % % % % % % % % % % % % % % % % % % % % % % % % % % % % % % % % % % %
\par {\sc  ╥хюЁхьр 1.3.}\: {\it ╧єёЄ№ $K(t)\in M_1$,
$\int\limits_0^\infty K(t)t^{i\lambda}dt\neq 0$ яЁш
$\lambda\in(-\infty,\infty)$, яєёЄ№ $\mu$  -- Ёр\-фю\-эю\-тр ьхЁр эр
яюыє\-юёш $(0,\infty)$ Єр\-ър , ўЄю фы  эх\-ъюЄю\-Ёюую $B>0$ ш
ы■сюую Ўхыюую  $n$ т√\-яюыэ \-хЄё  эх\-Ёртхэ\-ёЄтю
$|\mu|([e^n,e^{n+1}])\leq Be^n$. ╧єёЄ№
$$
\lim\limits_{r\to\infty} \frac{1}{r} \int\limits_0^\infty
K\left(\frac{t}{r}\right)d\mu(t)=A\int\limits_0^\infty K(t)dt.
$$
╥юуфр фы  ы■сющ Їєэъ\-Ўшш $Q\in M_1$ т√\-яюыэ \-хЄё  Ёртхэ\-ёЄ\-тю
$$
\lim\limits_{r\to\infty} \frac{1}{r} \int\limits_0^\infty
Q\left(\frac{t}{r}\right)d\mu(t)=A\int\limits_0^\infty Q(t)dt.
$$}
\par ╘юЁьєыш\-Ёютър Єхю\-Ёх\-ь√  1.3 эр\-ёых\-фєхЄ тшэх\-Ёют\-ёъє■
ЇюЁ\-ьє\-ыш\-Ёютъє Єрє\-сх\-Ёю\-тющ Єхю\-Ёх\-ь√. ╬фэръю, ¤Єр
Єхю\-Ёх\-ьр ¤ътш\-тр\-ыхэЄ\-эр Єхю\-Ёхьх ё сюыхх ёЄрэ\-фрЁЄ\-эющ
ЇюЁ\-ьє\-ыш\-Ёют\-ъющ фы  Єрє\-сх\-Ёю\-т√ї Єхю\-Ёхь.
% % % % % % % % % % % % % % % % % % % % % % % % % % % % % % % % % % % % %
\par {\sc  ╥хюЁхьр 1.4.} \; {\it ╧єёЄ№ $K(t)\in M_1$,
$\int\limits_0^\infty K(t)t^{i\lambda}dt\neq 0$  яЁш
$\lambda\in(-\infty,\infty)$, яєёЄ№ $\mu$  -- Ёр\-фю\-эю\-тр ьхЁр эр
яюыє\-юёш $(0,\infty)$ Єр\-ър , ўЄю фы  эх\-ъюЄю\-Ёюую $B>0$ ш
ы■сюую Ўхыюую  $n$ т√\-яюыэ \-хЄё  эх\-Ёртхэ\-ёЄтю
$|\mu|([e^n,e^{n+1}])\leq Be^n$. ╧єёЄ№ ёє\-∙хёЄ\-тєхЄ яЁх\-фхы
$$
\lim\limits_{r\to\infty} \frac{1}{r} \int\limits_0^\infty
K\left(\frac{t}{r}\right)d\mu(t)=c.
$$
╥юуфр яЁх\-фхы№\-эюх ьэю\-цх\-ёЄтю $Fr[\mu]$ ьхЁ√ $\mu$ ёюёЄю\-шЄ шч
хфшэ\-ёЄ\-тхэ\-эющ ьхЁ√ $\nu$, $d\nu(x)=\frac{c}{c_1}dx$, $c_1$ $=$
$\int\limits_0^\infty K(t)dt$.}

\par ╥хюЁхьр  1.4 хёЄ№ ЄЁштшры№\-эюх ёыхф\-ёЄтшх Єхю\-Ёх\-ь√  1.3,
т ў╕ь ыхуъю єсхфшЄ№\-ё , хёыш яю\-ёьюЄЁхЄ№ т ЄхъёЄх ёЄрЄ№ш эр
Єюў\-эюх юяЁх\-фхых\-эшх яЁх\-фхы№\-эю\-ую ьэю\-цх\-ёЄ\-тр, уфх
Ёрё\-°шЇЁю\-т√\-тр\-хЄё  ёююЄ\-эю°х\-эшх $\mu_{t_n}\to \nu$. ┬ ётю■
юўхЁхф№, Єхю\-Ёх\-ьр  1.3 хёЄ№ ёыхф\-ёЄтшх Єхю\-Ёх\-ь√  1.4 ш
ёююЄ\-тхЄ\-ёЄ\-тє■\-∙хщ рсхых\-тющ Єхю\-Ёх\-ь√.

\par ╨рё\-ёьюЄЁшь Єрє\-сх\-Ёю\-тє Єхю\-Ёх\-ьє, фю\-ърчрэ\-эє■ т Ёрсю\-Єх.

\par {\sc  ╥хюЁхьр 5.8.}\; {\it ╧єёЄ№ $\rho(r)$ -- яЁюшч\-тюы№\-э√щ
єЄюў\-э╕э\-э√щ яю\-Ё фюъ, $\mu$  -- Ёр\-фю\-эю\-тр ьхЁр эр яюыє\-юёш
$(0,\infty)$  шч ъырёёр $\frak{M}(\rho(r))$, $K$  -- сюЁхыхт\-ёър 
Їєэъ\-Ўш  эр яюыє\-юёш $(0,\infty)$ Єр\-ър , ўЄю
$t^{\rho-1}\gamma(t)K(t)\in L_1(0,\infty)$. ╧єёЄ№ Їєэъ\-Ўш 
$\int\limits_0^\infty K(t)t^{\rho-1+i\lambda}dt$ эх юсЁр∙р\-хЄ\-ё  т
эюы№ эр тх∙хёЄ\-тхэ\-эющ юёш ш яєёЄ№ Їєэъ\-Ўш  $\Psi(r)$
юяЁх\-фхы \-хЄ\-ё  Ёртхэ\-ёЄ\-тюь (\ref{introduction_Psi(r)}). ╥юуфр
хёыш ьхЁр $s$, $ds(t)=\Psi(t)dt$,  ты \-хЄё  Ёхуєы Ё\-эющ ьхЁющ
юЄ\-эюёш\-Єхы№\-эю єЄюў\-э╕э\-эю\-ую яю\-Ё ф\-ър $\rho(r)+1$, Єю
ьхЁр $\mu$  ты \-хЄё  Ёхуєы Ё\-эющ ьхЁющ юЄ\-эюёш\-Єхы№эю
єЄюў\-э╕э\-эю\-ую яю\-Ё ф\-ър  $\rho(r)$, яЁш\-ў╕ь хёыш $Fr[s]$
ёюёЄю\-шЄ шч ьхЁ√ ё яыюЄ\-эюёЄ№■ $ct^\rho$, Єю $Fr[\mu]$ ёюёЄю\-шЄ
шч ьхЁ√ ё яыюЄ\-эюёЄ№■ $\frac{c}{c_1} t^{\rho-1}$, уфх $c_1$ $=$
$\int\limits_0^\infty K(t)t^{\rho-1}dt$.}

\par {\bf 1.}  ┬ Єхю\-Ёхьх 1.4 Ёрё\-ёьрЄЁш\-трхЄё  ёыєўрщ  $\rho(r)\equiv 1$,
т Єхю\-Ёх\-ьх  5.8 Ёрё\-ёьрЄ\-Ёш\-тр\-хЄ\-ё  яЁю\-шч\-тюы№\-э√щ
єЄюў\-э╕э\-э√щ яю\-Ё фюъ. ╬фэръю ш ўрёЄ\-э√щ ёыєўрщ Єхю\-Ёх\-ь√ 5.8,
ъюуфр $\rho(r)\equiv 1$, чэрўш\-Єхы№\-эю ёшы№\-эхх Єхю\-Ёх\-ь√  1.4.
╨рё\-ёьюЄЁшь сюыхх яю\-фЁюс\-эю ¤ЄюЄ ёыєўрщ.

\par {\bf 2.} ╟эрўшЄхы№эю юёырсы \-хЄё  ЄЁхсю\-тр\-эшх эр  фЁю  $K$. ╧Ёш ¤Єюь
фю\-яєёър\-■Єё  Ёрч\-Ё√т√ є  фЁр  $K$ ш ыю\-ъры№\-эр 
эх\-юуЁрэш\-ўхэ\-эюёЄ№  фЁр  $K$. ┬ Єхю\-Ёх\-ьх  1.4 ЄЁхсє\-хЄё ,
ўЄюс√  фЁю  $K$ с√ыю эх\-яЁх\-Ё√т\-э√ь ш ёїюфшы\-ё  Ё ф
$\sum\limits_{n=-\infty}^\infty K_n e^n$. ┬ Ёрё\-ёьрЄЁш\-трх\-ьюь
ўрёЄ\-эюь ёыєўрх Єхю\-Ёх\-ь√  5.8 ¤Єю ЄЁхсю\-тр\-эшх юёырсы \-хЄ\-ё 
фю ЄЁхсю\-тр\-эш  $K(t)\in L_1(0,\infty)$.

\par {\bf 3.} ╥Ёхсютр\-эшх ёє\-∙хёЄ\-тю\-тр\-эш  яЁх\-фхыр
$\lim\limits_{r\to\infty} \frac{1}{r}\Psi(r)$ т
Ёрё\-ёьрЄЁш\-трх\-ьюь трЁш\-рэ\-Єх Єхю\-Ёх\-ь√  5.8 юёырсы \-хЄ\-ё ,
уЁєсю уютюЁ  (Єюў\-э√х ЇюЁ\-ьє\-ыш\-Ёют\-ъш ёь. т Єхю\-Ёх\-ьрї 3.19,
3.20), фю ЄЁхсю\-тр\-эш  ёє\-∙хёЄ\-тю\-тр\-эш  яЁх\-фхыр
$\lim\limits_{r\to\infty} \frac{1}{r^2} \int\limits_1^r \Psi(t)dt$.

\par ┬ьхёЄх ё Єхь чр\-ъы■ўх\-эш  т Єхю\-Ёх\-ьх  1.4 ш Ёрё\-ёьрЄЁш\-трх\-ьюь
трЁш\-рэ\-Єх Єхю\-Ёх\-ь√  5.8 юфшэр\-ъю\-т√х.

\par ╥хю\-Ёхьр 5.8  ты \-хЄё  Єръ\-цх єёшых\-эшхь Єхю\-Ёх\-ь√ ┴шэуїхьр,
├юыфш ш ╥хєухы\-ёр  \cite{Bingham}, Ёрч\-фхы  4.9, Єхю\-Ёх\-ьр
4.9.1. ╧юёыхф\-э   Єхю\-Ёх\-ьр Єръ\-цх ЎшЄш\-Ёє\-хЄ\-ё  т
\cite{Korevaar}, уыртр  4, Єхю\-Ёх\-ьр  9.3. ┬ Єхю\-Ёх\-ьх
┴шэу\-їхьр, ├юыфш ш ╥хєухыёр Єръ цх, ъръ ш т Єхю\-Ёх\-ьх  5.8,
Ёрё\-ёьрЄЁш\-тр\-хЄё  яЁю\-шч\-тюы№\-э√щ єЄюў\-э╕э\-э√щ яю\-Ё фюъ.
╬фэръю т єърчрэ\-эющ Єхю\-Ёх\-ьх Ёрё\-ёьрЄЁш\-тр\-■Є\-ё  Єюы№ъю
яюыю\-цш\-Єхы№\-э√х ьхЁ√ ш ЄЁхсє\-хЄё  ъю\-эхў\-эюёЄ№
рьры№\-урь-эюЁь√  фЁр  $K$.

\par ═хёъюы№ъю ёыют ю фю\-ърчр\-Єхы№\-ёЄтх Єхю\-Ёх\-ь√  5.8. ╥хю\-Ёх\-ьр  4.8
ётюфшЄ фю\-ърчр\-Єхы№\-ёЄтю Єхю\-Ёх\-ь√  5.8 ъ фю\-ърчр\-Єхы№\-ёЄтє
хфшэ\-ёЄтхэ\-эюё\-Єш Ёх°х\-эш  шэ\-Єху\-Ёры№\-эю\-ую єЁртэх\-эш 
$$
\int\limits_0^\infty K\left(\frac{t}{r}\right)d\nu(t)=cr^\rho
$$
ё эх\-шчтхёЄ\-эющ ьхЁющ  $\nu$. ╥ръюх фю\-ърчр\-Єхы№\-ёЄтю
яЁю\-тюфшЄ\-ё  ьхЄю\-фюь ╩рЁых\-ьрэр. ╧Ёш ¤Єюь яЁш\-їюфшЄ\-ё 
єёшыш\-трЄ№ ыхььє ╩рЁых\-ьрэр юс рэрыш\-Єш\-ўхё\-ъюь
яЁю\-фюы\-цх\-эшш.

\par ┬ Ёрсю\-Єх Єръ\-цх фю\-ърч√\-тр\-хЄё  трЁш\-рэЄ Єхю\-Ёх\-ь√  5.8, уфх
Їєэъ\-Ўшш $\int\limits_0^\infty K(t)t^{\rho-1+i\lambda}dt$
Ёрч\-Ёх°р\-хЄё  юсЁр∙рЄ№\-ё  т эюы№ эр тх∙хёЄ\-тхэ\-эющ юёш т
ъю\-эхў\-эюь ўшёых Єюўхъ.

 %%%%%%%%%%%%%%%%%%%%%%%%%%%%%%%%%%%%%%%%%%%%%%%%%%%%%%%%%%%%%%%%%%%%%%%%%%%%%%%%%%%%%%%%%
\section{╬с єЄюў\-э╕ээюь яю\-Ё фъх}
%%%%%%%%%%%%%%%%%%%%%%%%%%%%%%%%%%%%%%%%%%%%%%%%%%%%%%%%%%%%%%%%%%%%%%%%%%%%%%%%%%%%%%%%%
\qquad ╧єёЄ№ $f(r)$ -- яюыю\-цш\-Єхы№\-эр  Їєэъ\-Ўш  эр яюыє\-юёш
$(0,\infty)$ ш яхЁхф эрьш ёЄюшЄ чрфрўр юяшёрЄ№
рёшья\-Єю\-Єш\-ўхё\-ъюх яю\-тхфх\-эшх Їєэъ\-Ўшш $f$ эр
схё\-ъю\-эхў\-эюё\-Єш. ┬рцэющ ўшёыю\-тющ їрЁръ\-Єх\-Ёшё\-Єш\-ъющ
Їєэъ\-Ўшш  $f$  ты \-хЄё  х╕ {\it яю\-Ё фюъ} $\rho$, ъюЄю\-Ё√щ
юяЁх\-фхы \-хЄ\-ё  яю ЇюЁ\-ьєых
$$
\rho=\mathop{\overline{\lim}}\limits_{r\to\infty} \frac{\ln
f(r)}{\ln r}.
$$
┬ юс∙хь ёыєўрх тхыш\-ўшэр $\rho$ хёЄ№ ¤ыхьхэЄ Ёрё\-°шЁхэ\-эющ
ўшёыю\-тющ яЁ ьющ $[-\infty,\infty]$. ╤ююЄ\-эю°х\-эшхь
$\rho\in(-\infty,\infty)$ т√\-фхы \-хЄё  трц\-э√щ ъырёё Їєэъ\-Ўшщ --
Їєэъ\-Ўшщ ъю\-эхў\-эю\-ую яю\-Ё ф\-ър. ┬ фры№\-эхщ\-°хь ь√ сєфхь
Ёрё\-ёьрЄЁш\-трЄ№ Їєэъ\-Ўшш ъю\-эхў\-эю\-ую яю\-Ё ф\-ър.

\par ┼ёыш $\rho$ -- яю\-Ё фюъ Їєэъ\-Ўшш $f$, р $\varepsilon$ --
яЁю\-шч\-тюы№\-эюх ёЄЁюую яю\-ыюцш\-Єхы№\-эюх ўшёыю, Єю
т√\-яюы\-э \-хЄ\-ё  ёых\-фє\-■\-∙р  ёшёЄх\-ьр эх\-Ёр\-тхэёЄт
\begin{equation}\label{order-ineq-f}
\begin{split}
& f(r)<r^{\rho+\varepsilon},\qquad r\geq R(\varepsilon),\\
& f(r)>r^{\rho-\varepsilon},\qquad r\in E,
\end{split}
\end{equation}
уфх $E$ -- эх\-ъюЄю\-Ёюх, чр\-тшё \-∙хх юЄ  $\varepsilon$ ш  $f$,
эх\-юуЁрэш\-ўхэ\-эюх ьэю\-цх\-ёЄтю. ┼ёыш т ъюэ\-ъЁхЄ\-эющ чрфрўх
ёшёЄхьр эх\-Ёр\-тхэёЄт  (\ref{order-ineq-f})  ты \-хЄё  ёыш°ъюь
уЁєсющ, Єю яЁш\-їюфшЄ\-ё  ттюфшЄ№ сюыхх Єюэ\-ъшх їрЁръ\-ЄхЁшё\-Єшъш
ЁюёЄр Їєэъ\-Ўшш  $f$, ўхь яю\-Ё фюъ.

\par {\it ╥шяюь} Їєэъ\-Ўшш $f$ яЁш яю\-Ё фъх $\rho$ эр\-ч√тр\-хЄ\-ё 
тхыш\-ўшэр
$$
\sigma=\mathop{\overline{\lim}}\limits_{r\to\infty}
\frac{f(r)}{r^\rho}.
$$
\par ╩ръ яю\-ърч√\-трхЄ яЁш\-ьхЁ $f(r)=Ar^\rho (\ln(e+r))^\beta$,
$A>0$, $\beta\in(-\infty,\infty)$, фы  Їєэъ\-Ўшщ яю\-Ё ф\-ър $\rho$
тхыш\-ўшэр $\sigma$ ью\-цхЄ с√Є№ яЁю\-шч\-тюы№\-э√ь ¤ыхьхэ\-Єюь
ьэю\-цх\-ёЄтр $[0,\infty]$.

\par ┬ чр\-тшёш\-ьюёЄш юЄ Єюую, ъръюх шч ёююЄ\-эю°х\-эшщ
$\sigma=0$, $\sigma\in(0,\infty)$, $\sigma=\infty$
т√\-яюыэ \-хЄ\-ё , Їєэъ\-Ўш  $f$ эр\-ч√тр\-хЄё  Їєэъ\-Ўшхщ {\it
ьшэш\-ьры№\-эю\-ую, эюЁьры№\-эю\-ую} шыш {\it ьръёш\-ьры№\-эю\-ую}
Єшяр яЁш яю\-Ё фъх  $\rho$. ┼ёыш $\sigma<\infty$, Єю фы  ы■сю\-ую
$\varepsilon>0$ т√\-яюыэ \-■Єё  эх\-Ёртхэ\-ёЄ\-тр
\begin{equation}\label{order-f-sigma}
\begin{split}
& f(r)<(\sigma+\varepsilon)r^\rho,\qquad r\geq R(\varepsilon),\\
& f(r)>(\sigma-\varepsilon)r^\rho,\qquad r\in E,
\end{split}
\end{equation}
уфх $E$ -- эх\-ъюЄю\-Ёюх эх\-юуЁрэш\-ўхэ\-эюх ьэю\-цх\-ёЄтю,
чр\-тшё \-∙хх юЄ  $f$ ш  $\varepsilon$. ═х\-Ёртхэ\-ёЄ\-тр
(\ref{order-f-sigma}) -- ¤Єю чэрўш\-Єхы№\-эю сюыхх Єюўэ√х
эх\-Ёртхэ\-ёЄ\-тр, ўхь эх\-Ёртхэ\-ёЄ\-тр  (\ref{order-ineq-f}).

\par ┴єфхь уютюЁшЄ№, ўЄю Їєэъ\-Ўш  $f(r)$ Ёрё\-Є╕Є эр схё\-ъю\-эхў\-эюё\-Єш
ъръ Їєэъ\-Ўш  $\varphi(r)$, хёыш ёє\-∙хёЄ\-тєхЄ яю\-ёЄю э\-э√х $a$ ш
$b$, $0<a<b$, Єръшх, ўЄю т√\-яюы\-э \-■Є\-ё  эх\-Ёртхэ\-ёЄ\-тр
$$
\begin{array}{lcl}
f(r)<b\varphi(r),\qquad r> R,\\
f(r)>a\varphi(r),\qquad r\in E,
\end{array}
$$
уфх $E$  -- эх\-ъюЄюЁюх эх\-юуЁрэш\-ўхээюх ьэю\-цх\-ёЄтю.

\par ┼ёыш Їєэъ\-Ўш  $f(r)$  ты \-хЄё  Їєэъ\-Ўшхщ эюЁ\-ьры№\-эю\-ую Єшяр
яЁш яю\-Ё ф\-ъх  $\rho$, Єю, ъръ ёых\-фє\-хЄ шч эх\-Ёр\-тхэёЄт
(\ref{order-f-sigma}), Їєэъ\-Ўш  $f(r)$ Ёрё\-Є╕Є эр
схё\-ъю\-эхў\-эюё\-Єш ъръ Їєэъ\-Ўш  $r^\rho$.

\par ┼ёыш Їєэъ\-Ўш  $f(r)$  ты \-хЄё  Їєэъ\-Ўшхщ ьшэш\-ьры№\-эю\-ую
шыш ьръёш\-ьры№\-эю\-ую Єшяр яЁш яю\-Ё фъх $\rho$, Єю ¤Єющ
шэ\-ЇюЁьр\-Ўшш эх\-фю\-ёЄр\-Єюў\-эю, ўЄюс√ єърчрЄ№ фю\-ёЄр\-Єюў\-эю
яЁюё\-Єє■ Їєэъ\-Ўш■ $\varphi(r)$ Єръє■, ўЄюс√ Їєэъ\-Ўш  $f(r)$ Ёюёыр
эр схё\-ъю\-эхў\-эюё\-Єш ъръ Їєэъ\-Ўш  $\varphi(r)$.

\par ┬ ёт чш ёю ёърчрэ\-э√ь т√°х хёЄхёЄ\-тхэ\-эю тюч\-эшър\-хЄ ёых\-фє\-■\-∙р 
чрфрўр, эр\-чю\-т╕ь х╕ чрфр\-ўхщ $A$. ═хюс\-їюфшью єърчрЄ№ ъырёё
Їєэъ\-Ўшщ  $\mathfrak{A}$, ёюёЄю \-∙шщ шч фю\-ёЄр\-Єюўэю яЁюё\-Є√ї
Їєэъ\-Ўшщ, яю ётюшь ётющ\-ёЄ\-трь яюїю\-цшї эр Їєэъ\-Ўшш $r^\rho$ ш
Єр\-ъющ, ўЄюс√ фы  ы■сющ Їєэъ\-Ўшш  $f$ ъю\-эхў\-эю\-ую яю\-Ё ф\-ър
т ъырё\-ёх  $\mathfrak{A}$ эр°\-ырё№ с√ Єр\-ър  Їєэъ\-Ўш 
$\varphi(r)$, ўЄюс√ Їєэъ\-Ўш  $f(r)$ Ёюёыр эр схё\-ъю\-эхў\-эюё\-Єш
ъръ Їєэъ\-Ўш  $\varphi(r)$.

\par ╩ръ ь√ єцх чэрхь, ъырёё Їєэъ\-Ўшщ, ёюёЄю \-∙шщ шч Їєэъ\-Ўшщ
$r^\rho$, эх  ты \-хЄё  Єръшь ъырё\-ёюь. ╧єёЄ№ $\ln_k r$ -- $k$-Єр 
шЄхЁр\-Ўш  ыюур\-ЁшЇьр, эр\-яЁш\-ьхЁ, $\ln_2 r=\ln\ln r$.
╬сю\-чэрўшь ўхЁхч $e_k$ Єръє■ яю\-ёыхфю\-тр\-Єхы№\-эюёЄ№: $e_1=e$,
$e_{k+1}=e^{e_k}$. ╠юцэю яЁю\-тхЁшЄ№, ўЄю ъырёё, ёю\-ёЄю \-∙шщ шч
Їєэъ\-Ўшщ тшфр
\begin{equation}\label{order-varphi}
\varphi(r)=r^\rho (\ln(r+e))^{\alpha_1}...(\ln_k(r+e_k))^{\alpha_k},
\end{equation}
 уфх $k$-- яЁю\-шч\-тюы№\-эюх эрЄєЁры№\-эюх ўшёыю, р $\rho$, $\alpha_1$,...,
 $\alpha_k$  -- яЁю\-шч\-тюы№\-э√х тх∙хёЄ\-тхэ\-э√х ўшёыр, Єръ\-цх эх  ты \-хЄ\-ё 
 Єръшь ъырёёюь. ┬шфэю, ўЄю ёЇюЁ\-ьє\-ыш\-Ёю\-трэ\-эр  чрфрўр эх  ты \-хЄ\-ё 
 ЄЁштш\-ры№\-эющ. ╧єЄ№ ъ Ёх°х\-эш■ ¤Єющ чрфрўш єърчры ┬рышЁюэ  \cite{Valiron}.

\par  └сёюы■Єэю эх\-яЁх\-Ё√тэр  Їєэъ\-Ўш   $\rho(r)$ эр яюыє\-юёш $(0,\:\infty)$
  эр\-ч√тр\-хЄ\-ё  {\it єЄюў\-э╕э\-э√ь яю\-Ё ф\-ъюь} (т ёь√ёых ┬рыш\-Ёюэр),
  хёыш т√\-яюыэ \-■Єё  фтр єёыю\-тш  :\\
   1) $\lim\limits_{r\to+\infty}\rho(r)
   =\rho(\infty)= \rho\in(-\infty,\infty)$,\\
   2) $\lim\limits_{r\to+\infty}r\ln{r}\rho'(r)=0$.\\
╟рьхЄшь, ўЄю яюф $\rho'(r)$ ёых\-фєхЄ яюэш\-ьрЄ№ ьръёш\-ьры№\-эюх яю
ьюфєы■ яЁю\-шч\-тюф\-эюх ўшёыю.

\par  ╠√ сєфхь шё\-яюы№\-чю\-трЄ№ ёых\-фє\-■\-∙хх юсю\-чэрўх\-эшх
$V(r)=r^{\rho(r)}$. ╟рьхЄшь, ўЄю $V(1)$ $=1$.

\par ╬ЄьхЄшь ёых\-фє\-■\-∙хх ётющ\-ёЄ\-тю єЄюў\-э╕э\-эю\-ую яю\-Ё ф\-ър
(ёьюЄЁш, эр\-яЁш\-ьхЁ,  \cite{Levin} уыртр  1, \S  12).

%%%%%%%%%%%%%%%%%%%%%%%%%%%%%%%%%%%%%%%%%%%%%%%%%%%%%%%%%%%%%%%%%%%%%%%%%%%%%%%%%%%%%%%
\begin{theorem}\label{order_th_V:V}\hskip-2mm{.}\:
%%%%%%%%%%%%%%%%%%%%%%%%%%%%%%%%%%%%%%%%%%%%%%%%%%%%%%%%%%%%%%%%%%%%%%%%%%%%%%%%%%%%%%%
╧єёЄ№ $\rho(r)$ -- яЁю\-шч\-тюы№\-э√щ єЄюў\-э╕э\-э√щ яю\-Ё фюъ,
$\rho=\rho(\infty)$. ╥юуфр фы  ы■сюую $t>0$ ёє\-∙хёЄ\-тєхЄ яЁх\-фхы
$$
\lim\limits_{r\to\infty} \frac{V(tr)}{V(r)}=t^\rho,
$$
ш ¤ЄюЄ яЁх\-фхы  ты хЄё  Ёртэю\-ьхЁ\-э√ь эр ы■сюь ёху\-ьхэЄх
$[a,b]\subset(0,\infty)$.
\end{theorem}
%%%%%%%%%%%%%%%%%%%%%%%%%%%%%%%%%%%%%%%%%%%%%%%%%%%%%%%%%%%%%%%%%%%%%%%%%%%%%%%%%%%%%%%
\par ╙Єюў\-э╕э\-э√щ яю\-Ё фюъ $\rho(r)$ эр\-ч√тр\-хЄ\-ё  єЄюў\-э╕э\-э√ь
яю\-Ё ф\-ъюь Їєэъ\-Ўшш $f(r)$, хёыш т√\-яюыэ \-хЄ\-ё 
ёююЄ\-эю°х\-эшх
$$
\mathop{\overline{\lim}}\limits_{r\to\infty}
\frac{f(r)}{V(r)}=\sigma\in(0,\infty).
$$
\par ┬√\-яюыэх\-эшх ¤Єюую ёююЄ\-эю°х\-эш  ¤ътш\-трыхэЄ\-эю єЄ\-тхЁ\-ц\-фх\-эш■,
ўЄю Їєэъ\-Ўш  $f(r)$ Ёрё\-Є╕Є эр схё\-ъю\-эхў\-эюё\-Єш ъръ Їєэъ\-Ўш 
$V(r)$. ╟рьхЄшь Єръ\-цх, ўЄю хёыш $\rho(r)$  -- єЄюў\-э╕э\-э√щ
яю\-Ё фюъ Їєэъ\-Ўшш $f(r)$, Єю ўшёыю $\rho=\rho(\infty)$
 ты \-хЄ\-ё  яю\-Ё ф\-ъюь Їєэъ\-Ўшш $f(r)$.

\par ┬рцэюёЄ№ яюэ Єш  єЄюў\-э╕э\-эю\-ую яю\-Ё ф\-ър ёых\-фєхЄ шч Єюую, ўЄю
ъырёё  $\frak{A}$, ёю\-ёЄрт\-ыхэ\-э√щ шч Їєэъ\-Ўшщ тшфр
$V(r)=r^{\rho(r)}$  ты \-хЄ\-ё  Ёх°х\-эшхь яю\-ёЄрт\-ыхэ\-эющ т√°х
чрфрўш  $A$. ╧юёыхф\-эхх т√\-Єхър\-хЄ шч ёых\-фє\-■\-∙хщ
Єхю\-Ёх\-ь√.
%%%%%%%%%%%%%%%%%%%%%%%%%%%%%%%%%%%%%%%%%%%%%%%%%%%%%%%%%%%%%%%%%%%%%%%%%%%%%%%%%%%%%%%
\begin{theorem}\label{order_rho(r)}\hskip-2mm{.}\:
%%%%%%%%%%%%%%%%%%%%%%%%%%%%%%%%%%%%%%%%%%%%%%%%%%%%%%%%%%%%%%%%%%%%%%%%%%%%%%%%%%%%%%%
╧єёЄ№ $f(r)$ -- Їєэъ\-Ўш  ъю\-эхў\-эю\-ую яю\-Ё ф\-ър $\rho$. ╥юуфр
ёє\-∙хёЄ\-тєхЄ єЄюў\-э╕э\-э√щ яю\-Ё фюъ $\rho(r)$ Єр\-ъющ, ўЄю
т√\-яюыэ \-■Єё  єёыю\-тш :\\
1) $\lim\limits_{r\to\infty}\rho(r)=\rho$,\\
2) Їєэъ\-Ўш  $\rho(r)$  ты \-хЄё  ьюэю\-Єюэ\-эющ Їєэъ\-Ўшхщ
эр яюыє\-юёш $[1,\infty)$, \\
3) т√\-яюыэ \-хЄё  эх\-Ёртхэ\-ёЄтю
$$
(r+e)\ln(r+e)|\rho'(r)|\leq |\rho(r)-\rho|,\qquad r\geq 1,
$$
 4) т√\-яюыэ \-хЄё  ёююЄ\-эю°х\-эшх
$$
\mathop{\overline{\lim}}\limits_{r\to\infty}
\frac{f(r)}{V(r)}=\sigma\in(0,\infty).
$$
\end{theorem}
%%%%%%%%%%%%%%%%%%%%%%%%%%%%%%%%%%%%%%%%%%%%%%%%%%%%%%%%%%%%%%%%%%%%%%%%%%%%%%%%%%%%%%%
\par ─юърчр\-Єхы№ёЄтю ¤Єющ Єхю\-Ёх\-ь√ ьюц\-эю эрщЄш
т  \cite{Grishin_Malyutina}. ╧Ёртфр, т  \cite{Grishin_Malyutina}
фю\-яюыэш\-Єхы№\-эю ЄЁхсє\-хЄ\-ё , ўЄюс√ Їєэъ\-Ўш  $f(r)$ с√ыр
эх\-яЁх\-Ё√т\-эющ. ╤хщўрё ь√ єтшфшь, ўЄю ¤Єю эх\-ёє\-∙хёЄ\-тхэ\-эюх
ЄЁхсю\-тр\-эшх.

\par ╧єёЄ№ $f(r)$ -- яЁю\-шч\-тюы№\-эр  Їєэъ\-Ўш  яю\-Ё ф\-ър $\rho$.
═х юуЁрэш\-ўштр  юс∙\-эюё\-Єш, ьюц\-эю ёўш\-ЄрЄ№, ўЄю Їєэъ\-Ўш 
$f(r)$  ты \-хЄё  юуЁрэш\-ўхэ\-эющ эр ы■сюь ёхуьхэ\-Єх $[0,N]$.
╧єёЄ№ $n\geq 0$ -- яЁю\-шч\-тюы№\-эюх Ўхыюх ўшёыю. ╬сючэрўшь
$m_n=\inf\{f(x):\: x\in[n,n+1]\}$, $M_n=\sup\{f(x):\:x\in[n,n+1]\}$,
$\alpha_n=n+\frac{1}{3}$, $\beta_n=n+\frac{2}{3}$. ╤ЄЁюшь Їєэъ\-Ўш■
$f_1(r)$ ёых\-фє\-■\-∙шь юсЁр\-чюь. ═р ърцфюь шч ёхуьхэ\-Єют
$[n,\alpha_n]$, $[\alpha_n,\beta_n]$, $[\beta_n,n+1]$ юэр ышэхщ\-эр
ш, ъЁюьх Єюую, т√\-яюыэ \-■Єё  Ёртхэ\-ёЄ\-тр $f_1(n)=f(n)$,
$f_1(\alpha_n)=m_n$, $f_1(\beta_n)=M_n$. ╘єэъ\-Ўш  $f_1(r)$
эх\-яЁх\-Ё√т\-эр эр яюыє\-юёш $[0,\infty)$. ╬ўхтшфэю, ўЄю ы■сющ
єЄюў\-э╕э\-э√щ яю\-Ё фюъ Їєэъ\-Ўшш $f_1(r)$  ты \-хЄё  Єръ\-цх
єЄюў\-э╕э\-э√ь яю\-Ё ф\-ъюь Їєэъ\-Ўшш $f(r)$.

\par ╬ЄьхЄшь, ўЄю Їєэъ\-Ўш  $\rho(r)$, ёє\-∙хёЄ\-тю\-тр\-эшх ъю\-Єю\-Ёющ
єЄ\-тхЁ\-ц\-фр\-хЄ\-ё  т Єхю\-Ёх\-ьх \ref{order_rho(r)},
юсыр\-фр\-хЄ эх\-ъюЄю\-Ё√\-ьш фю\-яюыэш\-Єхы№\-э√\-ьш
ётющёЄ\-тр\-ьш, ъюЄю\-Ё√ї эхЄ є яЁю\-шч\-тюы№\-э√ї єЄюў\-э╕э\-э√ї
яю\-Ё ф\-ъют. ▌Єю, тю-яхЁ\-т√ї, ётющ\-ёЄ\-тю ьюэю\-Єюэ\-эюё\-Єш
Їєэъ\-Ўшш $\rho(r)$, р тю-тЄюЁ√ї, єёыю\-тшх  3) шч ЄхъёЄр
Єхю\-Ёх\-ь√ -- ¤Єю сюыхх ёшы№\-эюх юуЁрэш\-ўхэшх эр Їєэъ\-Ўш■
$\rho(r)$, ўхь ЄЁхсю\-тр\-эшх $\lim\limits_{r\to+\infty}
r\ln{r}\rho'(r)$ $=0$ шч юяЁх\-фхых\-эш  єЄюў\-э╕э\-эю\-ую
яю\-Ё ф\-ър. ┬ёыхфёЄ\-тшх ¤Єюую, Єхю\-Ёх\-ьр \ref{order_rho(r)} эх
т√\-Єхър\-хЄ шч рэрыю\-ушў\-э√ї єЄ\-тхЁ\-ц\-фх\-эшщ шч \cite{Levin},
\cite{Bingham}.

\par ╬ЄьхЄшь х∙╕ ёт ч№ єЄюў\-э╕э\-эю\-ую яю\-Ё ф\-ър ё Їєэъ\-Ўш \-ьш
Ёхує\-ы Ё\-эю ьхэ ■\-∙ш\-ьш\-ё  т ёь√ёых ╩рЁр\-ьр\-Є√.

\par ╧юыюцш\-Єхы№\-эр  эр яюыє\-юёш $(0,\infty)$ Їєэъ\-Ўш  $f(r)$
эр\-ч√тр\-хЄ\-ё  {\it Ёхуєы Ё\-эю ьхэ ■\-∙хщё  т ёь√ёых
╩рЁр\-ьр\-Є√}, хёыш фы  ы■сю\-ую $\lambda>0$ ёє\-∙хёЄ\-тє\-хЄ
ъю\-эхў\-э√щ яЁх\-фхы
$$
\lim\limits_{r\to\infty} \frac{f(\lambda r)}{f(r)}.
$$
\par ╤яЁр\-тхфыштр ёых\-фє\-■∙р  Єхю\-Ёхьр.
%%%%%%%%%%%%%%%%%%%%%%%%%%%%%%%%%%%%%%%%%%%%%%%%%%%%%%%%%%%%%%%%%%%%%%%%%%%%%%%%%%%%%%%
\begin{theorem}\label{order_Karamata}\hskip-2mm{.}\:
%%%%%%%%%%%%%%%%%%%%%%%%%%%%%%%%%%%%%%%%%%%%%%%%%%%%%%%%%%%%%%%%%%%%%%%%%%%%%%%%%%%%%%%
┼ёыш $f(r)$ -- шчьхЁш\-ьр  Їєэъ\-Ўш , Ёхує\-ы Ё\-эю ьхэ ■\-∙р \-ё  т
ёь√ёых ╩р\-Ёр\-ьр\-Є√, Єю ёє\-∙хёЄ\-тє\-хЄ Їєэъ\-Ўш  $C(r)\to 1$
$(r\to\infty)$ ш єЄюў\-э╕э\-э√щ яю\-Ё фюъ $\rho(r)$ Єръшх, ўЄю
$f(r)=C(r)V(r)$.
\end{theorem}
%%%%%%%%%%%%%%%%%%%%%%%%%%%%%%%%%%%%%%%%%%%%%%%%%%%%%%%%%%%%%%%%%%%%%%%%%%%%%%%%%%%%%%%
\par ▌Єю їюЁю°ю шчтхёЄ\-эр  Єхю\-Ёхьр ю яЁхф\-ёЄрт\-ых\-эшш фы  Ёхує\-ы Ё\-эю
ьхэ ■\-∙шї\-ё  \\Їєэъ\-Ўшщ (ёьюЄЁш, эр\-яЁш\-ьхЁ,  \cite{Bingham},
Єхю\-Ёх\-ьр  1.3.1).

\par ╙Єюў\-э╕ээ√щ яю\-Ё фюъ $\rho(r)$ эр\-ч√тр\-хЄё  {\it эєых\-т√ь
єЄюў\-э╕э\-э√ь яю\-Ё ф\-ъюь}, хёыш\\
$\lim\limits_{r\to\infty}\rho(r)$ $=0$.

\par ┼ёыш $\rho(r)$ -- яЁю\-шч\-тюы№\-э√щ єЄюў\-э╕э\-э√щ яю\-Ё фюъ,
Єю $\rho(r)=$ $\rho+$ $\hat{\rho}(r)$, уфх $\rho=$
$\lim\limits_{r\to\infty} \rho(r)$, р $\hat{\rho}(r)$ -- эєых\-тющ
єЄюў\-э╕э\-э√щ яю\-Ё фюъ.

\par ╤ь√ёы ттхфх\-эш  єЄюў\-э╕э\-эю\-ую яю\-Ё ф\-ър ёюёЄю\-шЄ т Єюь, ўЄюс√
фы  ы■сющ Їєэъ\-Ўшш $f(r)$ ъю\-эхў\-эю\-ую яю\-Ё ф\-ър эрщЄш
Їєэъ\-Ўш■ $V(r)=r^{\rho(r)}$ Єръє■, ўЄюс√ Їєэъ\-Ўш  $f(r)$ Ёюёыр эр
схё\-ъю\-эхў\-эюё\-Єш ъръ Їєэъ\-Ўш  $V(r)$. ╧Ёш Ёх°х\-эшш Єр\-ъющ
чрфрўш яю\-тхфх\-эшх Їєэъ\-Ўшш $\rho(r)$ т юъ\-ЁхёЄ\-эюё\-Єш эєы  эх
шуЁрхЄ эшър\-ъющ Ёюыш. ┬ьхёЄх ё Єхь, яЁш шё\-ёыхфю\-тр\-эшш
ъюэ\-ъЁхЄ\-э√ї чрфрў т Ёрч\-ышў\-э√ї юсырё\-Є ї ьрЄхьр\-Єш\-ъш ўрёЄю
тюч\-эшър\-■Є шэ\-Єху\-Ёры√ тшфр $\int\limits_0^\infty
K(t,r)V(t)dt$. ╧Ёш шчєўх\-эшш ётющёЄт Єръюую шэ\-Єху\-Ёр\-ыр
яю\-тхфх\-эшх Їєэъ\-Ўшш $\rho(t)$ т юъ\-ЁхёЄ\-эюё\-Єш эєы  Єръ цх
трцэю, ъръ ш яю\-тхфх\-эшх $\rho(t)$ т юъ\-ЁхёЄ\-эюё\-Єш
схё\-ъю\-эхў\-эюё\-Єш. ╧ю¤Єю\-ьє т ¤Єющ Ёрсю\-Єх ь√ сєфхь
яЁхф\-яюыр\-урЄ№, ўЄю фы  эєых\-тюую єЄюў\-э╕э\-эю\-ую яю\-Ё ф\-ър
$\rho(r)$ т√\-яюыэ \-хЄ\-ё  фю\-яюыэш\-Єхы№\-эюх єёыю\-тшх
$\rho\left(\frac{1}{r}\right)$ $=$ $-\rho(r)$. ▌Єю Ёртхэ\-ёЄ\-тю
¤ътш\-тр\-ыхэЄ\-эю ёых\-фє\-■\-∙х\-ьє $V\left(\frac{1}{r}\right)$
$=V(r)$. ▌Єюьє Ёртхэ\-ёЄтє, эр\-яЁш\-ьхЁ, єфют\-ыхЄ\-тю\-Ё \-хЄ
Їєэъ\-Ўш  $V(r)=1+|\ln r|^\alpha$, $\alpha>1$.
╤ююЄ\-тхЄ\-ёЄ\-тє■\-∙шщ єЄюў\-э╕э\-э√щ яю\-Ё фюъ чр\-фр╕Є\-ё 
ЇюЁ\-ьє\-ыющ
$$
\rho(r)=\frac{\ln(1+|\ln r|^\alpha)}{\ln r}.
$$
╟рьхЄшь, ўЄю яЁш $\alpha\leq 1$ Їєэъ\-Ўш , чр\-фртрх\-ьр  ¤Єшь
Ёртхэ\-ёЄтюь, эх сєфхЄ єЄюў\-э╕э\-э√ь яю\-Ё ф\-ъюь, Єръ ъръ т
юяЁх\-фхых\-эшш єЄюў\-э╕э\-эю\-ую яю\-Ё ф\-ър ЄЁхсє\-хЄ\-ё , ўЄюс√
Їєэъ\-Ўш  $\rho(r)$ с√ыр рсёю\-ы■Є\-эю эх\-яЁх\-Ё√т\-эющ эр
яюыє\-юёш $(0,\infty)$.

\par ╬ЄьхЄшь х∙╕, ўЄю эрышўшх Ёртхэ\-ёЄ\-тр
$\rho\left(\frac{1}{r}\right)$ $=$ $-\rho(r)$ т√\-фхы \-хЄ Єюўъє $1$
ёЁхфш фЁєушї Єюўхъ яюыє\-юёш $(0,\infty)$. ┬ ўрёЄ\-эюё\-Єш, шч
эр\-яшёрэ\-эю\-ую Ёртхэ\-ёЄ\-тр ёых\-фєхЄ, ўЄю $\rho(1)=0$.

\par ╧Ёш Ёрё\-ёьюЄЁх\-эшш єЄюў\-э╕э\-эю\-ую яю\-Ё ф\-ър эрЁ фє ё
Їєэъ\-Ўшхщ $\rho(r)$ яюыхч\-эю шё\-яюы№чю\-трЄ№ Їєэъ\-Ўш■ $\eta(r)$
$=$ $\rho(r)+r\ln r \rho'(r)$. ┬√\-яюыэ \-хЄ\-ё  Ёртхэ\-ёЄ\-тю
\begin{equation}\label{order_V_eta}
V(r)=e^{\int\limits_1^r \frac{\eta(t)}{t}dt}.
\end{equation}
╤яЁр\-тхф\-ыш\-тюёЄ№ ¤Єюую Ёртхэ\-ёЄ\-тр ыхуъю яЁю\-тхЁшЄ№ ё яюью∙№■
ыюурЁшЇ\-ьш\-Ёю\-тр\-эш  ш яю\-ёых\-фє\-■\-∙хую
фшЇ\-ЇхЁхэ\-Ўш\-Ёю\-тр\-эш .

\par ┼ёыш $\rho(r)$ -- эєых\-тющ єЄюў\-э╕э\-э√щ яю\-Ё фюъ, Єю $\eta(r)\to 0$
яЁш $r\to\infty$. ┼ёыш, ъЁюьх Єюую, $\rho\left(\frac{1}{r}\right)$
$=$ $-\rho(r)$, Єю $\eta\left(\frac{1}{r}\right)$ $=$ $-\eta(r)$ ш т
¤Єюь ёыєўрх шч ЇюЁ\-ьє\-ы√  (\ref{order_V_eta}) ыхуъю ёых\-фє\-хЄ,
ўЄю фы  ы■сюую $\varepsilon>0$ эр яюыє\-юёш $(0,\infty)$
т√\-яюыэ \-хЄ\-ё  эх\-Ёртхэ\-ёЄ\-тю
\begin{equation}\label{order-V_leq}
V(r)\leq M_\varepsilon(r^\varepsilon+r^{-\varepsilon}).
\end{equation}
▌Єю їюЄ  ш уЁєсюх эх\-Ёртхэ\-ёЄтю, эю т эх\-ъюЄю\-Ё√ї ёыєўр\- ї юэю
юърч√\-тр\-хЄ\-ё  яюыхч\-э√ь.

\par ┼ёыш $\eta(t)$ -- ыю\-ъры№\-эю шэ\-Єху\-Ёш\-Ёє\-хьр  Їєэъ\-Ўш  эр
яюыє\-юёш $[1,\infty)$, шьх■\-∙р  эєых\-тющ яЁх\-фхы эр
схё\-ъю\-эхў\-эюё\-Єш, Єю ёє\-∙хёЄ\-тє\-хЄ схё\-ъю\-эхў\-эю
фшЇ\-ЇхЁхэ\-Ўш\-Ёєх\-ьр  эр яюыє\-юёш $[1,\infty)$ Їєэъ\-Ўш 
$\eta_1(t)$, шьх■\-∙р  эєых\-тющ яЁх\-фхы эр схё\-ъю\-эхў\-эюё\-Єш ш
Єр\-ър , ўЄю сєфхЄ ёїю\-фшЄ\-ё  шэ\-Єху\-Ёры $\int\limits_1^\infty
\frac{|\eta(t)-\eta_1(t)|}{t}dt$. ─юяюыэш\-Єхы№\-эю ьюц\-эю
яю\-ЄЁхсю\-трЄ№, ўЄюс√ т√\-яюыэ \-ыюё№ Ёртхэ\-ёЄтю
$\int\limits_1^\infty \frac{\eta(t)-\eta_1(t)}{t}dt$ $=0$. ┬  ¤Єюь
ёыєўрх фы  Їєэъ\-Ўшш $V(r)$, юяЁх\-фхы \-хьющ ЇюЁ\-ьє\-ыющ
(\ref{order_V_eta}) сєфхЄ ёяЁртхф\-ыш\-тю Ёртхэ\-ёЄтю
$$
V(r)=C(r)V_1(r),
$$
уфх
$$
V_1(r)=e^{\int\limits_1^r \frac{\eta_1(t)}{t}dt}, \qquad
C(r)=e^{\int\limits_r^\infty \frac{\eta_1(t)-\eta(t)}{t}dt}.
$$
\par ╚ч ёърчрэ\-эю\-ую ыхуъю ёых\-фєхЄ Єр\-ър  Єхю\-Ёхьр.
%%%%%%%%%%%%%%%%%%%%%%%%%%%%%%%%%%%%%%%%%%%%%%%%%%%%%%%%%%%%%%%%%%%%%%%%%%%%%%%%%%%%%%%
\begin{theorem}\label{order_th_V_V_1}\hskip-2mm{.}\:
%%%%%%%%%%%%%%%%%%%%%%%%%%%%%%%%%%%%%%%%%%%%%%%%%%%%%%%%%%%%%%%%%%%%%%%%%%%%%%%%%%%%%%%
╧єёЄ№ $\rho(r)$ -- яЁю\-шч\-тюы№\-э√щ эєыхтющ єЄюў\-э╕э\-э√щ
яю\-Ё фюъ Єр\-ъющ, ўЄю т√\-яюыэ \-хЄё  Ёртхэ\-ёЄтю
$\rho\left(\frac{1}{r}\right)$ $=$ $-\rho(r)$. ╥юуфр т√\-яюыэ \-хЄё 
Ёртхэ\-ёЄтю
\begin{equation}\label{order_V=cV_1}
V(r)=C(r)V_1(r),
\end{equation}
уфх $V_1(r)=r^{\rho_1(r)}$, $\rho_1(r)$ -- эєых\-тющ єЄюў\-э╕э\-э√щ
яю\-Ё фюъ, схё\-ъю\-эхў\-эю фшЇ\-Їх\-Ёхэ\-Ўш\-Ёєх\-ь√щ эр
ьэю\-цх\-ёЄ\-тх $(0,\infty)\setminus \{1\}$, єфют\-ыхЄ\-тю\-Ё ■\-∙шщ
Ёртхэ\-ёЄтє $\rho_1\left(\frac{1}{r}\right)$ $=$ $-\rho_1(r)$, р
$C(r)$ -- эх\-яЁх\-Ё√т\-эр  эр яюыє\-юёш $(0,\infty)$ Їєэъ\-Ўш 
Єр\-ър , ўЄю $C(r)\to 1$ яЁш $r\to\infty$ ш $r\to 0$.
─юяюыэш\-Єхы№\-эю ьюц\-эю чр\-фртрЄ№ яЁю\-шч\-тюы№\-эю яю\-Ё фюъ
ьрыюё\-Єш Їєэъ\-Ўшш $C(r)-1$ яЁш $r\to\infty$ ш  $r\to 0$.
\end{theorem}
%%%%%%%%%%%%%%%%%%%%%%%%%%%%%%%%%%%%%%%%%%%%%%%%%%%%%%%%%%%%%%%%%%%%%%%%%%%%%%%%%%%%%%%
\par ─рыхх сєфхЄ шё\-ёыхфю\-трэр Їєэъ\-Ўш  $\gamma(t)$, ю ъю\-Єю\-Ёющ єцх
уютюЁш\-ыюё№ тю тёЄєяых\-эшш. ┬ ёых\-фє\-■\-∙хщ ыхььх яЁш\-тюф Є\-ё 
ыхуъю фю\-ърч√\-трх\-ь√х ётющёЄ\-тр Їєэъ\-Ўшш $\gamma(t)$.
%%%%%%%%%%%%%%%%%%%%%%%%%%%%%%%%%%%%%%%%%%%%%%%%%%%%%%%%%%%%%%%%%%%%%%%%%%%%%%%%%%%%%%%
\begin{lemma}\label{order_lemma_gamma_properties}\hskip-2mm{.}\:
%%%%%%%%%%%%%%%%%%%%%%%%%%%%%%%%%%%%%%%%%%%%%%%%%%%%%%%%%%%%%%%%%%%%%%%%%%%%%%%%%%%%%%%
╧єёЄ№ $\rho(r)$ -- эєыхтющ єЄюў\-э╕э\-э√щ яю\-Ё фюъ,
єфютыхЄ\-тю\-Ё ■\-∙шщ єёыю\-тш■ $\rho\left(\frac{1}{r}\right)$ $=$
$-\rho(r)$ ш яєёЄ№
$$
\gamma(t)=\gamma(\rho(\cdot),t)=\sup\limits_{r>0}\frac{V(rt)}{V(r)},
\qquad \underline{\gamma}(t)= \underline{\gamma}(\rho(\cdot),t)
=\inf\limits_{r>0} \frac{V(rt)}{V(r)}.
$$
╥юуфр ёяЁртхф\-ышт√ ёых\-фє\-■\-∙шх єЄ\-тхЁ\-ц\-фх\-эш :\\
\vspace{3mm}
% % % % % % % % % % % % % % % % % % % % % % % % % % % % % % % % % % % % % % % % %
1) $\gamma(t), \underline{\gamma}(t)\in(0,\infty)$,\\ \vspace{3mm}
% % % % % % % % % % % % % % % % % % % % % % % % % % % % % % % % % % % % % % % % %
2) $\underline{\gamma}(t)\leq \gamma(t),\qquad
\underline{\gamma}(1)=\gamma(1)=1,$ \\ \vspace{3mm}
% % % % % % % % % % % % % % % % % % % % % % % % % % % % % % % % % % % % % % % % %
3) $\gamma\left(\frac{1}{t}\right)=
\frac{1}{\underline{\gamma}(t)},\qquad
\gamma\left(\frac{1}{t}\right)= \gamma\left(\rho(\cdot),
\frac{1}{t}\right) =\gamma(-\rho(\cdot),t),$ \\ \vspace{3mm}
% % % % % % % % % % % % % % % % % % % % % % % % % % % % % % % % % % % % % % % % %
4) $\gamma(t_1t_2)\leq \gamma(t_1)\gamma(t_2),\qquad
\underline{\gamma}(t_1t_2) \geq \underline{\gamma}(t_1)
\underline{\gamma}(t_2),$  \\ \vspace{3mm}
% % % % % % % % % % % % % % % % % % % % % % % % % % % % % % % % % % % % % % % % %
5) $\gamma(t)\geq V(t), \qquad \underline{\gamma}(t)\leq V(t)$,\\
\vspace{3mm}
% % % % % % % % % % % % % % % % % % % % % % % % % % % % % % % % % % % % % % % % %
6) Їєэъ\-Ўшш $\gamma(t)$ ш  $\underline{\gamma}(t)$  ты \-■Єё 
эхяЁх\-Ё√т\-э√ьш Їєэъ\-Ўш ьш эр яюыє\-юёш $(0,\infty)$.
\end{lemma}
%%%%%%%%%%%%%%%%%%%%%%%%%%%%%%%%%%%%%%%%%%%%%%%%%%%%%%%%%%%%%%%%%%%%%%%%%%%%%%%%%%%%%%%
\par {\sc  ─юърчрЄхы№ёЄтю.}\:  ╙ЄтхЁцфхэшх 2) юўхтшфэю. ─рыхх шьххь
$$
\gamma\left(\frac{1}{t}\right)=\sup\limits_{r>0}
\frac{V\left(\frac{r}{t}\right)}{V(r)}=\sup\limits_{R>0}
\frac{V(R)}{V(tR)}= \frac{1}{\inf\limits_{R>0} \frac{V(tR)}{V(R)}}
=\frac{1}{\underline{\gamma}(t)},
$$
$$
\gamma\left(\frac{1}{t}\right)= \sup\limits_{R>0}
\frac{V(R)}{V(tR)}= \sup\limits_{R>0}
\frac{R^{\rho(R)}}{(tR)^{\rho(tR)}}= \sup\limits_{R>0}
\frac{(tR)^{-\rho(tR)}}{R^{-\rho(R)}}=\gamma(-\rho(\cdot),t).
$$
╥хь ёрь√ь єЄ\-тхЁ\-ц\-фх\-эшх  3) фюърчрэю.
\par ╚ч ЁртхэёЄтр $V(1)=1$ ёых\-фєхЄ єЄ\-тхЁ\-ц\-фхэшх 5).
\par ╥ръ ъръ $\lim\limits_{r\to\infty} \frac{V(rt)}{V(r)}=1$,
$\lim\limits_{r\to 0} \frac{V(rt)}{V(r)}=1$, Єю ёє\-∙хёЄ\-тє■Є ўшёыр
$r_1$, $r_2$, $0<r_1<r_2$ Єръшх, ўЄю яЁш
$r\in(0,r_1)\cup(r_2,\infty)$ сєфхЄ т√\-яюыэ Є№\-ё  эх\-Ёртхэ\-ёЄтю
$\frac{V(rt)}{V(r)}\leq 2$. ╥ръ ъръ Їєэъ\-Ўш  $\frac{V(rt)}{V(r)}$
эх\-яЁх\-Ё√тэр яЁш $r\in[r_1,r_2]$, Єю ёє\-∙хёЄ\-тєхЄ ўшёыю $M>0$
Єръюх, ўЄю яЁш $r\in[r_1,r_2]$ сєфхЄ т√\-яюыэ Є№\-ё 
эх\-Ёртхэ\-ёЄ\-тю $\frac{V(rt)}{V(r)}\leq M$. ╚ч ¤Єюую ёых\-фєхЄ,
ўЄю $\gamma(t)\leq\max(M,2)$. ┬ьхёЄх ё эх\-Ёртхэ\-ёЄ\-тюь
$\gamma(t)\geq V(t)$ ¤Єю фр╕Є $\gamma(t)\in(0,\infty)$. ─рыхх шч  3)
ёых\-фєхЄ, ўЄю $\underline{\gamma}(t)\in(0,\infty)$. ╥хь ёрь√ь
єЄ\-тхЁ\-ц\-фх\-эшх 1) фюърчрэю.

\par ╙ЄтхЁц\-фх\-эшх 4) юўхтшфэю.

\par ╬сючэрўшь
$$
\gamma(r,t)=\frac{V(tr)}{V(r)}.
$$

\par ╧єёЄ№ $[a,b]$ -- яЁю\-шч\-тюы№\-э√щ ёху\-ьхэЄ эр яюыє\-юёш $(0,\infty)$,
$\varepsilon$  -- яЁю\-шч\-тюы№\-эюх ёЄЁюую яюыю\-цш\-Єхы№\-эюх
ўшёыю. ╥ръ ъръ яЁхфхы√
$$
\lim\limits_{r\to\infty} \frac{V(rt)}{V(r)}=1,\qquad
\lim\limits_{r\to 0} \frac{V(rt)}{V(r)}=1
$$
яю Єхю\-Ёхьх  \ref{order_th_V:V} Ёртэю\-ьхЁ\-э√х эр ёху\-ьхэЄх
$[a,b]$, Єю ёє\-∙хёЄ\-тє■Є ўшёыр $r_3$ ш  $r_4$, $0<r_3<r_4$, Єръшх,
ўЄю яЁш ы■с√ї $t_1$ ш  $t_2\in[a,b]$ яЁш ы■с√ї
$r\in(0,r_3)\cup(r_4,\infty)$ сєфєЄ т√\-яюыэ Є№\-ё 
эх\-Ёртхэ\-ёЄ\-тр
\begin{equation}\label{order_gamma-0}
-\varepsilon\leq \gamma(r,t_2)-\gamma(r,t_1)\leq \varepsilon.
\end{equation}
╥ръ ъръ ЇєэъЎш  $\gamma(r,t)$ эх\-яЁх\-Ё√т\-эр т яЁ ью\-єуюы№\-эшъх
$[r_3,r_4]\times[a,b]$, Єю яю Єхю\-Ёхьх ╩рэЄюЁр юэр Ёртэю\-ьхЁ\-эю
эх\-яЁх\-Ё√т\-эр т  ¤Єюь яЁ ью\-єуюы№\-эшъх. ╧ю¤Єю\-ьє
ёє\-∙хёЄ\-тє\-хЄ $\delta>0$ Єръюх, ўЄю яЁш $t_1,t_2\in[a,b]$,
$|t_1-t_2|<\delta$, $r\in[r_3,r_4]$ сєфхЄ т√\-яюыэ Є№\-ё 
эх\-Ёртхэ\-ёЄ\-тю  (\ref{order_gamma-0}). ┬ьхёЄх ё Ёрэхх
фю\-ърчрэ\-э√ь ¤Єю фр╕Є, ўЄю эх\-Ёртхэ\-ёЄ\-тр (\ref{order_gamma-0})
сєфєЄ т√\-яюыэ Є№ё  фы  ы■с√ї $t_1,t_2\in [a,b]$, $|t_2-t_1|<\delta$
ш фы  ы■с√ї $r>0$.

\par ╧єёЄ№ $t_1,t_2\in [a,b]$, $|t_2-t_1|<\delta$. ─ы  ы■сюую $r>0$
т√\-яюыэ \-хЄё  эх\-Ёртхэ\-ёЄтю
$$
\gamma(t_2)-\gamma(t_1)\leq \gamma(t_2)-\gamma(r,t_1).
$$
╤є∙хёЄтєхЄ $r>0$ Єръюх, ўЄю
${\gamma(t_2)<\gamma(r,t_2)+\varepsilon}$. ┬ьхёЄх ё
(\ref{order_gamma-0}) ¤Єю фр╕Є
$\gamma(t_2)-\gamma(t_1)<2\varepsilon$. ┬  ¤Єюь эх\-Ёртхэ\-ёЄтх
ьюц\-эю яю\-ьхэ Є№ ьхёЄр\-ьш $t_1$ ш  $t_2$. ╧ю¤Єю\-ьє
$|\gamma(t_2)-\gamma(t_1)|<2\varepsilon$. ╚ч ¤Єюую ёых\-фє\-хЄ
эх\-яЁх\-Ё√т\-эюёЄ№ Їєэъ\-Ўшш $\gamma(t)$ эр ёху\-ьхэЄх $[a,b]$, р
чэрўшЄ ш эр яюыє\-юёш $(0,\infty)$. ╙ЄтхЁцфх\-эшх  6), р тьхёЄх ё
эшь ш ыхььр, фю\-ърчр\-э√.
%%%%%%%%%%%%%%%%%%%%%%%%%%%%%%%%%%%%%%%%%%%%%%%%%%%%%%%%%%%%%%%%%%%%%%%%%%%%%%%%%%%%%%%
\begin{theorem}\label{order_th_ln_gamma}\hskip-2mm{.}\:
%%%%%%%%%%%%%%%%%%%%%%%%%%%%%%%%%%%%%%%%%%%%%%%%%%%%%%%%%%%%%%%%%%%%%%%%%%%%%%%%%%%%%%%
╧єёЄ№ $\rho(r)$ -- яЁю\-шч\-тюы№\-э√щ эєыхтющ єЄюў\-э╕э\-э√щ
яю\-Ё фюъ Єр\-ъющ, ўЄю $\rho\left(\frac{1}{r}\right)$ $=$ $-\rho(r)$
ш яєёЄ№
\begin{equation}\label{order_ln_gamma_1}
\gamma(t)=\sup\limits_{r>0} \frac{V(tr)}{V(t)}.
\end{equation}
╥юуфр т√\-яюыэ \-■Єё  Ёртхэ\-ёЄ\-тр
$$
\lim\limits_{t\to\infty}\frac{\ln\gamma(t)}{\ln{t}}=0, \qquad
\lim\limits_{t\to\infty}
\frac{\ln\gamma\left(\frac{1}{t}\right)}{\ln{t}}=0.
$$
\end{theorem}
%%%%%%%%%%%%%%%%%%%%%%%%%%%%%%%%%%%%%%%%%%%%%%%%%%%%%%%%%%%%%%%%%%%%%%%%%%%%%%%%%%%%%%%
\par {\sc  ─юърчрЄхы№ёЄтю.}\: ╧єёЄ№ ЇєэъЎш  $V_1(r)$ юяЁх\-фхы \-хЄё 
Ёртхэ\-ёЄ\-тюь  (\ref{order_V=cV_1}). ╤яЁр\-тхф\-ышт√
эх\-Ёртхэ\-ёЄ\-тр
$$
\frac{1}{M}\leq \frac{V_1(r)}{V(r)}\leq M.
$$
\par ─рыхх шьххь
$$
\gamma(t)= \sup\limits_{r>0} \frac{V(rt)}{V(r)}= \sup\limits_{r>0}
\frac{V_1(rt)}{V_1(r)}\frac{V(rt)}{V_1(rt)}\frac{V_1(r)}{V(r)}\leq
M^2\sup\limits_{r>0}
\frac{V_1(rt)}{V_1(r)}=M^2\gamma(\rho_1(\cdot),t).
$$
╚ч ¤Єюую ёых\-фєхЄ, ўЄю Єхю\-Ёхьє фю\-ёЄрЄюў\-эю фю\-ърч√\-трЄ№ фы 
ёыєўр , ъюуфр єЄюў\-э╕э\-э√щ яю\-Ё фюъ $\rho(r)$  ты \-хЄё 
фшЇ\-Їх\-Ёхэ\-Ўш\-Ёєх\-ьющ Їєэъ\-Ўшхщ эр ьэю\-цх\-ёЄтх
$(0,\infty)\setminus\{1\}$. ┬ фры№\-эхщ\-°хь фю\-ърчр\-Єхы№\-ёЄтх
сєфхь ёўш\-ЄрЄ№, ўЄю ¤Єю єёыю\-тшх т√\-яюыэ \-хЄё .

\par ╬сючэрўшь $h(x)=\ln{V(e^x)}$. ╘єэъ\-Ўш  $h(x)$ сєфхЄ эх\-яЁх\-Ё√т\-эющ
ў╕Єэющ Їєэъ\-Ўшхщ фшЇ\-Їх\-Ёхэ\-Ўш\-Ёєх\-ьющ тё■фє чр тюч\-ьюц\-э√ь
шё\-ъы■ўх\-эшхь эєы . ╥ю, ўЄю $\rho(r)$  -- эєых\-тющ єЄюў\-э╕э\-э√щ
яю\-Ё фюъ, яЁш\-тюфшЄ ъ ёююЄ\-эю°х\-эш■
\begin{equation}\label{order_lim_h(x)}
\lim\limits_{x\to +\infty} \frac{h(x)}{x}=0.
\end{equation}
 ╙ёыютшх $\lim\limits_{r\to\infty} r\ln r\rho'(r)=0$, ъюЄюЁюх
¤ътш\-тр\-ыхэЄ\-эю єёыю\-тш■ $\lim\limits_{r\to\infty}
\frac{rV'(r)}{V(r)}=0$, яЁш\-тюфшЄ ъ ёююЄ\-эю°х\-эш■
\begin{equation}\label{order_lim_h'(x)}
\lim\limits_{x\to+\infty} h'(x)=0.
\end{equation}
┬√\-яюыэ \-хЄё  Ёртхэ\-ёЄтю
$$
\varphi(y)=\ln\gamma(e^y)=\sup\limits_{x\in(-\infty,\infty)}(h(x+y)-h(x)).
$$
╧хЁтюх єЄ\-тхЁ\-ц\-фх\-эшх Єхю\-Ёх\-ь√ ¤ътш\-тр\-ыхэЄэю Ёртхэ\-ёЄтє
\begin{equation}\label{order_lim_varphi_1}
\lim\limits_{y\to+\infty} \frac{\varphi(y)}{y}=0.
\end{equation}
\par ┼ёыш ¤Єю Ёртхэ\-ёЄтю эх тхЁэю, Єю ёє\-∙хёЄ\-тє■Є ўшёыю $a>0$,
яю\-ёыхфю\-тр\-Єхы№\-эюёЄ№ $y_n\to +\infty$ ш
яю\-ёыхфю\-тр\-Єхы№\-эюёЄ№ $x_n\in(-\infty,\infty)$ Єръшх, ўЄю
\begin{equation}\label{order_h-h}
|h(x_n+y_n)-h(x_n)|\geq ay_n.
\end{equation}
─юяюыэш\-Єхы№\-эю ьюц\-эю ёўш\-ЄрЄ№, ўЄю ёє\-∙хёЄ\-тєхЄ ¤ыхьхэЄ
$\alpha\in[-\infty,\infty]$ Єр\-ъющ, ўЄю $x_n\to\alpha$ яЁш
$n\to\infty$. ╧єёЄ№ $\varepsilon$ -- яЁю\-шч\-тюы№\-эюх ўшёыю шч
шэ\-ЄхЁ\-трыр $\left(0,\frac{1}{2}a\right)$.

\par ─юяєёЄшь, ўЄю $\alpha\in(-\infty,\infty)$. ╥юуфр фы  тёхї
фю\-ёЄрЄюў\-эю сюы№°шї $n$ сєфхЄ т√\-яюыэ Є№\-ё  эх\-Ёртхэ\-ёЄтю
$$
|h(x_n+y_n)-h(x_n)|\leq \varepsilon(x_n+y_n)+|h(\alpha)|+1 \leq
\varepsilon y_n+\varepsilon(|\alpha|+1)+|h(\alpha)|+1.
$$
▌Єю эх\-Ёртхэ\-ёЄтю яЁю\-Єш\-тю\-ЁхўшЄ эх\-Ёртхэ\-ёЄтє
(\ref{order_h-h}).

\par ─юяєёЄшь, ўЄю $\alpha=+\infty$. ╥юуфр ёє\-∙хёЄ\-тєхЄ ўшёыю
$\xi_n\in(x_n,x_n+y_n)$ Єръюх, ўЄю т√\-яюыэ \-хЄ\-ё  Ёртхэ\-ёЄтю
$h(x_n+y_n)-h(x_n)=h'(\xi_n)y_n$. ╥хяхЁ№ шч Ёртхэ\-ёЄ\-тр
(\ref{order_lim_h'(x)}) ёых\-фє\-хЄ, ўЄю фы  тёхї фю\-ёЄрЄюў\-эю
сюы№\-°шї $n$ т√\-яюыэ \-хЄ\-ё  эх\-Ёртхэ\-ёЄтю
$|h(x_n+y_n)-h(x_n)|\leq \varepsilon y_n$. ▌Єю эх\-Ёртхэ\-ёЄтю
яЁю\-Єштю\-ЁхўшЄ эх\-Ёртхэ\-ёЄ\-тє  (\ref{order_h-h}).

\par ╥хяхЁ№ фю\-яєёЄшь, ўЄю $\alpha=-\infty$. ╠юцэю фю\-яюыэш\-Єхы№\-эю
яЁхф\-яюыю\-цшЄ№, ўЄю ёє\-∙хёЄ\-тє\-хЄ ¤ыхьхэЄ
$\beta\in[-\infty,\infty]$ Єр\-ъющ, ўЄю $x_n+y_n\to\beta$ яЁш
$n\to\infty$.

\par ─юяєёЄшь, ўЄю $\beta\in(-\infty,\infty)$. ╥юуфр фы  тёхї фю\-ёЄр\-Єюўэю
сюы№°шї $n$ сєфєЄ т√\-яюыэ Є№\-ё  эх\-Ёртхэ\-ёЄ\-тр
$$
|h(x_n+y_n)-h(x_n)|\leq |h(\beta)|+1+\varepsilon|x_n|\leq
|h(\beta)|+1+\varepsilon|x_n+y_n|+\varepsilon y_n
$$
$$
\leq |h(\beta)|+1+\varepsilon(|\beta|+1)+\varepsilon y_n.
$$
▌Єю эх\-Ёртхэ\-ёЄтю яЁю\-Єштю\-ЁхўшЄ эх\-Ёртхэ\-ёЄтє
(\ref{order_h-h}).

\par ─юяєёЄшь, ўЄю $\beta=-\infty$. ╥юуфр ёє\-∙хёЄ\-тє\-хЄ ўшёыю
$\xi_n\in(x_n,x_n+y_n)$ Єръюх, ўЄю т√\-яюыэ \-хЄ\-ё  Ёртхэ\-ёЄ\-тю
$h(x_n+y_n)-h(x_n)=h'(\xi_n)y_n$. ╥хяхЁ№ шч Ёртхэ\-ёЄ\-тр
(\ref{order_lim_h'(x)}) ш ў╕Є\-эюё\-Єш Їєэъ\-Ўшш  $h$ ёых\-фє\-хЄ,
ўЄю фы  тёхї фю\-ёЄр\-Єюўэю сюы№°шї $n$ сєфхЄ т√\-яюыэ Є№\-ё 
эх\-Ёртхэ\-ёЄтю $|h(x_n+y_n)-h(x_n)|\leq \varepsilon y_n$. ▌Єю
эх\-Ёртхэ\-ёЄ\-тю яЁю\-Єш\-тю\-ЁхўшЄ эх\-Ёртхэ\-ёЄ\-тє
(\ref{order_h-h}).

\par ─юяєёЄшь ЄхяхЁ№, ўЄю $\beta=+\infty$. ╥юуфр фы  тёхї фю\-ёЄр\-Єюўэю
сюы№\-°шї $n$ сєфхЄ т√\-яюыэ Є№\-ё  эх\-Ёртхэ\-ёЄтю
$$
|h(x_n+y_n)-h(x_n)|\leq \varepsilon (x_n+y_n)+\varepsilon|x_n|
=\varepsilon y_n.
$$
▌Єю эх\-Ёртхэ\-ёЄтю яЁю\-Єш\-тю\-ЁхўшЄ эх\-Ёртхэ\-ёЄтє
(\ref{order_h-h}).

\par ╧юыєўхэ\-э√х яЁю\-Єш\-тю\-Ёхўш  фю\-ърч√\-тр■Є
Ёртхэ\-ёЄ\-тю  (\ref{order_lim_varphi_1}), р чэрўшЄ ш Ёртхэ\-ёЄтю
\begin{equation}\label{order_lim_gamma}
\lim\limits_{t\to+\infty} \frac{\ln\gamma(t)}{\ln t}=0.
\end{equation}
╚ч ¤Єюую Ёртхэ\-ёЄ\-тр ш єЄ\-тхЁ\-ц\-фх\-эш   3) ыхьь√
\ref{order_lemma_gamma_properties} ёых\-фє\-хЄ Ёртхэ\-ёЄтю
\begin{equation}\label{order_lim_gamma_1}
\lim\limits_{t\to+\infty}
\frac{\ln\gamma\left(\frac{1}{t}\right)}{\ln t}=0.
\end{equation}
╥хюЁхьр фюърчрэр.

\par ╧єёЄ№ ЄхяхЁ№ $\rho(r)$ -- єЄюў\-э╕э\-э√щ яю\-Ё фюъ, яЁхф\-ёЄртш\-ь√щ
т тшфх $\rho(r)=$ $\rho+\hat{\rho}(r)$, уфх $\hat{\rho}(r)$  --
эєыхтющ єЄюў\-э╕э\-э√щ яю\-Ё фюъ, єфютыхЄ\-тюЁ ■\-∙шщ єёыю\-тш■
$\hat{\rho}\left(\frac{1}{r}\right)$ $=$
$-\rho\left(\frac{1}{r}\right)$. ┬  ¤Єюь ёыєўрх яЁш $r>0$ ш  $t>0$
т√\-яюыэ \-хЄ\-ё  эх\-Ёртхэ\-ёЄ\-тю
\begin{equation}\label{order_V_V}
V(rt)\leq t^\rho\gamma(t)V(r),
\end{equation}
уфх эхяЁх\-Ё√т\-эр  Їєэъ\-Ўш  $\gamma(t)$ єфют\-ыхЄ\-тю\-Ё хЄ
ёююЄ\-эю°х\-эш ь  (\ref{order_lim_gamma}),
(\ref{order_lim_gamma_1}). ╩ръ єцх юЄ\-ьхўр\-ыюё№ тю тёЄєя\-ых\-эшш,
єЄ\-тхЁ\-ц\-фх\-эшх ю ёяЁртхф\-ыш\-тюёЄш ¤Єюую эх\-Ёртхэ\-ёЄ\-тр
 ты \-хЄ\-ё  эх\-ъюЄю\-Ё√ь єёшых\-эш\-хь Ёх\-чєы№\-ЄрЄр ╧юЄЄхЁр
\cite{Potter} ш чэрўш\-Єхы№\-э√ь єяЁю∙х\-эшхь хую
ЇюЁ\-ьє\-ыш\-Ёют\-ъш. ╨хчєы№ЄрЄ ╧юЄЄхЁр ЎшЄш\-Ёє\-хЄё  ш
шё\-яюы№чє\-хЄё  т  \cite{Bingham}. ═хЁртхэ\-ёЄтю  (\ref{order_V_V})
фю\-тюы№эю ўрёЄю шё\-яюы№чє\-хЄ\-ё  т яЁхф\-ёЄрт\-ыхэ\-эющ Ёрсю\-Єх.
╧юыхч\-эюёЄ№ ¤Єюую эх\-Ёртхэ\-ёЄ\-тр ьюц\-эю Єръ\-цх єтшфхЄ№,
рэрышчш\-Ёє  фю\-ърчр\-Єхы№\-ёЄтю ёых\-фє\-■\-∙хщ ыхьь√.

%%%%%%%%%%%%%%%%%%%%%%%%%%%%%%%%%%%%%%%%%%%%%%%%%%%%%%%%%%%%%%%%%%%%%%%%%%%%%%%%%%%%%%%
\begin{lemma}\label{order_lemma_V_1}\hskip-2mm{.}\:
%%%%%%%%%%%%%%%%%%%%%%%%%%%%%%%%%%%%%%%%%%%%%%%%%%%%%%%%%%%%%%%%%%%%%%%%%%%%%%%%%%%%%%%
╧єёЄ№ $\rho(r)$ -- яЁю\-шч\-тюы№\-э√щ эєыхтющ єЄюў\-э╕э\-э√щ
яю\-Ё фюъ, єфют\-ыхЄ\-тю\-Ё ■\-∙шщ єёыю\-тш■
$\rho\left(\frac{1}{r}\right)$ $=$ $-\rho(r)$ ш яєёЄ№
$$
V_1(r)=\frac{2r}{\pi}\int\limits_0^\infty \frac{V(t)}{t^2+r^2}dt.
$$
╥юуфр ёяЁртхф\-ышт√ ёых\-фє\-■\-∙шх єЄ\-тхЁ\-ц\-фх\-эш :\\
1) Їєэъ\-Ўш  $V_1(r)$ фю\-яєёър\-хЄ уюыю\-ьюЁЇ\-эюх
яЁю\-фюы\-цх\-эшх
ё яюыє\-юёш $(0,\infty)$ т яюыє\-яыюё\-ъюёЄ№ $\Re z>0$,\\
2) $V_1(r)=r^{\rho_1(r)}$, уфх $\rho_1(r)$  -- эєыхтющ
єЄюў\-э╕э\-э√щ яю\-Ё фюъ, єфют\-ыхЄ\-тю\-Ё ■\-∙шщ єёыю\-тш■
$\rho_1\left(\frac{1}{r}\right)$
$=$ $-\rho_1(r)$,\\
3) т√яюыэ \-хЄё  Ёртхэ\-ёЄтю
$$
\lim\limits_{r\to\infty} \frac{V_1(r)}{V(r)}=1.
$$
\end{lemma}
%%%%%%%%%%%%%%%%%%%%%%%%%%%%%%%%%%%%%%%%%%%%%%%%%%%%%%%%%%%%%%%%%%%%%%%%%%%%%%%%%%%%%%%
\par {\sc  ─юърчрЄхы№ёЄтю.}\: ╥ю, ўЄю Їєэъ\-Ўш 
$$
V_1(z)=\frac{2z}{\pi}\int\limits_0^\infty \frac{V(t)}{t^2+z^2}dt
$$
 ты \-хЄё  уюыю\-ьюЁЇ\-эющ т яюыє\-яыюё\-ъюёЄш $\Re z>0$, ыхуъю
фю\-ърчрЄ№, шё\-їюф  шч эх\-Ёртхэ\-ёЄ\-тр $V(t)\leq
M(t^{\frac{1}{2}}+t^{-\frac{1}{2}})$ (ёьюЄЁш эхЁртхэ\-ёЄтю
(\ref{order-V_leq})). ╥ръшь ёяюёюсюь юсюёэю\-т√\-тр\-хЄ\-ё 
єЄ\-тхЁ\-ц\-фх\-эшх  1) ыхьь√.

\par ╨рч\-сштр  яюыє\-юё№ $[0,\infty)$ эр ўрёЄш $[0,\varepsilon r]$,
$[\varepsilon r,Nr]$, $[Nr,\infty)$, яЁхф\-ёЄртшь $V_1(r)$ т тшфх
ёєь\-ь√ ЄЁ╕ї шэ\-ЄхуЁр\-ыют
$$
V_1(r)=I_1(r)+I_2(r)+I_3(r).
$$
╚ьххь
$$
\frac{I_2(r)}{V(r)}=\frac{2}{\pi}\int\limits_\varepsilon^N
\frac{V(ur)}{V(r)}\frac{du}{u^2+1}.
$$
╚ч Єхю\-Ёх\-ь√  \ref{order_th_V:V} ёых\-фєхЄ, ўЄю
\begin{equation}\label{order_I_2}
\lim\limits_{r\to\infty} \frac{I_2(r)}{V(r)}=
\frac{2}{\pi}\int\limits_\varepsilon^N \frac{du}{u^2+1}.
\end{equation}
\par ─рыхх эрїюфшь, ўЄю
$$
I_1(r)=\frac{2}{\pi}\int\limits_0^\varepsilon
\frac{V(ur)}{u^2+1}du\leq \frac{2V(r)}{\pi}
\int\limits_0^\varepsilon \frac{\gamma(u)}{u^2+1}du.
$$
╬Єё■фр ёых\-фєхЄ эхЁртхэёЄтю

\begin{equation}\label{order_I_1}
\mathop{\overline{\lim}} \limits_{r\to\infty}
\frac{I_1(r)}{V(r)}\leq \frac{2}{\pi}\int\limits_0^\varepsilon
\frac{\gamma(u)}{u^2+1}du.
\end{equation}
└эрыюушўэю яюыєўрхь, ўЄю
\begin{equation}\label{order_I_3}
\mathop{\overline{\lim}} \limits_{r\to\infty}
\frac{I_3(r)}{V(r)}\leq \frac{2}{\pi}\int\limits_N^\infty
\frac{\gamma(u)}{u^2+1}du.
\end{equation}

\par ╚ч ёююЄ\-эю°х\-эшщ  (\ref{order_I_2})-(\ref{order_I_3}) ёых\-фєхЄ
$$
\mathop{\overline{\lim}} \limits_{r\to\infty}
\left|\frac{V_1(r)}{V(r)}-1\right| \leq \frac{2}{\pi}
\int\limits_0^\varepsilon \frac{\gamma(u)}{u^2+1}du+1- \frac{2}{\pi}
\int\limits_\varepsilon^N \frac{du}{u^2+1} +\frac{2}{\pi}
\int\limits_N^\infty \frac{\gamma(u)}{u^2+1}du.
$$
╚ч ¤Єюую эх\-Ёртхэ\-ёЄ\-тр яЁх\-фхы№\-э√ь яхЁх\-їюфюь яЁш
$\varepsilon\to 0$, $N\to\infty$ яю\-ыєўрхь єЄ\-тхЁ\-ц\-фх\-эшх  3)
ыхьь√.

\par ┼ёыш Їєэъ\-Ўш■ $\rho_1(r)$ юяЁх\-фхышЄ№ Ёртхэ\-ёЄ\-тюь
$V_1(r)=r^{\rho_1(r)}$, Єю шч єЄ\-тхЁ\-ц\-фх\-эш   3) ёых\-фє\-хЄ
Ёртхэ\-ёЄтю
\begin{equation}\label{order_lim_rho_1}
\lim\limits_{r\to\infty} \rho_1(r)=0.
\end{equation}
─рыхх шьххь
$$
rV'_1(r)=\frac{2r}{\pi} \int\limits_0^\infty
\frac{t^2-r^2}{(t^2+r^2)^2}V(t)dt.
$$
╧ютЄюЁ   Ёрё\-ёєцфх\-эш , юсюёэю\-т√\-тр■\-∙шх єЄ\-тхЁ\-ц\-фх\-эшх
3), яю\-ыєўшь, ўЄю
$$
\lim\limits_{r\to\infty} \frac{rV_1'(r)}{V(r)}= -\frac{2}{\pi}
\int\limits_0^\infty \frac{d}{dt}\frac{t}{t^2+1}dt=0,
$$
\begin{equation}\label{order_lim_frac}
\lim\limits_{r\to\infty} \frac{rV_1'(r)}{V_1(r)}=0.
\end{equation}
╚ч ЁртхэёЄт  (\ref{order_lim_rho_1}), (\ref{order_lim_frac}),
ёых\-фєхЄ, ўЄю Їєэъ\-Ўш  $\rho_1(r)$  ты \-хЄ\-ё  эєыхт√ь
єЄюў\-э╕э\-э√ь яю\-Ё ф\-ъюь.

─рыхх эрїюфшь
$$
V_1\left(\frac{1}{r}\right)= \frac{2}{\pi r} \int\limits_0^\infty
\frac{V(t)}{t^2+\frac{1}{r^2}}dt= \frac{2}{\pi r}
\int\limits_0^\infty \frac{V\left(\frac{1}{u}\right)}
{\frac{1}{u^2}+\frac{1}{r^2}}\frac{du}{u^2}= \frac{2r}{\pi}
\int\limits_0^\infty \frac{V(u)}{u^2+r^2} du=V_1(r).
$$
╚ч ¤Єюую ёых\-фєхЄ, ўЄю $\rho_1\left(\frac{1}{r}\right)$ $=$
$-\rho_1(r)$. ╙ЄтхЁц\-фх\-эшх  2), р тьхёЄх ё эшь ш ыхььр,
фю\-ърчрэ√.

\par ╚ч ыхьь√  \ref{order_lemma_V_1}, т ўрёЄ\-эюё\-Єш, ёых\-фє\-хЄ, ўЄю яЁш
Ёх°х\-эшш чрфрўш  $A$, ёЇюЁ\-ьє\-ыш\-Ёю\-трэ\-эющ т эрўрых
Ёрч\-фхыр, т ърўхёЄ\-тх ъырёёр  $\frak{A}$ ьюц\-эю тч Є№ эх
ьэю\-цх\-ёЄ\-тю Їєэъ\-Ўшщ тшфр $V(r)=r^{\rho(r)}$, уфх $\rho(r)$  --
яЁю\-шч\-тюы№\-э√щ єЄюў\-э╕э\-э√щ яю\-Ё фюъ, р чэрўш\-Єхы№\-эю сюыхх
єчъшщ ъырёё, ёюёЄю \-∙шщ шч Єхї Їєэъ\-Ўшщ $V(r)$, фы  ъюЄю\-Ё√ї
Їєэъ\-Ўш  $\rho(r)$  ты \-хЄё  рэрыш\-Єш\-ўхё\-ъющ Їєэъ\-Ўшхщ эр
яюыє\-юёш $(0,\infty)$, ш яЁхф\-ёЄрт\-ы \-хЄ\-ё  т тшфх $\rho(r)$
$=$ $\rho+\hat{\rho}(r)$, уфх $\hat{\rho}(r)$  -- эєых\-тющ
єЄюў\-э╕э\-э√щ яю\-Ё фюъ, єфют\-ыхЄ\-тю\-Ё ■\-∙шщ єёыю\-тш■
$\hat{\rho}\left(\frac{1}{r}\right)$ $=$ $-\hat{\rho}(r)$.

\par ┬ Єхю\-Ёхьх  \ref{order_th_ln_gamma} юЄ\-эюёш\-Єхы№\-эю Їєэъ\-Ўшш
$\gamma(t)$ єЄ\-тхЁ\-ц\-фр\-хЄ\-ё , ўЄю фы  эх╕ ёяЁр\-тхф\-ышт√
Ёртхэ\-ёЄ\-тр (\ref{order_lim_gamma}), (\ref{order_lim_gamma_1}).
─ы  фю\-ёЄрЄюў\-эю °ш\-Ёю\-ъюую ъырёёр Їєэъ\-Ўшщ тшфр
$V(r)=r^{\rho(r)}$ ёяЁртхф\-ышт√ сюыхх Єюўэ√х юЎхэъш. ┬ ыхььх
\ref{order_lemma_gamma_properties} єЄ\-тхЁ\-ц\-фр\-хЄё , ўЄю
$\gamma(t)\geq V(t)$. ┬ ёых\-фє\-■\-∙хщ Єхю\-Ёхьх яЁш\-тюфшЄ\-ё 
фю\-ёЄрЄюў\-эю °ш\-Ёюъшщ ъырёё єЄюў\-э╕э\-э√ї яю\-Ё ф\-ъют Єръшї,
ўЄю фы  Їєэъ\-Ўшщ $\gamma(t)$ юяЁх\-фхы╕э\-э√ї ё яюью∙№■ ¤Єшї
єЄюў\-э╕э\-э√ї яю\-Ё ф\-ъют яЁш $t\geq 1$ т√\-яюыэ \-хЄё 
эх\-Ёртхэ\-ёЄтю $\gamma(t)\leq MV(t)$ ё эх\-ъю\-Єю\-Ёющ
яю\-ёЄю э\-эющ $M$.

%%%%%%%%%%%%%%%%%%%%%%%%%%%%%%%%%%%%%%%%%%%%%%%%%%%%%%%%%%%%%%%%%%%%%%%%%%%%%%%%%%%%%%%
\begin{theorem}\label{order_th_gamma_V}\hskip-2mm{.}\:
%%%%%%%%%%%%%%%%%%%%%%%%%%%%%%%%%%%%%%%%%%%%%%%%%%%%%%%%%%%%%%%%%%%%%%%%%%%%%%%%%%%%%%%
╧єёЄ№ $\rho(r)$ -- эєыхтющ єЄюў\-э╕э\-э√щ яю\-Ё фюъ, фтрцф√
фшЇ\-Їх\-Ёхэ\-Ўш\-Ёєх\-ь√щ эр ьэю\-цх\-ёЄ\-тх
$(0,\infty)\setminus\{1\}$ ш єфют\-ыхЄ\-тю\-Ё ■\-∙шщ єёыю\-тш■
$\rho\left(\frac{1}{r}\right)$ $=$ $-\rho(r)$. ╧єёЄ№
$V(r)=r^{\rho(r)}$, $V(r)\to\infty$ $(r\to\infty)$, ш Їєэъ\-Ўш 
$h(x)=\ln{V(e^x)}$  ты \-хЄ\-ё  тюуэє\-Єющ т эх\-ъю\-Єю\-Ёющ
юъ\-ЁхёЄ\-эюё\-Єш схё\-ъю\-эхў\-эюё\-Єш. ╧єёЄ№ Їєэъ\-Ўш  $\gamma(t)$
юяЁх\-фхы \-хЄ\-ё  Ёртхэ\-ёЄ\-тюь  (\ref{order_ln_gamma_1}). ╥юуфр
ёє\-∙хёЄ\-тєхЄ яю\-ёЄю э\-эр  $M$ Єр\-ър , ўЄю яЁш $t\geq 1$
т√\-яюыэ \-хЄ\-ё  эх\-Ёртхэ\-ёЄ\-тю $\gamma(t)\leq MV(t)$. ┼ёыш
Їєэъ\-Ўш  $h(x)$ тюуэєЄр эр яюыє\-юёш $(0,\infty)$, Єю
$\gamma(t)=V(t)$ яЁш $t\geq 1$.
\end{theorem}
%%%%%%%%%%%%%%%%%%%%%%%%%%%%%%%%%%%%%%%%%%%%%%%%%%%%%%%%%%%%%%%%%%%%%%%%%%%%%%%%%%%%%%%
\par {\sc  ─юърчрЄхы№ёЄтю.}\: ╨рё\-ёьюЄ\-Ёшь
Їєэъ\-Ўш■ $a(x)=h(x)-xh'(x)$. ╚ьххь $a'(x)$ $=$  $-xh''(x)$. ┬ ёшыє
єёыю\-тшщ Єхю\-Ёх\-ь√ Їєэъ\-Ўш  $a(x)$  ты \-хЄё  тюч\-ЁрёЄр■\-∙хщ т
эх\-ъю\-Єю\-Ёющ юъ\-ЁхёЄ\-эюё\-Єш схё\-ъю\-эхў\-эюё\-Єш. ╧Ёш $x\geq
1$ т√\-яюыэ \-хЄ\-ё  Ёртхэ\-ёЄ\-тю
\begin{equation}\label{order_h_c}
h(x)=-x\left(\int\limits_1^x \frac{a(t)}{t^2}dt+c\right),\qquad
c=-h(1).
\end{equation}
─юърцхь, ўЄю шэ\-ЄхуЁры
\begin{equation}\label{order_int_a}
\int\limits_1^\infty \frac{a(t)}{t^2}dt
\end{equation}
 ты \-хЄё  ёїюф \-∙шь\-ё . ┬ яЁюЄшт\-эюь ёыєўрх, хёыш  т
юъ\-ЁхёЄ\-эюё\-Єш схё\-ъю\-эхў\-эюё\-Єш т√\-яюыэ \-хЄё 
эх\-Ёртхэ\-ёЄтю $a(t)\leq 0$, Єю яюыєўшь яЁю\-Єштю\-Ёхўшх ё
Ёртхэ\-ёЄ\-тюь (\ref{order_lim_h(x)}), р т ёыєўрх, хёыш $a(t)\geq 0$
т эх\-ъю\-Єю\-Ёющ юъ\-ЁхёЄ\-эюё\-Єш схё\-ъю\-эхў\-эюё\-Єш, Єю
яю\-ыєўшь яЁю\-Єштю\-Ёхўшх ё ёююЄ\-эю°х\-эшхь $V(r)\to\infty$
$(r\to\infty)$ шч єёыю\-тш  Єхю\-Ёх\-ь√. ╥хь ёрь√ь ёїюфш\-ьюёЄ№
шэ\-ЄхуЁр\-ыр (\ref{order_int_a}) фю\-ърчрэр. ╥хяхЁ№ Ёртхэ\-ёЄ\-тю
(\ref{order_h_c}) ьюц\-эю яхЁх\-яшёрЄ№ т тшфх
\begin{equation}\label{order_h_a}
h(x)=x\int\limits_x^\infty \frac{a(t)}{t^2}dt.
\end{equation}
┬ ёшыє Ёртхэ\-ёЄ\-тр  (\ref{order_lim_h(x)}) фю\-яюыэш\-Єхы№\-эюх
ёырурх\-ьюх $c_1x$ т яЁртющ ўрёЄш ¤Єюую Ёртхэ\-ёЄ\-тр, ъюЄю\-Ёюх
яю\- ты \-хЄ\-ё  яЁш яхЁх\-їюфх юЄ  (\ref{order_h_c}) ъ
(\ref{order_h_a}) юЄ\-ёєЄ\-ёЄ\-тєхЄ.

\par ╧юърцхь, ўЄю
\begin{equation}\label{order_lim_a(t)}
\lim\limits_{t\to\infty} a(t)=+\infty.
\end{equation}

\par ┬ яЁюЄшт\-эюь ёыєўрх Їєэъ\-Ўш  $h(x)$, р тьхёЄх ё эхщ ш
Їєэъ\-Ўш   $V(r)$, с√ыш с√ юуЁрэш\-ўхэ√ т эх\-ъю\-Єю\-Ёющ
юъ\-ЁхёЄ\-эюё\-Єш $+\infty$, р ¤Єю яЁю\-Єштю\-ЁхўшЄ єёыю\-тш ь
Єхю\-Ёх\-ь√. ╥хь ёрь√ь Ёртхэ\-ёЄтю  (\ref{order_lim_a(t)})
фю\-ърчр\-эю.

\par ╧єёЄ№ ўшёыю $x_0>0$ Єръютю, ўЄю Їєэъ\-Ўш  $h(x)$  ты \-хЄё  тюуэє\-Єющ
яЁш $x\geq x_0$. ╙Ёртэх\-эшх ърёр\-Єхы№\-эющ ъ ъЁштющ $y=h(x)$ т
Єюўъх $x_0$ шьххЄ тшф
$$
Y(x)=h(x_0)-x_0h'(x_0)+h'(x_0)x.
$$
┼ёыш $x_0$ фю\-ёЄрЄюў\-эю тхышъю, Єю Єюуфр сєфхЄ т√\-яюы\-э Є№\-ё 
эх\-Ёртхэ\-ёЄтю $Y(0)>0$. ┬  ¤Єюь ёыєўрх ёє\-∙хёЄ\-тєхЄ Їєэъ\-Ўш 
$h_1(x)$, юсырфр\-■\-∙р  ётющёЄ\-трьш:\\
1) $h_1(x)$ -- ў╕Єэр  Їєэъ\-Ўш  эх\-яЁх\-Ё√т\-эр  эр юёш
$(-\infty,\infty)$ тюуэє\-Єр  ш фшЇ\-ЇхЁхэ\-Ўш\-Ёєх\-ьр 
эр яюыє\-юёш $(0,\infty)$,\\
2) $h_1(x)\geq 0$ яЁш $x\geq 0$, $h_1(0)=0$,\\
3) $h_1(x)=h(x)$ яЁш $x\geq x_0$.\\
╚ч ёЇюЁ\-ьєыш\-Ёю\-трэ\-э√ї ётющёЄт Їєэъ\-Ўшш $h_1(x)$ ёых\-фє\-■Є
х∙╕ Єръшх
ётющ\-ёЄ\-тр:\\
4) ёє\-∙хёЄ\-тєхЄ яю\-ёЄю э\-эр  $M_1$ Єр\-ър , ўЄю эр тёхщ
тх∙хёЄ\-тхэ\-эющ юёш т√\-яюыэ \-хЄ\-ё  эх\-Ёртхэ\-ёЄтю
$|h(x)-h_1(x)|\leq M_1$,\\
5) Їєэъ\-Ўш  $h_1(x)$  ты \-хЄ\-ё  тюч\-ЁрёЄр■\-∙хщ Їєэъ\-Ўшхщ эр
яюыє\-юёш $[0,\infty)$,\\
6) яЁш $x\geq 0$, $y\geq 0$, т√\-яюыэ \-хЄё  эх\-Ёртхэ\-ёЄтю
$h_1(x+y)\leq h_1(x)+h_1(y)$.

\par ─рыхх эрїюфшь яЁш $y>0$
$$
\ln\gamma(e^y)=\sup\limits_{x\in(-\infty,\infty)} (h(x+y)-h(x))\leq
2M_1+\sup\limits_{x\in(-\infty,\infty)} (h_1(x+y)-h_1(x))
$$
$$
=2M_1+ \sup\limits_{x\in(-\infty,\infty)} (h_1(|x+y|)-h_1(|x|))\leq
2M_1+\sup\limits_{x\in(-\infty,\infty)} (h_1(|x|+y)-h_1(|x|))
$$
$$
\leq 2M_1+h_1(y)\leq 3M_1+h(y)=3M_1+\ln V(e^y).
$$
╚ч ¤Єюую ёых\-фєхЄ, ўЄю яЁш $t>1$ т√\-яюыэ \-хЄё  эх\-Ёртхэ\-ёЄтю
$\gamma(t)\leq e^{3M_1}V(t)$. ┬ ёыєўрх, хёыш Їєэъ\-Ўш  $h(x)$
тюуэєЄр эр яюыє\-юёш $[0,\infty)$, Єю Єюуфр $h_1(x)=h(x)$, $M_1=0$,
$\gamma(t)\leq V(t)$. ┬ьхёЄх ё эх\-Ёртхэ\-ёЄ\-тюь $\gamma(t)\geq
V(t)$ ¤Єю фр╕Є $\gamma(t)=V(t)$. ╥хю\-Ёхьр фю\-ърчрэр.

\par ┬ю ьэюушї тюяЁюёрї юёэют\-эє■ Ёюы№ шуЁрхЄ эх єЄюў\-э╕э\-э√щ яю\-Ё фюъ
$\rho(r)$, р Їєэъ\-Ўш  $V(r)=r^{\rho(r)}$. ╚эюуфр трцэ√ ётющ\-ёЄтр
Їєэъ\-Ўшш $V(r)$ эх Єюы№\-ъю т юъ\-ЁхёЄ\-эюё\-Єш
схё\-ъю\-эхў\-эюё\-Єш, эю ш т юъ\-ЁхёЄ\-эюё\-Єш эєы , эр\-яЁш\-ьхЁ,
яЁш шё\-ёыхфю\-тр\-эшш ётющёЄт шэ\-ЄхуЁр\-ыр $\int\limits_0^\infty
K(t,r)$ $V(t)dt$. ╧ю¤Єю\-ьє, ъръ єцх с√ыю ёърчр\-эю, т
фры№\-эхщ\-°хь ь√ яЁхф\-яю\-ырур\-хь, ўЄю эєых\-тющ єЄюў\-э╕э\-э√щ
яю\-Ё фюъ єфют\-ыхЄ\-тю\-Ё хЄ єёыю\-тш■
$\rho\left(\frac{1}{r}\right)$ $=-\rho(r)$. ▌Єю ¤ътш\-тр\-ыхэЄ\-эю
ёююЄ\-эю°х\-эш■ $V\left(\frac{1}{r}\right)$ $=V(r)$. ╚ч Єхю\-Ёх\-ь√
\ref{order_th_V_V_1} ёых\-фє\-хЄ, ўЄю ЄЁхсю\-тр\-эшх
фшЇ\-ЇхЁхэ\-Ўш\-Ёєх\-ьюёЄш Їєэъ\-Ўшш $\rho(r)$ эр ьэю\-цх\-ёЄ\-тх
$(0,\infty)\setminus\{1\}$ тю ьэюушї ёыєўр\- ї эх  ты \-хЄ\-ё 
ёє\-∙хёЄ\-тхэ\-э√ь юуЁрэш\-ўх\-эш\-хь. ╨ртхэ\-ёЄ\-тю
$$
\rho(r)=\frac{\ln V(r)}{\ln r}
$$
яю\-ърч√\-трхЄ, ўЄю ЄЁхсю\-тр\-эшх фшЇ\-Їх\-Ёхэ\-Ўш\-Ёєх\-ьюё\-Єш
Їєэъ\-Ўшш $\rho(r)$ т Єюўъх  $1$ фю\-ёЄрЄюў\-эю ёЄхёэш\-Єхы№\-эю.
▌Єю ЄЁхсю\-тр\-эшх шё\-ъы■ўр\-хЄ шч Ёрё\-ёьюЄ\-Ёх\-эш  Їєэъ\-Ўшш
$\rho(r)$ $=$ $\frac{A|\ln r|^\alpha}{\ln r}$, $\alpha\in(0,1)$, фы 
ъюЄю\-Ё√ї Їєэъ\-Ўш  $V(r)$ шьххЄ тшф $V(r)$ $=$ $\exp(A|\ln
r|^\alpha)$. ┬тшфє яЁюё\-Єю\-Є√ ¤Єшї Їєэъ\-Ўшщ шї ўрёЄю
шё\-яюы№\-чє\-■Є т ърўхёЄ\-тх Їєэъ\-Ўшщ ёЁрт\-эх\-эш . ─ы  Єръшї
Їєэъ\-Ўшщ Єръ\-цх ёяЁр\-тхф\-ыштр ЇюЁ\-ьєыр $\gamma(t)$ $=$ $V(t)$
яЁш $t\geq 1$.

\par ╩юэхўэю, шч ыхьь√  \ref{order_lemma_V_1}  ёых\-фє\-хЄ, ўЄю ЄЁхсю\-тр\-эшх
ёє\-∙хёЄ\-тю\-тр\-эш  яЁю\-шч\-тюф\-эющ $\rho'(1)$ Єръ\-цх эх
 ты \-хЄё  ёє\-∙хёЄ\-тхэ\-э√ь. ─ы  ¤Єюую эєцэю чр\-ьх\-эшЄ№
Їєэъ\-Ўш■ $V(r)$ эр Їєэъ\-Ўш■ $V_1(r)$. ╬фэръю шё\-яюы№\-чю\-тр\-эшх
Їєэъ\-Ўшш $V_1(r)$ эх\-ёъюы№\-ъю чр\-ЄЁєф\-э \-хЄё 
фю\-ёЄр\-Єюў\-эющ ёыюц\-эюё\-Є№■ ¤Єющ Їєэъ\-Ўшш.

\par ┬ фры№\-эхщ\-°хь ЄхъёЄх ёЄрЄ№ш яюф {\it эєых\-т√ь єЄюў\-э╕э\-э√ь
яю\-Ё ф\-ъюь} яю\-эшьр\-хЄ\-ё  Їєэъ\-Ўш  $\rho(r)$, ъюЄю\-Ёр 
єфютыхЄ\-тю\-Ё хЄ єёыю\-тш\- ь:\\
1) $\lim\limits_{r\to\infty} \rho(r)=0$,\\
2) $\rho\left(\frac{1}{r}\right)=-\rho(r)$,\\
3) Їєэъ\-Ўш  $\rho(r)$  ты \-хЄё  эх\-яЁх\-Ё√т\-эю
фшЇ\-ЇхЁхэ\-Ўш\-Ёєх\-ьющ эр ьэю\-цх\-ёЄтх $(0,\infty)\setminus\{1\}$,\\
4) $\lim\limits_{r\to\infty} r\ln r \rho'(r)=0$,\\
5) Їєэъ\-Ўш  $V(r)=r^{\rho(r)}$ эх\-яЁх\-Ё√т\-эю яЁю\-фюыцр\-хЄ\-ё 
т Єюў\-ъє  {$1$}, яЁш\-ў╕ь $V(1)$ $=1$ (ёрьр Їєэъ\-Ўш  $\rho(r)$ т
Єюўъх  $1$ ью\-цхЄ ш эх с√Є№ юяЁх\-фхы╕э\-эющ).

\par ─Ёєушх єЄюў\-э╕э\-э√х яю\-Ё фъш $\rho(r)$ шьх■Є тшф
$\rho(r)=$ $\rho+$ $\hat{\rho}(r)$, уфх $\hat{\rho}(r)$  -- эєыхтющ
єЄюў\-э╕э\-э√щ яю\-Ё фюъ, єфютыхЄ\-тю\-Ё ■\-∙шщ эр\-яшёрэ\-э√ь т√°х
єёыю\-тш ь.

%%%%%%%%%%%%%%%%%%%%%%%%%%%%%%%%%%%%%%%%%%%%%%%%%%%%%%%%%%%%%%%%%%%%%%%%%%%%%%%%%%%%%%
\section{╠хЁ√, яЁх\-фхы№э√х ьэю\-цх\-ёЄтр ьхЁ}
%%%%%%%%%%%%%%%%%%%%%%%%%%%%%%%%%%%%%%%%%%%%%%%%%%%%%%%%%%%%%%%%%%%%%%%%%%%%%%%%%%%%%%
\qquad ┬  ¤Єюь Ёрч\-фхых сєфхЄ шч\-ыюцх\-эр ЄхюЁш  яЁх\-фхы№\-э√ї
ьэю\-цх\-ёЄт Ёр\-фю\-эю\-т√ї ьхЁ ъю\-эхў\-эю\-ую яю\-Ё ф\-ър эр
яюыє\-юёш $(0,\infty)$. ─ы  єфюс\-ёЄтр ўшЄр\-Єх\-ыхщ сєфєЄ
яЁш\-тхфх\-э√ юяЁх\-фхых\-эш  шё\-яюы№чє\-хь√ї ЄхЁьш\-эют ш
ЇюЁ\-ьє\-ыш\-Ёютъш шё\-яюы№чє\-хь√ї Ёх\-чєы№\-Єр\-Єют.
─юърчр\-Єхы№\-ёЄтр яЁш\-тхфхэ\-э√ї схч фю\-ърчр\-Єхы№\-ёЄтр
єЄ\-тхЁ\-ц\-фх\-эшщ ьюц\-эю эрщЄш т  \cite{Burbaki}, \cite{Danford},
\cite{Landkof}.

\par ╠√ эрўэ╕ь ёю ётхфхэшщ юЄэюё \-∙шїё  ъ "рсёЄЁръЄ\-эющ" \, ЄхюЁшш ьхЁ√.

\par {\it ╚чьхЁш\-ьюх яЁюёЄ\-Ёрэ\-ёЄ\-тю ё тх∙хёЄ\-тхэ\-эющ ьхЁющ }$\mu$  --
¤Єю ЄЁющър $(X, \cal{A},\mu)$, уфх $X$  -- ьэю\-цх\-ёЄтю, $\cal{A}$
-- ёшуьр-рыухсЁр яюф\-ьэю\-цх\-ёЄт ьэю\-цх\-ёЄ\-тр  $X$, $\mu$  --
Їєэъ\-Ўш , юяЁх\-фхы╕э\-эр  эр ьэю\-цх\-ёЄ\-трї  $E$, тїюф \-∙шї т
ёшуьр-рыухсЁє  $\cal{A}$, яЁш\-эшьр■\-∙р  чэрўх\-эш  шч
Ёрё\-°шЁхэ\-эющ ўшёыю\-тющ яЁ ьющ $[-\infty,\infty]$ ш
юсыр\-фр■\-∙р  ётющёЄ\-тюь ёў╕Є\-эющ рффш\-Єшт\-эюё\-Єш:
$\mu\left(\bigcup\limits_{k=1}^\infty E_k\right)=
\sum\limits_{k=1}^\infty \mu(E_k)$, уфх $E_k\in
 \cal{A}$ ш ьэю\-цх\-ёЄ\-тр $E_k$ юсЁрчє\-■Є фшч·■э\-ъЄ\-эє■
яю\-ёыхфю\-тр\-Єхы№\-эюёЄ№ ьэю\-цх\-ёЄт ($E_k\cap E_j=\emptyset$,
хёыш $k\neq j$).

\par ╠хЁр $\mu$ эр\-ч√тр\-хЄ\-ё  {\it яюыю\-цш\-Єхы№\-эющ}, хёыш $\mu(E)\geq 0$
фы  ы■сюую $E\in \cal{A}$.

\par ╧юыюцш\-Єхы№\-эр  ьхЁр эр\-ч√тр\-хЄ\-ё  {\it ъю\-эхў\-эющ}, хёыш
$\mu(X)<\infty$.

\par ╧єёЄ№ $A\in\cal{A}$. {\it ╬уЁрэш\-ўх\-эшхь ьхЁ√} $\mu$ эр ьэю\-цх\-ёЄ\-тю $A$
(юсю\-чэрўх\-эшх $\mu_A$) эр\-ч√тр\-хЄ\-ё  ьхЁр, юяЁх\-фхы х\-ьр 
Ёртхэ\-ёЄ\-тюь $\mu_A(E)=\mu(A\cap E)$.

\par ├ютюЁ Є, ўЄю ьхЁр $\mu$ {\it ёю\-ёЁхфю\-Єюўхэр эр ьэю\-цх\-ёЄ\-тх  $A$},
хёыш $\mu_A=\mu$.

\par ╧єёЄ№ $\mu_1$ ш $\mu_2$ фтх ьхЁ√, юяЁх\-фхы╕э\-э√х эр юфэющ
ёшуьр-рыухсЁх  $\cal{A}$. ╥ръшх ьхЁ√ эр\-ч√тр\-■Є\-ё  {\it тчршь\-эю
ёшэує\-ы Ё\-э√ьш}, хёыш юэш ёю\-cЁхфю\-Єюўх\-э√ эр
эх\-яхЁх\-ёхър■\-∙шї\-ё  ьэю\-цх\-ёЄ\-трї $A_1$ ш  $A_2$.
╥хю\-Ёх\-ьр ╒рэр єЄ\-тхЁ\-ц\-фрхЄ, ўЄю хёыш $\mu$  --
тх∙хёЄ\-тхэ\-эр  ьхЁр, Єю ёє\-∙хёЄ\-тє■Є ьэю\-цх\-ёЄ\-тр $A_1$ ш
$A_2$ шч рыухсЁ√ $\cal{A}$ Єръшх, ўЄю $A_1\cap A_2=\emptyset$,
$X=A_1\cup A_2$, юуЁрэш\-ўх\-эшх $\mu_{A_1}$  -- яюыю\-цш\-Єхы№\-эр 
ьхЁр, р юуЁрэш\-ўх\-эшх $\mu_{A_2}$  -- юЄЁшЎр\-Єхы№\-эр  ьхЁр.

\par ╤ююЄтхЄ\-ёЄ\-тє■∙р  ярЁр ьэю\-цх\-ёЄт $A_1$, $A_2$ эр\-ч√тр\-хЄ\-ё 
{\it Ёрч\-ыюцх\-эшхь ╒рэр} фы  ьхЁ√ $\mu$. ╠хЁр $\mu_+=\mu_{A_1}$
эр\-ч√тр\-хЄ\-ё  {\it яюыю\-цш\-Єхы№\-эющ ёю\-ёЄрт\-ы ■\-∙хщ} ьхЁ√
$\mu$, ьхЁр $\mu_-=-\mu_{A_2}$ эр\-ч√тр\-хЄ\-ё  {\it
юЄЁшЎр\-Єхы№\-эющ ёю\-ёЄрт\-ы ■\-∙хщ} ьхЁ√ $\mu$. ╥ръшь юсЁр\-чюь,
ы■ср  тх∙хёЄ\-тхэ\-эр  ьхЁр  $\mu$ хёЄ№ Ёрч\-эюёЄ№ $\mu=\mu_+-\mu_-$
фтєї тчршьэю ёшэує\-ы Ё\-э√ї яюыю\-цш\-Єхы№\-э√ї ьхЁ.

\par ╧ЁхфёЄртых\-эшх ьхЁ√ $\mu=\mu_+-\mu_-$ т тшфх Ёрч\-эюё\-Єш фтєї
тчршь\-эю ёшэує\-ы Ё\-э√ї яюыю\-цш\-Єхы№\-э√ї ьхЁ эр\-ч√тр\-хЄ\-ё 
{\it Ёрч\-ыюцх\-эшхь ╞юЁфрэр} ьхЁ√ $\mu$.

\par ╒юЄ  Ёрч\-ыюцх\-эшх ╒рэр $X=A_1\cup A_2$ эх\-юфэю\-чэрў\-эю,
Ёрч\-ыюцх\-эшх ╞юЁфрэр $\mu=\mu_+-\mu_-$ юфэю\-чэрў\-эю. ╠хЁ√
$\mu_+$ ш  $\mu_-$ юфэю\-чэрўэю юяЁх\-фхы \-хЄё   ьхЁющ  $\mu$.

\par ┼ёыш $\mu=\mu_+-\mu_-$ хёЄ№ Ёрч\-ыюцх\-эшх ╞юЁфрэр ьхЁ√ $\mu$,
Єю їюЄ  с√ юфэр шч ьхЁ $\mu_+$, $\mu_-$  ты \-хЄё  ъю\-эхў\-эющ,
шэрўх Ёртхэ\-ёЄ\-тю $\mu(X)=\mu_+(X)-\mu_-(X)$ эх шьхыю с√ ёь√ёыр.

\par ╠хЁр $|\mu|=\mu_++\mu_-$ эр\-ч√тр\-хЄё  {\it ьюфєыхь} шыш {\it яюыэющ
трЁшр\-Ўшхщ} ьхЁ√ $\mu$.

\par ┬рцэ√ь  ты \-хЄё  ёыєўрщ, ъюуфр $X$  -- Єюяю\-ыюуш\-ўхё\-ъюх
яЁюёЄ\-Ёрэ\-ёЄ\-тю, р ёююЄ\-тхЄ\-ёЄтє■\-∙р  ёшуьр-рыухсЁр  ты \-хЄё 
ёшуьр-рыухсЁющ $\frak{B}$ сюЁхыхт\-ёъшї ьэю\-цх\-ёЄт
яЁюёЄ\-Ёрэ\-ёЄтр $X$. ┬  ¤Єюь ёыєўрх ьхЁр $\mu$ эр\-ч√тр\-хЄ\-ё 
{\it сюЁх\-ыхт\-ёъющ} ьхЁющ.

\par ╧єёЄ№ $X=K$ -- ъюь\-яръЄ. ╥юуфр ьэю\-цх\-ёЄтю ъю\-эхў\-э√ї
тх∙хёЄ\-тхэ\-э√ї сюЁх\-ыхт\-ёъшї ьхЁ эр  $K$ юсЁрчє\-■Є срэрїю\-тю
яЁюёЄ\-Ёрэ\-ёЄтю. ▌Єю яЁюёЄ\-Ёрэ\-ёЄтю ьюц\-эю юЄюцфхёЄ\-тшЄ№ ё
яЁюёЄ\-Ёрэ\-ёЄ\-тюь ёюяЁ \-ц╕э\-э√ь ъ срэрїю\-тюьє яЁюёЄ\-Ёрэ\-ёЄтє
$C(K)$  -- яЁюёЄ\-Ёрэ\-ёЄтє эх\-яЁх\-Ё√т\-э√ї тх∙хёЄ\-тхэ\-э√ї
Їєэъ\-Ўшщ  $f$ эр  $K$ ё эюЁьющ $\|f\|=\max\limits_{x\in K}|f(x)|$.
▌Єю ёых\-фєхЄ шч Єхю\-Ёх\-ь√ ╨шёёр, ъюЄю\-Ёр  єЄ\-тхЁ\-ц\-фрхЄ, ўЄю
тё ъшщ ышэхщ\-э√щ эх\-яЁх\-Ё√т\-э√щ Їєэъ\-Ўшю\-эры $T$ т
яЁюёЄ\-Ёрэ\-ёЄтх $C(K)$ шьххЄ тшф
$$
(T,f)=\int\limits_K f(x)d\mu(x),
$$
уфх $\mu$ -- эх\-ъюЄю\-Ёр  ъю\-эхў\-эр  сюЁхыхт\-ёър  ьхЁр эр  $K$.

\par ─рыхх тьхёЄю $(T,f)$ ь√ сєфхь яшёрЄ№ $(\mu,f)$.
╚ьххЄ ьхёЄю Ёртхэ\-ёЄ\-тю $\|\mu\|=|\mu|(K)$ ($\|\mu\|$  -- ¤Єю {\it
эюЁьр ышэхщ\-эю\-ую Їєэъ\-Ўшю\-эрыр} $\mu$,
$\|\mu\|=\sup\limits_{\|f\|\leq 1} (\mu,f)$).

\par ┴рэрїютю яЁюёЄ\-Ёрэ\-ёЄтю ъю\-эхў\-э√ї тх∙хёЄ\-тхэ\-э√ї сюЁх\-ыхт\-ёъшї
ьхЁ эр ъюь\-яръЄх $K$ ь√ юсю\-чэрўшь ${\cal{M}}_r(K)$. ┬
яЁюёЄ\-Ёрэ\-ёЄтх ${\cal{M}}_r(K)$ сюыхх трцэє■ Ёюы№, ўхь
ёїюфш\-ьюёЄ№ яю эюЁьх, шуЁрхЄ ёырср  ёїюфш\-ьюёЄ№. ╤юуырёэю
юс∙х\-яЁш\-э Єющ ЄхЁьш\-эю\-ыюушш яю\-ёыхфю\-тр\-Єхы№\-эюёЄ№ $\mu_n$
{\it ёырсю ёїю\-фшЄ\-ё } ъ  $\mu$ (юсю\-чэрўх\-эшх
$\mu=W\lim\limits_{n\to\infty}\mu_n$), хёыш фы  ы■сющ Їєэъ\-Ўшш
$f\in C(K)$ ўшёыю\-тр  яю\-ёыхфю\-тр\-Єхы№\-эюёЄ№ $(\mu_n,f)$
ёїю\-фшЄ\-ё  ъ  $(\mu,f)$.

\par ╚ч Єхю\-Ёх\-ь√ └ырюуыє ёых\-фє\-хЄ, ўЄю хёыш $H$ $\subset$ ${\cal{M}}_r(K)$ ш
$\sup\{|\mu|(K):\:\mu\in H\}$ $<\infty$, Єю є ы■сющ
яю\-ёыхфю\-тр\-Єхы№\-эюё\-Єш $\mu_n\in H$ хёЄ№ ёырсю ёїюф \-∙р \-ё 
яюф\-яю\-ёыхфю\-тр\-Єхы№\-эюёЄ№.

\par ═рЁ фє ё тх∙хёЄ\-тхэ\-э√ьш срэрїю\-т√ь яЁюёЄ\-Ёрэ\-ёЄ\-тюь
${\cal{M}}_r(K)$ Ёрё\-ёьрЄ\-Ёш\-тр\-хЄ\-ё  ъюья\-ыхъё\-эюх
срэр\-їю\-тю яЁюёЄ\-Ёрэ\-ёЄтю ${\cal{M}}_c(K)$, ёю\-ёЄю \-∙хх шч
ъюья\-ыхъё\-э√ї сюЁх\-ыхт\-ёъшї ьхЁ $\mu=\mu_1+i\mu_2$, уфх $\mu_1$
ш  $\mu_2$  -- ъю\-эхў\-э√х тх∙хёЄ\-тхэ\-э√х сюЁх\-ыхт\-ёъшх ьхЁ√.
╧ЁюёЄ\-Ёрэ\-ёЄ\-тю${\cal{M}}_c(K)$ ьюц\-эю юЄюц\-фхёЄ\-тшЄ№ ё
яЁюёЄ\-Ёрэ\-ёЄ\-тюь ёю\-яЁ ц╕э\-э√ь срэр\-їю\-тє яЁюёЄ\-Ёрэ\-ёЄтє
$C(K)$, ёю\-ёЄю \-∙хьє шч эх\-яЁх\-Ё√т\-э√ї ъюья\-ыхъё\-э√ї
Їєэъ\-Ўшщ  $f$ эр  $K$, яЁш\-ў╕ь $\|f\|=\max\limits_{x\in K}
|f(x)|$. ┬ яЁюёЄ\-Ёрэ\-ёЄтх ${\cal{M}}_c(K)$ Єръ\-цх
т√\-яюыэ \-хЄ\-ё  Ёртхэ\-ёЄ\-тю $\|\mu\|=\sup\limits_{\|f\|\leq 1}
|(\mu,f)|=|\mu|(K)$.

\par ┬ ъюья\-ыхъё\-эюь ёыєўрх ьхЁр  $|\mu|$ юяЁх\-фхы \-хЄё  ёыюц\-эхх,
эю ёяЁр\-тхф\-ыш\-т√ ёыхфє\-■\-∙шх яЁюёЄ√х эх\-Ёртхэ\-ёЄ\-тр. ┼ёыш
$\mu=\mu_1+i\mu_2$, Єю $|\mu_1|\leq |\mu|$, $|\mu_2|\leq |\mu|$,
$|\mu|\leq |\mu_1|+|\mu_2|$.

\par ╥хяхЁ№ ь√ яхЁх\-їюфшь ъ Ёрё\-ёьюЄЁх\-эш■ юёэют\-эю\-ую юс·хъ\-Єр
эр°хую шё\-ёыхфю\-тр\-эш   -- яЁюёЄ\-Ёрэ\-ёЄтє Ёр\-фю\-эю\-т√ї ьхЁ
эр яюыє\-юёш $(0,\infty)$.

\par ┬эрўрых юяЁх\-фхышь {\it яЁюёЄ\-Ёрэ\-ёЄтю ЄхёЄю\-т√ї Їєэъ\-Ўшщ} $\Phi$
эр яюыє\-юёш $(0,\infty)$ (ь√ яю\-чршь\-ёЄ\-тю\-трыш ЄхЁьшэ шч
ЄхюЁшш юсюс∙╕э\-э√ї Їєэъ\-Ўшщ, р юсю\-чэрўх\-эшх  -- шч ъэшуш
╦рэф\-ъюЇр \cite{Landkof}). ╘єэъ\-Ўш  $f\in\Phi$, хёыш $f$  --
эх\-яЁх\-Ё√т\-эр  Їєэъ\-Ўш  эр яюыє\-юёш $(0,\infty)$ ш
ёє\-∙хёЄ\-тєхЄ Єр\-ъющ ёху\-ьхэЄ $[a,b]$, ўЄю т√\-яюыэ \-хЄ\-ё 
ёююЄ\-эю°х\-эшх $\supp f$ $\subset$ $[a,b]$ $\subset(0,\infty)$. ╠√
сєфхь Ёрё\-ёьрЄ\-Ёш\-трЄ№ ъръ тх∙хёЄ\-тхэ\-э√х, Єръ ш
ъюья\-ыхъё\-э√х яЁюёЄ\-Ёрэ\-ёЄ\-тр  $\Phi$.

\par ┬ яЁюёЄ\-Ёрэ\-ёЄ\-тх $\Phi$ шчтхёЄ\-э√ь ёяюёю\-сюь ттюфшЄ\-ё 
яюэ \-Єшх ёїюфш\-ьюё\-Єш. ╧юёыхфю\-тр\-Єхы№\-эюёЄ№ Їєэъ\-Ўшщ $f_n$
{\it ёїю\-фшЄ\-ё } ъ Їєэъ\-Ўшш  $f$ т яЁюёЄ\-Ёрэ\-ёЄ\-тх $\Phi$,
хёыш ёє\-∙хёЄ\-тєхЄ ёху\-ьхэЄ $[a,b]\subset(0,\infty)$ Єр\-ъющ, ўЄю
фы  ы■сюую $n$ $\supp f_n$ $\subset$ $[a,b]$ ш
яю\-ёыхфю\-тр\-Єхы№\-эюёЄ№ $f_n(x)$ Ёртэю\-ьхЁэю ёїю\-фшЄ\-ё  ъ
$f(x)$ эр яюыє\-юёш $(0,\infty)$.

\par ─ы  Їєэъ\-Ўшщ $f\in\Phi$ ь√ сєфхь шё\-яюы№чю\-трЄ№ юсю\-чэрўх\-эшх
$\|f\|=\max|f(x)|$.

\par ┴юЁхыхтёър  ьхЁр $\mu$ эр яюыє\-юёш $(0,\infty)$ эр\-ч√тр\-хЄ\-ё 
{\it ыю\-ъры№\-эю ъю\-эхў\-эющ}, хёыш фы  ы■сюую ёху\-ьхэЄр $[a,b]$
$\subset$ $(0,\infty)$ т√\-яюыэ \-хЄ\-ё  эх\-Ёртхэ\-ёЄ\-тю
$|\mu|([a,b])<\infty$.

\par ╘єэъЎш  ьэю\-цх\-ёЄт $\mu$ эр\-ч√тр\-хЄё  {\it тх∙хёЄ\-тхэ\-эющ
Ёр\-фю\-эю\-тющ} ьхЁющ эр яюыє\-юёш $(0,\infty)$, хёыш юэр
яЁхф\-ёЄртш\-ьр т тшфх $\mu$ $=$ $\mu_1-\mu_2$, уфх $\mu_1$, $\mu_2$
-- тчршьэю ёшэує\-ы Ё\-э√х ыю\-ъры№эю ъю\-эхў\-э√х
яюыю\-цш\-Єхы№\-э√х сюЁх\-ыхт\-ёъшх ьхЁ√ эр яюыє\-юёш $(0,\infty)$.

\par ─ы  ы■сюую сюЁхыхт\-ёъюую ьэю\-цх\-ёЄтр $A\subset(0,\infty)$
Ёртхэ\-ёЄ\-тюь $\mu_A=(\mu_1)_A-(\mu_2)_A$ юяЁх\-фхы \-хЄ\-ё  {\it
юуЁрэш\-ўх\-эшх $\mu$ эр ьэю\-цх\-ёЄтю $A$}. ╒юЄ  ттхф╕э\-эр 
Їєэъ\-Ўш  ьэю\-цх\-ёЄт $\mu$ эх  ты \-хЄё  сюЁх\-ыхт\-ёъющ ьхЁющ эр
яюыє\-юёш $(0,\infty)$ (юэр эх юяЁх\-фхых\-эр эр Єхї сюЁх\-ыхт\-ёъшї
ьэю\-цх\-ёЄ\-трї  $E$, фы  ъюЄюЁ√ї т√\-яюыэ \-хЄё  юфэю шч
ёююЄ\-эю°х\-эшщ $\mu_1(E)=\mu_2(E)=+\infty$,
$\mu_1(E)=\mu_2(E)=-\infty$), Єхь эх ьхэхх, юуЁрэш\-ўх\-эшх
$\mu_{[a,b]}$ фы  ы■сюую ёхуьхэ\-Єр $[a,b]$ $\subset$ $(0,\infty)$
 ты \-хЄё  ъю\-эхў\-эющ сюЁх\-ыхт\-ёъющ ьхЁющ.

\par ╥хь ёрь√ь фы  ы■сющ Їєэъ\-Ўшш $f\in\Phi$ шьххЄ ёь√ёы
$(\mu,f)=\int\limits_0^\infty f(x)d\mu(x)$. ╬ўхтшф\-эю, ўЄю Єръ
юяЁх\-фхы╕э\-эр  Їєэъ\-Ўш  $(\mu,f)$  ты \-хЄё  ышэхщ\-э√ь
Їєэъ\-Ўшю\-эр\-ыюь эр яЁюёЄ\-Ёрэ\-ёЄ\-тх $\Phi$. ╬ўхтшф\-эю Єръ\-цх,
ўЄю ¤ЄюЄ Їєэъ\-Ўшю\-эры эх\-яЁх\-Ё√тхэ эр  $\Phi$, Єю хёЄ№ шч
ёїюфш\-ьюё\-Єш $f_n$ ъ  $f$ т яЁюёЄ\-Ёрэ\-ёЄ\-тх  $\Phi$
ёых\-фє\-хЄ, ўЄю $(\mu,f_n)\to(\mu,f)$.

\par ╚чтхёЄэю ш юсЁрЄэюх. ┼ёыш $T$  -- эх\-яЁх\-Ё√т\-э√щ ышэхщ\-э√щ
Їєэъ\-Ўшю\-эры эр яЁюёЄ\-Ёрэ\-ёЄ\-тх  $\Phi$, Єю ёє\-∙хёЄ\-тє\-хЄ
Ёр\-фю\-эю\-тр ьхЁр  $\mu$ Єр\-ър , ўЄю $(T,f)=(\mu,f)$ фы  ы■сющ
Їєэъ\-Ўшш $f\in\Phi$. ╚чтхёЄэю Єръ\-цх, ўЄю ёє\-∙хёЄ\-тєхЄ
хфшэ\-ёЄ\-тхэ\-эр  Ёр\-фю\-эютр ьхЁр  $\mu$, юсырфр■\-∙р  ¤Єшь
ётющ\-ёЄ\-тюь.

\par ┴єЁсръш яхЁтю\-эрўры№\-эю юяЁх\-фхы хЄ {\it Ёр\-фюэют√} ьхЁ√ т
яЁю\-шч\-тюы№\-эюь  ыю\-ъры№\-эю ъюь\-яръЄ\-эюь Єюяю\-ыюуш\-ўхёъюь
яЁюёЄ\-Ёрэ\-ёЄтх ъръ эх\-яЁх\-Ё√т\-э√х ышэхщ\-э√х Їєэъ\-Ўшю\-эры√ эр
яЁюёЄ\-Ёрэ\-ёЄ\-тх  эх\-яЁх\-Ё√т\-э√ї Їєэъ\-Ўшщ  $f$ ё
ъюь\-яръЄ\-э√ь эюёш\-Єх\-ыхь. ╘єэъ\-Ўшю\-эры $T$ эр\-ч√тр\-хЄ\-ё 
{\it яюыю\-цш\-Єхы№\-э√ь}, хёыш $(T,f)\geq 0$ яЁш $f\geq 0$. ╬э
фю\-ърч√\-тр\-хЄ, ўЄю тё ъшщ эх\-яЁх\-Ё√т\-э√щ ышэхщ\-э√щ
яюыю\-цш\-Єхы№\-э√щ Їєэъ\-Ўшю\-эры ёют\-ярфр\-хЄ ё
яюыю\-цш\-Єхы№\-эющ ыю\-ъры№\-эю ъю\-эхў\-эющ сюЁхыхт\-ёъющ ьхЁющ эр
Єюяю\-ыюуш\-ўхё\-ъюь яЁюёЄ\-Ёрэ\-ёЄ\-тх. ╟рЄхь юэ фю\-ърч√\-трхЄ,
ўЄю  ы■сющ эх\-яЁх\-Ё√т\-э√щ  ышэхщ\-э√щ Їєэъ\-Ўшю\-эры
яЁхф\-ёЄрты \-хЄё  т тшфх  Ёрч\-эюё\-Єш яюыю\-цш\-Єхы№\-э√ї
ышэхщ\-э√ї эх\-яЁх\-Ё√т\-э√ї Їєэъ\-Ўшю\-эр\-ыют.

\par ╥ръшь юсЁр\-чюь, т ёыєўрх, ъюуфр ыю\-ъры№\-эю ъюь\-яръЄ\-эюх
Єюяю\-ыюуш\-ўхё\-ъюх яЁюёЄ\-Ёрэ\-ёЄтю хёЄ№ яюыє\-юё№ $(0,\infty)$,
яЁхф\-ыю\-цхэ\-эюх т ёЄрЄ№х юяЁх\-фхых\-эшх Ёр\-фю\-эю\-тющ ьхЁ√
¤ътш\-трыхэЄ\-эю юяЁх\-фхых\-эш■ ┴єЁ\-сръш.

\par ╠эю\-цх\-ёЄтю тх∙хёЄ\-тхэ\-э√ї Ёр\-фю\-эю\-т√ї ьхЁ эр яюыє\-юёш $(0,\infty)$
юсЁрчє\-хЄ тх∙хёЄ\-тхэ\-эюх ышэхщ\-эюх яЁюёЄ\-Ёрэ\-ёЄтю. ╠√ сєфхь
юсю\-чэр\-ўрЄ№ хую ўхЁхч  $\frak{R}$.

 \par ╠√ сєфхь Єръ\-цх Ёрё\-ёьрЄЁш\-трЄ№ ышэхщ\-эюх ъюья\-ыхъё\-эюх
 яЁюёЄ\-Ёрэ\-ёЄ\-тю ${\frak{R}}_c$, ёю\-ёЄю \-∙хх шч ъюья\-ыхъё\-э√ї
 Ёр\-фю\-эю\-т√ї ьхЁ $\mu=\mu_1+i\mu_2$, уфх  $\mu_1$ ш  $\mu_2$  --
 тх∙хёЄ\-тхэ\-э√х Ёр\-фю\-эю\-т√ ьхЁ√.

 \par ┬ яЁюёЄ\-Ёрэ\-ёЄтрї $\frak{R}$, ${\frak{R}}_c$ ттюфшЄ\-ё  яюэ \-Єшх
 °ш\-Ёю\-ъющ ёїюфш\-ьюёЄш. ▌ЄюЄ ЄхЁьшэ яЁш\-эрф\-ыхцшЄ ┴єЁсръш.
 ╧юёыхфю\-тр\-Єхы№\-эюёЄ№ Ёр\-фю\-эю\-т√ї ьхЁ $\mu_n$ {\it °ш\-Ёю\-ъю ёїю\-фшЄ\-ё }
 ъ Ёр\-фю\-эю\-тющ ьхЁх $\mu$ (юсю\-чэрўх\-эшх $\mu=\lim\limits_{n\to\infty}\mu_n$
 шыш $\mu_n\to\mu$), хёыш фы  ы■сющ Їєэъ\-Ўшш $f\in\Phi$ ўшёыю\-тр 
 яю\-ёыхфю\-тр\-Єхы№\-эюёЄ№ $(\mu_n,f)$ ёїю\-фшЄ\-ё  ъ  $(\mu,f)$.
 ╥хЁьшэ °ш\-Ёюър  ёїюфш\-ьюёЄ№ чр\-эшьрхЄ т эр°хщ ёЄрЄ№х яЁш\-тхыш\-уш\-Ёю\-трэ\-эюх
 яю\-ыюцх\-эшх. ┼ёыш яш°хЄ\-ё  $\mu_n\to\mu$ ш эх єЄюў\-э \-хЄ\-ё ,
 ю ъръюь тшфх ёїюфш\-ьюё\-Єш шф╕Є Ёхў№, Єю шьххЄ\-ё  т тшфє °ш\-Ёюър 
 ёїюфш\-ьюёЄ№. ┬ яЁхф\-ёЄрт\-ыхэ\-эющ ёЄрЄ№х Ёрё\-ёьрЄЁш\-тр\-■Є\-ё 
 Ёрч\-ышў\-э√х тшф√ ёїюфш\-ьюёЄш яю\-ёыхфю\-тр\-Єхы№\-эюёЄхщ ьхЁ $\mu_n$.

\par ╠эю\-цх\-ёЄтю $E$ эрч√тр\-хЄё  {\it шч\-ьхЁш\-ь√ь яю ╞юЁ\-фр\-эє}
юЄ\-эюёш\-Єхы№\-эю Ёр\-фю\-эю\-тющ ьхЁ√ $\mu$, хёыш $|\mu|(\partial
E)=0$.

\par ╒юЁю°ю шчтхёЄэ√ ёых\-фє\-■\-∙шх єЄ\-тхЁ\-цфх\-эш .
%%%%%%%%%%%%%%%%%%%%%%%%%%%%%%%%%%%%%%%%%%%%%%%%%%%%%%%%%%%%%%%%%%%%%%%%%%%%%%%%%%%%%
\begin{theorem}\label{measures_th_1}\hskip-2mm{.}\:
%%%%%%%%%%%%%%%%%%%%%%%%%%%%%%%%%%%%%%%%%%%%%%%%%%%%%%%%%%%%%%%%%%%%%%%%%%%%%%%%%%%%%
╧єёЄ№ $\mu_n\to \mu$, $|\mu_n|\to\hat{\mu}$. ┼ёыш  ьэю\-цх\-ёЄтю $E$
шч\-ьхЁш\-ью яю ╞юЁфр\-эє юЄ\-эюёш\-Єхы№\-эю ьхЁ√ $\hat{\mu}$, Єю
$(\mu_n)_E\to\mu_E$.
\end{theorem}
%%%%%%%%%%%%%%%%%%%%%%%%%%%%%%%%%%%%%%%%%%%%%%%%%%%%%%%%%%%%%%%%%%%%%%%%%%%%%%%%%%%%%
\begin{theorem}\label{measures_th_2}\hskip-2mm{.}\:
%%%%%%%%%%%%%%%%%%%%%%%%%%%%%%%%%%%%%%%%%%%%%%%%%%%%%%%%%%%%%%%%%%%%%%%%%%%%%%%%%%%%%
╧єёЄ№  $\mu_n\to \mu$, $|\mu_n|\to\hat{\mu}$. ╧єёЄ№ $K\subset$
$(0,\infty)$ -- ъюь\-яръЄ шч\-ьхЁш\-ь√щ яю ╞юЁфр\-эє
юЄ\-эюёш\-Єхы№эю ьхЁ√ $\hat{\mu}$. ╥юуфр яю\-ёыхфю\-тр\-Єхы№\-эюёЄ№
$(\mu_n)_K$ ёырсю ёїю\-фшЄ\-ё  ъ  $\mu_K$.
\end{theorem}
%%%%%%%%%%%%%%%%%%%%%%%%%%%%%%%%%%%%%%%%%%%%%%%%%%%%%%%%%%%%%%%%%%%%%%%%%%%%%%%%%%%%%
\par ─ы  яю\-ёыхфю\-тр\-Єхы№\-эюё\-Єш яюыю\-цш\-Єхы№\-э√ї ьхЁ $\mu_n$
фю\-ърчр\-Єхы№\-ёЄтр ёЇюЁ\-ьє\-ыш\-Ёю\-трэ\-э√ї єЄ\-тхЁ\-ц\-фх\-эшщ
ьюц\-эю эрщЄш т  \cite{Landkof}, ттхфх\-эшх, \S1. ┬ юс∙хь ёыєўрх
эєцэю яЁш\-ьхэ Є№ ¤Єш єЄ\-тхЁ\-ц\-фх\-эш  юЄ\-фхы№эю ъ ьхЁрь
$(\mu_n)_+$ ш $(\mu_n)_-$. ═рь сєфхЄ эєцэр ёых\-фє\-■\-∙р 
Єхю\-Ёх\-ьр, ъюЄю\-Ёє■ ьюц\-эю яю\-ыєўшЄ№ ъръ ёыхф\-ёЄ\-тшх
Єхю\-Ёх\-ь√  \ref{measures_th_1}.

%%%%%%%%%%%%%%%%%%%%%%%%%%%%%%%%%%%%%%%%%%%%%%%%%%%%%%%%%%%%%%%%%%%%%%%%%%%%%%%%%%%%%
\begin{theorem}\label{measures_th_3}\hskip-2mm{.}\:
%%%%%%%%%%%%%%%%%%%%%%%%%%%%%%%%%%%%%%%%%%%%%%%%%%%%%%%%%%%%%%%%%%%%%%%%%%%%%%%%%%%%%
╧єёЄ№ яю\-ёыхфю\-тр\-Єхы№\-эюёЄ№ $\mu_n$ Ёр\-фю\-эю\-т√ї ьхЁ эр
яюыє\-юёш $(0,\infty)$ °ш\-Ёю\-ъю ёїю\-фшЄ\-ё  ъ Ёр\-фю\-эю\-тющ
ьхЁх $\nu$. ╧єёЄ№, ъЁюьх Єюую, $|\mu_n|\to \hat{\nu}$ ш
$\hat{\nu}(\{\xi\})=0$. ╥юуфр фы  ы■сющ Їєэъ\-Ўшш $\varphi\in\Phi$
т√\-яюыэ \-■Є\-ё  Ёртхэ\-ёЄ\-тр
$$
\lim\limits_{n\to\infty} \int\limits_0^\infty
\varphi(t)\chi_{(0,\xi]}(t) d\mu_n(t)= \int\limits_0^\infty
\varphi(t)\chi_{(0,\xi]}(t) d\nu(t),
$$
$$
\lim\limits_{n\to\infty} \int\limits_0^\infty
\varphi(t)\chi_{(\xi,\infty)}(t) d\mu_n(t)= \int\limits_0^\infty
\varphi(t)\chi_{(\xi,\infty)}(t) d\nu(t).
$$
\end{theorem}
%%%%%%%%%%%%%%%%%%%%%%%%%%%%%%%%%%%%%%%%%%%%%%%%%%%%%%%%%%%%%%%%%%%%%%%%%%%%%%%%%%%%%
\par\quad {\sc  ─юърчрЄхы№ёЄтю.}\: ┼ёыш т Єхю\-Ёхьх  \ref{measures_th_1} т
ърўхёЄ\-тх ьэю\-цх\-ёЄ\-тр  $E$ т√\-сшЁрЄ№ ьэю\-цх\-ёЄтр $(0,\xi]$,
$(\xi,\infty)$, Єю яюыєўшь єЄ\-тхЁ\-ц\-фх\-эш  Єхю\-Ёх\-ь√.
╥хю\-Ёх\-ьр фю\-ърчр\-эр.
%%%%%%%%%%%%%%%%%%%%%%%%%%%%%%%%%%%%%%%%%%%%%%%%%%%%%%%%%%%%%%%%%%%%%%%%%%%%%%%%%%%%%
\par ╠эю\-цхёЄ\-тю  $H\subset\frak{R}_c(\frak{R})$ эр\-ч√тр\-хЄ\-ё  {\it °ш\-Ёю\-ъю
юуЁрэш\-ўхэ\-э√ь}, хёыш фы  ы■сющ Їєэъ\-Ўшш $\varphi\in\Phi$
ьэю\-цх\-ёЄтю $\{(\mu,f):\,\mu\in H\}$  ты \-хЄё  юуЁрэш\-ўхэ\-э√ь.

\par ╠эю\-цхёЄ\-тю $H$ эр\-ч√тр\-хЄ\-ё  {\it ёшы№эю
юуЁрэш\-ўхэ\-э√ь}, хёыш фы  ы■сю\-ую ёхуьхэ\-Єр $[a,b]$ $\subset$
$(0,\infty)$ ьэю\-цх\-ёЄтю $\{|\mu|([a,b]):\,\mu\in H\}$  ты \-хЄё 
юуЁрэш\-ўхэ\-э√ь.

\par ╠эю\-цхёЄ\-тю $H$ эрч√тр\-хЄ\-ё  {\it ъюь\-яръЄ\-э√ь}, хёыш
шч ы■сющ яюёых\-фю\-тр\-Єхы№\-эюё\-Єш $\mu_n\in H$ ьюц\-эю
т√\-фхышЄ№ °ш\-Ёю\-ъю ёїюф \-∙є■\-ё 
яюф\-яюёых\-фю\-тр\-Єхы№\-эюёЄ№.

\par Cыхфє\-■∙хх єЄ\-тхЁ\-ц\-фхэшх -- ¤Єю уырт\-э√щ Ёх\-чєы№\-ЄрЄ шч ЄхюЁшш ьхЁ√,
ъюЄю\-Ё√щ ь√ яЁш\-ьхэ \-хь т эр°хщ Ёрсю\-Єх.

%%%%%%%%%%%%%%%%%%%%%%%%%%%%%%%%%%%%%%%%%%%%%%%%%%%%%%%%%%%%%%%%%%%%%%%%%%%%%%%%%%%%%
\begin{theorem}\label{measures_th_4}\hskip-2mm{.}\:
%%%%%%%%%%%%%%%%%%%%%%%%%%%%%%%%%%%%%%%%%%%%%%%%%%%%%%%%%%%%%%%%%%%%%%%%%%%%%%%%%%%%%
╩ырёё√ °ш\-Ёю\-ъю юуЁрэш\-ўхэ\-э√ї, ёшы№эю юуЁрэш\-ўхэ\-э√ї ш
ъюь\-яръЄ\-э√ї ьэю\-цх\-ёЄт т яЁюёЄ\-Ёрэ\-ёЄтх $\frak{R}$
$(\frak{R}_c)$ ёют\-ярфр\-■Є.
\end{theorem}

\par ─юърчр\-Єхы№\-ёЄтю ьюц\-эю эрщЄш т \cite{Burbaki}, уыртр  3, \S1,
яЁхф\-ыюцх\-эшх  15  ш яЁш\-ьхўр\-эшх ъ эхьє.

\par ╟рьхЄшь, ўЄю ы╕уъшь ёыхфёЄ\-тшхь ёЇюЁ\-ьє\-ыш\-Ёю\-трэ\-эющ Єхю\-Ёх\-ь√
 ты \-хЄё  ёых\-фє\-■\-∙хх єЄ\-тхЁ\-ц\-фх\-эшх.
%%%%%%%%%%%%%%%%%%%%%%%%%%%%%%%%%%%%%%%%%%%%%%%%%%%%%%%%%%%%%%%%%%%%%%%%%%%%%%%%%%%%%
\begin{theorem}\label{measures_th_4-1}\hskip-2mm{.}\:
%%%%%%%%%%%%%%%%%%%%%%%%%%%%%%%%%%%%%%%%%%%%%%%%%%%%%%%%%%%%%%%%%%%%%%%%%%%%%%%%%%%%%
╬ЄюсЁрцхэшх $(\mu,\varphi):\: \frak{R}_c\times \Phi\to \mathbb{C}$
 ты \-хЄё  эх\-яЁх\-Ё√т\-э√ь яю ёютю\-ъєя\-эюё\-Єш яхЁх\-ьхэ\-э√ї.
\end{theorem}
%%%%%%%%%%%%%%%%%%%%%%%%%%%%%%%%%%%%%%%%%%%%%%%%%%%%%%%%%%%%%%%%%%%%%%%%%%%%%%%%%%%%%
\par {\sc  ─юърчрЄхы№ёЄтю.}\: ╧єёЄ№ $\mu_n\to \mu$,
$\varphi_n\to\varphi$. ╧ю\-ёъюы№ъє $\varphi_n\to \varphi$, Єю
ёє\-∙хёЄ\-тєхЄ ёху\-ьхэЄ $[a,b]$ $\subset$ $(0,\infty)$ Єр\-ъющ, ўЄю
$\supp\varphi_n$ $\subset$ $[a,b]$, $\supp\varphi$ $\subset$
$[a,b]$. ╧ю\-ёъюы№ъє $\mu_n\to \mu$, Єю яю\-ёыхфю\-тр\-Єхы№\-эюёЄ№
$\mu_n$  ты \-хЄё  °ш\-Ёю\-ъю юуЁрэш\-ўхэ\-эющ, р чэрўшЄ ш ёшы№эю
юуЁрэш\-ўхэ\-эющ. ╧ю¤Єю\-ьє ёє\-∙хёЄ\-тєхЄ яюёЄю э\-эр   $M$
Єр\-ър , ўЄю т√\-яюыэ \-■Є\-ё  эх\-Ёртхэ\-ёЄ\-тр $|\mu_n|([a,b])$
$\leq$ $M$. ─рыхх шьххь
$$
|(\mu_n,\varphi_n)-(\mu,\varphi)|\leq
|\mu_n|([a,b])\|\varphi_n-\varphi\|+ |((\mu-\mu_n),\varphi)| \leq
M\|\varphi_n-\varphi\|+|(\mu_n-\mu),\varphi|.
$$
╚ч ¤Єюую ёых\-фєхЄ єЄ\-тхЁ\-ц\-фхэшх Єхю\-Ёх\-ь√. ╥хю\-Ёхьр
фюърчрэр.

%%%%%%%%%%%%%%%%%%%%%%%%%%%%%%%%%%%%%%%%%%%%%%%%%%%%%%%%%%%%%%%%%%%%%%%%%%%%%%%%%%%%%
\begin{theorem}\label{measures_th_5}\hskip-2mm{.}\:
%%%%%%%%%%%%%%%%%%%%%%%%%%%%%%%%%%%%%%%%%%%%%%%%%%%%%%%%%%%%%%%%%%%%%%%%%%%%%%%%%%%%%
─ы  Єюую, ўЄюс√ яюёыхфю\-тр\-Єхы№\-эюёЄ№ яюыю\-цш\-Єхы№\-э√ї
Ёр\-фю\-эю\-т√ї ьхЁ $\mu_n$ °ш\-Ёю\-ъю ёїюфш\-ырё№ ъ Ёр\-фю\-эю\-тющ
ьхЁх  $\mu$, фю\-ёЄр\-Єюў\-эю, ўЄюс√ ёююЄ\-эю°х\-эшх $(\mu_n,f)$
$\to$ $(\mu,f)$ т√\-яюыэ \-ыюё№ эр тё■фє яыюЄ\-эюь ьэю\-цх\-ёЄ\-тх т
яЁюёЄ\-Ёрэ\-ёЄ\-тх  $\Phi$.
\end{theorem}
%%%%%%%%%%%%%%%%%%%%%%%%%%%%%%%%%%%%%%%%%%%%%%%%%%%%%%%%%%%%%%%%%%%%%%%%%%%%%%%%%%%%%
\begin{theorem}\label{measures_th_6}\hskip-2mm{.}\:
%%%%%%%%%%%%%%%%%%%%%%%%%%%%%%%%%%%%%%%%%%%%%%%%%%%%%%%%%%%%%%%%%%%%%%%%%%%%%%%%%%%%%
╧єёЄ№ Ёр\-фюэют√ ьхЁ√ $\mu$ ш $\nu$ Єръют√, ўЄю Ёртхэ\-ёЄ\-тю
$(\mu,f)=(\nu,f)$ т√\-яюыэ \-хЄ\-ё  эр тё■фє яыюЄ\-эюь
ьэю\-цх\-ёЄ\-тх т яЁюёЄ\-Ёрэ\-ёЄ\-тх  $\Phi$. ╥юуфр $\mu=\nu$.
\end{theorem}
%%%%%%%%%%%%%%%%%%%%%%%%%%%%%%%%%%%%%%%%%%%%%%%%%%%%%%%%%%%%%%%%%%%%%%%%%%%%%%%%%%%%%
\par ─юърчр\-Єхы№\-ёЄтр Єхю\-Ёхь  \ref{measures_th_5} ш  \ref{measures_th_6}
ьюц\-эю эрщЄш т  \cite{Landkof}, ттхфх\-эшх, \S1.

\par ╬с√ўэю ь√ сєфхь шё\-яюы№чю\-трЄ№ °ш\-Ёю\-ъє■
ёїю\-фш\-ьюё\-Є№ Ёр\-фю\-эю\-т√ї ьхЁ. ╠хцфє Єхь, ўрёЄю
Ёрё\-ёєц\-фх\-эш  єфюс\-эхх яЁю\-тю\-фшЄ№ т ьхЄЁш\-ўхё\-ъшї
яЁюёЄ\-Ёрэ\-ёЄтрї. ╧ю\-¤Єюьє ь√ сє\-фхь шёяюы№\-чю\-трЄ№ т
яЁюёЄ\-Ёрэ\-ёЄтх $\frak{R}_c$ ёых\-фє\-■\-∙шх шчтхёЄ\-э√х ьхЄЁш\-ъш.
╧єёЄ№ $\{\varphi_n: n=1,2,...\}$ -- ёў╕Є\-эюх тё■фє яыюЄ\-эюх
ьэю\-цхёЄ\-тю т яЁюёЄ\-Ёрэ\-ёЄтх  $\Phi$. ▌Єю ючэр\-ўр\-хЄ, ўЄю фы 
ы■сющ Їєэъ\-Ўшш $\varphi\in\Phi$ ёє\-∙хёЄ\-тє\-хЄ
яюф\-яю\-ёыхфю\-тр\-Єхы№\-эюёЄ№ $\varphi_{n_k}$
яю\-ёыхфю\-тр\-Єхы№\-эюё\-Єш $\varphi_n$ Єр\-ър , ўЄю
$\varphi_{n_k}\to \varphi$ т яЁюёЄ\-Ёрэ\-ёЄтх  $\Phi$. ─рыхх яю
яю\-ёыхфю\-тр\-Єхы№\-эюё\-Єш $\varphi_n$ юяЁх\-фхы \-хь Їєэъ\-Ўш■

\begin{equation}\label{measures_def_metrika_d}
d(\mu_1,\mu_2)=\sum\limits_{n=1}^{\infty}
\frac{|(\mu_1-\mu_2)(\varphi_n)|}{2^n(1+|(\mu_1-\mu_2)(\varphi_n)|)},
\end{equation}
$\mu_1,\mu_2\in\frak{R}_c$. ╦хуъю яЁю\-тхЁ \-хЄ\-ё  ўЄю  $d$ --
ьхЄЁш\-ър т яЁюёЄ\-Ёрэ\-ёЄтх  $\frak{R}_c$.

\par ┴юыхх Єюую, ыхуъю тшфхЄ№, ўЄю шч ёююЄ\-эю°х\-эш  $\mu_k\to \mu$
ёых\-фє\-хЄ, ўЄю $d(\mu_k,\mu)\to 0$. ╬фэръю, юсЁрЄ\-эюх
єЄ\-тхЁ\-ц\-фх\-эшх эх\-тхЁ\-эю. ╧єёЄ№, эр\-яЁш\-ьхЁ, тёх Їєэъ\-Ўшш
$\varphi_n$ эх\-яЁх\-Ё√т\-эю фшЇЇх\-Ёхэ\-Ўш\-Ёєх\-ь√, р
$$
\mu_k=k\left(\delta\Bigl(x-1-\frac{1}{2k}\Bigr)-
\delta\Bigl(x-1+\frac{1}{2k}\Bigr)\right)
$$
$$
-\sqrt{2}k\left(\delta\Bigl(x-1-\frac{1}{2\sqrt{2}k}\Bigr)-
\delta\Bigl(x-1+\frac{1}{2\sqrt{2}k}\Bigr)\right),
$$
уфх $\delta(x-a)$  -- ьхЁр ─шЁрър, ёю\-ёЁх\-фю\-Єю\-ўхэ\-эр  т
Єюў\-ъх  $a$. ╥юу\-фр $d(\mu_k,0)\to 0$, т Єю тЁхь  ъръ
ёююЄ\-эю\-°х\-эшх $\mu_k\to 0$ эх т√яюы\-э \-хЄ\-ё . ╬фэръю,
ёяЁртхф\-ыш\-тр ёых\-фє\-■\-∙р  Єхю\-Ёх\-ьр.
%%%%%%%%%%%%%%%%%%%%%%%%%%%%%%%%%%%%%%%%%%%%%%%%%%%%%%%%%%%%%%%%%%%%%%%%%%%%%%%%%%%%%%%
\begin{theorem}\label{measures_th_7}\hskip-2mm{.}\:
%%%%%%%%%%%%%%%%%%%%%%%%%%%%%%%%%%%%%%%%%%%%%%%%%%%%%%%%%%%%%%%%%%%%%%%%%%%%%%%%%%%%%%%
┼ёыш $\mu_k$ -- ъюь\-яръЄ\-эр  яю\-ёыхфю\-тр\-Єхы№\-эюёЄ№ т
$\frak{R}_c$ ш $d(\mu_k,\mu)\to 0$, Єю $\mu_k$ °ш\-Ёю\-ъю
ёїю\-фшЄ\-ё  ъ  $\mu$.
\end{theorem}
%%%%%%%%%%%%%%%%%%%%%%%%%%%%%%%%%%%%%%%%%%%%%%%%%%%%%%%%%%%%%%%%%%%%%%%%%%%%%%%%%%%%%%%
\par {\sc  ─юърчрЄхы№ёЄтю.}\: ┼ёыш єЄ\-тхЁ\-ц\-фх\-эшх Єхю\-Ёх\-ь√
эх\-тхЁэю, Єю ёє\-∙хёЄ\-тє\-хЄ Їєэъ\-Ўш  $\varphi\in \Phi$ ш фтх
яюф\-яю\-ёыхфю\-тр\-Єхы№\-эюё\-Єш $\mu_{k_p^1}$, $\mu_{k_p^2}$
яю\-ёыхфю\-тр\-Єхы№\-эюё\-Єш $\mu_k$ Єръшх, ўЄю
$$
\lim\limits_{p\to\infty}(\mu_{k_p^1},\varphi) \neq
\lim\limits_{p\to\infty}(\mu_{k_p^2},\varphi).
$$
╧єёЄ№ $\nu_p=\mu_{k_p^1}-\mu_{k_p^2}$, р $\varphi_n$ -- Єр
яю\-ёыхфю\-тр\-Єхы№\-эюёЄ№ Їєэъ\-Ўшщ шч  $\Phi$, ъюЄю\-Ёр 
юяЁхфх\-ы \-хЄ ьхЄЁш\-ъє  $d$. ╧ю\-ёъюы№\-ъє $\varphi_n$ тё■фє
яыюЄ\-эр  яю\-ёыхфю\-тр\-Єхы№\-эюёЄ№ т  $\Phi$, Єю ёє\-∙хёЄ\-тє\-хЄ
яюф\-яю\-ёыхфю\-тр\-Єхы№\-эюёЄ№ $\psi_n$
яю\-ёыхфю\-тр\-Єхы№\-эюё\-Єш $\varphi_n$, ъюЄю\-Ёр  ёїю\-фшЄ\-ё  ъ
$\varphi$ т яЁюёЄ\-Ёрэ\-ёЄ\-тх  $\Phi$. ╤є∙хёЄ\-тє\-хЄ ёху\-ьхэЄ
$[a,b]$ эр яюыє\-юёш $(0,\infty)$, Єр\-ъющ ўЄю $\supp\psi_n$
$\subset$ $[a,b]$ фы  ы■сюую  $n$. ╧ю¤Єю\-ьє ё эх\-ъю\-Єю\-Ёющ
ъюэ\-ёЄрэ\-Єющ $M$, эх чртшё \-∙хщ юЄ  $n$, т√\-яюыэ \-хЄ\-ё 
эх\-Ёртхэ\-ёЄ\-тю
$\mathop{\overline{\lim}}\limits_{p\to\infty}|\nu_p(\varphi)|$
$\leq$ $ M\|\varphi-\psi_n\|$. ╚ч ¤Єюую ёых\-фє\-хЄ, ўЄю
$\nu_p(\varphi)\to 0$. ▌Єю яЁюЄш\-тю\-ЁхўшЄ т√\-сюЁє  $\varphi$.
╥хю\-Ёх\-ьр фю\-ърчр\-эр.

\par ╥ръшь юсЁр\-чюь, т юс∙хь ёыєўрх ёїю\-фш\-ьюёЄ№ т
ьхЄЁш\-ъх  $d$ ёырсхх °ш\-Ёю\-ъющ ёїю\-фш\-ьюё\-Єш т  $\frak{R}_c$,
т Єю тЁхь  ъръ эр ъюь\-яръЄ\-э√ї ьэю\-цхёЄ\-трї юср тшфр
ёїю\-фш\-ьюё\-Єш ¤ътш\-тр\-ыхэЄ\-э√. ╠хЄЁшър  $d$ юяЁх\-фхы \-хЄ\-ё 
ёў╕Є\-эющ тё■фє яыюЄ\-эющ яю\-ёыхфю\-тр\-Єхы№\-эюёЄ№■ $\varphi_n$.
╤ыхфю\-тр\-Єхы№\-эю ёє\-∙хёЄ\-тє\-хЄ схё\-ъю\-эхў\-эюх ўшёыю ьхЄЁшъ
Єръюую Єшяр. ┬ юс∙хь ёыєўрх шч ёїюфш\-ьюё\-Єш
яю\-ёыхфю\-тр\-Єхы№\-эюё\-Єш $\mu_n$ т юфэющ ьхЄЁш\-ъх эх
ёых\-фє\-хЄ ёїюфш\-ьюёЄ№ т фЁєующ ьхЄЁш\-ъх. ╬фэръю эр
ъюь\-яръЄ\-э√ї ьэю\-цх\-ёЄ\-трї ёїюфш\-ьюёЄ№
яю\-ёыхфю\-тр\-Єхы№\-эюё\-Єш $\mu_n$ т юфэющ шч ьхЄЁшъ тыхў╕Є
°ш\-Ёю\-ъє■ ёїюфш\-ьюёЄ№ ¤Єющ яю\-ёыхфю\-тр\-Єхы№\-эюё\-Єш ш,
чэрўшЄ, ёїюфш\-ьюёЄ№ $\mu_n$ т ы■сющ ьхЄЁш\-ъх
Ёрё\-ёьрЄЁш\-трх\-ьюую Єшяр.

\par ─рыхх фы  ьхЁ $\mu\in \frak{R}_c$ ь√ сєфхь Ёрё\-ёьрЄ\-Ёш\-трЄ№
шэ\-Єху\-Ёр\-ы√ $\int\limits_0^\infty f(x)d\mu(x)$ эх Єюы№ъю фы 
Їєэъ\-Ўшщ $f\in\Phi$. ┬ ёт чш ё ¤Єшь юЄ\-ьхЄшь ёых\-фє\-■\-∙хх. ┼ёыш
$f\in\Phi$, Єю эр\-яшёрэ\-э√щ шэ\-ЄхуЁры ьюц\-эю Ёрё\-ёьрЄЁш\-трЄ№
ъръ шэ\-Єху\-Ёры ╨шьр\-эр-╤Єшы\-Є№хёр Їєэъ\-Ўшш $f(x)$ яю Їєэъ\-Ўшш
$\mu(x)$, уфх $\mu(x)$  -- Їєэъ\-Ўш  Ёрё\-яЁх\-фхых\-эш  ьхЁ√ $\mu$.
╬фэръю, эрь яЁш\-ф╕Єё  шьхЄ№ фхыю ёю ёыєўрхь, ъюуфр $f$  --
яЁю\-шч\-тюы№\-эр  сюЁх\-ыхт\-ёър  Їєэъ\-Ўш  эр яюыє\-юёш
$(0,\infty)$. ┬эрўрых яЁхф\-яюыю\-цшь, ўЄю $\supp f$ $\subset$ $
[a,b]$ $\subset$ $(0,\infty)$. ╬уЁрэш\-ўх\-эшх $\mu$ эр ёху\-ьхэЄ
$[a,b]$ хёЄ№ ъю\-эхў\-эр  сюЁхыхт\-ёър  ьхЁр. ┬ ¤Єюь ёыєўрх
$\int\limits_a^b f(x)d\mu(x)$ Ёрё\-ёьрЄЁш\-тр\-хЄё  ъръ шэ\-Єху\-Ёры
╦хсхур Їєэъ\-Ўшш $f$ яю ьхЁх  $\mu$. ╩ю\-эхўэю, эх тё ър  Єрър 
Їєэъ\-Ўш   ты \-хЄ\-ё  шэ\-Єху\-Ёш\-Ёєх\-ьющ яю ьхЁх $\mu$.

\par ╧Ёхфяюыю\-цшь ЄхяхЁ№, ўЄю $f(x)$ -- яЁюшч\-тюы№\-эр  сюЁх\-ыхт\-ёър 
Їєэъ\-Ўш  эр яюыє\-юёш $(0,\infty)$, ъюЄю\-Ёр   ты \-хЄё 
ыю\-ъры№\-эю шэ\-Єху\-Ёш\-Ёєх\-ьющ яю ьхЁх  $\mu$. ┬ ¤Єюь ёыєўрх
яю\-ырур\-хь
$$
\int\limits_0^\infty f(x)d\mu(x)=\lim\limits_{{a\to +0} \atop{b\to
+\infty}} \int\limits_a^b f(x)d\mu(x).
$$
▌Єю ючэрўрхЄ, ўЄю ь√ Ёрё\-ёьрЄЁш\-трхь шэ\-ЄхуЁры
$\int\limits_0^\infty f(x)d\mu(x)$ ъръ эх\-ёюсёЄ\-тхэ\-э√щ
шэ\-Єху\-Ёры ё юёюс√\-ьш Єюўър\-ьш эюы№ ш схё\-ъю\-эхў\-эюёЄ№. ─ы 
эх\-ъюЄю\-Ё√ї Їєэъ\-Ўшщ  $f$ Єр\-ъющ шэ\-Єху\-Ёры ёїю\-фшЄ\-ё .

%%%%%%%%%%%%%%%%%%%%%%%%%%%%%%%%%%%%%%%%%%%%%%%%%%%%%%%%%%%%%%%%%%%%%%%%%%%%%%%%%%%%%%%%%
%         ╧ ╨ ┼ ─ ┼ ╦ ▄ ═ █ ┼    ╠ ═ ╬ ╞ ┼ ╤ ╥ ┬ └    ╠ ┼ ╨
%%%%%%%%%%%%%%%%%%%%%%%%%%%%%%%%%%%%%%%%%%%%%%%%%%%%%%%%%%%%%%%%%%%%%%%%%%%%%%%%%%%%%%%%%
\par ╧єёЄ№ $\rho(t)$ -- яЁю\-шч\-тюы№\-э√щ єЄюў\-э╕э\-э√щ яю\-Ё фюъ. ═р
яЁюёЄ\-Ёрэ\-ёЄтх  $\frak{R}_c$ юяЁхфх\-ышь {\it
юфэю\-ярЁр\-ьхЄЁш\-ўхё\-ъюх ёхьхщ\-ёЄтю яЁх\-юсЁр\-чю\-тр\-эшщ
└чр\-Ёш\-эр}  $A_t:\frak{R}_c\to\frak{R}_c$, $t\in(0,\infty)$,
ёюуырё\-эю ЇюЁ\-ьє\-ырь
$$
\mu_t=A_t\mu, \quad \mu_t(E)=\frac{\mu(tE)}{V(t)}.
$$
\par ┼ёыш ьхЁр  $\mu$ шьххЄ ётюхщ Їєэъ\-Ўшхщ ЁрёяЁх\-фхых\-эш   $\mu(x)$,
Єю Їєэъ\-Ўш  ЁрёяЁх\-фхых\-эш  ьхЁ√ $\mu_t$ сєфхЄ
$\frac{1}{V(t)}\mu(tx)$.

\par ╧єёЄ№ $f\in\Phi$. ╘юЁьє\-ыр чрьх\-э√ яхЁх\-ьхэ\-э√ї фр╕Є
$$
\int\limits_E f(x)d\mu_t(x)=\frac{1}{V(t)}\int\limits_{tE}
f\left(\frac{y}{t}\right)d\mu(y),
$$
ш, т ўрёЄ\-эюё\-Єш,
\begin{equation}{\label{azarin_2}}
\int\limits_0^\infty
f(x)d\mu_t(x)=\frac{1}{V(t)}\int\limits_0^\infty
f\left(\frac{y}{t}\right)d\mu(y).
\end{equation}
╚ч ЇюЁ\-ьє\-ы√  (\ref{azarin_2}) ш Єхю\-Ёх\-ь√ \ref{measures_th_4-1}
ыхуъю ёых\-фє\-хЄ, ўЄю Їєэъ\-Ўш  $\mu_t:\frak{R}_c \times(0,\infty)
\to\frak{R}_c$ эх\-яЁх\-Ё√т\-эр яю ёютю\-ъєя\-эюё\-Єш
яхЁх\-ьхэ\-э√ї, Єю хёЄ№, хёыш $t_n\to\tau$, $\mu_n\to\mu$, Єю
$(\mu_n)_{t_n}\to\mu_\tau$.

\par  ╩ырёёш\-ўхё\-ъшх {\it фшэрьш\-ўхё\-ъшх ёшёЄх\-ь√} т ьхЄЁш\-ўхё\-ъюь
яЁюёЄ\-Ёрэ\-ёЄтх  $X$ юяЁх\-фхы \-■Є\-ё   {\cite{Nemytskii}} ъръ
юфэю\-ярЁр\-ьхЄЁш\-ўхё\-ъшх ёхьхщ\-ёЄтр юЄюсЁр\-цх\-эшщ $B_t: X\to
X$, $t\in(-\infty,\infty)$, єфютыхЄ\-тю\-Ё ■\-∙шх єёыю\-тш\- ь:
\\ 1) $B_0x=x$ (эрўры№\-эюх єёыю\-тшх);
\\ 2) юЄюсЁр\-цх\-эшх $B_t: X\times(-\infty,\infty)\to X$
эх\-яЁх\-Ё√т\-эю яю ёютюъєя\-эюё\-Єш яхЁх\-ьхэ\-э√ї (єёыю\-тшх
эх\-яЁх\-Ё√т\-эюё\-Єш);
\\ 3) $B_{t_1}B_{t_2}=B_{t_1+t_2}$ (єёыю\-тшх уЁєя\-я√).

\par ┬ ёыєўрх $\rho(r)\equiv \rho$ ёшёЄх\-ьр └чр\-Ёш\-эр $A_t$
 ты \-хЄ\-ё  фшэрьш\-ўхё\-ъющ ёшёЄх\-ьющ т яЁюёЄ\-Ёрэ\-ёЄтх
$\frak{R}_c$, уфх ёїюфш\-ьюёЄ№ яю\-эшьр\-хЄ\-ё  ъръ °ш\-Ёюър 
ёїю\-фш\-ьюёЄ№, рффшЄшт\-эр  уЁєя\-яр тх∙хёЄ\-тхэ\-э√ї ўшёхы
чрьх\-э \-хЄ\-ё  эр ьєы№Єш\-яыш\-ър\-Єшт\-эє■ уЁєя\-яє ёЄЁю\-ую
яюыю\-цш\-Єхы№\-э√ї ўшёхы, яЁш\-ў╕ь эрўры№\-эюх єёыю\-тшх т√уы \-фшЄ
Єръ: $A_1\mu=\mu$, р єёыю\-тшх уЁєя\-я√ яЁш\-эшьр\-хЄ тшф
$A_{t_1}A_{t_2}\mu=A_{t_1t_2}\mu$.

\par ╤шёЄх\-ьє  $A_t\mu$ фы  яюыю\-цш\-Єхы№\-э√ї сюЁхыхт\-ёъшї ьхЁ $\mu$ т
яЁюёЄ\-Ёрэ\-ёЄтх $\mathbb{R}^m$ тт╕ы ш ё єёяхїюь яЁшьх\-э ы т ЄхюЁшш
ёєс\-урЁью\-эш\-ўхё\-ъшї Їєэъ\-Ўшщ └чр\-Ёшэ  ({\cite{Azarin_79}},
\cite{Azarin_2009}).

\par ┬ ъырёёш\-ўхёъющ ЄхюЁшш фшэрьш\-ўхё\-ъшї ёшёЄхь ьэю\-цхёЄ\-тю
$$
\{y\in X:\, y=\lim\limits_{n\to\infty}B_{t_n}x,
\lim\limits_{n\to\infty}t_n=+\infty\}
$$
эрч√тр\-хЄ\-ё  $\omega$-{\it яЁх\-фхы№\-э√ь ьэю\-цхёЄ\-тюь
ЄЁрхъ\-Єю\-Ёшш} $B_tx$.

\par ╤шёЄх\-ьє  $A_t\mu$ ь√, шчсхур  єёыюц\-эх\-эш 
т ЄхЁьш\-эю\-ыю\-ушш, Єръ\-цх сєфхь эрч√\-трЄ№ {\it
фшэр\-ьш\-ўхё\-ъющ ёшёЄх\-ьющ └чрЁш\-эр} фрцх т ёыєўрх
яЁю\-шч\-тюы№\-эю\-ую єЄюў\-э╕э\-эю\-ую яю\-Ё ф\-ър. ╠эю\-цхёЄ\-тю
$$
\{\nu\in\frak{R}_c:\,\nu=\lim\limits_{n\to\infty}A_{t_n}\mu,
\,\lim\limits_{n\to\infty}t_n=+\infty\}
$$
ь√ сєфхь эрч√\-трЄ№ {\it яЁхфхы№\-э√ь ьэю\-цхёЄ\-тюь └чр\-Ёш\-эр фы 
ьхЁ√} $\mu$ ш юсю\-чэр\-ўрЄ№, ёых\-фє  └чр\-Ёш\-эє, $Fr[\mu]$. ╧Ёш
юЄёєЄ\-ёЄ\-тшш ёт чш ьхцфє ьхЁющ  $\mu$ ш єЄюў\-э╕э\-э√ь
яю\-Ё ф\-ъюь $\rho(r)$ эхы№ч  ўЄю-ышсю ёърчрЄ№ ю ётющёЄ\-трї
ьэю\-цх\-ёЄтр $Fr[\mu]$. ╧ю¤Єю\-ьє ь√ сєфхь яЁхф\-яюыр\-урЄ№, ўЄю
т√\-яюыэ \-хЄ\-ё  ёююЄ\-эю°х\-эшх $\mu\in$ $\frak{M}_\infty$
$(\rho(r))$. ╬яЁх\-фхых\-эшх ьэю\-цх\-ёЄ\-тр $\frak{M}_\infty$
$(\rho(r))$ фрэю тю ттхфх\-эшш. ╚ч ыю\-ъры№\-эющ ъю\-эхў\-эюё\-Єш
Ёр\-фю\-эю\-тющ ьхЁ√ ёых\-фє\-хЄ, ўЄю ёююЄ\-эю∙х\-эшх $\mu\in$
$\frak{M}_\infty$ $(\rho(r))$ ¤ътш\-трыхэЄ\-эю ёююЄ\-эю°х\-эш■
$$
\mathop{\overline{\lim}}\limits_{r\to\infty}
\frac{|\mu|([r,er])}{V(r)}<\infty.
$$
┼ёыш $\mu\in$ $\frak{M}_\infty (\rho(r))$, Єю ъръ ь√ єтшфшь фры№°х,
ьэю\-цх\-ёЄ\-тю  $Fr[\mu]$ юсырфр\-хЄ ётющ\-ёЄ\-трьш,
рэрыю\-ушў\-э√ьш ётющёЄ\-трь яЁх\-фхы№\-э√ї ьэю\-цх\-ёЄт т ЄхюЁшш
фшэрьш\-ўхё\-ъшї ёшёЄхь.

\par ╚ч ёююЄ\-эю°х\-эш  $\mu\in$ $\frak{M}_\infty$ $(\rho(r))$ ыхуъю ёых\-фє\-хЄ,
ўЄю яюыю\-цш\-Єхы№\-эр  яюыє\-ЄЁрхъ\-ЄюЁш  $\mu_t$ (ьэю\-цх\-ёЄтю
$\{\mu_t:\:t\in[1,\infty)\}$)  ты \-хЄё  ъюь\-яръЄ\-э√ь
ьэю\-цх\-ёЄ\-тюь. ╤ЇюЁ\-ьєыш\-Ёєхь ¤Єю єЄ\-тхЁ\-ц\-фх\-эшх т тшфх
юЄфхы№\-эющ ыхьь√.
%%%%%%%%%%%%%%%%%%%%%%%%%%%%%%%%%%%%%%%%%%%%%%%%%%%%%%%%%%%%%%%%%%%%%%%%%%%%%%%%%%%%%%%
\begin{lemma}\label{azarin_lemma_1}\hskip-2mm{.}\:
%%%%%%%%%%%%%%%%%%%%%%%%%%%%%%%%%%%%%%%%%%%%%%%%%%%%%%%%%%%%%%%%%%%%%%%%%%%%%%%%%%%%%%%
╧єёЄ№ $\rho(r)$ -- яЁю\-шч\-тюы№\-э√щ єЄюў\-э╕э\-э√щ яю\-Ё фюъ,
$\mu\in \frak{M}_\infty$ $(\rho(r))$. ╥юуфр яюыє\-ЄЁрхъ\-Єю\-Ёш 
$\mu_t$, $t\geq 1$, хёЄ№ ъюь\-яръЄ\-эюх ьэю\-цх\-ёЄ\-тю
Ёр\-фю\-эю\-т√ї ьхЁ.
\end{lemma}
%%%%%%%%%%%%%%%%%%%%%%%%%%%%%%%%%%%%%%%%%%%%%%%%%%%%%%%%%%%%%%%%%%%%%%%%%%%%%%%%%%%%%%%
\par {\sc  ─юърчрЄхы№ёЄтю.}\: ╧єёЄ№ $0<a<b<\infty$. ╚ёяюы№чє 
ЇюЁ\-ьє\-ы√  (\ref{order_ln_gamma_1}) ш  (\ref{order_V_V}),
яю\-ыєўр\-хь юЎхэъє
$$
|\mu_t|([a,b])= \frac{|\mu|\left(\left[at,bt\right]\right)}{V(at)}
\frac{V(at)}{V(t)}\leq \gamma(a)a^\rho
\frac{|\mu|\left(\left[at,\frac{b}{a}at\right]\right)}{V(at)}\leq
M(a,b).
$$
╧юёыхф\-эхх эх\-Ёртхэ\-ёЄтю ёых\-фє\-хЄ шч ыю\-ъры№\-эющ
ъю\-эхў\-эюё\-Єш $\mu$ ш юяЁх\-фхых\-эш  ьэю\-цх\-ёЄ\-тр
$\frak{M}_\infty$ $(\rho(r))$. ╥хяхЁ№ єЄ\-тхЁ\-ц\-фх\-эшх ыхьь√
ёых\-фє\-хЄ шч Єхю\-Ёх\-ь√  \ref{measures_th_4}. ╦хььр фю\-ърчрэр.

\par ─рыхх эр°р чрфрўр ёюёЄюшЄ т юяшёр\-эшш ётющёЄт яЁх\-фхы№\-э√ї
ьэю\-цх\-ёЄт ьхЁ шч ъырёёр $\frak{M}_\infty$ $(\rho(r))$. ┬эрўрых ь√
ёЇюЁ\-ьє\-ыш\-Ёєхь Єхю\-Ёхьє, ъюЄю\-Ёр  яю\-чтюы хЄ ётхёЄш ¤Єє
чрфрўє ъ сюыхх яЁюёЄю\-ьє ёыєўр■, ъюуфр $\rho(r)$ $\equiv \rho$.
%%%%%%%%%%%%%%%%%%%%%%%%%%%%%%%%%%%%%%%%%%%%%%%%%%%%%%%%%%%%%%%%%%%%%%%%%%%%%%%%%%%%%%%
\begin{theorem}\label{azarin_th_1}\hskip-2mm{.}\:
%%%%%%%%%%%%%%%%%%%%%%%%%%%%%%%%%%%%%%%%%%%%%%%%%%%%%%%%%%%%%%%%%%%%%%%%%%%%%%%%%%%%%%%
╧єёЄ№ $\rho_1(r)$ ш $\rho_2(r)$ -- єЄюў\-э╕э\-э√х яю\-Ё ф\-ъш Єръшх,
ўЄю $\lim\rho_1(r)$ $=$ $\rho_1$, $\lim\rho_2(r)$ $=$ $\rho_2$
$(r\to\infty)$. ╧єёЄ№ $\mu\in$ $\frak{M}_\infty(\rho_1(r))$, р ьхЁр
$\lambda$ юяЁх\-фхы \-хЄё  Ёртхэ\-ёЄ\-тюь
$$
d\lambda(t)=\frac{V_2(t)}{V_1(t)}d\mu(t).
$$
╥юуфр $\lambda\in\frak{M}_\infty(\rho_2(r))$ ш ёююЄ\-эю°х\-эшх
$\mu_{t_n}\to\nu$ ¤ътш\-тр\-ыхэЄ\-эю ёююЄ\-эю°х\-эш■
$\lambda_{t_n}\to\nu_1$, уфх $d\nu_1(t)= t^{\rho_2-\rho_1} d\nu(t)$.
╟фхё№
$$
\mu_t(E)=\frac{\mu(tE)}{V_1(t)},\quad
\lambda_t(E)=\frac{\lambda(tE)}{V_2(t)}, \quad
V_1(r)=r^{\rho_1(r)},\quad V_2(r)=r^{\rho_2(r)}, \; t_n\to\infty.
$$
\end{theorem}
%%%%%%%%%%%%%%%%%%%%%%%%%%%%%%%%%%%%%%%%%%%%%%%%%%%%%%%%%%%%%%%%%%%%%%%%%%%%%%%%%%%%%%%
\par \quad {\sc ─юърчрЄхы№ёЄтю.}\: ╚ьххь
$$
|\lambda|([r,er])=\int\limits_r^{er} \frac{V_2(t)}{V_1(t)}d|\mu|(t).
$$
╚ч Єхю\-Ёх\-ь√  \ref{order_th_V:V} ёых\-фєхЄ, ўЄю
$$
\mathop{\overline{\lim}}\limits_{r\to\infty}
\frac{|\lambda|([r,er])}{V_2(r)}\leq
\left(e^{(\rho_2-\rho_1)_+}\right)
\mathop{\overline{\lim}}\limits_{r\to\infty}
\frac{|\mu|([r,er])}{V_1(r)} \qquad (a_+=\max\{a,0\}).
$$
яю¤Єю\-ьє $\lambda\in\frak{M}_\infty(\rho_2(r))$. ╧єёЄ№
$\mu_{t_n}\to\nu$, $\varphi\in\Phi$. ╚ьххь
$$
\lim\limits_{n\to\infty} \int\limits_0^\infty
\varphi(t)d\lambda_{t_n}(t)= \lim\limits_{n\to\infty}
\frac{1}{V_2(t_n)} \int\limits_0^\infty
\varphi\left(\frac{u}{t_n}\right)d\lambda(u)
$$
$$
=\lim\limits_{n\to\infty} \frac{1}{V_2(t_n)} \int\limits_0^\infty
\varphi\left(\frac{u}{t_n}\right) \frac{V_2(u)}{V_1(u)} d\mu(u)=
\lim\limits_{n\to\infty} \int\limits_0^\infty \varphi(\tau)
\frac{V_2(\tau t_n)}{V_1(\tau t_n)} \frac{V_1(t_n)}{V_2(t_n)}
d\mu_{t_n}(\tau).
$$
╥ръ ъръ яюёыхфю\-тр\-Єхы№\-эюёЄ№
$$
\varphi(\tau)\frac{V_2(\tau t_n)}{V_1(\tau t_n)}
\frac{V_1(t_n)}{V_2(t_n)}
$$
ёїю\-фшЄё  т яЁюёЄ\-Ёрэ\-ёЄтх $\Phi$ ъ Їєэъ\-Ўшш
$\tau^{\rho_2-\rho_1}\varphi(\tau)$ ш $\mu_{t_k}\to\nu$, Єю шч
Єхю\-Ёх\-ь√  \ref{measures_th_4-1} ёых\-фєхЄ, ўЄю
$$
\lim\limits_{n\to\infty} \int\limits_0^\infty
\varphi(t)d\lambda_{t_n}(t)=\int\limits_0^\infty
\tau^{\rho_2-\rho_1} \varphi(\tau)d\nu(\tau).
$$
╥хь ёрь√ь ь√ фюърчрыш, ўЄю $\lambda_{t_n}\to\nu_1$. ┬ фюърчрэ\-эюь
єЄ\-тхЁ\-ц\-фх\-эшш ьюц\-эю яю\-ьхэ Є№ ьхёЄрьш $\mu$ ш  $\nu$.
╥хю\-Ёхьр фю\-ърчрэр.

\par ╧юїюцр  Єхю\-Ёхьр фы  ьхЁ т яыюёъюёЄш хёЄ№ т  \cite{Azarin_Giner_84},
Єхю\-Ёх\-ьр  4.
%%%%%%%%%%%%%%%%%%%%%%%%%%%%%%%%%%%%%%%%%%%%%%%%%%%%%%%%%%%%%%%%%%%%%%%%%%%%%%%%%%%%%%%
\begin{remark}\label{rho1=rho2}\hskip-2mm{.}\:
%%%%%%%%%%%%%%%%%%%%%%%%%%%%%%%%%%%%%%%%%%%%%%%%%%%%%%%%%%%%%%%%%%%%%%%%%%%%%%%%%%%%%%%
┼cыш $\rho_1=\rho_2$, Єю ьхЁ√ $\nu$ ш $\nu_1$ ёют\-ярфр■Є.
\end{remark}
\par ╬ЄюсЁрцхэшх $A_t$, $A_t\mu=\mu_t$, т ёыєўрх, ъюуфр
$\rho(r)$ $\equiv \rho$ ь√ сєфхь юсю\-чэрўрЄ№ $F_t$ шыш $F_t(\rho)$.

\par ┬ ёых\-фє\-■\-∙хщ Єхю\-Ёхьх яЁш\-тюфшЄ\-ё  Ё ф ётющёЄт яЁх\-фхы№\-эю\-ую
ьэю\-цх\-ёЄ\-тр $Fr[\mu]$.
%%%%%%%%%%%%%%%%%%%%%%%%%%%%%%%%%%%%%%%%%%%%%%%%%%%%%%%%%%%%%%%%%%%%%%%%%%%%%%%%%%%%%%%
\begin{theorem}\label{azarin_th_2}\hskip-2mm{.}\:
%%%%%%%%%%%%%%%%%%%%%%%%%%%%%%%%%%%%%%%%%%%%%%%%%%%%%%%%%%%%%%%%%%%%%%%%%%%%%%%%%%%%%%%
╧єёЄ№ $\rho(r)$ -- яЁю\-шч\-тюы№\-э√щ єЄюў\-э╕э\-э√щ яю\-Ё фюъ,
$\mu\in$ $\frak{M}_\infty$ $(\rho(r))$.
╥юуфр ёяЁртхф\-ыш\-т√ ёых\-фє\-■\-∙шх єЄ\-тхЁ\-ц\-фх\-эш :\\
1) $Fr[\mu]$ хёЄ№ эх\-яєёЄющ ъюь\-яръЄ, \\
2) $Fr[\mu]$ хёЄ№ ёт чэюх ьэю\-цхёЄ\-тю т ьхЄЁш\-ўхё\-ъюь
яЁюёЄ\-Ёрэ\-ёЄ\-тх $(\frak{R}_c,d)$,\\
3) ьэю\-цхёЄ\-тю $Fr[\mu]$ шэтрЁш\-рэЄ\-эю юЄ\-эюёш\-Єхы№\-эю
яЁх\-юсЁрчю\-тр\-эш  $F_t$, сюыхх Єюую $F_t$ хёЄ№ тчршь\-эю
юфэю\-чэрў\-эюх юЄюсЁрцх\-эшх ьэю\-цхёЄ\-тр $Fr[\mu]$ эр ёхс .
\end{theorem}
%%%%%%%%%%%%%%%%%%%%%%%%%%%%%%%%%%%%%%%%%%%%%%%%%%%%%%%%%%%%%%%%%%%%%%%%%%%%%%%%%%%%%%%
\par \quad {\sc ─юърчрЄхы№ёЄтю.}\: ╩ръ ёых\-фє\-хЄ шч чр\-ьхўр\-эш  ъ
Єхю\-Ёх\-ьх  \ref{azarin_th_1}, эх юуЁрэш\-ўш\-тр  юс∙\-эюё\-Єш,
ьюц\-эю ёўш\-ЄрЄ№, ўЄю $\rho(r)$ $\equiv$ $\rho$. ╧ю ыхььх
\ref{azarin_lemma_1} яюыє\-ЄЁрхъ\-ЄюЁш  $\mu_t$, $t\geq 1$,
 ты \-хЄё  ъюь\-яръЄ\-э√ь ьэю\-цх\-ёЄ\-тюь. ═р  ч√ъх ЄхюЁшш
фшэрьш\-ўхё\-ъшї ёшёЄхь \cite{Nemytskii} ¤Єю єЄ\-тхЁ\-ц\-фх\-эшх
яхЁх\-ЇюЁ\-ьє\-ыш\-Ёє\-хЄё  Єръ. ─тшцхэшх $\mu_t$ яюыю\-цш\-Єхы№\-эю
єёЄющ\-ўштю яю ╦руЁрэцє. ═р ъюь\-яръЄ\-э√ї ьэю\-цх\-ёЄ\-трї т
$\frak{R}_c$ °ш\-Ёюър  ёїюфш\-ьюёЄ№ яю\-ёыхфю\-тр\-Єхы№\-эюё\-Єш
$\mu_n$ ¤ътш\-трыхэЄ\-эр ёїюфш\-ьюё\-Єш т ьхЄЁш\-ўхё\-ъюь
яЁюёЄ\-Ёрэ\-ёЄтх $(\frak{R}_c,d)$. ╧ю¤Єю\-ьє яЁх\-фхы№\-эюх
ьэю\-цх\-ёЄ\-тю └чрЁш\-эр $Fr[\mu]$ ёют\-ярфр\-хЄ ё
$\omega$-яЁх\-фхы№\-э√ь ьэю\-цх\-ёЄ\-тюь ЄЁрхъ\-Єю\-Ёшш $\mu_t$ т
ьхЄЁш\-ўхёъюь яЁюёЄ\-Ёрэ\-ёЄ\-тх $(\frak{R}_c,d)$. ╥хяхЁ№ ь√ ьюцхь
шё\-яюы№чю\-трЄ№ фю\-ёЄр\-Єюўэю Ёрч\-тшЄє■ ЄхюЁш■ фшэрьш\-ўхё\-ъшї
ёшёЄхь т ьхЄЁш\-ўхё\-ъшї яЁюёЄ\-Ёрэ\-ёЄ\-трї. ╥хю\-Ёх\-ьр  10, уыртр
5, \S  3 шч ъэшуш \cite{Nemytskii} єЄ\-тхЁ\-ц\-фр\-хЄ, ўЄю $Fr[\mu]$
хёЄ№ эх\-яєёЄющ ъюь\-яръЄ, ъюЄю\-Ё√щ юЄюсЁрцх\-эшх $F_t$ тчршь\-эю
юфэю\-чэрў\-эю юЄю\-сЁрцр\-хЄ эр ёхс , р Єхю\-Ёх\-ьр  14 шч ¤Єющ
ъэшуш єЄ\-тхЁ\-ц\-фр\-хЄ, ўЄю $Fr[\mu]$ хёЄ№ ёт ч\-эюх
ьэю\-цх\-ёЄ\-тю. ╥хю\-Ёхьр фю\-ърчрэр.

\par ╤ыхфє\-■∙шх Єхю\-Ёх\-ь√  ты ■Єё  яюыхч\-э√ь фю\-яюыэх\-эш\-хь
ъ яЁх\-ф√\-фє∙хщ.
%%%%%%%%%%%%%%%%%%%%%%%%%%%%%%%%%%%%%%%%%%%%%%%%%%%%%%%%%%%%%%%%%%%%%%%%%%%%%%%%%%%%%%%
\begin{theorem}\label{azarin_th_3}\hskip-2mm{.}\:
%%%%%%%%%%%%%%%%%%%%%%%%%%%%%%%%%%%%%%%%%%%%%%%%%%%%%%%%%%%%%%%%%%%%%%%%%%%%%%%%%%%%%%%
╧єёЄ№ $\mu\in$ $\frak{M}_\infty$ $(\rho(r))$,
яю\-ёыхфю\-тр\-Єхы№\-эюёЄ№ $t_n\to\infty$ Єръютр, ўЄю
$\mu_{t_n}\to\nu$, р яю\-ёыхфю\-тр\-Єхы№\-эюёЄ№ $\tau_n$
ёїю\-фшЄ\-ё  ъ  $\tau>0$. ╥юуфр яю\-ёыхфю\-тр\-Єхы№\-эюёЄ№
$\mu_{\tau_n t_n}$ ёїю\-фшЄ\-ё  ъ  $\nu_\tau$.
\end{theorem}
%%%%%%%%%%%%%%%%%%%%%%%%%%%%%%%%%%%%%%%%%%%%%%%%%%%%%%%%%%%%%%%%%%%%%%%%%%%%%%%%%%%%%%%
\par {\sc  ─юърчрЄхы№ёЄтю.}\: ╧єёЄ№ $\varphi$ -- яЁю\-шч\-тюы№\-эр 
Їєэъ\-Ўш  шч яЁюёЄ\-Ёрэ\-ёЄтр $\Phi$. ╚ьххь
$$
(\varphi,\mu_{t_n\tau_n})=\frac{1}{V(t_n\tau_n)}
\int\limits_0^\infty \varphi\left(\frac{u}{t_n\tau_n}\right)d\mu(u)=
\frac{V(t_n)}{V(t_n\tau_n)} \int\limits_0^\infty
\varphi\left(\frac{\xi}{\tau_n}\right)d\mu_{t_n}(\xi).
$$
╥ръ ъръ $\mu_{t_n}\to\nu$, р яю\-ёыхфю\-тр\-Єхы№\-эюёЄ№
$\frac{V(t_n)}{V(t_n\tau_n)} \varphi\left(\frac{\xi}{\tau_n}\right)$
ёїю\-фшЄ\-ё  т яЁюёЄ\-Ёрэ\-ёЄтх $\Phi$ ъ Їєэъ\-Ўшш
$\frac{1}{\tau^\rho}\varphi\left(\frac{\xi}{\tau}\right)$, Єю шч
Єхю\-Ёх\-ь√  \ref{measures_th_4-1} ёых\-фє\-хЄ, ўЄю
$$
\lim\limits_{n\to\infty}
(\varphi,\mu_{t_n\tau_n})=\frac{1}{\tau^\rho} \int\limits_0^\infty
\varphi\left(\frac{\xi}{\tau}\right)d\nu(\xi)=(\varphi,\nu_\tau).
$$
╥ръшь юсЁр\-чюь $\mu_{t_n\tau_n}\to\nu_\tau$, ш Єхю\-Ёхьр
фю\-ърчр\-эр.

%%%%%%%%%%%%%%%%%%%%%%%%%%%%%%%%%%%%%%%%%%%%%%%%%%%%%%%%%%%%%%%%%%%%%%%%%%%%%%%%%%%%%%%
\begin{theorem}\label{azarin_th_4}\hskip-2mm{.}\:
%%%%%%%%%%%%%%%%%%%%%%%%%%%%%%%%%%%%%%%%%%%%%%%%%%%%%%%%%%%%%%%%%%%%%%%%%%%%%%%%%%%%%%%
╧єёЄ№ $\mu\in$ $\frak{M}_\infty$ $(\rho(r))$ ш яєёЄ№ $t_n\to\infty$
Єр\-ър  яю\-ёыхфю\-тр\-Єхы№\-эюёЄ№, ўЄю $\mu_{t_n}\to\nu$,
$|\mu|_{t_n}\to \hat{\nu}$. ╥юуфр $|\nu|\leq \hat{\nu}$.
\end{theorem}
%%%%%%%%%%%%%%%%%%%%%%%%%%%%%%%%%%%%%%%%%%%%%%%%%%%%%%%%%%%%%%%%%%%%%%%%%%%%%%%%%%%%%%%
\par \quad {\sc ─юърчрЄхы№ёЄтю.}\: ╧єёЄ№ $\varphi$ -- яЁюшч\-тюы№\-эр 
Їєэъ\-Ўш  шч яЁюёЄ\-Ёрэ\-ёЄтр $\Phi$. ╥юуфр
$$
|(\nu,\varphi)|=\lim\limits_{n\to\infty}|(\mu_{t_n},\varphi)| \leq
\lim\limits_{n\to\infty} (|\mu|_{t_n},|\varphi|)=
(\hat{\nu},|\varphi|).
$$
╚ч ¤Єюую ыхуъю ёых\-фєхЄ эх\-Ёртхэ\-ёЄтю $|\nu|(E)\leq
\hat{\nu}(E)$. ╥хю\-Ёхьр фюърчрэр.

\par ┬ ёых\-фє\-■\-∙хщ Єхю\-Ёх\-ьх ь√ яю\-ърч√\-трхь, ўЄю эх\-ъюЄюЁ√х
рёшья\-Єю\-Єш\-ўхёъшх юЎхэъш фы  ьхЁ√ $\mu$ яю\-Ёюцфр\-■Є
уыюсры№\-э√х юЎхэъш фы  ьхЁ $\nu$ шч яЁх\-фхы№\-эю\-ую
ьэю\-цх\-ёЄ\-тр $Fr[\mu]$. ╧Ёхф\-трЁш\-Єхы№\-эю ттхф╕ь эх\-ъюЄю\-Ё√х
эют√х яюэ Єш .

\par ┬ ёыєўрх тх∙хёЄ\-тхэ\-э√ї Ёр\-фю\-эю\-т√ї ьхЁ эрЁ фє ё яЁх\-фхы№\-э√ь
ьэю\-цх\-ёЄ\-тюь └чрЁш\-эр $Fr[\mu]$ трцэ√ьш
рёшья\-Єю\-Єш\-ўхё\-ъшьш їрЁръЄх\-Ёшё\-Єш\-ърьш ьхЁ√ $\mu$
 ты \-■Є\-ё  х╕ тхЁї\-э   ш эшц\-э   яыюЄ\-эюё\-Єш $N(\alpha)$ ш
$\underline{N}(\alpha)$. ╘єэъ\-Ўшш яыюЄ\-эюё\-Єш $N(\alpha)$,
$\underline{N}(\alpha)$ -- ¤Єю Їєэъ\-Ўшш эр яюыє\-юёш $[0,\infty)$ ш
яю¤Єю\-ьє  ты \-■Єё  сюыхх яЁюёЄ√\-ьш ьрЄхьр\-Єш\-ўхё\-ъшьш
юс·хъ\-Єр\-ьш, ўхь ьэю\-цх\-ёЄ\-тю ьхЁ $Fr[\mu]$. ┬  ¤Єюь
яЁх\-шьє∙хёЄтю Їєэъ\-Ўшщ яыюЄ\-эюё\-Єш. ╤ фЁєующ ёЄюЁюэ√, ъръ ь√
єтшфшь т фры№\-эхщ\-°хь, ьэю\-цх\-ёЄ\-тю $Fr[\mu]$ фр╕Є, тююс∙х
уютюЁ , чэрўш\-Єхы№\-эю сюы№\-°х шэЇюЁ\-ьр\-Ўшш ю ьхЁх  $\mu$, ўхь
Їєэъ\-Ўшш $N(\alpha)$ ш  $\underline{N}(\alpha)$.

\par ╧єёЄ№ $\mu$ -- тх∙хёЄ\-тхэ\-эр  ьхЁр ╨рфюэр эр яюыє\-юёш $(0,\infty)$
ё Їєэъ\-Ўшхщ Ёрё\-яЁх\-фхых\-эш  $\mu(r)$, $\rho(r)$ --
єЄюў\-э╕э\-э√щ яю\-Ё фюъ. {\it ┬хЁїэхщ яыюЄ\-эюёЄ№■} ьхЁ√ $\mu$
юЄ\-эюёш\-Єхы№\-эю єЄюў\-э╕э\-эю\-ую яю\-Ё ф\-ър $\rho(r)$
эр\-ч√тр\-хЄ\-ё  тхыш\-ўш\-эр
\begin{equation}\label{def_N(alpha)}
N(\alpha)=\mathop{\overline{\lim}}\limits_{r\to\infty}
\frac{\mu(r+\alpha r)-\mu(r)}{V(r)}.
\end{equation}
┬ ёыєўрх $\alpha>0$ ьюц\-эю Єръ\-цх яшёрЄ№
\begin{equation}\label{def_N(alpha>0)}
N(\alpha)=\mathop{\overline{\lim}}\limits_{r\to\infty}
\frac{\mu(r,r+\alpha r]}{V(r)}.
\end{equation}
╬ЄьхЄшь, ўЄю ЇюЁ\-ьєыющ  (\ref{def_N(alpha)}) тхышўшэр $N(\alpha)$
юяЁх\-фхы \-хЄё  яЁш $\alpha>-1$, р ЇюЁ\-ьє\-ыющ
(\ref{def_N(alpha>0)}) яЁш $\alpha>0$. ┬ фры№\-эхщ\-°хь сєфхь
ёўш\-ЄрЄ№, ўЄю $N(0)=0$.

\par └эрыюушўэю юяЁх\-фхы \-хЄё  {\it эшцэ   яыюЄ\-эюёЄ№} ьхЁ√ $\mu$
юЄ\-эюёш\-Єхы№\-эю єЄюў\-э╕э\-эю\-ую яю\-Ё ф\-ър $\rho(r)$:
$$
\underline{N}(\alpha)=\mathop{\underline{\lim}}\limits_{r\to\infty}
\frac{\mu(r+\alpha r)-\mu(r)}{V(r)}.
$$
╚ч ётющёЄт яЁх\-фхыют ш єЄюў\-э╕э\-эю\-ую яю\-Ё ф\-ър $\rho(r)$
ыхуъю яюыєўр\-хЄё , ўЄю Їєэъ\-Ўшш $N(\alpha)$ ш
$\underline{N}(\alpha)$ єфют\-ыхЄтю\-Ё ■Є эх\-Ёртхэ\-ёЄ\-трь
\begin{equation}\label{N(alpha+beta)leq}
N(\alpha+\beta)\leq N(\alpha)+(1+\alpha)^\rho
N\left(\frac{\beta}{1+\alpha}\right),
\end{equation}
\begin{equation}\label{N(alpha+beta)geq}
N(\alpha+\beta)\geq N(\alpha)+(1+\alpha)^\rho
\underline{N}\left(\frac{\beta}{1+\alpha}\right),
\end{equation}
\begin{equation}\label{underline_N(alpha+beta)geq}
\underline{N}(\alpha+\beta)\geq \underline{N}(\alpha)
+(1+\alpha)^\rho \underline{N}\left(\frac{\beta}{1+\alpha}\right),
\end{equation}
\begin{equation}\label{underline_N(alpha+beta)leq}
\underline{N}(\alpha+\beta)\leq \underline{N}(\alpha)
+(1+\alpha)^\rho N\left(\frac{\beta}{1+\alpha}\right),
\end{equation}
уфх $\rho=\rho(\infty)=\lim\limits_{r\to\infty} \rho(r)$. ╟рьхЄшь,
ўЄю хёыш яЁртр  ўрёЄ№ т ъръюь-Єю шч эр\-яшёрэ\-э√ї эх\-ЁртхэёЄт
шьххЄ тшф $\infty-\infty$, Єю ¤Єю эх\-Ёртхэ\-ёЄ\-тю эєцэю ёўш\-ЄрЄ№
яєёЄ√ь єЄ\-тхЁ\-ц\-фх\-эшхь. ┬ ёыєўрх, хёыш $\mu\in \frak{M}_\infty
(\rho(r))$, Єю т√\-яюыэ \-■Єё  эх\-Ёртхэ\-ёЄ\-тр
$-\infty<\underline{N}(\alpha) \leq N(\alpha)<\infty$. ╥ръцх
чр\-ьхЄшь, ўЄю тхышўш\-э√ $\underline{N}(\alpha)$ ш $N(\alpha)$
ьюуєЄ с√Є№ ъю\-эхў\-э√\-ьш фрцх фы  ьхЁ эх тїюф \-∙шї т
ьэю\-цх\-ёЄ\-тю $\frak{M}_\infty (\rho(r))$. ▌Єю єърч√\-тр\-хЄ эр
трц\-эюёЄ№ Їєэъ\-Ўшщ $N(\alpha)$ ш $\underline{N}(\alpha)$ яЁш
шё\-ёыхфю\-трэшш ётющёЄт ьхЁ√ $\mu$. ═ряюь\-эшь, ўЄю Єхю\-Ёхьр ю
ётющ\-ёЄтрї ьэю\-цх\-ёЄтр $Fr[\mu]$ фю\-ърчр\-эр т
яЁхф\-яюыю\-цх\-эшш, ўЄю $\mu\in$ $\frak{M}_\infty$ $(\rho(r))$.
%%%%%%%%%%%%%%%%%%%%%%%%%%%%%%%%%%%%%%%%%%%%%%%%%%%%%%%%%%%%%%%%%%%%%%%%%%%%%%%%%%%%%%%
\begin{theorem}\label{azarin_th_nu_N}\hskip-2mm{.}\:
%%%%%%%%%%%%%%%%%%%%%%%%%%%%%%%%%%%%%%%%%%%%%%%%%%%%%%%%%%%%%%%%%%%%%%%%%%%%%%%%%%%%%%%
╧єёЄ№ тх∙хёЄ\-тхэ\-эр  ьхЁр $\mu\in$ $\frak{M}_\infty$ $(\rho(r))$,
$Fr[\mu]$ -- х╕ яЁх\-фхы№\-эюх ьэю\-цх\-ёЄ\-тю, $N(\alpha)$ ш
$\underline{N}(\alpha)$ -- х╕ тхЁї\-э   ш эшц\-э   яыюЄ\-эюё\-Єш.
╥юуфр фы  ы■сющ ьхЁ√ $\nu\in Fr[\mu]$ ёє\-∙хёЄ\-тє\-хЄ эх сюыхх ўхь
ёў╕Є\-эюх ьэю\-цх\-ёЄ\-тю $E(\nu)$ Єръюх, ўЄю яЁш $a,b\notin
E(\nu)$, $0<a<b<\infty$ т√\-яюыэ \-■Єё  эх\-Ёртхэ\-ёЄ\-тр
$$
\nu([a,b])\leq a^\rho N\left(\frac{b}{a}-1\right),\qquad
\nu([a,b])\geq a^\rho \underline{N}\left(\frac{b}{a}-1\right).
$$
─ы  ы■с√ї $a$ ш $b$, $0<a<b<\infty$, т√\-яюыэ \-■Є\-ё 
эх\-Ёртхэ\-ёЄ\-тр
$$
\nu([a,b])\leq a^\rho
\mathop{\underline{\lim}}\limits_{\varepsilon\to +0}
N\left(\frac{b}{a}-1+\varepsilon\right),
$$
$$
\nu([a,b])\geq a^\rho
\mathop{\overline{\lim}}\limits_{\varepsilon\to +0}
\underline{N}\left(\frac{b}{a}-1+\varepsilon\right).
$$
\end{theorem}
%%%%%%%%%%%%%%%%%%%%%%%%%%%%%%%%%%%%%%%%%%%%%%%%%%%%%%%%%%%%%%%%%%%%%%%%%%%%%%%%%%%%%%%
\par\quad {\sc  ─юърчрЄхы№ёЄтю.}\: ╧єёЄ№
$\nu= \lim\limits_{n\to\infty} \mu_{t_n}$. ─юяюыэш\-Єхы№\-эю сєфхь
ёўш\-ЄрЄ№, ўЄю ёє\-∙хёЄ\-тє\-хЄ яЁх\-фхы
$\hat{\nu}=\lim\limits_{n\to\infty} |\mu|_{t_n}$. ╬сючэрўшь
$E(\nu)=\{x\in(0,\infty):\: \hat{\nu}(\{x\})>0\}$. ╠эю\-цх\-ёЄтю
$E(\nu)$ эх сюыхх ўхь ёў╕Єэю. ╧єёЄ№ ЄхяхЁ№ $[a,b]$ $\subset$
$(0,\infty)$, $a,b\notin E(\nu)$. ╚ч Єхю\-Ёх\-ь√ \ref{measures_th_2}
ёых\-фє\-хЄ, ўЄю
$$
\nu([a,b])=\lim\limits_{n\to\infty} \mu_{t_n}([a,b]).
$$
╚ч ЁртхэёЄтр $\nu(\{a\})=0$ ш Єющ цх Єхю\-Ёх\-ь√ ёых\-фєхЄ, ўЄю
$\lim\limits_{n\to\infty}\mu_{t_n}(\{a\})=0$. ╧ю¤Єю\-ьє
$$
\nu([a,b])=\lim\limits_{n\to\infty} \mu_{t_n}((a,b]) \leq
\mathop{\overline{\lim}} \limits_{r\to\infty} \mu_r((a,b])= a^\rho
N\left(\frac{b}{a}-1\right).
$$
└эрыюушўэю фю\-ърч√\-трхЄё  эх\-Ёртхэ\-ёЄтю $\nu([a,b])\geq a^\rho
\underline{N}\left(\frac{b}{a}-1\right)$ яЁш $a,b\notin E(\nu)$.
╨рё\-ёьюЄ\-Ёшь юс∙шщ ёыєўрщ. ╧єёЄ№ $a_k\to a$, $a_k<a$, $b_k\to b$,
$b_k>b$, яЁш\-ў╕ь $a_k, b_k\notin E(\nu)$ ш т√\-яюыэ \-хЄё 
Ёртхэ\-ёЄтю
$$
\lim\limits_{k\to\infty}
N\left(\frac{b_k}{a_k}-1\right)=\mathop{\underline{\lim}}
\limits_{\varepsilon\to +0} N\left(\frac{b}{a}-1+\varepsilon\right).
$$
╚ьххь
$$
\nu([a,b])=\lim\limits_{k\to\infty} \nu([a_k,b_k])\leq
\lim\limits_{k\to\infty} a_k^\rho N\left(\frac{b_k}{a_k}-1\right)
=a^\rho \mathop{\underline{\lim}} \limits_{\varepsilon\to +0}
N\left(\frac{b}{a}-1+\varepsilon\right).
$$
└эрыюушўэю яю\-ыєўр\-хь эх\-ЁртхэёЄтю
$$
\nu([a,b])\geq a^\rho \mathop{\overline{\lim}}
\limits_{\varepsilon\to +0}
\underline{N}\left(\frac{b}{a}-1+\varepsilon\right).
$$
╥хю\-Ёхьр фюърчрэр.

\par ╬ЄьхЄшь Єръюх ёыхфёЄ\-тшх фю\-ърчрэ\-эющ Єхю\-Ёх\-ь√.
%%%%%%%%%%%%%%%%%%%%%%%%%%%%%%%%%%%%%%%%%%%%%%%%%%%%%%%%%%%%%%%%%%%%%%%%%%%%%%%%%%%%%%%
\begin{remark}\label{azarin_remark}\hskip-2mm{.}\:
%%%%%%%%%%%%%%%%%%%%%%%%%%%%%%%%%%%%%%%%%%%%%%%%%%%%%%%%%%%%%%%%%%%%%%%%%%%%%%%%%%%%%%%
╧єёЄ№ $\mu\in$ $\frak{M}_\infty$ $(\rho(r))$ ш яєёЄ№
$N(\alpha)\equiv \underline{N}(\alpha)\equiv 0$. ╥юуфр
$Fr[\mu]=\{0\}$.
\end{remark}
%%%%%%%%%%%%%%%%%%%%%%%%%%%%%%%%%%%%%%%%%%%%%%%%%%%%%%%%%%%%%%%%%%%%%%%%%%%%%%%%%%%%%%%

\par ╠юцэю юЄьхЄшЄ№, ўЄю ёююЄ\-эю°х\-эшх
$\underline{N}(\alpha)\equiv N(\alpha)\equiv 0$ ш
$$
\mathop{\underline{\lim}} \limits_{r\to \infty}
\frac{|\mu|([ar,br])}{V(r)}>0
$$
эх яЁюЄш\-тюЁхўрЄ фЁєу фЁєує.

\par  ╠хЁр $\nu$ эр\-ч√тр\-хЄё  {\it эх\-яЁх\-Ё√т\-эющ},
хёыш фы  ы■сю\-ую  $x\;$ $\nu(\{x\})=0$. ╚эюуфр трцэю, ўЄюс√ ьхЁ√ шч
яЁх\-фхы№\-эю\-ую ьэю\-цх\-ёЄтр $Fr[\mu]$ с√ыш эх\-яЁх\-Ё√т\-э√ьш.
╤яЁр\-тхф\-ышт ёых\-фє\-■\-∙шщ яЁш\-чэръ эх\-яЁх\-Ё√т\-эюё\-Єш.
%%%%%%%%%%%%%%%%%%%%%%%%%%%%%%%%%%%%%%%%%%%%%%%%%%%%%%%%%%%%%%%%%%%%%%%%%%%%%%%%%%%%%%%
\begin{theorem}\label{azarin_th_5}\hskip-2mm{.}\:
%%%%%%%%%%%%%%%%%%%%%%%%%%%%%%%%%%%%%%%%%%%%%%%%%%%%%%%%%%%%%%%%%%%%%%%%%%%%%%%%%%%%%%%
╧єёЄ№ $\mu\in$ $\frak{M}_\infty$ $(\rho(r))$, р $N(\alpha)$ ш
$\underline{N}(\alpha)$ -- х╕ Їєэъ\-Ўшш яыюЄ\-эюё\-Єш
юЄ\-эюёш\-Єхы№\-эю єЄюў\-э╕э\-эю\-ую яю\-Ё ф\-ър $\rho(r)$. ┼ёыш
Їєэъ\-Ўшш $N(\alpha)$ ш  $\underline{N}(\alpha)$ эх\-яЁх\-Ё√т\-э√ эр
яюыє\-юёш $[0,\infty)$, Єю ы■ср  ьхЁр  $\nu$ шч ьэю\-цх\-ёЄ\-тр
$Fr[\mu]$ эх\-яЁх\-Ё√т\-эр.
\end{theorem}
%%%%%%%%%%%%%%%%%%%%%%%%%%%%%%%%%%%%%%%%%%%%%%%%%%%%%%%%%%%%%%%%%%%%%%%%%%%%%%%%%%%%%%%
\par \quad {\sc ─юърчрЄхы№ёЄтю.}\: ╧єёЄ№ $\nu\in Fr[\mu]$,
$x\in(0,\infty)$, $a_n$ -- ёЄЁюую тюч\-ЁрёЄр■\-∙р , $b_n$ -- ёЄЁюую
єс√тр■\-∙р  яю\-ёыхфю\-тр\-Єхы№\-эюё\-Єш, ёїюф \-∙шх\-ё  ъ  $x$,
яЁш\-ў╕ь $a_n, b_n\notin E(\nu)$. ╧Ёш\-ьхэх\-эшх Єхю\-Ёх\-ь√
\ref{azarin_th_nu_N} фр╕Є
$$
\nu(\{x\})=\lim\limits_{n\to\infty} \nu([a_n,b_n])\leq
\lim\limits_{n\to\infty} a_n^\rho N\left(\frac{b_n}{a_n}-1\right)=0,
$$
$$
\nu(\{x\})=\lim\limits_{n\to\infty} \nu([a_n,b_n])\geq
\lim\limits_{n\to\infty} a_n^\rho
\underline{N}\left(\frac{b_n}{a_n}-1\right)=0.
$$
╥хь ёрь√ь Єхю\-Ёхьр фюърчрэр.

\par ─ы  яюыю\-цш\-Єхы№\-э√ї ьхЁ $\mu$ шьххЄ ьхёЄю ёых\-фє\-■\-∙шщ ъЁшЄх\-Ёшщ
эх\-яЁх\-Ё√т\-эюё\-Єш ьхЁ шч ьэю\-цх\-ёЄ\-тр $Fr[\mu]$.
%%%%%%%%%%%%%%%%%%%%%%%%%%%%%%%%%%%%%%%%%%%%%%%%%%%%%%%%%%%%%%%%%%%%%%%%%%%%%%%%%%%%%%%
\begin{theorem}\label{azarin_th_6}\hskip-2mm{.}\:
%%%%%%%%%%%%%%%%%%%%%%%%%%%%%%%%%%%%%%%%%%%%%%%%%%%%%%%%%%%%%%%%%%%%%%%%%%%%%%%%%%%%%%%
╧єёЄ№ $\mu$ -- яюыю\-цш\-Єхы№\-эр  ьхЁр шч ьэю\-цх\-ёЄтр $\mu\in$
$\frak{M}_\infty$ $(\rho(r))$, $N(\alpha)$ -- тхЁї\-э   яыюЄ\-эюёЄ№
ьхЁ√ $\mu$ юЄ\-эюёш\-Єхы№\-эю єЄюў\-э╕э\-эю\-ую яю\-Ё ф\-ър
$\rho(r)$. ─ы  Єюую, ўЄюс√ тёх ьхЁ√ шч яЁх\-фхы№\-эю\-ую
ьэю\-цх\-ёЄ\-тр $Fr[\mu]$ с√ыш эх\-яЁх\-Ё√т\-э√ьш, эх\-юс\-їюфш\-ью
ш фю\-ёЄр\-Єюўэю, ўЄюс√ $N(\alpha)\to 0$ $(\alpha\to +0)$.
\end{theorem}
%%%%%%%%%%%%%%%%%%%%%%%%%%%%%%%%%%%%%%%%%%%%%%%%%%%%%%%%%%%%%%%%%%%%%%%%%%%%%%%%%%%%%%%
\par \quad {\sc ─юърчрЄхы№ёЄтю.}\: ╧єёЄ№ $N(\alpha)\to 0$ яЁш $\alpha\to
+0$. ╥юуфр ш $\underline{N}(\alpha)\to 0$ яЁш $\alpha\to +0$. ┬ ¤Єюь
ёыєўрх Їєэъ\-Ўшш $N(\alpha)$ ш $\underline{N}(\alpha)$ сєфєЄ
эх\-яЁх\-Ё√т\-э√ьш. ╧ю яЁх\-ф√фє\-∙хщ Єхю\-Ёхьх ы■ср  ьхЁр шч
ьэю\-цх\-ёЄ\-тр $Fr[\mu]$ сєфхЄ эх\-яЁх\-Ё√т\-эющ.

\par ─юърцхь ЄхяхЁ№ эхюс\-їюфш\-ьюёЄ№ єёыю\-тш  $N(\alpha)\to 0$ $(\alpha\to
+0)$. ╧Ёхф\-яюыю\-цшь яЁюЄшт\-эюх. ╥юуфр $N(\alpha)\to 2a$
$(\alpha\to +0)$, уфх $a>0$. ┬  ¤Єюь ёыєўрх ёє\-∙хёЄ\-тє■Є
яю\-ёыхфю\-тр\-Єхы№\-эюё\-Єш $r_n\to\infty$,
$\varepsilon_n\downarrow 0$ Єръшх, ўЄю $\mu((r_n,$
$(1+\varepsilon_n)r_n])>aV(r_n)$. ─юяюыэш\-Єхы№\-эю ьюц\-эю
ёўш\-ЄрЄ№, ўЄю $\mu_{r_n}\to\nu$. ╧єёЄ№ $0<a_k<1<b_k$, $a_k\to 1$,
$b_k\to 1$, $\nu(\{a_k\})=0$, $\nu(\{b_k\})=0$. ╥юуфр, шёяюы№\-чє 
Єхю\-Ёх\-ьє \ref{measures_th_2}, эр\-їюфшь, ўЄю
$$
\nu(\{1\})=\lim\limits_{k\to\infty} \nu([a_k,b_k])=
\lim\limits_{k\to\infty} \lim\limits_{n\to\infty}
\mu_{r_n}((a_k,b_k])= \lim\limits_{k\to\infty}
\lim\limits_{n\to\infty} \frac{\mu((a_kr_n,b_kr_n])}{V(r_n)}.
$$
╧Ёш ЇшъёшЁю\-трээюь $k$ ш фюёЄр\-Єюў\-эю сюы№°шї $n$ сєфхЄ
т√\-яюыэ Є№\-ё  ёююЄ\-эю°х\-эшх $(r_n,$ $(1+\varepsilon_n)r_n]$
$\subset$ $(a_kr_n,b_kr_n)$. ╧ю¤Єю\-ьє
$$
\lim\limits_{n\to\infty} \frac{\mu((a_kr_n,b_kr_n])}{V(r_n)}\geq
\mathop{\overline{\lim}}\limits_{n\to\infty}
\frac{\mu((r_n,(1+\varepsilon_n)r_n])}{V(r_n)}\geq a.
$$
╚ч ¤Єюую ёых\-фє\-хЄ, ўЄю $\nu(\{1\})\geq a$. ▌Єю яЁю\-Єштю\-ЁхўшЄ
Єюьє, ўЄю ьхЁ√ шч ьэю\-цх\-ёЄ\-тр $Fr[\mu]$ эх\-яЁх\-Ё√т\-э√.
╥хю\-Ёхьр фю\-ърчр\-эр.

\par ╬сючэрўшь ўхЁхч $M(\rho,\sigma)$ ьэю\-цх\-ёЄтю тх∙хёЄ\-тхэ\-э√ї
шыш ъюья\-ыхъё\-э√ї Ёрфю\-эю\-т√ї ьхЁ, ёюёЄю \-∙хх шч тёхї ьхЁ
$\mu$, ъюЄю\-Ё√х єфют\-ыхЄ\-тю\-Ё ■Є эх\-Ёртхэ\-ёЄ\-трь:
$$
|\mu|((0,r])\leq \sigma r^\rho \quad \text{яЁш} \quad 0<r<\infty,
\quad \rho>0,
$$
$$
|\mu|([a,b])\leq \sigma \ln \frac{b}{a},\quad 0<a<b<\infty,\;
\rho=0,
$$
$$
|\mu|((r,\infty))\leq \sigma r^\rho \quad \text{яЁш} \quad
0<r<\infty, \quad \rho<0.
$$
%%%%%%%%%%%%%%%%%%%%%%%%%%%%%%%%%%%%%%%%%%%%%%%%%%%%%%%%%%%%%%%%%%%%%%%%%%%%%%%%%%%%%%%
\begin{theorem}\label{azarin_th_7}\hskip-2mm{.}\:
%%%%%%%%%%%%%%%%%%%%%%%%%%%%%%%%%%%%%%%%%%%%%%%%%%%%%%%%%%%%%%%%%%%%%%%%%%%%%%%%%%%%%%%
╧єёЄ№ $\mu\in$ $\frak{M}_\infty$ $(\rho(r))$. ╥юуфр ёє\-∙хёЄ\-тєхЄ
Єръюх $\sigma>0$, ўЄю $Fr[\mu]$ $\subset$ $M(\rho,\sigma)$.
\end{theorem}
%%%%%%%%%%%%%%%%%%%%%%%%%%%%%%%%%%%%%%%%%%%%%%%%%%%%%%%%%%%%%%%%%%%%%%%%%%%%%%%%%%%%%%%
\par\quad {\sc  ─юърчрЄхы№ёЄтю.}\: ╬сючэрўшь ўхЁхч $N_1(\alpha)$
тхЁї\-э■■ яыюЄ\-эюёЄ№ ьхЁ√ $|\mu|$ юЄ\-эюёш\-Єхы№\-эю
єЄюў\-э╕э\-эю\-ую яю\-Ё ф\-ър $\rho(r)$. ╥ръ ъръ $\mu\in$
$\frak{M}_\infty$ $(\rho(r))$, Єю $N_1(\alpha)<\infty$. ╧єёЄ№ $q>1$
-- яЁю\-шч\-тюы№\-эюх ўшёыю. ╧Ёш\-ьхэ   Єхю\-Ёх\-ь√
\ref{azarin_th_4}, \ref{azarin_th_nu_N}, эрїюфшь, ўЄю
$$
|\nu|((0,r])\leq \hat{\nu}((0,r])=\sum\limits_{n=0}^\infty
\hat{\nu}\left(\left(\frac{r}{q^{n+1}},\frac{r}{q^n}\right]\right)
$$
$$
\leq r^\rho N_1(q-1+0) \sum\limits_{n=0}^\infty
\frac{1}{q^{(n+1)\rho}}=\frac{N_1(q-1+0)}{q^\rho-1}r^\rho,
$$
хёыш $\rho>0$. └эрыюушў\-эю яЁш $\rho<0$ яюыєўр\-хь, ўЄю
$$
|\nu|((r,\infty))\leq \hat{\nu}((r,\infty))=
\sum\limits_{n=0}^\infty \hat{\nu}((q^nr,q^{n+1}r])
$$
$$
\leq r^\rho N_1(q-1+0) \sum\limits_{n=0}^\infty
q^{n\rho}=\frac{N_1(q-1+0)}{1-q^\rho}r^\rho.
$$
╤ыєўрщ, ъюуфр $\rho=0$, Ёрё\-ёьрЄЁш\-тр\-хЄё  рэрыюушў\-эю. ╥хь
ёрь√ь, Єхю\-Ёхьр фю\-ърчр\-эр.

%%%%%%%%%%%%%%%%%%%%%%%%%%%%%%%%%%%%%%%%%%%%%%%%%%%%%%%%%%%%%%%%%%%%%%%%%%
\begin{remark}\label{remark_th_M}\hskip-2mm{.}\:
%%%%%%%%%%%%%%%%%%%%%%%%%%%%%%%%%%%%%%%%%%%%%%%%%%%%%%%%%%%%%%%%%%%%%%%%%%
╚ч яЁш\-тхфхэ\-эю\-ую фю\-ърчр\-Єхы№\-ёЄтр ёых\-фєхЄ, ўЄю т
ёююЄ\-эю°х\-эшш $Fr[\mu]$ $\subset$ $M(\rho,\sigma)$ т ърўхёЄ\-тх
$\sigma$ ьюц\-эю тч Є№ ўшёыю
$$
\hat{\sigma}=\inf\limits_{q>1} \frac{N_1(q-1+0)}{|q^\rho-1|}.
$$
\end{remark}
%%%%%%%%%%%%%%%%%%%%%%%%%%%%%%%%%%%%%%%%%%%%%%%%%%%%%%%%%%%%%%%%%%%%%%%%%%
╬ЄьхЄшь, ўЄю, ъръ ёых\-фєхЄ шч Єхю\-Ёхь  6 ш  14 шч
\cite{Grishin_2000}, т ёыєўрх $\rho>0$ т√\-яюыэ \-хЄ\-ё 
Ёртхэ\-ёЄ\-тю
$$
\hat{\sigma}=\lim\limits_{q\to\infty} \frac{N_1(q-1+0)}{q^\rho-1}
=\mathop{\overline{\lim}}\limits_{r\to\infty}
\frac{|\mu|((1,r])}{V(r)}.
$$
╤ыєўрщ $\rho<0$ Єръ\-цх Ёрё\-ёьрЄЁш\-тр\-хЄ\-ё  т Ёрсю\-Єх
\cite{Grishin_2000}.

\par ╨рё\-ёьюЄЁшь яЁш\-ьхЁ√ эр т√\-ўшёых\-эшх яЁх\-фхы№\-эю\-ую
ьэю\-цх\-ёЄ\-тр $Fr[\mu]$.
%%%%%%%%%%%%%%%%%%%%%%%%%%%%%%%%%%%%%%%%%%%%%%%%%%%%%%%%%%%%%%%%%%%%%%%%%%
\begin{lemma}\label{azarin_example_1}\hskip-2mm{.}\:
%%%%%%%%%%%%%%%%%%%%%%%%%%%%%%%%%%%%%%%%%%%%%%%%%%%%%%%%%%%%%%%%%%%%%%%%%%
╧єёЄ№ $\rho(r)$ -- єЄюў\-э╕э\-э√щ яю\-Ё фюъ, $\mu$  -- ьхЁр эр
яюыє\-юёш $(0,\infty)$ ё яыюЄ\-эюёЄ№■ $\frac{V(x)}{x}$. ╥юуфр
яЁх\-фхы№\-эюх ьэю\-цх\-ёЄ\-тю $Fr[\mu]$ ьхЁ√ $\mu$, яю\-ёЄЁюхэ\-эюх
яю єЄюў\-э╕э\-эюьє яю\-Ё ф\-ъє $\rho(r)$, ёюёЄю\-шЄ шч юфэющ ьхЁ√
$\nu$, $d\nu(x)=x^{\rho-1}dx$.
\end{lemma}
%%%%%%%%%%%%%%%%%%%%%%%%%%%%%%%%%%%%%%%%%%%%%%%%%%%%%%%%%%%%%%%%%%%%%%%%%%
\par \quad {\sc ─юърчрЄхы№ёЄтю.}\: ╧єёЄ№ $\varphi\in\Phi$. ╚ё\-яюы№\-чє 
Єхю\-Ёх\-ьє  \ref{order_th_V:V}, яю\-ыєўр\-хь
$$
(\varphi,\mu_t)=\int\limits_0^\infty \varphi(x)d\mu_t(x)=
\frac{1}{V(t)} \int\limits_0^\infty \varphi(x)\frac{V(xt)}{x}dx,
$$
$$
\lim\limits_{t\to\infty}(\varphi,\mu_t)=\int\limits_0^\infty
\varphi(x)x^{\rho-1}dx=(\varphi,\nu).
$$
╦хььр фюърчрэр.

%%%%%%%%%%%%%%%%%%%%%%%%%%%%%%%%%%%%%%%%%%%%%%%%%%%%%%%%%%%%%%%%%%%%%%%%%%
\begin{lemma}\label{azarin_example}\hskip-2mm{.}\:
%%%%%%%%%%%%%%%%%%%%%%%%%%%%%%%%%%%%%%%%%%%%%%%%%%%%%%%%%%%%%%%%%%%%%%%%%%
╧єёЄ№ $\rho(r)$ -- єЄюў\-э╕э\-э√щ яю\-Ё фюъ, $\mu$ -- ьхЁр эр
яюыє\-юёш $(0,\infty)$ ё яыюЄ\-эюёЄ№■
$x^{i\lambda_0}\frac{V(x)}{x}$. ╥юуфр
$$
Fr[\mu]=\{\nu: \; d\nu(u)=e^{i\lambda_0c}u^{i\lambda_0+\rho-1}du:
\quad c\in(-\infty,\infty)\}.
$$
\end{lemma}
%%%%%%%%%%%%%%%%%%%%%%%%%%%%%%%%%%%%%%%%%%%%%%%%%%%%%%%%%%%%%%%%%%%%%%%%%%
\par \quad {\sc ─юърчрЄхы№ёЄтю.}\: ╧єёЄ№ $\varphi\in \Phi$. ╚ьххь
$$
(\varphi,\mu_r)=\frac{1}{V(r)} \int\limits_0^\infty
\varphi\left(\frac{x}{r}\right)
x^{i\lambda_0}\frac{V(x)}{x}dx=r^{i\lambda_0}\int\limits_0^\infty
\varphi(u)u^{i\lambda_0} \frac{V(ur)}{V(r)}\frac{1}{u}du.
$$
╚ч ¤Єюую Ёртхэ\-ёЄ\-тр ш Єхю\-Ёх\-ь√  \ref{order_th_V:V} ыхуъю
ёых\-фє\-хЄ єЄ\-тхЁ\-ц\-фх\-эшх ыхьь√.

\par ╠хЁр $\mu$ эр\-ч√тр\-хЄ\-ё  {\it яхЁш\-ю\-фш\-ўхё\-ъющ} ьхЁющ
яю\-Ё ф\-ър  $\rho$ ё яхЁш\-ю\-фюь $T>1$, хёыш фы  ¤Єюую  $T$ ш
ы■сю\-ую сюЁх\-ыхт\-ёъю\-ую ьэю\-цх\-ёЄ\-тр  $E$ т√\-яюы\-э \-хЄ\-ё 
Ёртхэ\-ёЄ\-тю
\begin{equation}\label{for_example_1}
\mu(TE)=T^\rho \mu(E).
\end{equation}

%%%%%%%%%%%%%%%%%%%%%%%%%%%%%%%%%%%%%%%%%%%%%%%%%%%%%%%%%%%%%%%%%%%%%%%%%%
\begin{lemma}\label{azarin_example_2}\hskip-2mm{.}\:
%%%%%%%%%%%%%%%%%%%%%%%%%%%%%%%%%%%%%%%%%%%%%%%%%%%%%%%%%%%%%%%%%%%%%%%%%%
╧єёЄ№ $\mu$ -- ыю\-ъры№\-эю ъю\-эхў\-эр  яхЁш\-ю\-фш\-ўхёър  ьхЁр
яю\-Ё ф\-ър  $\rho$ ё яхЁш\-ю\-фюь $T>1$. ╥юуфр $\mu\in$
$\frak{M}_\infty(\rho)$, $Fr[\mu]=\{\mu_t:\: 1\leq t<T\}$.
\end{lemma}
%%%%%%%%%%%%%%%%%%%%%%%%%%%%%%%%%%%%%%%%%%%%%%%%%%%%%%%%%%%%%%%%%%%%%%%%%%
\par\quad {\sc  ─юърчрЄхы№ёЄтю.}\: ╤ююЄ\-эю°хэшх $\mu\in$ $\frak{M}_\infty(\rho)$
юўхтшф\-эю. ╧єёЄ№ $t\in[1,T)$, $t_n=tT^n$. ╚ч Ёртхэ\-ёЄ\-тр
(\ref{for_example_1}) ёых\-фє\-хЄ, ўЄю $\mu_{t_n}=\mu_t$. ╧ю¤Єю\-ьє
$\mu_{t_n}\to\mu_t$. ╥хь ёрь√ь, фю\-ърчр\-эю ёююЄ\-эю°х\-эшх
$\{\mu_t:\: t\in[1,T)\}$ $\subset$ $Fr[\mu]$.
\par ╧єёЄ№ ЄхяхЁ№ $t_n$ Єр\-ър  яю\-ёыхфю\-тр\-Єхы№\-эюёЄ№, ўЄю
$\mu_{t_n}\to\nu$. ╤є∙хёЄ\-тє\-хЄ Єръюх Ўхыюх ўшёыю $m(n)$, ўЄю
$t_n=\tau_n T^{m(n)}$, уфх $\tau_n\in[1,T)$. ╚ч Ёртхэ\-ёЄ\-тр
(\ref{for_example_1}) ёых\-фє\-хЄ, ўЄю $\mu_{t_n}=\mu_{\tau_n}$. ╧ю
Єхю\-Ёхьх ┴юы№Ўр\-эю-┬хщхЁ\-°Є\-Ёрёёр ёє\-∙хёЄ\-тєхЄ ёїюф \-∙р ё 
яюф\-яю\-ёыхфю\-тр\-Єхы№\-эюёЄ№ $\tau_{n_k}$,
$\lim\limits_{k\to\infty}\tau_{n_k}=\tau$, яЁш\-ў╕ь $\tau\in[1,T]$.
╚ьххь $\mu_{\tau_{n_k}}\to\mu_\tau$. ╥юуфр
$\nu=\lim\limits_{k\to\infty} \mu_{t_{n_k}}
=\lim\limits_{k\to\infty} \mu_{\tau_{n_k}}=\mu_\tau$. ┬ьхёЄх ё
Ёртхэ\-ёЄ\-тр\-ьш $\mu=\mu_1=\mu_T$ ¤Єю фр╕Є ёююЄ\-эю°х\-эшх $
Fr[\mu]$ $\subset$ $\{\mu_t:\: t\in[1,T)\}$. ╦хььр фю\-ърчр\-эр.

%%%%%%%%%%%%%%%%%%%%%%%%%%%%%%%%%%%%%%%%%%%%%%%%%%%%%%%%%%%%%%%%%%%%%%%%%%
\begin{remark}\label{remark_ex_2}\hskip-2mm{.}\:
%%%%%%%%%%%%%%%%%%%%%%%%%%%%%%%%%%%%%%%%%%%%%%%%%%%%%%%%%%%%%%%%%%%%%%%%%%
╧єёЄ№ $\alpha$ -- яЁюшч\-тюы№\-эюх ёЄЁюую яюыю\-цш\-Єхы№\-эюх ўшёыю,
р $n$ -- яЁюшч\-тюы№\-эюх Ўхыюх ўшёыю. ╧Ёш т√\-яюыэх\-эшш єёыю\-тшщ
ыхьь√  \ref{azarin_example_2} cяЁртхф\-ыштю ёююЄ\-эю°х\-эшх
$$
Fr[\mu]=\left\{\mu_t:\: t \in \left[\alpha T^n,\alpha
T^{n+1}\right)\right\}.
$$
╤ыєўрщ $\alpha=1$, $n=0$ ёююЄтхЄ\-ёЄ\-тєхЄ єЄ\-тхЁ\-ц\-фх\-эш■ ыхьь√
\ref{azarin_example_2}.
\end{remark}
%%%%%%%%%%%%%%%%%%%%%%%%%%%%%%%%%%%%%%%%%%%%%%%%%%%%%%%%%%%%%%%%%%%%%%%%%%
\begin{lemma}\label{azarin_example_3}\hskip-2mm{.}\:
%%%%%%%%%%%%%%%%%%%%%%%%%%%%%%%%%%%%%%%%%%%%%%%%%%%%%%%%%%%%%%%%%%%%%%%%%%
╧єёЄ№ $\rho>0$, р $R_n$ -- ёЄЁюую тюч\-ЁрёЄр■\-∙р 
яю\-ёыхфю\-тр\-Єхы№\-эюёЄ№ Єр\-ър , ўЄю $R_1\geq 1$,
$\lim\limits_{n\to\infty}\frac{R_n}{R_{n-1}}=\infty$. ╧єёЄ№ ьхЁр
$\mu$ чрфр╕Є\-ё  ЇюЁ\-ьє\-ыющ
$$
\mu=\sum\limits_{n=1}^\infty R^\rho_n\delta(x-R_n).
$$
╥юуфр яЁх\-фхы№\-эюх ьэю\-цх\-ёЄтю $Fr[\mu]$ ьхЁ√ $\mu$
юЄ\-эюёш\-Єхы№\-эю єЄюў\-э╕э\-эю\-ую яю\-Ё ф\-ър $\rho(r)$ $\equiv$
$\rho$ шьххЄ тшф
\begin{equation}\label{for_example_3}
Fr[\mu]=\left\{ t^\rho \delta(x-t):\: t\in(0,\infty) \right\} \cup
\{0\}.
\end{equation}
\end{lemma}
%%%%%%%%%%%%%%%%%%%%%%%%%%%%%%%%%%%%%%%%%%%%%%%%%%%%%%%%%%%%%%%%%%%%%%%%%%
\par \quad {\sc ─юърчрЄхы№ёЄтю.}\: ╚ч єёыю\-тш  $R_{n-1}=o(R_n)$ ё яюью∙№■
Єхю\-Ёх\-ь√ ╪Єюы№Ўр ыхуъю яю\-ыєўшЄ№, ўЄю
$$
\sum\limits_{k=1}^n R_k^\rho \sim R^\rho_n\qquad (n\to\infty).
$$
╧єёЄ№ $r\geq R_1$ ш яЁюшч\-тюы№\-эюх, р $n$  -- эрш\-сюы№\-°хх Ўхыюх
ўшёыю, фы  ъюЄю\-Ёюую т√\-яюыэ \-хЄ\-ё  эх\-Ёртхэ\-ёЄтю $R_n\leq r$.
╥юуфр
$$
\mu((0,r])=\sum\limits_{k=1}^n R^\rho_k= (1+o(1))R_n^\rho \leq
(1+o(1))r^\rho \quad (r\to\infty).
$$
╚ч ¤Єюую ёых\-фєхЄ, ўЄю $\mu\in$ $\frak{M}_\infty(\rho)$.

\par ╧єёЄ№ ЄхяхЁ№ $t>0$ ш яЁюшч\-тюы№эю, р $t_n=\frac{1}{t}R_n$.
╧єёЄ№ $\varphi$  -- яЁю\-шч\-тюы№\-эр  Їєэъ\-Ўш  шч яЁюёЄ\-Ёрэ\-ёЄтр
$\Phi$. ╥юуфр фы  тёхї фю\-ёЄрЄюў\-эю сюы№\-°шї  $n$
т√\-яюыэ \-хЄ\-ё  Ёртхэ\-ёЄтю
$$
(\mu_{t_n},\varphi)=\frac{t^\rho}{R_n^\rho} \int\limits_0^\infty
\varphi\left(\frac{ut}{R_n}\right)d\mu(u)=t^\rho\varphi(t)=(t^\rho
\delta(x-t),\varphi(x)).
$$
╬Єё■фр ёых\-фєхЄ, ўЄю $\mu_{t_n}\to t^\rho \delta(x-t)$. ╥юўэю Єръ
цх фю\-ърч√\-тр\-хЄё , ўЄю хёыш $t_n=\frac{1}{\tau_n}R_n$, уфх
$\tau_n$ -- яю\-ёыхфю\-тр\-Єхы№\-эюёЄ№, єфют\-ыхЄ\-тю\-Ё ■\-∙р 
єёыю\-тш ь $\tau_n\to\infty$, $\frac{\tau_n R_{n-1}}{R_n}\to 0$, Єю
$\mu_{t_n}\to 0$.

\par ╬сючэрўшь ўхЁхч $H$ яЁртє■ ўрёЄ№ Ёртхэ\-ёЄ\-тр  (\ref{for_example_3}).
╠√ фю\-ърчр\-ыш, ўЄю $H\subset$ $Fr[\mu]$.

\par ╧єёЄ№ ЄхяхЁ№ $\nu$ -- яЁюшч\-тюы№\-эр  ьхЁр шч $Fr[\mu]$. ╥юуфр
$\nu=\lim\mu_{t_n}$, уфх $t_n\to\infty$. ─юяюыэш\-Єхы№эю ьюц\-эю
ёўш\-ЄрЄ№, ўЄю ърцф√щ ёху\-ьхэЄ $[R_{k-1},R_{k+1}]$ ёю\-фхЁ\-цшЄ эх
сюыхх, ўхь юфэє Єюўъє  $t_n$. ╬сю\-чэрўшь ўхЁхч $m(n)$ Єръюх Ўхыюх
ўшёыю, ўЄю Єюўър $\ln R_{m(n)}$ хёЄ№ сышцрщ\-°р  Єюўър
яю\-ёыхфю\-тр\-Єхы№\-эюё\-Єш $\ln R_k$ ъ Єюўъх  $\ln t_n$ (хёыш
Єръшї ўшёхы фтр, Єю т ърўхёЄ\-тх $m(n)$ т√\-сшЁр\-хь ьхэ№\-°хх).
─юяюыэш\-Єхы№\-эюх юуЁрэш\-ўх\-эшх эр яю\-ёыхфю\-тр\-Єхы№\-эюёЄ№
$t_n$ юсхё\-ях\-ўш\-тр\-хЄ т√\-яюыэх\-эшх єёыю\-тш  $m(n_1)$ $\neq$
$m(n_2)$ яЁш $n_1$ $\neq$ $n_2$. ╥ръшь юсЁр\-чюь, ёє\-∙хёЄ\-тє\-хЄ
юфэю\-чэрў\-эюх яЁхф\-ёЄрт\-ых\-эшх ўшёыр  $t_n$ т тшфх
\begin{equation}\label{example_3_t_n}
t_n=\frac{1}{\tau_n}R_{m(n)}.
\end{equation}
─юяюыэш\-Єхы№эю ьюц\-эю яЁхф\-яюыю\-цшЄ№, ўЄю
яю\-ёыхфю\-тр\-Єхы№\-эюёЄ№ $\tau_n$ ёїю\-фшЄ\-ё  ышсю ъ
ъю\-эхў\-эюьє яЁх\-фхыє, ышсю ъ схё\-ъю\-эхў\-эюё\-Єш.
╧Ёхф\-яюыю\-цшь, ўЄю
\begin{equation}\label{example_tau_n}
\lim\limits_{n\to\infty} \tau_n=\infty.
\end{equation}
┬√яюыэ \-■Єё  ёююЄ\-эю°х\-эш 
$$
\ln t_n-\ln R_{m(n)-1}\geq \ln R_{m(n)}-\ln t_n,
$$
$$
t_n\geq \sqrt{\frac{R_{m(n)}}{R_{m(n)-1}}},\qquad \tau_n\leq
\sqrt{R_{m(n)-1}R_{m(n)}},
$$
$$
\frac{\tau_n R_{m(n)-1}}{R_{m(n)}}\leq
\sqrt{\frac{R_{m(n)-1}}{R_{m(n)}}},
$$
\begin{equation}\label{example_3_lim_frac}
\lim\limits_{n\to\infty}\frac{\tau_n R_{m(n)-1}}{R_{m(n)}}=0.
\end{equation}
╚ч ёююЄ\-эю°х\-эшщ (\ref{example_3_t_n})-(\ref{example_3_lim_frac})
ёых\-фєхЄ, ўЄю $\mu_{t_n}\to 0$. ╟эрўшЄ, т Ёрё\-ёьрЄЁш\-трх\-ьюь
ёыєўрх $\nu=0$.

\par └эрыюушўэю фю\-ърч√\-тр\-хЄё , ўЄю т ёыєўрх, хёыш $\tau_n\to 0$ Єръ\-цх
т√\-яюыэ \-хЄ\-ё  ёююЄ\-эю°х\-эшх $\mu_{t_n}\to 0$. ▌Єю ёых\-фє\-хЄ
шч ыхуъю яЁю\-тхЁ х\-ьюую ёююЄ\-эю°х\-эш 
\begin{equation}\label{example_3_lim_infty}
\lim\limits_{n\to\infty} \frac{\tau_n R_{m(n)+1}}{R_{m(n)}}=\infty.
\end{equation}

\par ╧єёЄ№ ЄхяхЁ№ $\tau_n\to \tau\in(0,\infty)$. ╙цх фюърчрэю
ёююЄ\-эю°х\-эшх $\mu_{{\tilde{t}}_n}\to\tau^\rho\delta(x-\tau)$, уфх
${\tilde{t}}_n=\frac{1}{\tau}R_{m(n)}$. ╥хяхЁ№ шч Єхю\-Ёх\-ь√
\ref{azarin_th_3}, яЁш\-ьхэ╕э\-эющ ъ яю\-ёыхфю\-тр\-Єхы№\-эюё\-Єш
${\tilde{t}}_n$ ёых\-фє\-хЄ, ўЄю $\mu_{t_n}\to\tau^\rho
\delta(x-\tau)$. ╤ыхфю\-тр\-Єхы№эю, т  ¤Єюь ёыєўрх
$\nu=\tau^\rho\delta(x-\tau)$. ╠√ фю\-ърчр\-ыш, ўЄю т ы■сюь ёыєўрх
$\nu\in H$. ╥хь ёрь√ь фю\-ърчр\-эю ёююЄ\-эю°х\-эшх $Fr[\mu]$
$\subset H$. ╚ч яЁш\-тхфхэ\-э√ї Ёрё\-ёєцфх\-эшщ ёых\-фє\-хЄ
єЄ\-тхЁ\-ц\-фх\-эшх ыхьь√. ╦хььр фю\-ърчр\-эр.

\par ╬ЄьхЄшь, ўЄю фюърчр\-Єхы№\-ёЄтю ыхьь√ ьюц\-эю шч\-тыхў№ шч
Ёрсю\-Є√  \cite{Shouigi}.

%%%%%%%%%%%%%%%%%%%%%%%%%%%%%%%%%%%%%%%%%%%%%%%%%%%%%%%%%%%%%%%%%%%%%%%%%%%%%
\par ╤ Єюўъш чЁхэш  ЄхюЁшш яЁх\-фхы№\-э√ї ьэю\-цх\-ёЄт эрш\-сюыхх яЁюё\-Є√\-ьш
 ты \-■Є\-ё  Єх ьхЁ√ $\mu\in$ $\frak{M}_\infty$ $(\rho(r))$, фы 
ъюЄю\-Ё√ї яЁх\-фхы№\-эюх ьэю\-цхёЄ\-тю  $Fr[\mu]$ ёю\-ёЄюшЄ шч юфэющ
ьхЁ√ $\nu$. ╥ръшх ьхЁ√ ь√ эр\-ч√тр\-хь {\it Ёхує\-ы Ё\-э√\-ьш} шыш
{\it Ёхує\-ы Ё\-э√\-ьш т ёь√ёых └чрЁш\-эр}. ╚ёёых\-фє\-хь, ъръющ
ью\-цхЄ с√Є№ яЁх\-фхы№\-эр  ьхЁр фы  Ёхує\-ы Ё\-эющ ьхЁ√ $\mu$.
%%%%%%%%%%%%%%%%%%%%%%%%%%%%%%%%%%%%%%%%%%%%%%%%%%%%%%%%%%%%%%%%%%%%%%%%%%%%%%%%%%%%%%%
\begin{theorem}\label{azarin_th_8}\hskip-2mm{.}\:
%%%%%%%%%%%%%%%%%%%%%%%%%%%%%%%%%%%%%%%%%%%%%%%%%%%%%%%%%%%%%%%%%%%%%%%%%%%%%%%%%%%%%%%
╧єёЄ№ $\rho(r)$ -- єЄюў\-э╕э\-э√щ яю\-Ё фюъ,
$\rho=\lim\limits_{n\to\infty}\rho(r)$ ш $\mu\in$ $\frak{M}_\infty$
$(\rho(r))$. ╧єёЄ№ ьхЁр $\mu$ Ёхує\-ы Ё\-эр  ш $\{\nu\}=Fr[\mu]$.
╥юу\-фр ёє\-∙хёЄ\-тє\-хЄ ъюья\-ыхъё\-эюх ўшёыю $c$ Єръюх, ўЄю
$d\nu(r)=cr^{\rho-1}dr$.
\end{theorem}
%%%%%%%%%%%%%%%%%%%%%%%%%%%%%%%%%%%%%%%%%%%%%%%%%%%%%%%%%%%%%%%%%%%%%%%%%%%%%%%%%%%%%%%
\par \quad {\sc ─юърчрЄхы№ёЄтю.}\: ┼ёыш ъюья\-ыхъё\-эр  Ёр\-фю\-эю\-тр ьхЁр
$\mu=\mu_1+i\mu_2$ Ёхує\-ы Ё\-эр, Єю Ёхує\-ы Ё\-э√ьш сєфєЄ ш ьхЁ√
$\mu_1$ ш  $\mu_2$. ╧ю¤Єю\-ьє фюёЄр\-Єюў\-эю фю\-ърч√\-трЄ№
Єхю\-Ёхьє фы  тх∙хёЄ\-тхэ\-э√ї ьхЁ ╨рфюэр. ╧ю\-ёъюы№\-ъє
ьэю\-цхёЄ\-тю $Fr[\mu]$ шэтр\-ЁшрэЄ\-эю юЄэю\-ёш\-Єхы№\-эю
яЁх\-юсЁр\-чю\-тр\-эш  $F_t$, Єю т√\-яюы\-э \-хЄ\-ё  Ёртхэ\-ёЄ\-тю
$\nu=F_t\nu$. ╧єёЄ№ $0<a<b<\infty$. ╥юу\-фр, т ўрёЄ\-эюё\-Єш,
$\nu((a,b])= (F_t\nu)((a,b])= t^{-\rho}\nu((at,bt])$. ┬√сшЁр 
$t=\frac{1}{a}$, яюыє\-ўшь $\nu((a,b])=a^\rho\nu((1,\frac{b}{a}])$.
╬сю\-чэр\-ўшь $N(s)=\nu((1,1+s])$, $s>0$. ╥юу\-фр
$$
N(s_1+s_2)=\nu((1,1+s_1])+\nu((1+s_1,1+s_1+s_2])=N(s_1)
+(1+s_1)^\rho N\left(\frac{s_2}{1+s_1}\right).
$$
╘єэъ\-Ўшш  $N$, єфютыхЄ\-тю\-Ё ■\-∙шх Єръюьє Ёртхэ\-ёЄ\-тє т
{\cite{Grishin_2000}} эр\-ч√тр\-■Є\-ё  {\it
$\rho$-рффш\-Єшт\-э√\-ьш}. ╧ю\-ёъюы№\-ъє $\nu$ -- ыю\-ъры№\-эю
ъю\-эхў\-эр  ьхЁр эр яюыє\-юёш  $(0,\infty)$, Єю Їєэъ\-Ўш   $N$
юуЁрэш\-ўх\-эр эр яюыє\-шэ\-ЄхЁ\-тр\-ых  $(0,1]$. ╧ю Єхю\-Ёх\-ьх 4
{\cite{Grishin_2000}} ёє\-∙хёЄ\-тє\-хЄ тх∙хёЄ\-тхэ\-эюх  $c$ Єръюх,
ўЄю ёяЁр\-тхф\-ыш\-т√ ёююЄ\-эю\-°х\-эш 
$$
N(s)=\frac{c}{\rho}\left((1+s)^\rho-1\right),\:\rho\neq 0, \quad
N(s)=c\ln(1+s),\:\rho=0.
$$
╠√ яюыєўр\-хь, ўЄю яЁш  $b>1$ т√яюыэ \-■Є\-ё  Ёртхэ\-ёЄ\-тр
$$
\nu((1,b])=\frac{c}{\rho}\left(b^\rho-1\right),\: \rho\neq 0, \quad
\nu((1,b])=c\ln b,\:\rho=0.
$$
╚ч ¤Єшї ЁртхэёЄт ыхуъю ёых\-фє\-хЄ єЄ\-тхЁ\-ц\-фх\-эшх Єхю\-Ёх\-ь√.
╥хю\-Ёх\-ьр фю\-ърчр\-эр.

%%%%%%%%%%%%%%%%%%%%%%%%%%%%%%%%%%%%%%%%%%%%%%%%%%%%%%%%%%%%%%%%%%%%%%%%%%%%%%%%%%%%%%%
\par ╤ыхфє\-■\-∙хх ётющ\-ёЄтю Ёхуєы Ё\-э√ї ьхЁ хёЄ№ яЁюёЄюх
ёыхф\-ёЄ\-тшх юяЁхфх\-ых\-эшщ.

\par ┴єфхь уютюЁшЄ№, ўЄю ёхЄ№ ьхЁ $\mu_R$ $(R\in(0,\infty))$ °ш\-Ёю\-ъю
ёїю\-фшЄ\-ё  ъ ьхЁх  $\mu$ яЁш $R\to\infty$, хёыш фы  ы■сющ
Їєэъ\-Ўшш $\varphi\in\Phi$ т√\-яюыэ \-хЄё  Ёртхэ\-ёЄтю
$\lim\limits_{R\to\infty}(\mu_R,\varphi)$ $=(\mu,\varphi)$.
%%%%%%%%%%%%%%%%%%%%%%%%%%%%%%%%%%%%%%%%%%%%%%%%%%%%%%%%%%%%%%%%%%%%%%%%%%%%%%%%%%%%%%%
\begin{theorem}\label{azarin_th_9}\hskip-2mm{.}\:
%%%%%%%%%%%%%%%%%%%%%%%%%%%%%%%%%%%%%%%%%%%%%%%%%%%%%%%%%%%%%%%%%%%%%%%%%%%%%%%%%%%%%%%
╧єёЄ№ $\rho(r)$ -- єЄюў\-э╕э\-э√щ яюЁ \-фюъ, $\mu\in$
$\frak{M}_\infty$ $(\rho(r))$  -- Ёхує\-ы Ё\-эр  ьхЁр, яЁш\-ў╕ь
$\{\nu\}=Fr[\mu]$. ╥юу\-фр
\begin{equation}{\label{azarin_8}}
\nu=\lim\limits_{R\to\infty}\mu_R.
\end{equation}
\end{theorem}
%%%%%%%%%%%%%%%%%%%%%%%%%%%%%%%%%%%%%%%%%%%%%%%%%%%%%%%%%%%%%%%%%%%%%%%%%%%%%%%%%%%%%%%
\par \quad {\sc ─юърчрЄхы№\-ёЄтю.}\: ╧юърцхь тэрўр\-ых, ўЄю фы  ы■сющ
Їєэъ\-Ўшш $\varphi\in\Phi$ ёє\-∙хёЄ\-тєхЄ яЁх\-фхы
\begin{equation}{\label{azarin_9}}
\lim\limits_{R\to\infty}\int\limits_0^\infty
\varphi(x)d\mu_R(x)=a(\varphi).
\end{equation}
┼ёыш ¤Єю єЄ\-тхЁ\-цфх\-эшх эхтхЁэю, Єю ёє\-∙хёЄ\-тє\-■Є Їєэъ\-Ўш 
$\varphi\in \Phi$ ш фтх яю\-ёых\-фю\-тр\-Єхы№\-эюё\-Єш
$r_n\to\infty$ ш $R_n\to\infty$ Єръшх, ўЄю
\begin{equation}{\label{azarin_10}}
\lim\limits_{n\to\infty}\int\limits_0^\infty
\varphi(x)d\mu_{r_n}(x)\neq \lim\limits_{n\to\infty}
\int\limits_0^\infty \varphi(x)d\mu_{R_n}(x).
\end{equation}
╥ръ ъръ ёхьхщ\-ёЄ\-тю ьхЁ $\mu_R$, $R\geq 1$, ъюь\-яръЄ\-эю, Єю, эх
юуЁрэш\-ўш\-тр  юс∙\-эюё\-Єш, ьюц\-эю ёўш\-ЄрЄ№, ўЄю
яю\-ёых\-фю\-тр\-Єхы№\-эюё\-Єш $\mu_{r_n}$ ш $\mu_{R_n}$ °ш\-Ёю\-ъю
ёїюф Є\-ё . ╧ю\-ёъюы№\-ъє яЁх\-фхы№\-эюх ьэю\-цхёЄ\-тю  $Fr[\mu]$
ёюёЄю\-шЄ шч юфэющ ьхЁ√ $\nu$, Єю $\mu_{r_n}\to\nu$ ш
$\mu_{R_n}\to\nu$. ╧ю юяЁх\-фх\-ых\-эш■ °ш\-Ёю\-ъющ ёїю\-фш\-ьюё\-Єш
$$
\lim\limits_{n\to\infty}\int\limits_0^\infty
\varphi(x)d\mu_{r_n}(x)=\int\limits_0^\infty \varphi(x)d\nu(x),
$$
\begin{equation}\label{azarin_lim_int}
\lim\limits_{n\to\infty}\int\limits_0^\infty
\varphi(x)d\mu_{R_n}(x)=\int\limits_0^\infty \varphi(x)d\nu(x).
\end{equation}
▌Єю яЁюЄш\-тю\-Ёх\-ўшЄ  (\ref{azarin_10}) ш, Єхь ёрь√ь, Ёртхэ\-ёЄтю
(\ref{azarin_9}) фюърчр\-эю. ╚ч Ёртхэ\-ёЄ\-тр (\ref{azarin_lim_int})
ёых\-фє\-хЄ, ўЄю $a(\varphi)=\int\limits_0^\infty
\varphi(x)d\nu(x)$. ╥хю\-Ёх\-ьр фюърчр\-эр.

%%%%%%%%%%%%%%%%%%%%%%%%%%%%%%%%%%%%%%%%%%%%%%%%%%%%%%%%%%%%%%%%%%%%%%%%%%%%%%%%%%%%%%%%%
\par ┬ ёых\-фє\-■\-∙хщ Єхю\-Ёх\-ьх юяшё√\-тр\-■Є\-ё  яюыю\-цш\-Єхы№\-э√х
Ёхує\-ы Ё\-э√х ьхЁ√.
%%%%%%%%%%%%%%%%%%%%%%%%%%%%%%%%%%%%%%%%%%%%%%%%%%%%%%%%%%%%%%%%%%%%%%%%%%%%%%%%%%%%%%%%%
\begin{theorem}\label{azarin_th_10}\hskip-2mm{.}\:
%%%%%%%%%%%%%%%%%%%%%%%%%%%%%%%%%%%%%%%%%%%%%%%%%%%%%%%%%%%%%%%%%%%%%%%%%%%%%%%%%%%%%%%%%
╧єёЄ№ $\rho(r)$ -- яЁюшч\-тюы№э√щ єЄюў\-э╕э\-э√щ яю\-Ё фюъ, $\mu$ --
яюыю\-цш\-Єхы№\-эр  ьхЁр эр яюыє\-юёш  $(0,\infty)$, $\mu\in$
$\frak{M}_\infty$ $(\rho(r))$. ─ы  Єюую, ўЄюс√ ьхЁр  $\mu$ с√ыр
Ёхуєы Ё\-эющ юЄ\-эюёш\-Єхы№\-эю єЄюў\-э╕э\-эю\-ую яюЁ \-фър
$\rho(r)$, эхюс\-їюфш\-ью ш фюёЄр\-Єюў\-эю, ўЄюс√ т√\-яюыэ \-ыюё№
ёых\-фє\-■\-∙хх єёыю\-тшх (т чр\-тшёш\-ьюё\-Єш юЄ чэрър  $\rho$ ¤Єю
єёыю\-тшх чр\-яшё√\-тр\-хЄ\-ё  т Ёрч\-ышўэющ ЇюЁьх):
\begin{equation}{\label{azarin_11}}
\lim\limits_{R\to\infty}\frac{\mu((1,R])}{V(R)}= c, \quad \rho>0,
\end{equation}\begin{equation}{\label{azarin_12}}
\lim\limits_{R\to\infty}\frac{\mu([R,\infty))}{V(R)}=c, \quad
\rho<0,
\end{equation}\begin{equation}{\label{azarin_13}}
\lim\limits_{R\to\infty}\frac{\mu((aR,bR])}{V(R)}= c\ln\frac{b}{a},
\quad \rho=0,
\end{equation}
фы  ы■с√ї  $a$ ш  $b$, $0<a<b<\infty$.
\end{theorem}
%%%%%%%%%%%%%%%%%%%%%%%%%%%%%%%%%%%%%%%%%%%%%%%%%%%%%%%%%%%%%%%%%%%%%%%%%%%%%%%%%%%%%%%%%
\par \quad {\sc ─юърчрЄхы№ёЄтю.}\: ╦хуъю єтшфхЄ№, чр\-яшё√\-тр 
ёююЄ\-тхЄ\-ёЄ\-тє■\-∙шх шэ\-ЄхуЁры№\-э√х ёєь\-ь√
╨ш\-ьр\-эр-╤Єшы\-Є№х\-ёр, ўЄю шч ЁртхэёЄт
(\ref{azarin_11})-(\ref{azarin_13}) ёых\-фє\-хЄ, ўЄю фы  ы■сющ
Їєэъ\-Ўшш $\varphi\in\Phi$ т√\-яюы\-э \-хЄ\-ё  Ёртхэ\-ёЄ\-тю
\begin{equation}{\label{azarin_14}}
\lim\limits_{R\to\infty}\int\limits_0^\infty \varphi(x)d\mu_R(x)=
\int\limits_0^\infty \varphi(x)d\nu(x),
\end{equation}
уфх $d\nu(x)=c|\rho|x^{\rho-1}dx$, хёыш $\rho\neq 0$ ш $d\nu(x)=$
$\frac{c}{x}dx$, хёыш $\rho=0$. ╬ўхтшф\-эю, ўЄю шч (\ref{azarin_14})
ёых\-фє\-хЄ Ёхує\-ы Ё\-эюёЄ№ ьхЁ√ $\mu$.

\par ╧хЁх\-їюфшь ЄхяхЁ№ ъ фюърчр\-Єхы№\-ёЄтє эхюс\-їю\-фш\-ьюё\-Єш.
╧єёЄ№ $\mu$ -- Ёхує\-ы Ё\-эр  ьхЁр. ╥юу\-фр яю Єхю\-Ёх\-ьх
\ref{azarin_th_9} $\lim\limits_{R\to\infty}\mu_R=\nu$. ╧ю
Єхю\-Ёх\-ьх  \ref{azarin_th_7} ьхЁр  $\nu$ эх\-яЁх\-Ё√т\-эр. ╚ч
Єхю\-Ёх\-ь√  \ref{measures_th_2} ёых\-фєхЄ, ўЄю
$\lim\limits_{R\to\infty} \mu_R([a,b])$ $=\nu([a,b])$. ▌Єю
¤ътш\-тр\-ыхэЄ\-эю Ёртхэ\-ёЄ\-трь
(\ref{azarin_11})-(\ref{azarin_13}). ╥хю\-Ёх\-ьр фю\-ърчр\-эр.
%%%%%%%%%%%%%%%%%%%%%%%%%%%%%%%%%%%%%%%%%%%%%%%%%%%%%%%%%%%%%%%%%%%%%%%%%%%%%%%%%%%%%%%%%
\par ╥хю\-Ёхьр  \ref{azarin_th_10} эх ёяЁртхф\-ыштр фы  тх∙хёЄ\-тхэ\-э√ї ьхЁ
╨рфюэр. ╚ч Єхю\-Ёх\-ь√ \ref{azarin_th_8} ш ыхьь√
\ref{azarin_example_2} ёых\-фє\-хЄ, ўЄю фы  Єюую, ўЄюс√ юяшёрЄ№
тх∙хёЄ\-тхэ\-э√х Ёхує\-ы Ё\-э√х ьхЁ√, фю\-ёЄр\-Єюў\-эю юяшёрЄ№
тх∙хёЄ\-тхэ\-э√х ьхЁ√ ё эєых\-т√ь яЁх\-фхы№\-э√ь ьэю\-цх\-ёЄ\-тюь.
▌Єю ёфхыр\-эю т ёых\-фє\-■\-∙хщ Єхю\-Ёхьх.
%%%%%%%%%%%%%%%%%%%%%%%%%%%%%%%%%%%%%%%%%%%%%%%%%%%%%%%%%%%%%%%%%%%%%%%%%%%%%%%%%%%%%%%%%
\begin{theorem}\label{azarin_th_11}\hskip-2mm{.}\:
%%%%%%%%%%%%%%%%%%%%%%%%%%%%%%%%%%%%%%%%%%%%%%%%%%%%%%%%%%%%%%%%%%%%%%%%%%%%%%%%%%%%%%%%%
╧єёЄ№ $\mu\in$ $\frak{M}_\infty$ $(\rho(r))$. ─ы  Єюую, ўЄюс√
т√\-яюыэ \-ыюё№ Ёртхэ\-ёЄтю $Fr[\mu]$ $=$ $\{0\}$, эхюс\-їюфш\-ью ш
фюёЄр\-Єюў\-эю, ўЄюс√ ёє\-∙хёЄ\-тю\-трыр тюч\-ЁрёЄр■\-∙р 
яю\-ёыхфю\-тр\-Єхы№\-эюёЄ№ $r_n$ ёїюф \-∙р \-ё  ъ
схё\-ъю\-эхў\-эюё\-Єш Єр\-ър , ўЄю
$$
\lim\limits_{n\to\infty} \frac{r_{n+1}}{r_n}=1,
$$
$$
\lim\limits_{R\to\infty} \frac{1}{V(R)}
\sum\limits_{(r_n,r_{n+1}]\cap[R,2R]\neq \emptyset}
|\mu((r_n,r_{n+1}])|=0.
$$
\end{theorem}
%%%%%%%%%%%%%%%%%%%%%%%%%%%%%%%%%%%%%%%%%%%%%%%%%%%%%%%%%%%%%%%%%%%%%%%%%%%%%%%%%%%%%%%%%
\par ─юърчр\-Єхы№ёЄтю Єхю\-Ёх\-ь√  \ref{azarin_th_11} эрь шчтхёЄ\-эю.
╬эю чэрўш\-Єхы№\-эю ёыюцэхх фю\-ърчр\-Єхы№\-ёЄтр Єхю\-Ёх\-ь√
\ref{azarin_th_10} ш т√\-їюфшЄ чр Ёрьъш эр°хщ ЁрсюЄ√. ╟рьхЄшь,
юфэръю, ўЄю фю\-ърчр\-Єхы№\-ёЄтю Єхю\-Ёх\-ь√  \ref{azarin_th_11} т
ёЄюЁю\-эє фю\-ёЄрЄюў\-эюё\-Єш фютюы№\-эю яЁюёЄю. ╥хю\-Ёхьр
\ref{azarin_th_11} яЁш\-тхфх\-эр фы  яюыэю\-Є√ шэЇюЁ\-ьр\-Ўшш.

 %%%%%%%%%%%%%%%%%%%%%%%%%%%%%%%%%%%%%%%%%%%%%%%%%%%%%%%%%%%%%%%%%%%%%%%%%%
\section{└схых\-т√ Єхю\-Ёх\-ь√ фы  шэ\-ЄхуЁрыют}
%%%%%%%%%%%%%%%%%%%%%%%%%%%%%%%%%%%%%%%%%%%%%%%%%%%%%%%%%%%%%%%%%%%%%%%%%%
\qquad ╠√ яЁш\-ёЄєяр\-хь ъ шчыюцх\-эш■ Ёх\-чєы№\-ЄрЄют, ъюЄю\-Ё√х
ьюуєЄ с√Є№ юяшёр\-э√ ъръ рсхых\-т√ Єхю\-Ёх\-ь√ фы  шэ\-ЄхуЁр\-ыют.
╬ёэют\-э√ьш юс·хъ\-Єр\-ьш эр°х\-ую шё\-ёыхфю\-тр\-эш   ты \-■Є\-ё 
Їєэъ\-Ўшш $\Psi(r)$, $J(r)$, ьхЁр  $s$, ьэю\-цх\-ёЄтю $L(J,\infty)$
(ёьюЄЁш ЇюЁ\-ьє\-ыє  (\ref{introduction_Psi(r)}) ш ёых\-фє\-■\-∙шщ
чр эхщ ЄхъёЄ).

\par ╬фэр шч Ўхыхщ ¤Єюую Ёрч\-фхыр -- юяшёрЄ№ ётющ\-ёЄтр Їєэъ\-Ўшш $J(r)$ т
чр\-тшёш\-ьюёЄш юЄ юуЁрэш\-ўх\-эшщ эр  фЁю  $K$ ш ьхЁє  $\mu$.

\par ╠√ эрўэ╕ь ёю ёых\-фє\-■\-∙хую яЁюёЄю фю\-ърч√\-трх\-ьюую
єЄ\-тхЁ\-ц\-фх\-эш . ╥хь эх ьхэхх ¤Єю єЄ\-тхЁ\-ц\-фх\-эшх
ёю\-фхЁ\-цшЄ юёэют\-эє■ шфх■ ьхЄюфр -- шё\-яюы№чю\-трЄ№ фы 
юяшёр\-эш  ётющёЄт Їєэъ\-Ўшш $J(r)$ яЁх\-фхы№\-эюх ьэю\-цх\-ёЄтю
$Fr[\mu]$ ьхЁ√ $\mu$.

%%%%%%%%%%%%%%%%%%%%%%%%%%%%%%%%%%%%%%%%%%%%%%%%%%%%%%%%%%%%%%%%%%%%%%%%%%%%%%%%%%%%%%%
\begin{theorem}\label{abel_th_1}\hskip-2mm{.}\:
%%%%%%%%%%%%%%%%%%%%%%%%%%%%%%%%%%%%%%%%%%%%%%%%%%%%%%%%%%%%%%%%%%%%%%%%%%%%%%%%%%%%%%%
╧єёЄ№ $\rho(r)$ -- яЁю\-шч\-тюы№\-э√щ єЄюў\-э╕э\-э√щ яюЁ \-фюъ,
$\mu\in$ $\frak{M}_\infty$ $(\rho(r))$, $K$  -- эх\-яЁх\-Ё√т\-эюх
ЇшэшЄ\-эюх  фЁю эр яюыє\-юёш $(0,\infty)$. ╥юуфр т√\-яюыэ \-хЄ\-ё 
Ёртхэ\-ёЄтю
\begin{equation}\label{abel_def_L(I)}
L(J,\infty)=\left\{\int\limits_0^\infty K(u)d\nu(u):\quad \nu\in
Fr[\mu]\right\}.
\end{equation}
\end{theorem}
%%%%%%%%%%%%%%%%%%%%%%%%%%%%%%%%%%%%%%%%%%%%%%%%%%%%%%%%%%%%%%%%%%%%%%%%%%%%%%%%%%%%%%%
\par \quad {\sc ─юърчрЄхы№ёЄтю.}\: ╬сючэрўшь ўхЁхч $H$ яЁртє■ ўрёЄ№
Ёртхэ\-ёЄ\-тр  (\ref{abel_def_L(I)}). ╧єёЄ№ $\nu$  --
яЁю\-шч\-тюы№\-эр  ьхЁр шч ьэю\-цх\-ёЄтр $Fr[\mu]$. ╤є∙хёЄ\-тєхЄ
яю\-ёыхфю\-тр\-Єхы№\-эюёЄ№ $r_n\to\infty$ Єр\-ър , ўЄю
$\mu_{r_n}\to\nu$. ╥юуфр яю юяЁх\-фхых\-эш■ °ш\-Ёю\-ъющ
ёїюфш\-ьюё\-Єш $J(r_n)\to\int\limits_0^\infty K(u)d\nu(u)$. ╠√
фюърчр\-ыш тъы■ўх\-эшх $H\subset L(J,\infty)$.

\par ╧єёЄ№ ЄхяхЁ№ $r_n\to\infty$ Єр\-ър  яю\-ёыхфю\-тр\-Єхы№\-эюёЄ№,
ўЄю яю\-ёыхфю\-тр\-Єхы№\-эюёЄ№ $J(r_n)$ ёїю\-фшЄ\-ё  т
ёюсёЄ\-тхэ\-эюь шыш эх\-ёюсёЄ\-тхээюь (ёїю\-фшЄ\-ё  ъ
схё\-ъю\-эхў\-эюё\-Єш) ёь√ёых. ╥ръ ъръ яю ыхььх \ref{azarin_lemma_1}
яюыє\-ЄЁрхъ\-Єю\-Ёш  $\mu_r$, $r\geq 1$, ъюь\-яръЄэр, Єю эх
юуЁрэш\-ўш\-тр  юс∙\-эюё\-Єш ьюц\-эю ёўш\-ЄрЄ№, ўЄю
яю\-ёыхфю\-тр\-Єхы№\-эюёЄ№ $\mu_{r_n}$ °ш\-Ёю\-ъю ёїю\-фшЄё  ъ
эх\-ъю\-Єю\-Ёющ ьхЁх $\nu$. ╧ю юяЁх\-фхых\-эш■ $Fr[\mu]\:$ $\nu\in
Fr[\mu]$. ╧ю фю\-ърчрэ\-эюьє $J(r_n)\to \int\limits_0^\infty
K(u)d\nu(u)$. ╠√ фю\-ърчр\-ыш ёяЁр\-тхф\-ыш\-тюёЄ№ тъы■ўх\-эш 
$L(J,\infty)\subset H$, р тьхёЄх ё ¤Єшь -- ш Єхю\-Ёхьє.

\par ╟эрўш\-Єхы№\-эр  ўрёЄ№ фры№\-эхщ\-°шї Ёх\-чєы№\-Єр\-Єют уырт√ сєфхЄ
трЁш\-рэ\-Єрьш Єхю\-Ёх\-ь√  \ref{abel_th_1}. ╠√ сєфхь
Ёрё\-ёьрЄЁш\-трЄ№ Ёрч\-ышў\-э√х юуЁрэш\-ўх\-эш  эр  фЁю  $K$ ш ьхЁє
$\mu$. ─ы  эрўрыр яю\-ёьюЄЁшь, ъ ўхьє тхф╕Є юЄърч юЄ
эх\-яЁх\-Ё√т\-эюё\-Єш  фЁр  $K$. ┬  ¤Єюь ёыєўрх ьэю\-цх\-ёЄ\-тю
$Fr[\mu]$ єцх эх юяЁх\-фхы \-хЄ ьэю\-цх\-ёЄ\-тю $L(J,\infty)$.

\par ╧єёЄ№ $\mu\in$ $\frak{M}_\infty$ $(\rho(r))$. {\it ╨рё°шЁхэ\-э√ь
яЁх\-фхы№\-э√ь ьэю\-цх\-ёЄ\-тюь} $\widehat{Fr}[\mu]$ ьхЁ√ $\mu$
эр\-чю\-т╕ь ьэю\-цх\-ёЄтю ярЁ ьхЁ $(\nu_1,\nu_2)$ Єръшї, ўЄю
ёє\-∙хёЄ\-тєхЄ яю\-ёыхфю\-тр\-Єхы№\-эюёЄ№ $r_n \to\infty$ Єр\-ър ,
ўЄю $\mu_{r_n} \to\nu$, ${\mu^1}_{r_n}\to\nu_1$,
${\mu^2}_{r_n}\to\nu_2$, уфх $\mu^1_{r_n}$  -- юуЁрэш\-ўх\-эшх ьхЁ√
$\mu_{r_n}$ эр яюыє\-шэ\-ЄхЁ\-тры $(0,1]$, р $\mu^2_{r_n}$  --
юуЁрэш\-ўх\-эшх ьхЁ√ $\mu_{r_n}$ эр яюыє\-юё№ $(1,\infty)$.

%%%%%%%%%%%%%%%%%%%%%%%%%%%%%%%%%%%%%%%%%%%%%%%%%%%%%%%%%%%%%%%%%%%%%%%%%%%%%%%%%%%%%%%
\begin{theorem}\label{abel_th_2}\hskip-2mm{.}\:
%%%%%%%%%%%%%%%%%%%%%%%%%%%%%%%%%%%%%%%%%%%%%%%%%%%%%%%%%%%%%%%%%%%%%%%%%%%%%%%%%%%%%%%
╧єёЄ№ $\rho(r)$ -- яЁю\-шч\-тюы№\-э√щ єЄюў\-э╕э\-э√щ яюЁ \-фюъ,
$\mu\in$ $\frak{M}_\infty$ $(\rho(r))$, $K$  -- Їєэъ\-Ўш  ЇшэшЄ\-эр 
эр яюыє\-юёш $(0,\infty)$ эх\-яЁх\-Ё√т\-эр  тё■фє, ъЁюьх Єюўъш
$t=1$, ъюЄю\-Ёр   ты \-хЄё  Єюўъющ Ёрч\-Ё√тр яхЁтю\-ую Ёюфр фы 
Їєэъ\-Ўшш $K(t)$, яЁш\-ў╕ь Їєэъ\-Ўш  $K(t)$ эх\-яЁх\-Ё√т\-эр ёыхтр т
Єюўъх $t=1$. ╥юуфр т√\-яюыэ \-хЄ\-ё  Ёртхэ\-ёЄ\-тю
$$
L(J,\infty)=\left\{\int\limits_{(0,1]} K(t)d\nu_1(t)+
\int\limits_{[1,\infty)} \widetilde{K}(t) d\nu_2(t)
:\;(\nu_1,\nu_2)\in\widehat{Fr}[\mu]\right\},
$$
уфх $\widetilde{K}$  -- эх\-яЁх\-Ё√т\-эюх яЁю\-фюы\-цх\-эшх
Їєэъ\-Ўшш $K(t)$ ё яюыє\-юёш $(1,\infty)$ эр яюыє\-юё№ $[1,\infty)$.
\end{theorem}

\par ╬ЄьхЄшь, ўЄю їюЄ  ьхЁр $\mu^2_{r_n}$ эх эр\-уЁєцр\-хЄ Єюўъє  $1$,
ьхЁр  $\nu_2$ ью\-цхЄ эр\-уЁєцрЄ№ ¤Єє Єюўъє. ╧ю¤Єю\-ьє т
ЇюЁ\-ьє\-ыш\-Ёют\-ъх Єхю\-Ёх\-ь√ Їєэъ\-Ўш■ $\widetilde{K}(t)$ эхы№ч 
чр\-ьхэшЄ№ эр Їєэъ\-Ўш■ $K(t)$.

%%%%%%%%%%%%%%%%%%%%%%%%%%%%%%%%%%%%%%%%%%%%%%%%%%%%%%%%%%%%%%%%%%%%%%%%%%%%%%%%%%%%%%%
\par {\sc  ─юърчрЄхы№ёЄтю.}\: ─юърчр\-Єхы№\-ёЄтю
Єхю\-Ёх\-ь√  \ref{abel_th_2} яЁю\-їюфшЄ яю ёїхьх
фю\-ърчр\-Єхы№\-ёЄ\-тр Єхю\-Ёх\-ь√  \ref{abel_th_1}. ═єцэю Єюы№ъю
яюы№чю\-трЄ№\-ё  Ёртхэ\-ёЄ\-тр\-ьш Єшяр эшцх\-ёых\-фє\-■\-∙хую.
╧єёЄ№ $K_1(t)$  -- эх\-яЁх\-Ё√т\-эюх ЇшэшЄ\-эюх яЁю\-фюы\-цх\-эшх
Їєэъ\-Ўшш $K(t)$ ё яюыє\-шэ\-ЄхЁ\-трыр $(0,1]$ эр яюыє\-юё№
$(0,\infty)$. ╥юуфр
$$
\int\limits_0^1 K(t)d\nu_1(t)= \int\limits_0^1 K_1(t)d\nu_1(t)=
\int\limits_0^\infty K_1(t)d\nu_1(t)=\lim\limits_{n\to\infty}
\int\limits_0^\infty K_1(t)d\mu^1_{r_n}(t).
$$

\par ╧Ёш ёююЄ\-тхЄ\-ёЄ\-тє■\-∙шї юуЁрэш\-ўх\-эш\- ї эр ьхЁє  $\mu$
ЄЁхсю\-тр\-эшх эх\-яЁх\-Ё√т\-эюё\-Єш  фЁр  $K$ ёЄрэю\-тшЄ\-ё 
шч\-ыш°\-эшь. ┬ ёыєўрх Ёрч\-Ё√т\-э√ї  фхЁ  $K$ ш рсёю\-ы■Є\-эю
эх\-яЁх\-Ё√т\-э√ї ьхЁ $\mu$ шэ\-ЄхуЁры $\int\limits_0^\infty
K(t)d\mu(t)$ ёых\-фє\-хЄ яюэш\-ьрЄ№ ъръ $\int\limits_0^\infty
K(t)\mu'(t)dt$.

%%%%%%%%%%%%%%%%%%%%%%%%%%%%%%%%%%%%%%%%%%%%%%%%%%%%%%%%%%%%%%%%%%%%%%%%%%%%%%%%%%%%%%%%%%%%%%%%%%%
\begin{theorem}\label{abel_th_3}\hskip-2mm{.}\:
%%%%%%%%%%%%%%%%%%%%%%%%%%%%%%%%%%%%%%%%%%%%%%%%%%%%%%%%%%%%%%%%%%%%%%%%%%%%%%%%%%%%%%%%%%%%%%%%%%%
╧єёЄ№ $\rho(r)$  -- яЁю\-шч\-тюы№\-э√щ єЄюў\-э╕э\-э√щ яю\-Ё фюъ,
$\mu$  -- Ёр\-фю\-эютр ьхЁр эр яюыє\-юёш $(0,\infty)$ ё яыюЄ\-эюёЄ№■
$\mu'(r)$, єфют\-ыхЄ\-тю\-Ё ■\-∙хщ эх\-Ёртхэ\-ёЄ\-тє $|\mu'(r)|$
$\leq$ $ M\frac{V(r)}{r}$, $r\in[1,\infty)$, $K$  -- ЇшэшЄ\-эюх  фЁю
шч яЁюёЄ\-Ёрэ\-ёЄтр $L_1(0,\infty)$. ╥юуфр
\begin{equation}\label{abel_L(I)}
L(J,\infty)=\left\{\int\limits_a^b K(u)d\nu(u):\; \nu\in
Fr[\mu]\right\}.
\end{equation}
\end{theorem}
%%%%%%%%%%%%%%%%%%%%%%%%%%%%%%%%%%%%%%%%%%%%%%%%%%%%%%%%%%%%%%%%%%%%%%%%%%%%%%%%%%%%%%%%%%%%%%%%%%%%
\par {\sc  ─юърчрЄхы№ёЄтю.}\: ╧єёЄ№ $N_1(t)$ -- тхЁї\-э   яыюЄ\-эюёЄ№ ьхЁ√
$|\mu|$. ╚ьххь
$$
N_1(\alpha)=\mathop{\overline\lim}\limits_{r\to\infty}
\frac{|\mu|((r,(1+\alpha)r])}{V(r)}\leq \lim\limits_{r\to\infty}
M\int\limits_r^{(1+\alpha)r} \frac{1}{t}\frac{V(t)}{V(r)}dt
$$
$$
= M\lim\limits_{r\to\infty} \int\limits_1^{1+\alpha}
\frac{1}{u}\frac{V(ur)}{V(r)}du =M\frac{(1+\alpha)^\rho-1}{\rho}.
$$
╚ч Єхю\-Ёх\-ь√  \ref{azarin_th_6} ёых\-фєхЄ, ўЄю тёх ьхЁ√ тю
ьэю\-цх\-ёЄтх $Fr[|\mu|]$ эх\-яЁх\-Ё√т\-э√. ┬ёх ьхЁ√ тю
ьэю\-цх\-ёЄтх $Fr[\mu]$ Єръ\-цх сєфєЄ эх\-яЁх\-Ё√т\-э√. ╚ч
ёърчрэ\-эю\-ую ш Єхю\-Ёх\-ь√ \ref{measures_th_2} ёых\-фєхЄ, ўЄю хёыш
$\nu=\lim\limits_{n\to\infty}\mu_{t_n}$, Єю т√\-яюы\-э \-хЄё 
Ёртхэ\-ёЄ\-тю
$$
\nu([a,b])=\lim\limits_{n\to\infty} \mu_{t_n}([a,b])
$$
фы  ы■сюую ёху\-ьхэЄр $[a,b]\subset(0,\infty)$. ╚ч ¤Єюую ёых\-фєхЄ,
ўЄю
$$
|\nu([a,b])|<\mathop{\overline\lim}\limits_{r\to\infty}
\int\limits_a^b d|\mu|_{t_n}(t)\leq M\lim\limits_{n\to\infty}
\int\limits_a^b \frac{1}{t} \frac{V(t_n
t)}{V(t_n)}dt=M\frac{b^\rho-a^\rho}{\rho}.
$$
╚ч ¤Єюую эхЁртхэ\-ёЄ\-тр ёых\-фєхЄ, ўЄю ы■ср  ьхЁр $\nu$ шч
ьэю\-цх\-ёЄтр $Fr[\mu]$  ты \-хЄё  рсёю\-ы■Є\-эю эхяЁх\-Ё√т\-эющ ш
т√яюыэ \-хЄё  эх\-Ёртхэ\-ёЄтю $|\nu'(x)|\leq  Mx^{\rho-1}$. ╧єёЄ№
$r_n\to\infty$ Єр\-ър  яю\-ёыхфю\-тр\-Єхы№\-эюёЄ№, ўЄю
$\mu_{r_n}\to\nu$. ╧єёЄ№ $\varepsilon$  -- яЁю\-шч\-тюы№\-эюх ёЄЁюую
яюыю\-цш\-Єхы№\-эюх ўшёыю, р $K_1$  -- ЇшэшЄ\-эр  эх\-яЁх\-Ё√т\-эр 
Їєэъ\-Ўш  эр яюыє\-юёш $(0,\infty)$ ш Єр\-ър , ўЄю
$$
\int\limits_0^\infty |K(t)-K_1(t)|dt<\varepsilon.
$$
╥юуфр
$$
\left| \int\limits_0^\infty K(t)d\mu_{r_n}(t)- \int\limits_0^\infty
K(t)d\nu(t) \right| \leq \left| \int\limits_0^\infty (K(t)-K_1(t))
d\mu_{r_n}(t) \right|
$$
\begin{equation}\label{abel_for_th_3}
+\left|\int\limits_0^\infty (K(t)-K_1(t)) d\nu(t) \right|+\left|
\int\limits_0^\infty K_1(t)d\mu_{r_n}(t)- \int\limits_0^\infty
K_1(t) d\nu(t) \right|.
\end{equation}
╚ьххь
$$
\left|\int\limits_0^\infty \left(K(t)-K_1(t)\right)
d\mu_{r_n}(t)\right| \leq M\int\limits_0^\infty
\left|K(t)-K_1(t)\right| \frac{V(r_n t)}{tV(r_n)}dt
$$
$$
\leq 2M \int\limits_0^\infty\left|K(t)-K_1(t)\right| t^{\rho-1}dt
\leq 2M\max\{t^{\rho-1}:\: t\in[a,b]\}\varepsilon,
$$
уфх $[a,b]\subset(0,\infty)$ Єр\-ъющ ёху\-ьхэЄ, ўЄю $\supp K$,
$\supp K_1 \subset[a,b]$. ┬ЄюЁюх ёырурх\-ьюх т яЁртющ ўрёЄш
эх\-Ёртхэ\-ёЄ\-тр  (\ref{abel_for_th_3}) фю\-яєёър\-хЄ Єръє■ цх
юЎхэ\-ъє. ╥ЁхЄ№х ёырур\-х\-ьюх ёЄЁх\-ьшЄ\-ё  ъ эєы■ яЁш
$n\to\infty$. ╚ч ёърчрэ\-эю\-ую ёых\-фє\-хЄ, ўЄю
$$
\lim\limits_{n\to\infty} J(r_n)=\int\limits_0^\infty K(u)d\nu(u).
$$
╥хь ёрь√ь фю\-ърчр\-эю, ўЄю ёяЁртхф\-ыштю тъы■ўх\-эшх $H\subset
L(J,\infty)$, уфх  $H$  -- яЁртр  ўрёЄ№ Ёртхэ\-ёЄ\-тр
(\ref{abel_L(I)}). ─юърчр\-Єхы№\-ёЄ\-тю чр\-ърэўш\-тр\-хЄ\-ё 
Ёрё\-ёєц\-фх\-эшхь рэрыю\-ушў\-э√ь Єюьє, ъюЄю\-Ёюх яЁш\-тхфх\-эю яЁш
фю\-ърчр\-Єхы№\-ёЄ\-тх Єхю\-Ёх\-ь√  \ref{abel_th_1}.

\par ─рыхх т єёыю\-тш ї Єхю\-Ёхь  \ref{abel_th_1}-\ref{abel_th_3}
ь√ чр\-ьхэшь ЄЁхсютр\-эшх ЇшэшЄ\-эюё\-Єш  фЁр  $K$ сюыхх ёырс√ь
юуЁрэш\-ўх\-эшхь эр ¤Єю  фЁю ё ёю\-їЁрэх\-эшхь єЄ\-тхЁ\-ц\-фх\-эшщ
¤Єшї Єхю\-Ёхь. ╧Ёш ¤Єюь єёшыш\-тр\-хЄё  юуЁрэш\-ўх\-эшх эр ьхЁє
$\mu$.

\par ═рўэ╕ь ё эхюс\-їюфш\-ь√ї юяЁх\-фхых\-эшщ.

\par ╠√ сєфхь уютюЁшЄ№, ўЄю ЄЁющър $(K,\rho(r),\mu)$ єфют\-ыхЄ\-тю\-Ё \-хЄ
 {\it єёыю\-тш■ эхщ\-ЄЁр\-ыш\-чр\-Ўшш эєы }, хёыш
$$
\lim\limits_{\varepsilon\to +0}\mathop{\overline{\lim}}_{r\to\infty}
\left|\int\limits_0^\varepsilon K(t)d\mu_r(t)\right|=0.
$$
\par ╠√ сєфхь уютюЁшЄ№, ўЄю ЄЁющър $(K,\rho(r),\mu)$ єфют\-ыхЄ\-тю\-Ё \-хЄ
{\it єёыю\-тш■ эхщ\-ЄЁр\-ыш\-чр\-Ўшш схё\-ъю\-эхў\-эюё\-Єш}, хёыш
$$
\lim\limits_{N\to\infty}\mathop{\overline{\lim}}_{r\to\infty}
\left|\int\limits_N^\infty K(t)d\mu_r(t)\right|=0.
$$
\par ╟рьхЄшь, ўЄю хёыш  $K$  -- ЇшэшЄ\-эюх  фЁю эр яюыє\-юёш
$(0,\infty)$, Єю ЄЁющър $(K,\rho(r),\mu)$ єфютыхЄ\-тю\-Ё хЄ
єёыю\-тш ь эхщ\-ЄЁр\-ыш\-чр\-Ўшш эєы  ш схё\-ъю\-эхў\-эюё\-Єш яЁш
ы■сюь єЄюў\-э╕э\-эюь яю\-Ё фъх $\rho(r)$ ш ы■сющ Ёр\-фю\-эю\-тющ
ьхЁх $\mu$ эр яюыє\-юёш  $(0,\infty)$.

\par ┬тхф╕ээ√х юяЁх\-фхых\-эш  фр■Є тючьюц\-эюёЄ№ єяЁюё\-ЄшЄ№
ЇюЁ\-ьє\-ыш\-Ёют\-ъш Ё фр Єхю\-Ёхь, яюч\-тюы   эх
фхЄрыш\-чш\-Ёю\-трЄ№ юуЁрэш\-ўх\-эш  эр  фЁю  $K$ ш ьхЁє  $\mu$,
ъюЄю\-Ё√х юсхё\-ях\-ўш\-тр\-■Є т√\-яюы\-эш\-ьюёЄ№ єёыю\-тшщ
эхщ\-ЄЁр\-ыш\-чр\-Ўшш эєы  ш схё\-ъю\-эхў\-эюё\-Єш. ─хЄрыш\-чр\-Ўш■
юуЁрэш\-ўх\-эшщ эр  $K$ ш  $\mu$ ьюц\-эю юЇюЁь\-ы Є№ ъръ
ёрью\-ёЄю \-Єхы№\-э√х єЄ\-тхЁ\-ц\-фх\-эш .

%%%%%%%%%%%%%%%%%%%%%%%%%%%%%%%%%%%%%%%%%%%%%%%%%%%%%%%%%%%%%%%%%%%%%%%%%%%%%%%%%%%%%%%
\begin{lemma}\label{abel_lemma_1}\hskip-2mm{.}\:
%%%%%%%%%%%%%%%%%%%%%%%%%%%%%%%%%%%%%%%%%%%%%%%%%%%%%%%%%%%%%%%%%%%%%%%%%%%%%%%%%%%%%%%
┼ёыш ЄЁющър $(K,\rho(r),\mu)$ єфют\-ыхЄ\-тю\-Ё хЄ єёыю\-тш■
эхщ\-ЄЁр\-ыш\-чр\-Ўшш эєы , Єю
$$
\lim\limits_{{\varepsilon_1\to +0}\atop{\varepsilon_2\to +0}}
\mathop{\overline{\lim}}_{r\to\infty}
\left|\:\int\limits_{\varepsilon_1}^{\varepsilon_2}
K(t)d\mu_r(t)\right|=0.
$$
\end{lemma}
%%%%%%%%%%%%%%%%%%%%%%%%%%%%%%%%%%%%%%%%%%%%%%%%%%%%%%%%%%%%%%%%%%%%%%%%%%%%%%%%%%%%%%%
\begin{lemma}\label{abel_lemma_2}\hskip-2mm{.}\:
%%%%%%%%%%%%%%%%%%%%%%%%%%%%%%%%%%%%%%%%%%%%%%%%%%%%%%%%%%%%%%%%%%%%%%%%%%%%%%%%%%%%%%%
┼ёыш ЄЁющър $(K,\rho(r),\mu)$ єфют\-ыхЄ\-тю\-Ё хЄ єёыю\-тш■
эхщ\-ЄЁр\-ыш\-чр\-Ўшш схё\-ъю\-эхў\-эюё\-Єш, Єю
$$
\lim\limits_{{N_1\to \infty}\atop{N_2\to \infty}}
\mathop{\overline{\lim}}_{r\to\infty}
\left|\,\int\limits_{N_1}^{N_2} K(t)d\mu_r(t)\right|=0.
$$
\end{lemma}
╦хьь√ \ref{abel_lemma_1} ш \ref{abel_lemma_2} юўхтшфэ√.

%%%%%%%%%%%%%%%%%%%%%%%%%%%%%%%%%%%%%%%%%%%%%%%%%%%%%%%%%%%%%%%%%%%%%%%%%%%%%%%%%%%%%%%%%%%%%%%%%%%
\begin{theorem}\label{abel_th_4}\hskip-2mm{.}\:
%%%%%%%%%%%%%%%%%%%%%%%%%%%%%%%%%%%%%%%%%%%%%%%%%%%%%%%%%%%%%%%%%%%%%%%%%%%%%%%%%%%%%%%%%%%%%%%%%%%
╧єёЄ№ $\rho(r)$  -- яЁю\-шч\-тюы№\-э√щ єЄюў\-э╕э\-э√щ яю\-Ё фюъ,
$\mu\in$ $\frak{M}(\rho(r))$, $K$  -- эх\-яЁх\-Ё√т\-эюх  фЁю эр
яюыє\-юёш $(0,\infty)$. ┼ёыш ЄЁющ\-ър $(K,\rho(r),\mu)$
єфют\-ыхЄ\-тю\-Ё хЄ єёыю\-тш\- ь эхщ\-ЄЁр\-ыш\-чр\-Ўшш эєы  ш
схё\-ъю\-эхў\-эюё\-Єш, Єю фы  ы■сющ ьхЁ√ $\nu\in Fr[\mu]$ шэ\-ЄхуЁры
$\int\limits_0^\infty K(t)d\nu(t)$ ёє\-∙хёЄ\-тє\-хЄ ъръ
эх\-ёюсёЄ\-тхэ\-э√щ шэ\-ЄхуЁры ё юёюс√\-ьш Єюўър\-ьш  $0$ ш
$\infty$.
\end{theorem}
%%%%%%%%%%%%%%%%%%%%%%%%%%%%%%%%%%%%%%%%%%%%%%%%%%%%%%%%%%%%%%%%%%%%%%%%%%%%%%%%%%%%%%%%%%%%%%%%%%%%
\par {\sc  ─юърчрЄхы№ёЄтю.}\: ╧єёЄ№ $\nu$  -- яЁю\-шч\-тюы№\-эр  ьхЁр шч
яЁх\-фхы№\-эю\-ую ьэю\-цх\-ёЄ\-тр $Fr[\mu]$. ╤є∙хёЄ\-тєхЄ
яю\-ёыхфю\-тр\-Єхы№\-эюёЄ№ $r_n\to\infty$ Єр\-ър , ўЄю $\mu_{r_n}\to
\nu$. ─юяюыэш\-Єхы№\-эю ьюц\-эю ёўш\-ЄрЄ№, ўЄю
$|\mu|_{r_n}\to\hat{\nu}$. ╥ръ ъръ ЄЁющър $(K,\rho(r),\mu)$
єфют\-ыхЄ\-тю\-Ё \-хЄ єёыю\-тш ь эхщ\-ЄЁр\-ыш\-чр\-Ўшш эєы , Єю шч
ыхьь√ \ref{abel_lemma_1} ёых\-фє\-хЄ, ўЄю фы  ы■сюую $\delta>0$
ёє\-∙хёЄ\-тєхЄ ўшёыю $\varepsilon_0>0$ Єръюх, ўЄю хёыш
$0<\varepsilon_1<\varepsilon_2\leq\varepsilon_0$, Єю сєфхЄ
т√\-яюыэ Є№\-ё  эх\-Ёртхэ\-ёЄ\-тю
$$
\mathop{\overline{\lim}}\limits_{r\to\infty}
\left|\,\int\limits_{\varepsilon_1}^{\varepsilon_2}
K(t)d\mu_r(t)\right|\leq \delta.
$$
┬ ўрёЄэюё\-Єш,
\begin{equation}\label{abel_lim_int_K}
\mathop{\overline{\lim}}\limits_{n\to\infty}
\left|\,\int\limits_{\varepsilon_1}^{\varepsilon_2}
K(t)d\mu_{r_n}(t)\right|\leq \delta.
\end{equation}
╧єёЄ№ $E$  -- ьэю\-цх\-ёЄ\-тю Єюўхъ, ъюЄю\-Ё√х эр\-уЁєцр\-хЄ ьхЁр
$\hat{\nu}$. ╠эю\-цх\-ёЄ\-тю  $E$  -- ¤Єю эх сюыхх ўхь ёў╕Є\-эюх
ьэю\-цх\-ёЄ\-тю. ─юяюыэш\-Єхы№\-эю яЁхф\-яюыю\-цшь, ўЄю
$\varepsilon_1,\varepsilon_2$ $\notin E$. ╚ёяюы№\-чє  Єхю\-Ёх\-ьє
\ref{measures_th_2}, эх\-Ёртхэ\-ёЄ\-тю (\ref{abel_lim_int_K})
ьюц\-эю яхЁх\-яшёрЄ№ т тшфх
\begin{equation}\label{abel_int_K}
\left|\,\int\limits_{\varepsilon_1}^{\varepsilon_2}
K(t)d\nu(t)\right|\leq \delta.
\end{equation}
┬  ¤Єюь эх\-Ёртхэ\-ёЄтх яЁхф\-яюырур\-хЄё , ўЄю ўшёыр
$\varepsilon_1$, $\varepsilon_2$ єфютыхЄ\-тю\-Ё ■Є юуЁрэш\-ўх\-эш■
$\varepsilon_1,\varepsilon_2$ $\notin E$. ╬фэръю шч Ёртхэ\-ёЄ\-тр
$$
\left|\:\int\limits_{\varepsilon_1}^{\varepsilon_2}
K(t)d\nu(t)\right|= \lim\limits_{h\to +0} \left|\:
\int\limits_{\varepsilon_1-h}^{\varepsilon_2+h} K(t)d\nu(t)\right|
$$
ёых\-фєхЄ, ўЄю ¤Єю юуЁрэш\-ўх\-эшх ьюц\-эю юЄ\-сЁюёшЄ№. ╥юуфр
эх\-Ёртхэ\-ёЄ\-тю  (\ref{abel_int_K}) ючэрўр\-хЄ, ўЄю фы 
эх\-ёюс\-ёЄ\-тхэ\-эю\-ую шэ\-ЄхуЁр\-ыр $\int\limits_0^\infty
K(t)d\nu(t)$ т√\-яюы\-э \-хЄ\-ё  єёыю\-тшх ╩ю°ш ёїюфш\-ьюё\-Єш
шэ\-Єху\-Ёрыр т юёюсющ Єюўъх $0$. ╤ыхфютр\-Єхы№\-эю, шэ\-Єху\-Ёры
ёїю\-фшЄ\-ё  т эєых. └эрыю\-ушў\-эю фю\-ърч√\-тр\-хЄ\-ё 
ёїюфш\-ьюёЄ№ ¤Єюую шэ\-Єху\-Ёрыр т схё\-ъю\-эхў\-эюё\-Єш.
╥хю\-Ёх\-ьр фю\-ърчр\-эр.

%%%%%%%%%%%%%%%%%%%%%%%%%%%%%%%%%%%%%%%%%%%%%%%%%%%%%%%%%%%%%%%%%%%%%%%%%%%%%%%%%%%%%%%
\par ╥хяхЁ№ ь√ ьюцхь фю\-ърчрЄ№ трЁш\-рэЄ Єхю\-Ёх\-ь√  \ref{abel_th_1},
т ъюЄю\-Ёюь ЄЁхсю\-тр\-эшх ЇшэшЄ\-эюё\-Єш  фЁр  $K$ юЄёєЄ\-ёЄ\-тєхЄ.
%%%%%%%%%%%%%%%%%%%%%%%%%%%%%%%%%%%%%%%%%%%%%%%%%%%%%%%%%%%%%%%%%%%%%%%%%%%%%%%%%%%%%%%%%%%%%%%%%%%
\begin{theorem}\label{abel_th_5}\hskip-2mm{.}\:
%%%%%%%%%%%%%%%%%%%%%%%%%%%%%%%%%%%%%%%%%%%%%%%%%%%%%%%%%%%%%%%%%%%%%%%%%%%%%%%%%%%%%%%%%%%%%%%%%%%
╧єёЄ№ $\rho(r)$  -- яЁю\-шч\-тюы№\-э√щ єЄюў\-э╕э\-э√щ яю\-Ё фюъ,
$\mu\in$ $\frak{M}_\infty$ $(\rho(r))$, $K$  -- эх\-яЁх\-Ё√т\-эюх
 фЁю эр яюыє\-юёш $(0,\infty)$. ┼ёыш ЄЁющър $(K,\rho(r),\mu)$
єфют\-ыхЄ\-тю\-Ё хЄ єёыю\-тш\- ь эхщ\-ЄЁр\-ыш\-чр\-Ўшш эєы  ш
схё\-ъю\-эхў\-эюё\-Єш, Єю
\begin{equation}\label{abel_L(J,infty)}
L(J,\infty)=\left\{\int\limits_0^\infty K(u)d\nu(u):\: \nu\in
Fr[\mu]\right\}.
\end{equation}
\end{theorem}
%%%%%%%%%%%%%%%%%%%%%%%%%%%%%%%%%%%%%%%%%%%%%%%%%%%%%%%%%%%%%%%%%%%%%%%%%%%%%%%%%%%%%%%%%%%%%%%%%%%%
\par {\sc  ─юърчрЄхы№ёЄтю.}\: ╧єёЄ№ $\nu$  -- яЁю\-шч\-тюы№\-эр  ьхЁр шч
яЁх\-фхы№\-эю\-ую ьэю\-цх\-ёЄ\-тр  $Fr[\mu]$. ╤є∙хёЄ\-тєхЄ
яю\-ёыхфю\-тр\-Єхы№\-эюёЄ№ $r_n\to\infty$ Єр\-ър , ўЄю
$\mu_{r_n}\to\nu$, $|\mu|_{r_n}\to\hat{\nu}$. ╧єёЄ№  $\delta$  --
яЁю\-шч\-тюы№\-эюх ёЄЁюую яюыю\-цш\-Єхы№\-эюх ўшёыю. ╚ч Єюую, ўЄю
ЄЁющър $(K,\rho(r),\mu)$ єфют\-ыхЄ\-тю\-Ё \-хЄ єёыю\-тш\- ь
эхщ\-ЄЁр\-ыш\-чр\-Ўшш эєы  ш схё\-ъю\-эхў\-эюё\-Єш ш шч Єхю\-Ёх\-ь√
\ref{abel_th_4} ёых\-фє\-хЄ, ўЄю ёє\-∙хёЄ\-тє■Є ўшёыр
$\varepsilon_0>0$ ш  $N_0>0$ Єръшх, ўЄю фы 
$\varepsilon\in(0,\varepsilon_0)$ ш $N>N_0$ сєфєЄ т√\-яюыэ Є№\-ё 
эх\-Ёртхэ\-ёЄ\-тр

\begin{equation}{\label{abel_for_th_5}}
\begin{split}
\mathop{\overline{\lim}} \limits_{r\to\infty}
\left|\int\limits_{0}^{\varepsilon} K(t)d\mu_r(t)\right|<
\delta,&\qquad \mathop{\overline{\lim}} \limits_{r\to\infty}
\left|\int\limits_{N}^{\infty} K(t)d\mu_r(t)\right|< \delta,\\
\left|\int\limits_{0}^{\varepsilon} K(t)d\nu(t)\right| < \delta,
&\qquad \left|\int\limits_{N}^{\infty} K(t)d\nu(t)\right| < \delta.
\end{split}
\end{equation}

\par ┴єфхь ёўш\-ЄрЄ№, ўЄю Єюўъш $\varepsilon$ ш  $N$
єфют\-ыхЄ\-тю\-Ё \-■Є т√\-яш\-ёрэ\-э√ь фы  эшї ёююЄ\-эю°х\-эш ь ш,
ъЁюьх Єюую, эх эр\-уЁєцр\-■Єё  ьхЁющ  $\hat{\nu}$. ╥юуфр
$$
\left|\int\limits_0^{\infty} K(t)d\mu_{r_n}(t)-
\int\limits_0^{\infty} K(t)d\nu(t)\right|\leq
\left|\int\limits_0^{\varepsilon} K(t)d\mu_{r_n}(t)\right| +
\left|\int\limits_0^{\varepsilon} K(t)d\nu(t)\right|
$$
\begin{equation}\label{abel_ineq_th_5}
+\left|\int\limits_N^\infty K(t)d\mu_{r_n}(t)\right| +
\left|\int\limits_N^\infty K(t)d\nu(t)\right|+
\left|\int\limits_{\varepsilon}^N K(t)d\mu_{r_n}(t)-
\int\limits_{\varepsilon}^N K(t)d\nu(t)\right| .
\end{equation}
╧ю Єхю\-Ёхьх  \ref{measures_th_2} яю\-ёыхф\-эхх шч т√\-яшёрэ\-э√ї
ёырурх\-ь√ї ёЄЁхьшЄ\-ё  ъ эєы■ яЁш $n\to\infty$. ╚ч ¤Єюую, р Єръ\-цх
шч эх\-ЁртхэёЄт  (\ref{abel_for_th_5}) ш  (\ref{abel_ineq_th_5})
ёых\-фє\-хЄ, ўЄю
$$
\mathop{\overline{\lim}} \limits_{n\to\infty}
\left|\int\limits_0^\infty K(t)d\mu_{r_n}(t)- \int\limits_0^{\infty}
K(t)d\nu(t)\right|\leq 4\delta.
$$
╚ч ¤Єюую, т ётю■ юўхЁхф№ ёых\-фє\-хЄ, ўЄю $H\subset L(J,\infty)$,
уфх $H$  -- яЁртр  ўрёЄ№ Ёртхэ\-ёЄ\-тр (\ref{abel_L(J,infty)}).
╬ъюэўр\-эшх фю\-ърчр\-Єхы№\-ёЄ\-тр Єръюх цх, ъръ ш т Єхю\-Ёх\-ьх
\ref{abel_th_1}. ╥хю\-Ёх\-ьр фю\-ърчр\-эр.

\par ─рыхх ь√ яЁш\-тюфшь трЁшрэ\-Є√ Єхю\-Ёхь  \ref{abel_th_2} ш
\ref{abel_th_3}, т ъюЄю\-Ё√ї ЄЁхсю\-тр\-эшх ЇшэшЄ\-эюё\-Єш  фЁр $K$
юЄ\-ёєЄ\-ёЄ\-тє\-хЄ. ╬эш фю\-ърч√\-тр\-■Є\-ё  ё яюью∙№■
Ёрё\-ёєцфх\-эшщ рэрыю\-ушў\-э√ї Єхь, ъюЄю\-Ё√х яЁш\-ьхэ \-ышё№ яЁш
фю\-ърчр\-Єхы№\-ёЄ\-тх ЄхюЁх\-ь√ \ref{abel_th_5}.

%%%%%%%%%%%%%%%%%%%%%%%%%%%%%%%%%%%%%%%%%%%%%%%%%%%%%%%%%%%%%%%%%%%%%%%%%%%%%%%%%%%%%%%%%%%%%%%%%%%
\begin{theorem}\label{abel_th_6}\hskip-2mm{.}\:
%%%%%%%%%%%%%%%%%%%%%%%%%%%%%%%%%%%%%%%%%%%%%%%%%%%%%%%%%%%%%%%%%%%%%%%%%%%%%%%%%%%%%%%%%%%%%%%%%%%
╧єёЄ№ $\rho(r)$  -- яЁю\-шч\-тюы№\-э√щ єЄюў\-э╕э\-э√щ яю\-Ё фюъ,
$\mu\in$ $\frak{M}_\infty$ $(\rho(r))$. ╧єёЄ№ $K(t)$  --  фЁю эр
яюыє\-юёш $(0,\infty)$ эх\-яЁх\-Ё√т\-эюх тё■фє ъЁюьх Єюўъш  $1$,
ъюЄю\-Ёр   ты \-хЄё  Єюўъющ Ёрч\-Ё√тр яхЁтю\-ую Ёюфр, яЁш\-ў╕ь
Їєэъ\-Ўш  $K(t)$ эх\-яЁх\-Ё√т\-эр ёыхтр т Єюўъх  $1$. ╧єёЄ№ ЄЁющър
$(K,\rho(r),\mu)$ єфют\-ыхЄ\-тю\-Ё \-хЄ єёыю\-тш\- ь
эхщ\-ЄЁр\-ыш\-чр\-Ўшш эєы  ш схё\-ъю\-эхў\-эюё\-Єш. ╥юуфр
$$
L(J,\infty)=\left\{\int\limits_0^1 K(t)d\nu_1(t)+
\int\limits_1^\infty \widetilde{K}(t)d\nu_2(t):\: (\nu_1,\nu_2)\in
\widehat{Fr}[\mu]\right\},
$$
уфх $\widetilde{K}(t)$  -- эх\-яЁх\-Ё√т\-эюх яЁю\-фюы\-цх\-эшх  фЁр
ё яюыє\-юёш $(1,\infty)$ эр яюыє\-юё№ $[1,\infty)$.
\end{theorem}
%%%%%%%%%%%%%%%%%%%%%%%%%%%%%%%%%%%%%%%%%%%%%%%%%%%%%%%%%%%%%%%%%%%%%%%%%%%%%%%%%%%%%%%%%%%%%%%%%%%%
\begin{theorem}\label{abel_th_7}\hskip-2mm{.}\:
%%%%%%%%%%%%%%%%%%%%%%%%%%%%%%%%%%%%%%%%%%%%%%%%%%%%%%%%%%%%%%%%%%%%%%%%%%%%%%%%%%%%%%%%%%%%%%%%%%%
╧єёЄ№ $\rho(r)$  -- яЁю\-шч\-тюы№\-э√щ єЄюў\-э╕э\-э√щ яю\-Ё фюъ,
$\mu$  -- Ёр\-фю\-эютр ьхЁр эр яюыє\-юёш $(0,\infty)$ ё яыюЄ\-эюёЄ№■
$\mu'(r)$, єфют\-ыхЄ\-тюЁ ■\-∙хщ эх\-Ёртхэ\-ёЄ\-тє $|\mu'(r)|$
$\leq$ $M\frac{V(r)}{r}$. ╧єёЄ№ $K(t)$  -- ыю\-ъры№\-эю
шэ\-Єху\-Ёш\-Ёєх\-ьюх  фЁю эр яюыє\-юёш $(0,\infty)$ ш яєёЄ№ ЄЁющър
$(K,\rho(r),\mu)$ єфют\-ыхЄ\-тю\-Ё хЄ єёыю\-тш\- ь
эхщ\-ЄЁр\-ыш\-чр\-Ўшш эєы  ш схё\-ъю\-эхў\-эюё\-Єш. ╥юуфр
$$
L(J,\infty)=\left\{\int\limits_0^\infty K(t)d\nu(t):\: \nu\in
Fr[\mu]\right\}.
$$
\end{theorem}
%%%%%%%%%%%%%%%%%%%%%%%%%%%%%%%%%%%%%%%%%%%%%%%%%%%%%%%%%%%%%%%%%%%%%%%%%%%%%%%%%%%%%%%%%%%%%%%%%%%%
\par ┬ ёт чш ё Єхю\-Ёхьр\-ьш  \ref{abel_th_5}-\ref{abel_th_7}
трцэюх чэрўх\-эшх яЁш\-юсЁх\-ЄрхЄ тюяЁюё юс юуЁрэш\-ўх\-эш\- ї эр
 фЁю  $K$ ш ьхЁє  $\mu$, ъюЄю\-Ё√х с√ урЁрэ\-Єш\-Ёю\-трыш, ўЄю
ЄЁющър $(K,\rho(r),\mu)$ єфют\-ыхЄ\-тю\-Ё хЄ єёыю\-тш\- ь
эхщ\-ЄЁр\-ыш\-чр\-Ўшш эєы  ш схё\-ъю\-эхў\-эюё\-Єш. ╠√ фю\-ърцхь фтр
Ёх\-чєы№\-ЄрЄр эр ¤Єє Єхьє. ┬ ёт чш ё эшцх\-ёых\-фє\-■\-∙шьш
єЄ\-тхЁ\-ц\-фх\-эш \-ьш эр\-яюь\-эшь, ўЄю Їєэъ\-Ўш  $\gamma(t)$
юяЁх\-фхых\-эр Ёртхэ\-ёЄ\-тюь (\ref{def_gamma}).

%%%%%%%%%%%%%%%%%%%%%%%%%%%%%%%%%%%%%%%%%%%%%%%%%%%%%%%%%%%%%%%%%%%%%%%%%%%%%%%%%%%%%%%%%%%%%%%%%%%
\begin{lemma}\label{abel_lemma_3}\hskip-2mm{.}\:
%%%%%%%%%%%%%%%%%%%%%%%%%%%%%%%%%%%%%%%%%%%%%%%%%%%%%%%%%%%%%%%%%%%%%%%%%%%%%%%%%%%%%%%%%%%%%%%%%%%
╧єёЄ№ $\rho(r)$  -- яЁю\-шч\-тюы№\-э√щ єЄюў\-э╕э\-э√щ яю\-Ё фюъ,
$\mu$  -- Ёр\-фю\-эютр ьхЁр эр яюыє\-юёш $(0,\infty)$ ё
яыюЄ\-эюё\-Є№■ $\mu'(r)$, єфют\-ыхЄ\-тю\-Ё ■\-∙хщ эх\-Ёртхэ\-ёЄ\-тє
$|\mu'(r)|$ $\leq$ $M\frac{V(r)}{r}$ $(r\in(0,\infty))$. ╧єёЄ№
$t^{\rho-1}\gamma(t)K(t)$ $\in$ $L_1(0,\infty)$. ╥юуфр ЄЁющър
$(K,\rho(r),\mu)$ єфют\-ыхЄ\-тю\-Ё хЄ єёыю\-тш\- ь
эхщ\-ЄЁр\-ыш\-чр\-Ўшш эєы  ш схё\-ъю\-эхў\-эюё\-Єш.
\end{lemma}
%%%%%%%%%%%%%%%%%%%%%%%%%%%%%%%%%%%%%%%%%%%%%%%%%%%%%%%%%%%%%%%%%%%%%%%%%%%%%%%%%%%%%%%%%%%%%%%%%%%%
\par \quad {\sc ─юърчр\-Єхы№ёЄтю.}\: ╚ёяюы№чє 
эх\-Ёртхэ\-ёЄтю  (\ref{order_V_V}), эр\-їюфшь, ўЄю
$$
\left|\int\limits_0^\varepsilon K(t)d\mu_r(t)\right| \leq
M\int\limits_0^\varepsilon |K(t)|\frac{V(rt)}{tV(r)}dt\leq M
\int\limits_0^\varepsilon t^{\rho-1}\gamma(t)|K(t)|dt.
$$
╚ч ¤Єюую эхЁртхэёЄтр ш єёыю\-тш  $t^{\rho-1}\gamma(t)|K(t)|\in
L_1(0,\infty)$ ёых\-фєхЄ, ўЄю ЄЁющър $(K,\rho(r),\mu)$
єфют\-ыхЄ\-тю\-Ё хЄ єёыю\-тш■ эхщ\-ЄЁр\-ыш\-чр\-Ўшш эєы .
└эрыю\-ушў\-эю фю\-ърч√\-тр\-хЄё , ўЄю ¤Єр ЄЁющър $(K,\rho(r),\mu)$
єфют\-ыхЄ\-тю\-Ё хЄ єёыю\-тш■ эхщ\-ЄЁр\-ыш\-чр\-Ўшш
схё\-ъю\-эхў\-эюё\-Єш. ╦хььр фю\-ърчр\-эр.
%%%%%%%%%%%%%%%%%%%%%%%%%%%%%%%%%%%%%%%%%%%%%%%%%%%%%%%%%%%%%%%%%%%%%%%%%%%%%%%%%%%%%%%%%%%%%%%%%%%
\begin{remark}\label{abel_for_lemma_3}\hskip-2mm{.}\:
%%%%%%%%%%%%%%%%%%%%%%%%%%%%%%%%%%%%%%%%%%%%%%%%%%%%%%%%%%%%%%%%%%%%%%%%%%%%%%%%%%%%%%%%%%%%%%%%%%%
┬ ЇюЁ\-ьє\-ыш\-Ёютъх ыхьь√  \ref{abel_lemma_3} єўрёЄ\-тєхЄ Їєэъ\-Ўш 
$\gamma(t)$, шё\-ёыхфю\-тр\-эшх ъю\-Єю\-Ёющ яЁхф\-ёЄрт\-ы хЄ
юяЁх\-фх\-ы╕э\-э√х ЄЁєф\-эюё\-Єш. ▌Єш шё\-ёыхфю\-тр\-эш 
ёЄрэю\-т Є\-ё  шч\-ыш°\-эшьш, хёыш  фЁю  $K$ єфютыхЄ\-тю\-Ё хЄ
єёыю\-тш■ $t^{\rho-1}
\frac{1+t^{2\varepsilon}}{t^\varepsilon}K(t)\in L_1(0,\infty)$ фы 
эх\-ъюЄю\-Ёю\-ую  $\varepsilon>0$. ╚ч ¤Єюую ёююЄ\-эю°х\-эш  ш
Єхю\-Ёх\-ь√  \ref{order_th_ln_gamma} т√\-Єхър\-хЄ, ўЄю
$t^{\rho-1}\gamma(t)K(t)\in L_1(0,\infty)$.
\end{remark}
%%%%%%%%%%%%%%%%%%%%%%%%%%%%%%%%%%%%%%%%%%%%%%%%%%%%%%%%%%%%%%%%%%%%%%%%%%%%%%%%%%%%%%%%%%%%%%%%%%%
\par ╤ыхфє\-■\-∙шщ Ёх\-чєы№ЄрЄ тюёїю\-фшЄ ъ ┬шэхЁє.

\par ╧єёЄ№ $K(t)$  --  фЁю эр яюыє\-юёш $(0,\infty)$. ╬сючэрўшь
$$
K_n=\sup\{|K(t)|:\:t\in(e^n,e^{n+1}]\},\quad n\in(-\infty,\infty).
$$

\par ╠эю\-цх\-ёЄтю $\frak{M}(\rho(r))$, єўрёЄ\-тє■\-∙шх т ЇюЁ\-ьє\-ыш\-Ёют\-ъх
ёых\-фє\-■\-∙хщ ыхьь√, юяЁх\-фхых\-эю тю ттхфх\-эшш.

%%%%%%%%%%%%%%%%%%%%%%%%%%%%%%%%%%%%%%%%%%%%%%%%%%%%%%%%%%%%%%%%%%%%%%%%%%%%%%%%%%%%%%%%%%%%%%%%%%%
\begin{lemma}\label{abel_lemma_4}\hskip-2mm{.}\:
%%%%%%%%%%%%%%%%%%%%%%%%%%%%%%%%%%%%%%%%%%%%%%%%%%%%%%%%%%%%%%%%%%%%%%%%%%%%%%%%%%%%%%%%%%%%%%%%%%%
╧єёЄ№ $\rho(r)$  -- яЁю\-шч\-тюы№\-э√щ єЄюў\-э╕э\-э√щ яю\-Ё фюъ,
$\mu\in$ $\frak{M}(\rho(r))$. ╧єёЄ№ $K(t)$  -- сюЁх\-ыхт\-ёър 
Їєэъ\-Ўш  эр яюыє\-юёш $(0,\infty)$ Єр\-ър , ўЄю ёїю\-фшЄ\-ё  Ё ф
\begin{equation}\label{abel_sum}
\sum\limits_{n=-\infty}^\infty e^{n\rho}\gamma(e^n)K_n, \quad
\rho=\rho(\infty).
\end{equation}
╥юуфр ЄЁющър $(K,\rho(r),\mu)$ єфют\-ыхЄ\-тю\-Ё хЄ єёыю\-тш\- ь
эхщ\-ЄЁр\-ыш\-чр\-Ўшш эєы  ш схё\-ъю\-эхў\-эюё\-Єш.
\end{lemma}
%%%%%%%%%%%%%%%%%%%%%%%%%%%%%%%%%%%%%%%%%%%%%%%%%%%%%%%%%%%%%%%%%%%%%%%%%%%%%%%%%%%%%%%%%%%%%%%%%%%%
\par \quad {\sc ─юърчр\-Єхы№ёЄтю.}\: ╚ч єёыю\-тш  $\mu\in$ $\frak{M}(\rho(r))$
ёых\-фє\-хЄ, ўЄю ёє\-∙хёЄ\-тє\-хЄ Єръюх $A>0$, ўЄю фы  ы■сю\-ую
$r>0$ т√\-яюыэ \-хЄё  эх\-Ёртхэ\-ёЄ\-тю $|\mu|((r,er])\leq AV(r)$.
─рыхх эрїюфшь
$$
\left|\frac{1}{V(r)}\int\limits_0^{\varepsilon r}
K\left(\frac{t}{r}\right)d\mu(t)\right|\leq \frac{1}{V(r)}
\sum\limits_{n=-\infty}^{n_0} \int\limits_{e^nr}^{e^{n+1}r}
\left|K\left(\frac{u}{r}\right)\right|d|\mu|(u)
$$
$$
\leq \frac{1}{V(r)} \sum\limits_{n=-\infty}^{n_0}
K_n|\mu|((e^nr,e^{n+1}r]) \leq \frac{A}{V(r)}
\sum\limits_{n=-\infty}^{n_0} K_n V\left(e^n r\right) \leq
A\sum\limits_{n=-\infty}^{n_0} K_n e^{n\rho}\gamma(e^n),
$$
уфх $n_0=[\ln\varepsilon]$. ╚ч яюыєўхэ\-эю\-ую эх\-Ёртхэ\-ёЄ\-тр ш
ёїюфш\-ьюё\-Єш Ё фр  (\ref{abel_sum}) ёых\-фє\-хЄ, ўЄю ЄЁющър
$(K,\rho(r),\mu)$ єфют\-ыхЄ\-тю\-Ё \-хЄ єёыю\-тш■
эхщ\-ЄЁр\-ыш\-чр\-Ўшш эєы . └эрыюушў\-эю фю\-ърч√\-тр\-хЄ\-ё , ўЄю
¤Єр ЄЁющър єфют\-ыхЄ\-тю\-Ё хЄ єёыю\-тш■ эхщ\-ЄЁр\-ыш\-чр\-Ўшш
схё\-ъю\-эхў\-эюё\-Єш. ╦хььр фю\-ърчрэр.
%%%%%%%%%%%%%%%%%%%%%%%%%%%%%%%%%%%%%%%%%%%%%%%%%%%%%%%%%%%%%%%%%%%%%%%%%%%%%%%%%%%%%%%%%%%%%%%%%%%%
\begin{remark}\label{abel_for_lemma_4}\hskip-2mm{.}\:
%%%%%%%%%%%%%%%%%%%%%%%%%%%%%%%%%%%%%%%%%%%%%%%%%%%%%%%%%%%%%%%%%%%%%%%%%%%%%%%%%%%%%%%%%%%%%%%%%%%%
╙ёыютшх ёїюфш\-ьюёЄш Ё фр  (\ref{abel_sum})  ьюц\-эю чр\-ьхэшЄ№ эр
сюыхх ёшы№\-эюх юуЁрэш\-ўх\-эшх, яю\-ЄЁхсютрт, ўЄюс√ яЁш
эх\-ъюЄю\-Ёюь $\varepsilon>0$ ёїюфшы\-ё  Ё ф
$$
\sum\limits_{n=-\infty}^\infty \frac{1+e^{2\varepsilon
n}}{e^{\varepsilon n}} e^{n\rho}K_n.
$$
╧Ёш ¤Єюь эх тюч\-эшърхЄ эхюс\-їюфш\-ьюёЄ№ шё\-ёыхфю\-трЄ№ Їєэъ\-Ўш■
$\gamma(t)$.
\end{remark}
%%%%%%%%%%%%%%%%%%%%%%%%%%%%%%%%%%%%%%%%%%%%%%%%%%%%%%%%%%%%%%%%%%%%%%%%%%%%%%%%%%%%%%%%%%%%%%%%%%%%
\par ╨рё\-ёьюЄЁшь Їєэъ\-Ўш■
$$
v(z)=\int\limits_0^\infty \ln\left|1-\frac{z}{t}\right|d\mu(t),
$$
уфх $\mu$  -- яюыю\-цш\-Єхы№\-эр  ьхЁр эр яюыє\-юёш $(0,\infty)$ шч
ъырёёр $\frak{M}_\infty$ $(\rho(r))$, $\rho=\rho(\infty)\in(0,1)$, ш
Єр\-ър , ўЄю эр\-яшёрэ\-э√щ шэ\-ЄхуЁры ёїю\-фшЄ\-ё  т эєых. ▌Єр
Їєэъ\-Ўш  їюЁю°ю шчтхёЄ\-эр т ЄхюЁшш ЁюёЄр
ёєс\-урЁ\-ью\-эш\-ўхё\-ъшї Їєэъ\-Ўшщ ш  тшырё№ юс·хъ\-Єюь
ьэю\-ую\-ўшёыхэ\-э√ї шё\-ёыхфю\-тр\-эшщ. ┬ ёыєўрх, хёыш ьхЁр $\mu$
 ты \-хЄ\-ё  Ёхує\-ы Ё\-эющ, Їєэъ\-Ўш  $v(z)$ яЁш\-эрф\-ыхцшЄ
ёяхЎш\-ры№\-эюьє ъырёёє ёєс\-урЁью\-эш\-ўхёъшї Їєэъ\-Ўшщ тяюыэх
Ёхуєы Ё\-эю\-ую ЁюёЄр (юЄ\-эюёш\-Єхы№\-эю єЄюў\-э╕э\-эю\-ую
яю\-Ё ф\-ър  $\rho(r)$) т ёь√ёых ╦хтшэр-╧Їы■ухЁр. ┬ ¤Єюь ёыєўрх фы 
Їєэъ\-Ўшш $v(r)$ эх юс чр\-Єхы№\-эю ёє\-∙хёЄ\-тєхЄ яЁх\-фхы
ўрёЄ\-эю\-ую $v(r)/V(r)$ яЁш $r\to\infty$, эю яЁш ¤Єюь
ёє\-∙хёЄ\-тє\-хЄ яЁх\-фхы
\begin{equation}\label{abel_lim_v}
\lim\limits_{r\to\infty} \frac{1}{rV(r)} \int\limits_0^r v(t)dt.
\end{equation}
─юърчр\-Єхы№\-ёЄтр т√\-ёърчрэ\-э√ї чфхё№ єЄ\-тхЁ\-ц\-фх\-эшщ ю
Їєэъ\-Ўшш $v(z)$, ьюц\-эю эрщЄш т  \cite{Grishin_2005}.

\par ╙ЄтхЁц\-фх\-эшх ю ёє\-∙хёЄ\-тю\-тр\-эшш яЁх\-фхыр  (\ref{abel_lim_v})
фы  Їєэъ\-Ўшш $v(r)$ эх ёых\-фє\-хЄ шч фю\-ърчрэ\-э√ї эрьш Єхю\-Ёхь.
─рыхх сєфхЄ фю\-ърчр\-эр Єхю\-Ёх\-ьр, шч ъю\-Єю\-Ёющ сєфхЄ
ёыхфю\-трЄ№ ¤Єю єЄ\-тхЁ\-ц\-фх\-эшх. ╤ыхфю\-тр\-Єхы№\-эю, ь√
яю\-ыєўшь Ёх\-чєы№\-ЄрЄ ёЁртэш\-ь√щ яю ёшых ё Ёх\-чєы№\-Єр\-Єр\-ьш,
яю\-ыєўхэ\-э√ьш т ЄхюЁшш ёєс\-урЁью\-эш\-ўхё\-ъшї Їєэъ\-Ўшщ, эю
ёяЁртхф\-ыш\-т√ьш фы   фхЁ чэрўш\-Єхы№\-эю сюыхх юс∙шї, ўхь  фЁю
$\ln\left|1-\frac{r}{t}\right|$, ъюЄю\-Ёюх єўрёЄ\-тєхЄ т
юяЁх\-фхых\-эшш Їєэъ\-Ўшш $v(r)$. ╬ЄьхЄшь х∙╕, ўЄю, ъръ
яю\-ърч√\-тр■Є Єхю\-Ёх\-ьр  \ref{abel_th_2} ш
ёЇюЁ\-ьєыш\-Ёю\-трэ\-э√х ётющёЄ\-тр Їєэъ\-Ўшш $v(r)$, т ёыєўрх
Ёрч\-Ё√т\-э√ї  фхЁ яЁх\-фхы№\-эюх ьэю\-цх\-ёЄтю $Fr[\mu]$ ьхЁ√ $\mu$
эх юяЁх\-фхы \-хЄ рёшья\-Єю\-Єш\-ўхёъюх яю\-тхфх\-эшх Їєэъ\-Ўшш
$\Psi$. ╬фэръю, ъръ ь√ єтшфшь фрыхх, ¤Єю ьэю\-цх\-ёЄ\-тю
юяЁх\-фхы \-хЄ яЁх\-фхы№\-эюх ьэю\-цх\-ёЄтю ьхЁ√ $s$. ╩ръ єцх
юЄ\-ьхўр\-ыюё№ тю тёЄєя\-ых\-эшш, ЇюЁ\-ьєыш\-Ёєх\-ьр  эшцх Єхю\-Ёхьр
-- юфшэ шч уырт\-э√ї Ёх\-чєы№\-Єр\-Єют ЁрсюЄ√.
%%%%%%%%%%%%%%%%%%%%%%%%%%%%%%%%%%%%%%%%%%%%%%%%%%%%%%%%%%%%%%%%%%%%%%%%%%%%%%%%%%%%%%%%%%%%%%%%%%%
\begin{theorem}\label{abel_th_8}\hskip-2mm{.}\:
%%%%%%%%%%%%%%%%%%%%%%%%%%%%%%%%%%%%%%%%%%%%%%%%%%%%%%%%%%%%%%%%%%%%%%%%%%%%%%%%%%%%%%%%%%%%%%%%%%%
╧єёЄ№ $\rho(r)$  -- яЁю\-шч\-тюы№\-э√щ єЄюў\-э╕э\-э√щ яю\-Ё фюъ,
$\mu\in$ $\frak{M}(\rho(r))$.\\ ╧єёЄ№ $K(t)$  -- сюЁх\-ыхт\-ёър 
Їєэъ\-Ўш  эр яюыє\-юёш $(0,\infty)$ Єр\-ър , ўЄю
$t^{\rho-1}\gamma(t)K(t)\in L_1(0,\infty)$ $(\rho=\rho(\infty))$.
╥юуфр ьхЁр  $s$, $ds(u)=\Psi(u)du$, $\Psi$ -- Їєэъ\-Ўш ,
юяЁх\-фхы \-хьр  Ёртхэ\-ёЄ\-тюь  (\ref{introduction_Psi(r)}),
яЁш\-эрф\-ыхцшЄ ъырёёє $\frak{M}(\rho(r)+1)$ ш х╕ яЁх\-фхы№\-эюх
ьэю\-цх\-ёЄ\-тю $Fr[\rho(r)+1,s]$ ёю\-ёЄюшЄ шч рсёю\-ы■Є\-эю
эх\-яЁх\-Ё√т\-э√ї ьхЁ, ьэю\-цх\-ёЄ\-тю яыюЄ\-эюё\-Єхщ ъюЄю\-Ё√ї
ёют\-ярфр\-хЄ ё ьэю\-цх\-ёЄ\-тюь
$$
\left\{\int\limits_0^\infty K\left(\frac{t}{u}\right)d\nu(t):\quad
\nu\in Fr[\mu]\right\}.
$$
\end{theorem}
%%%%%%%%%%%%%%%%%%%%%%%%%%%%%%%%%%%%%%%%%%%%%%%%%%%%%%%%%%%%%%%%%%%%%%%%%%%%%%%%%%%%%%%%%%%%%%%%%%%%
\par \quad {\sc ─юърчр\-Єхы№ёЄтю.}\: ┬эрўрых фю\-ърцхь Єхю\-Ёхьє яЁш
фю\-яюыэш\-Єхы№\-эюь яЁхф\-яю\-ыюцх\-эшш, ўЄю $\rho>0$. ╧єёЄ№
$\mu(t)$, $\hat{\mu}(t)$  -- Їєэъ\-Ўшш Ёрё\-яЁх\-фхых\-эш  ьхЁ $\mu$
ш  $|\mu|$, эюЁьш\-Ёю\-трэ\-э√х єёыю\-тш \-ьш
$\mu(0)=\hat{\mu}(0)=0$. ╥ръ ъръ $\mu\in$ $\frak{M}(\rho(r))$, Єю
ёє\-∙хёЄ\-тєхЄ яю\-ёЄю э\-эр  $A_1$ Єр\-ър , ўЄю эр яюыє\-юёш
$(0,\infty)$ т√\-яюыэ \-хЄё  эх\-Ёртхэ\-ёЄтю $|\mu|([r,er])\leq A_1
V(r)$. ┬ ёыєўрх $\rho>0$ юЄё■фр ёых\-фє\-хЄ, ўЄю
$$
|\mu|((0,r])=\sum\limits_{n=-\infty}^0 |\mu|((e^{n-1}r,e^n r])\leq
A_1 \sum\limits_{n=-\infty}^0 V\left(e^{n-1}r\right)
$$
$$
= A_1 \sum\limits_{n=-\infty}^0 \int\limits_{e^{n-1}r}^{e^nr}
\frac{V\left(e^{n-1}r\right)}{V(t)}\frac{V(t)}{t}dt \leq A_1
\max\limits_{\left[\frac{1}{e},1\right]}\gamma(t)
\sum\limits_{n=-\infty}^0 \int\limits_{e^{n-1}r}^{e^nr}
\frac{V(t)}{t}dt=A_2 \int\limits_0^r \frac{V(t)}{t}dt
$$
\begin{equation}\label{abel_|mu|}
=A_2\int\limits_0^1 \frac{V(ur)}{u}du \leq A_2 \int\limits_0^1
\frac{u^\rho \gamma(u)V(r)}{u}du=AV(r).
\end{equation}
┬ яЁштхфхээ√ї юЎхэърї шё\-яюы№чю\-трыюё№ эх\-Ёртхэ\-ёЄтю
(\ref{order_V_V}). ─рыхх шьххь
$$
|s|([R,eR])=\int\limits_R^{eR} |\Psi(r)|dr\leq \int\limits_R^{eR}
\int\limits_0^\infty
\left|K\left(\frac{t}{r}\right)\right|d|\mu|(t)dr
$$
$$
= \int\limits_0^\infty \int\limits_R^{eR}
\left|K\left(\frac{t}{r}\right)\right|drd\hat{\mu}(t)
=\int\limits_0^\infty t \int\limits_\frac{t}{eR}^\frac{t}{R}
\frac{1}{u^2} |K(u)|dud\hat{\mu}(t).
$$
╘юЁьєыр шэ\-Єху\-Ёш\-Ёю\-тр\-эш  яю ўрёЄ ь фр╕Є
$$
|s|([R,eR])\leq \left\{\left. t \int\limits_\frac{t}{eR}^\frac{t}{R}
\frac{1}{u^2} |K(u)|du \hat{\mu}(t)\right\} \right|_0^\infty
-\int\limits_0^\infty \int\limits_\frac{t}{eR}^\frac{t}{R}
\frac{1}{u^2} |K(u)|du \hat{\mu}(t)dt
$$
\begin{equation}\label{abel_for_th_8}
-R\int\limits_0^\infty \frac{1}{t}
\left|K\left(\frac{t}{R}\right)\right| \hat{\mu}(t)dt +
eR\int\limits_0^\infty \frac{1}{t}
\left|K\left(\frac{t}{eR}\right)\right| \hat{\mu}(t)dt.
\end{equation}
═рь эєцэю юЎхэшЄ№ яюыю\-цш\-Єхы№\-э√х ёырурх\-ь√х, тїюф \-∙шх т
яЁртє■ ўрёЄ№ эх\-Ёртхэ\-ёЄ\-тр  (\ref{abel_for_th_8}). ╠√ сєфхь
шё\-яюы№чю\-трЄ№ эх\-Ёртхэ\-ёЄтю (\ref{order_V_V}). ╟рьхўрхь, ўЄю
$$
\int\limits_\frac{t}{eR}^\frac{t}{R} \frac{1}{u^2} |K(u)|du
=\frac{1}{\frac{t}{R}V\left(\frac{t}{R}\right)}
\int\limits_\frac{t}{eR}^\frac{t}{R}
\frac{\frac{t}{R}V\left(\frac{t}{R}\right)}{uV(u)}
|K(u)|\frac{V(u)}{u}du
$$
\begin{equation}\label{abel_int_th_8}
\leq \frac{1}{\frac{t}{R}V\left(\frac{t}{R}\right)}
\int\limits_\frac{t}{eR}^\frac{t}{R}\left(\frac{t}{Ru}\right)^{\rho+1}
\gamma\left(\frac{t}{Ru}\right) |K(u)| \frac{V(u)}{u}du \leq
\frac{e^{\rho+1}}{\frac{t}{R}V\left(\frac{t}{R}\right)} \hat{\gamma}
\int\limits_\frac{t}{eR}^\frac{t}{R} |K(u)| \frac{V(u)}{u} du,
\end{equation}
уфх $\hat{\gamma}=\max\{\gamma(x):\:x\in[1,e]\}$. ╚ч
эх\-Ёртхэ\-ёЄ\-тр $V(t)\leq t^\rho\gamma(t)$ ш єёыю\-тш 
$t^{\rho-1}\gamma(t)K(t)\in L_1(0,\infty)$ ёых\-фє\-хЄ, ўЄю
$\frac{V(t)}{t}K(t)\in L_1(0,\infty)$. ╥хяхЁ№ шч эх\-Ёртхэ\-ёЄ\-тр
(\ref{abel_int_th_8}) ёых\-фє\-хЄ, ўЄю тэх\-шэ\-Єху\-Ёры№\-э√щ ўыхэ
т эх\-Ёртхэ\-ёЄ\-тх  (\ref{abel_for_th_8}) Ёртхэ эєы■.

\par ╬Ўхэштрхь ўхЄт╕ЁЄ√щ ўыхэ шч яЁртющ ўрёЄш
эх\-Ёртхэ\-ёЄ\-тр  (\ref{abel_for_th_8}).
$$
eR\int\limits_0^\infty \frac{1}{t}
\left|K\left(\frac{t}{eR}\right)\right| \hat{\mu}(t)dt\leq
AeR\int\limits_0^\infty
\frac{1}{t}\left|K\left(\frac{t}{eR}\right)\right| V(t)dt
$$
$$
= AeR\int\limits_0^\infty \frac{1}{u}|K(u)|V(eRu)du\leq
AeRV(eR)\int\limits_0^\infty u^{\rho-1}\gamma(u)|K(u)|du.
$$
╥хь ёрь√ь фюърчрэю, ўЄю яЁш $R>0$ ёяЁртхф\-ыштр юЎхэър
$$
|s|([R,eR])\leq A\int\limits_0^\infty
u^{\rho-1}\gamma(u)|K(u)|du\,eRV(eR).
$$
▌Єю ючэрўрхЄ, ўЄю $s\in\frak{M}(\rho(r)+1)$.

\par ╧єёЄ№ $\varphi\in\Phi$. ╚ьххь
$$
\frac{1}{rV(r)}\int\limits_0^\infty \varphi\left(\frac{u}{r}\right)
ds(u)= \frac{1}{rV(r)}\int\limits_0^\infty
\varphi\left(\frac{u}{r}\right) \Psi(u)du = \frac{1}{rV(r)}
\int\limits_0^\infty \int\limits_0^\infty
\varphi\left(\frac{u}{r}\right)K\left(\frac{t}{u}\right)d\mu(t)du
$$
$$
= \frac{1}{rV(r)} \int\limits_0^\infty \int\limits_0^\infty
\varphi\left(\frac{u}{r}\right)K\left(\frac{t}{u}\right)dud\mu(t)
=\frac{1}{rV(r)} \int\limits_0^\infty \int\limits_0^\infty
\frac{1}{x^2} \varphi\left(\frac{t}{xr}\right) K(x)dxd\mu(t)
$$
\begin{equation}\label{abel_int_varphi}
= \frac{1}{rV(r)} \int\limits_0^\infty \frac{1}{x^2} K(x)
\int\limits_0^\infty t\varphi\left(\frac{t}{xr}\right) d\mu(t)dx.
\end{equation}
\par ╚ч яЁхф√фє\-∙шї Ёрё\-ёєцфх\-эшщ ыхуъю ёых\-фє\-хЄ, ўЄю тёх шэ\-Єху\-Ёры√,
тёЄЁх\-ўр\-■\-∙шх\-ё  т ЇюЁ\-ьє\-ых  (\ref{abel_int_varphi})
 ты \-■Єё  рсёю\-ы■Є\-эю ёїюф \-∙ш\-ьш\-ё . ▌Єю юяЁрт\-ф√\-трхЄ
яхЁх\-ёЄрэют\-ъш яю\-Ё ф\-ър шэ\-Єху\-Ёш\-Ёю\-тр\-эш , ъюЄю\-Ё√х
тёЄЁхўр\-ышё№ яЁш т√тюфх ЇюЁ\-ьє\-ы√  (\ref{abel_int_varphi}).

\par ╧єёЄ№ ЄхяхЁ№ $r_n\to\infty$ Єр\-ър  яю\-ёыхфю\-тр\-Єхы№\-эюёЄ№,
ўЄю $\mu_{r_n}\to\nu$. ╚ьххь
$$
\lim\limits_{n\to\infty}\frac{1}{r_nV(r_n)} \int\limits_0^\infty
t\varphi\left(\frac{t}{xr_n}\right) d\mu(t)=
\lim\limits_{n\to\infty} \int\limits_0^\infty
u\varphi\left(\frac{u}{x}\right) d\mu_{r_n}(u)= \int\limits_0^\infty
u\varphi\left(\frac{u}{x}\right)d\nu(u).
$$
\par ╧єёЄ№ $\supp\varphi\subset[a,b]\subset(0,\infty)$.
╤яЁр\-тхф\-ышт√ юЎхэъш
$$
\left| \frac{1}{r_nV(r_n)} \int\limits_0^\infty
t\varphi\left(\frac{t}{xr_n}\right) d\mu(t)\right|\leq
\frac{bx}{V(r_n)} \|\varphi\|\hat{\mu}(bxr_n) \leq
\frac{Ab\|\varphi\|x}{V(r_n)} V(bxr_n)
$$
$$
\leq A\|\varphi\|(bx)^{\rho+1}\gamma(bx)\leq
A\|\varphi\|b^{\rho+1}\gamma(b)x^{\rho+1}\gamma(x).
$$
╠√ тюё\-яюы№чю\-тр\-ышё№ юЎхэъющ $\gamma(bx)\leq
\gamma(b)\gamma(x)$. ╚ч яюыєўхэ\-эющ юЎхэъш ёых\-фє\-хЄ, ўЄю
$$
\left|\frac{1}{x^2}K(x)\frac{1}{r_nV(r_n)}\int\limits_0^\infty
t\varphi\left(\frac{t}{xr_n}\right)d\mu(t)\right|\leq
A\|\varphi\|b^{\rho+1}\gamma(b)x^{\rho-1}\gamma(x)|K(x)|.
$$
╥хяхЁ№ шч Єхю\-Ёх\-ь√ ╦хсхур ю ьрцюЁш\-Ёєх\-ьющ ёїюфш\-ьюё\-Єш
ёых\-фє\-хЄ, ўЄю
$$
\lim\limits_{n\to\infty}\frac{1}{r_nV(r_n)} \int\limits_0^\infty
\frac{1}{x^2}K(x)\int\limits_0^\infty
t\varphi\left(\frac{t}{xr_n}\right) d\mu(t)dx
$$
$$
= \int\limits_0^\infty \frac{1}{x^2}K(x)\int\limits_0^\infty
u\varphi\left(\frac{u}{x}\right) d\nu(u)dx =\int\limits_0^\infty
u\int\limits_0^\infty \frac{1}{x^2} \varphi\left(\frac{u}{x}\right)
K(x)dxd\nu(u)
$$
$$
= \int\limits_0^\infty \int\limits_0^\infty
K\left(\frac{u}{t}\right)\varphi(t)dtd\nu(u)= \int\limits_0^\infty
\int\limits_0^\infty K\left(\frac{u}{t}\right)d\nu(u)\varphi(t)dt.
$$
╧юфёЄрты   т Ёртхэ\-ёЄтю  (\ref{abel_int_varphi}) $r=r_n$ ш
яхЁх\-їюф  ъ яЁх\-фхыє яЁш $n\to\infty$, яюыєўшь
$$
\lim\limits_{n\to\infty} \int\limits_0^\infty \varphi(t)ds_{r_n}(t)=
\int\limits_0^\infty \int\limits_0^\infty
K\left(\frac{u}{t}\right)d\nu(u)\varphi(t)dt.
$$
▌Єю ючэрўрхЄ, ўЄю яю\-ёыхфю\-тр\-Єхы№\-эюёЄ№ ьхЁ $s_{r_n}$
 ты \-хЄё  ёїюф \-∙хщё , ш х╕ яЁх\-фхы хёЄ№ рсёюы■Є\-эю
эхяЁх\-Ё√т\-эр  ьхЁр ё яыюЄ\-эюёЄ№■ $\int\limits_0^\infty
K\left(\frac{u}{t}\right)d\nu(u)$. ▌Єшь фю\-ърчр\-эю, ўЄю
ьэю\-цх\-ёЄтю $Fr[s]$ ёюфхЁ\-цшЄ тёх рсёю\-ы■Є\-эю эх\-яЁх\-Ё√т\-э√х
ьхЁ√ ё яыюЄ\-эюёЄ \-ьш тшфр $\int\limits_0^\infty
K\left(\frac{u}{t}\right)d\nu(u)$, уфх $\nu\in Fr[\mu]$.

\par ╧єёЄ№ ЄхяхЁ№ $\nu_1$  -- яЁюшч\-тюы№\-эр  ьхЁр шч  $Fr[s]$
ш  $R_n\to\infty$ Єр\-ър  яю\-ёыхфю\-тр\-Єхы№\-эюёЄ№, ўЄю
$s_{R_n}\to\nu_1$. ┬юч№ь╕ь т Ёртхэ\-ёЄ\-тх  (\ref{abel_int_varphi})
$r=R_n$. ═х юуЁрэш\-ўш\-тр  юс∙\-эюё\-Єш, ьюц\-эю ёўш\-ЄрЄ№, ўЄю
$\mu_{R_n}\to\nu\in Fr[\mu]$. ╥юуфр яЁю\-тхфхэ\-эюх
фю\-ърчр\-Єхы№\-ёЄтю яю\-ърч√\-тр\-хЄ, ўЄю ьхЁр  $\nu_1$  ты \-хЄё 
рсёюы■Є\-эю эх\-яЁх\-Ё√т\-эющ ш х╕ яыюЄ\-эюёЄ№ Ёртэр
$\int\limits_0^\infty K\left(\frac{u}{t}\right)d\nu(u)$. ╥хю\-Ёхьр
фы  ёыєўр  $\rho>0$ фю\-ърчр\-эр.

\par ╧єёЄ№ ЄхяхЁ№ $\rho(r)$  -- яЁюшч\-тюы№\-э√щ єЄюў\-э╕э\-э√щ яю\-Ё фюъ,
$\rho=\rho(\infty)$. ╧єёЄ№ $p$  -- Єръюх тх∙хёЄ\-тхэ\-эюх ўшёыю, ўЄю
$\rho_1=\rho+p>0$. ╬сючэрўшь $\rho_1(r)=\rho(r)+p$. ═рЁ фє ё
Ёртхэ\-ёЄ\-тюь
$$
\Psi(r)=\int\limits_0^\infty K\left(\frac{t}{r}\right)d\mu(t)
$$
шьххЄ ьхёЄю ЁртхэёЄтю
$$
\Psi_1(r)=r^p \Psi(r)= \int\limits_0^\infty
\left(\frac{r}{t}\right)^p K\left(\frac{t}{r}\right)t^p d\mu(t)=
\int\limits_0^\infty K_1\left(\frac{t}{r}\right) d\mu_1(t),
$$
уфх $K_1(t)=t^{-p}K(t)$, $d\mu_1(t)=t^p d\mu(t)$. ╚ьххь
$$
t^{\rho_1-1}\gamma(t)K_1(t)=t^{\rho+p-1}\gamma(t)t^{-p}K(t)
=t^{\rho-1}\gamma(t)K(t).
$$
─рыхх ё яюью∙№■ Єхю\-Ёх\-ь√  \ref{azarin_th_1} ыхуъю яЁю\-тхЁшЄ№,
ўЄю фы  єЄюў\-э╕э\-эю\-ую яю\-Ё ф\-ър  $\rho_1(r)$,  фЁр  $K_1$ ш
ьхЁ $\mu_1$ ш  $s_1$, $ds_1(t)=\Psi_1(t)dt$, т√\-яюыэ \-■Є\-ё  тёх
єёыю\-тш  Єхю\-Ёх\-ь√  \ref{abel_th_8} ш, ъЁюьх Єюую,
т√\-яюыэ \-хЄё  єёыю\-тшх $\rho_1>0$. ╧ю фю\-ърчрэ\-эюьє
ьэю\-цх\-ёЄтю $Fr[s_1]$ ёюёЄюшЄ шч рсёюы■Є\-эю эх\-яЁх\-Ё√т\-э√ї
ьхЁ, ьэю\-цх\-ёЄтю яыюЄ\-эюё\-Єхщ ъюЄю\-Ё√ї шьххЄ тшф
$$
\left\{\int\limits_0^\infty
K_1\left(\frac{u}{t}\right)d\nu_1(t):\quad \nu_1\in
Fr[\mu_1]\right\}.
$$
╥хяхЁ№ шч Єхю\-Ёх\-ь√  \ref{azarin_th_1} ыхуъю ёых\-фє\-хЄ
єЄ\-тхЁ\-ц\-фх\-эшх Єхю\-Ёх\-ь√. ╥хю\-Ёхьр фю\-ърчр\-эр.

\par ┬ Єхю\-Ёхьх  \ref{abel_th_5}  юяшё√\-тр\-хЄё  яю\-тхфх\-эшх Їєэъ\-Ўшш
\begin{equation}\label{abel_def_Psi}
\Psi(r)=\int\limits_0^\infty K\left(\frac{t}{r}\right)d\mu(t)
\end{equation}
т юъ\-ЁхёЄ\-эюё\-Єш схё\-ъю\-эхў\-эюё\-Єш т ЄхЁьш\-эрї
яЁх\-фхы№\-эю\-ую ьэю\-цх\-ёЄтр $Fr[\mu]$ ьхЁ√ $\mu$.
┬ючьюц\-эюё\-Єш ¤Єющ Єхю\-Ёх\-ь√ юуЁрэш\-ўш\-тр\-■Є\-ё 
тюч\-ьюц\-эюё\-Є ьш фю\-ърчрЄ№, ўЄю ЄЁющър $(K,\rho(r),\mu)$
єфют\-ыхЄ\-тю\-Ё \-хЄ єёыю\-тш\- ь эхщ\-ЄЁр\-ыш\-чр\-Ўшш эєы  ш
схё\-ъю\-эхў\-эюё\-Єш.

\par ╬фэръю, ёє\-∙хёЄ\-тє■Є ш фЁєушх яЁхя Є\-ёЄ\-тш  фы  шчєўх\-эш  ётющёЄт
Їєэъ\-Ўшш $\Psi(r)$. ╠√ яю\-фЁюс\-эю Ёрё\-ёьюЄ\-Ёшь ёыєўрщ, ъюуфр
$K$  -- ЇшэшЄ\-эюх схё\-ъю\-эхў\-эю фшЇ\-ЇхЁхэ\-Ўш\-Ёєх\-ьюх  фЁю эр
яюыє\-юёш $(0,\infty)$. ▌Єю юўхэ№ ёшы№\-эюх юуЁрэш\-ўх\-эшх эр  фЁю
$K$. ┬ ўрёЄ\-эюё\-Єш, фы  Єръюую  фЁр  $K$ ЄЁющър $(K,\rho(r),\mu)$
єфют\-ыхЄ\-тю\-Ё \-хЄ єёыю\-тш\- ь эхщ\-ЄЁр\-ыш\-чр\-Ўшш эєы  ш
схё\-ъю\-эхў\-эюё\-Єш фы  ы■сюую єЄюў\-э╕э\-эю\-ую яю\-Ё ф\-ър
$\rho(r)$ ш ы■сющ Ёр\-фю\-эю\-тющ ьхЁ√ эр яюыє\-юёш $(0,\infty)$.

\par ╧єёЄ№ $K$  -- ЇшэшЄ\-эюх схё\-ъю\-эхў\-эю фшЇ\-ЇхЁхэ\-Ўш\-Ёєх\-ьюх
 фЁю эр яюыє\-юёш $(0,\infty)$, $\mu$  -- Ёр\-фю\-эю\-тр ьхЁр эр
¤Єющ яюыє\-юёш. ╟рьхЄшь, ўЄю схё\-ъю\-эхў\-эр 
фшЇ\-Їх\-Ёхэ\-Ўш\-Ёєх\-ьюёЄ№  фЁр $K$, хую ЇшэшЄ\-эюёЄ№ ш
ЇюЁ\-ьє\-ыр шэ\-Єху\-Ёш\-Ёю\-тр\-эш  яю ўрёЄ ь яюч\-тюы \-■Є
эр\-яшёрЄ№ схё\-ъю\-эхў\-эюх ўшёыю ЁртхэёЄт
\begin{equation}\label{abel_F_n}
(-1)^{n+1}r^{n+1}\Psi(r)=\int\limits_0^\infty
K^{(n+1)}\left(\frac{t}{r}\right)F_n(t)dt,\qquad n=0,1,...,
\end{equation}
уфх $F_0(t)=\mu(t)$  -- Їєэъ\-Ўш  Ёрё\-яЁх\-фхых\-эш  ьхЁ√ $\mu$,
$F'_{n+1}(t)=F_n(t)$, $n=0,1,...$.

\par ┬юяЁюё чръы■ўр\-хЄё  т ёых\-фє\-■\-∙хь. ╧ючтюы ■Є ыш
Ёртхэ\-ёЄ\-тр  (\ref{abel_F_n}) ш Єхю\-Ёх\-ьр  \ref{abel_th_1}
юяЁх\-фх\-ышЄ№ яю\-Ё фюъ ЁюёЄр Їєэъ\-Ўшш $\Psi(r)$ эр
схё\-ъю\-эхў\-эюё\-Єш? ╩ръ ь√ єтшфшь фрыхх, юЄтхЄ ёых\-фє\-■\-∙шщ.
╧Ёхф\-ыюцхэ\-э√щ ьхЄюф тю ьэюушї ёыєўр ї яю\-чтюы хЄ юяЁх\-фхышЄ№
яю\-Ё фюъ ЁюёЄр Їєэъ\-Ўшш $\Psi(r)$ эр схё\-ъю\-эхў\-эюё\-Єш.
╬фэръю, шьх■Є\-ё  шё\-ъы■ўх\-эш .

\par ╧єёЄ№ Ёр\-фю\-эютр ьхЁр $\mu_1$ эр яюыє\-юёш $(0,\infty)$
ёю\-ёЁхфю\-Єюўх\-эр эр яюыє\-шэ\-ЄхЁ\-трых $(0,1]$, $K$  --
ЇшэшЄ\-эюх эр яюыє\-юёш $(0,\infty)$  фЁю,
$$
u(r)=\int\limits_0^\infty K\left(\frac{t}{r}\right)d\mu_1(t).
$$
╥юуфр Їєэъ\-Ўш  $u(r)$ Ёртэр эєы■ т эх\-ъю\-Єю\-Ёющ
юъ\-ЁхёЄ\-эюё\-Єш схё\-ъю\-эхў\-эюё\-Єш. ╧ю¤Єю\-ьє хёыш ь√ Ёх°рхь
чрфрўє ю яю\-тхфх\-эшш Їєэъ\-Ўшш $\Psi(r)$, чр\-фртрх\-ьющ
Ёртхэ\-ёЄ\-тюь (\ref{abel_def_Psi}) т юъ\-ЁхёЄ\-эюё\-Єш
схё\-ъю\-эхў\-эюё\-Єш, Єю, схч юуЁрэш\-ўх\-эш  юс∙\-эюё\-Єш, ьюц\-эю
ёўш\-ЄрЄ№, ўЄю ьхЁр  $\mu$ т ЇюЁ\-ьє\-ых (\ref{abel_def_Psi})
ёю\-ёЁхфю\-Єюўх\-эр эр яюыє\-юёш $(1,\infty)$. ▌Єю єёыю\-тшх т
фры№\-эхщ\-°хь ь√ сєфхь ёўш\-ЄрЄ№ т√\-яюыэхэ\-э√ь. ┬  ¤Єюь ёыєўрх
Їєэъ\-Ўш   $F_0(t)$ сєфхЄ юуЁрэш\-ўхэр эр ы■сюь шэ\-ЄхЁ\-трых
$(0,N)$, р тёх Їєэъ\-Ўшш $F_k(t)$, $k\geq 1$, ьюц\-эю ёўш\-ЄрЄ№
эх\-яЁх\-Ё√т\-э√ьш эр яюыє\-юёш $[0,\infty)$. ┬ю тё ъюь ёыєўрх,
ьюц\-эю ёўш\-ЄрЄ№, ўЄю тёх Їєэъ\-Ўшш $F_k(t)$  ты \-■Є\-ё 
ыю\-ъры№\-эю шэ\-Єху\-Ёш\-Ёєх\-ь√ьш эр яюыє\-юёш $(0,\infty)$ ш фы 
ы■сюую $k\geq 0$ ёє\-∙хёЄ\-тєхЄ шэ\-ЄхуЁры $\int\limits_0^1
F_k(t)dt$.

\par ╧єёЄ№ Їєэъ\-Ўш  $f(t)$ ыю\-ъры№\-эю шэ\-Єху\-Ёш\-Ёєх\-ьр эр яюыє\-юёш
$(0,\infty)$ ш ёє\-∙хёЄ\-тєхЄ шэ\-ЄхуЁры $\int\limits_0^1 f(t)dt$
їюЄ  с√ т эх\-ёюсёЄ\-тхэ\-эюь ёь√ёых. ╘єэъ\-Ўш■ $F(t)$,
юяЁх\-фхы х\-ьє■ Ёртхэ\-ёЄ\-тюь $F(t)=-\int\limits_t^\infty f(x)dx$,
хёыш эр\-яшёрэ\-э√щ шэ\-ЄхуЁры ёїю\-фшЄ\-ё  ш Ёртхэ\-ёЄ\-тюь
$F(t)=\int\limits_0^t f(x)dx$ т яЁюЄшт\-эюь ёыєўрх, сєфхь
эр\-ч√\-трЄ№ {\it ърэю\-эш\-ўхё\-ъющ яхЁтю\-юсЁрч\-эющ} Їєэъ\-Ўшш
$f$ эр яюыє\-юёш $(0,\infty)$. ╩рэюэш\-ўхё\-ър  яхЁтю\-юсЁрч\-эр 
$F$ юяЁх\-фхы \-хЄ\-ё  яю Їєэъ\-Ўшш $f$ юфэю\-чэрў\-эю. ─ы  ы■сющ
ыю\-ъры№\-эю шэ\-Єху\-Ёш\-Ёєх\-ьющ Їєэъ\-Ўшш  $f$ эр яюыє\-юёш
$(0,\infty)$ ёє\-∙хёЄ\-тє\-хЄ яхЁтю\-юсЁрч\-эр  $F$, юфэръю
ърэю\-эш\-ўхё\-ър  яхЁтю\-юсЁрч\-эр  ёє\-∙хёЄ\-тєхЄ эх тёхуфр.
═ряЁш\-ьхЁ, є Їєэъ\-Ўшш $f(t)=\frac{1}{t}$ эхЄ ърэю\-эш\-ўхё\-ъющ
яхЁтю\-юсЁрч\-эющ эр яюыє\-юёш $(0,\infty)$.

\par ╘єэъЎш■ $F_0(t)$ т ЇюЁ\-ьє\-ых  (\ref{abel_F_n}) юфэю\-чэрў\-эю
юяЁх\-фхы \-хь ёых\-фє\-■\-∙шь юсЁр\-чюь. $F_0(t)$ $=$
$-\mu((t,\infty))$, хёыш $F_0(t)$ ъю\-эхў\-эр  тхыш\-ўшэр
(эр\-яюьэшь, ўЄю $\mu$  -- Ёр\-фю\-эютр ьхЁр эр яюыє\-юёш
$(0,\infty)$, ш яю¤Єю\-ьє тхыш\-ўшэр $\mu(E)$ юяЁх\-фхыхэр эх фы 
ы■сюую сюЁх\-ыхт\-ёъюую ьэю\-цх\-ёЄтр $E\subset(0,\infty)$). ┬
яЁю\-Єшт\-эюь ёыєўрх яю\-ырур\-хь, ўЄю $F_0(t)$ $=$ $\mu((0,t])$.
▌Єр тхышўш\-эр тёхуфр шьххЄ ёь√ёы шч-чр яЁхф\-яюыюцх\-эш , ўЄю ьхЁр
$\mu$ ёюёЁхфю\-Єюўхэр эр яюыє\-юёш $[1,\infty)$. ─рыхх Їєэъ\-Ўш■
$F_{n+1}(t)$ юяЁх\-фхы \-хь ъръ ърэю\-эш\-ўхё\-ъє■ яхЁтю\-юсЁрч\-эє■
Їєэъ\-Ўшш $F_n(t)$. ╥хяхЁ№ т ЇюЁ\-ьє\-ых  (\ref{abel_F_n})
яю\-ёыхфю\-тр\-Єхы№\-эюёЄ№ $F_n(t)$  ты \-хЄ\-ё  юфэю\-чэрў\-эю
юяЁх\-фхы╕э\-эющ.

\par ┬ Ёрчфхых 2 яю\-Ё фюъ ш єЄюў\-э╕э\-э√щ яю\-Ё фюъ юяЁх\-фхы \-ышё№ фы 
яюыю\-цш\-Єхы№\-э√ї Їєэъ\-Ўшщ. ╨рёяЁюёЄ\-Ёр\-эшь ¤Єш юяЁх\-фхых\-эш 
эр ъюья\-ыхъё\-эю\-чэрў\-э√х Їєэъ\-Ўшш.

\par ╧єёЄ№ $f(t)$  -- ъюья\-ыхъё\-эю\-чэрў\-эр  Їєэъ\-Ўш  эр яюыє\-юёш
$(0,\infty)$. ╧юЁ фюъ  $\rho$ Їєэъ\-Ўшш  $f$ юяЁх\-фхы \-хЄ\-ё  яю
ЇюЁ\-ьє\-ых
$$
\rho=\mathop{\overline{\lim}}\limits_{r\to\infty}
\frac{\ln|f(r)|}{\ln r}.
$$
┼ёыш $\rho$  -- тх∙хёЄ\-тхэ\-эюх ўшёыю, Єю Їєэъ\-Ўш   $f$
эр\-ч√тр\-хЄё  Їєэъ\-Ўшхщ {\it ъю\-эхў\-эю\-ую яю\-Ё ф\-ър} $\rho$.

\par ╙Єюў\-э╕э\-э√щ яю\-Ё фюъ $\rho(r)$ сєфхь эр\-ч√трЄ№ {\it єЄюў\-э╕э\-э√ь
яю\-Ё ф\-ъюь} Їєэъ\-Ўшш $f(r)$, хёыш
$\overline{\mathop{\lim}}_{r\to\infty} \frac{|f(r)|}{V(r)}$ $=$
$\sigma\in(0,\infty)$. ┼ёыш $\rho(r)$  -- єЄюў\-э╕э\-э√щ яю\-Ё фюъ
Їєэъ\-Ўшш $f(r)$, Єю сєфхь уютю\-ЁшЄ№, ўЄю Їєэъ\-Ўш  $f(r)$ Ёрё\-Є╕Є
эр схё\-ъю\-эхў\-эюё\-Єш ъръ Їєэъ\-Ўш  $V(r)$.

\par ╚ч Єхю\-Ёх\-ь√  \ref{order_rho(r)} ёых\-фєхЄ, ўЄю хёыш Їєэъ\-Ўш  $f(r)$
шьххЄ ъю\-эхў\-э√щ яю\-Ё фюъ, Єю є эх╕ хёЄ№ єЄюў\-э╕э\-э√щ яю\-Ё фюъ
$\rho(r)$ (ъю\-эхў\-эю, Єръшї яю\-Ё ф\-ъют ёє\-∙хёЄ\-тєхЄ
схё\-ъю\-эхў\-эюх ўшёыю).

\par ┬ю ьэюушї чрфрўрї, т ўрёЄ\-эюё\-Єш, т ёт чш ё ЇюЁ\-ьє\-ыющ
(\ref{abel_F_n}) тюч\-эшър\-хЄ тюяЁюё юс юЎхэ\-ърї яхЁтю\-юсЁрч\-эющ
$F$ Їєэъ\-Ўшш $f$. ╬ёЄрэю\-тшь\-ё  эр ¤Єюь яюфЁюс\-эхх.
%%%%%%%%%%%%%%%%%%%%%%%%%%%%%%%%%%%%%%%%%%%%%%%%%%%%%%%%%%%%%%%%%%%%%%%%%%%%%%%%%%%%%%%%%%%%%%%%%%%
\begin{lemma}\label{abel_lemma_5}\hskip-2mm{.}\:
%%%%%%%%%%%%%%%%%%%%%%%%%%%%%%%%%%%%%%%%%%%%%%%%%%%%%%%%%%%%%%%%%%%%%%%%%%%%%%%%%%%%%%%%%%%%%%%%%%%
╧єёЄ№ $f$  -- ыю\-ъры№эю шэ\-Єху\-Ёш\-Ёєх\-ьр  Їєэъ\-Ўш 
ъю\-эхў\-эю\-ую яю\-Ё ф\-ър эр яюыє\-юёш $(0,\infty)$ ш Єр\-ър , ўЄю
ёє\-∙хёЄ\-тєхЄ шэ\-ЄхуЁры $\int\limits_0^1 f(t)dt$. ╧єёЄ№ $\rho(r)$
-- єЄюў\-э╕э\-э√щ яю\-Ё фюъ Їєэъ\-Ўшш  $f$, $F$  --
ърэю\-эш\-ўхё\-ър  яхЁтю\-юсЁрч\-эр  Їєэъ\-Ўшш  $f$. ╥юуфр
ёє\-∙хёЄ\-тє■Є яю\-ёЄю\- э\-э√х $M_k$ ш  $r_0$ Єръшх,
ўЄю \\
1) $|F(r+\alpha r)-F(r)|\leq M_1\alpha rV(r)$, $\;r\geq r_0$, $\;\alpha\in [0,1]$;\\
2) $|F(r)|\leq M_2\left(1+\int\limits_1^r V(t)dt\right)$, $\;r\geq
1$;\\
3) $|F(r)|\leq M_3 \int\limits_r^\infty V(t)dt$, $\;r\geq 1$;\\
4) $|F(r)|\leq M_4 rV(r)$, $\;r\geq 1$, хёыш $\rho=\rho(\infty)\neq
-1$;\\
5) яю\-Ё фъш $\rho(f)$, $\rho(F)$ Їєэъ\-Ўшщ $f$ ш  $F$
єфют\-ыхЄ\-тю\-Ё ■Є эх\-Ёртхэ\-ёЄ\-тє $\rho(F)\leq \rho(f)+1$.
\end{lemma}

\par ╟рьхЄшь, ўЄю хёыш шэ\-ЄхуЁры $\int\limits_1^\infty V(t)dt$ Ёрё\-їюфшЄ\-ё ,
Єю эх\-Ёртхэ\-ёЄ\-тю  3) хёЄ№ ЄЁштш\-ры№\-эюх эх\-Ёртхэ\-ёЄ\-тю
$|F(r)|\leq \infty$.
%%%%%%%%%%%%%%%%%%%%%%%%%%%%%%%%%%%%%%%%%%%%%%%%%%%%%%%%%%%%%%%%%%%%%%%%%%%%%%%%%%%%%%%%%%%%%%%%%%%%
\par \quad {\sc ─юърчр\-Єхы№ёЄтю.}\: ╧єёЄ№ $\sigma$  -- Єшя Їєэъ\-Ўшш $f$
юЄ\-эюёш\-Єхы№\-эю єЄюў\-э╕э\-эю\-ую яю\-Ё ф\-ър $\rho(r)$. ╥юуфр
фы  тёхї фю\-ёЄрЄюў\-эю сюы№°шї  $r$ ш  $\alpha\in[0,1]$ сєфєЄ
т√\-яюыэ Є№\-ё  ёююЄ\-эю°х\-эш 
$$
|F(r+\alpha r)-F(r)|= \left|\int\limits_r^{(1+\alpha)r}f(t)dt\right|
\leq 2\sigma \int\limits_r^{(1+\alpha)r}V(t)dt = 2\sigma
r\int\limits_1^{1+\alpha} V(ur)du
$$
$$
\leq 2\sigma rV(r)\int\limits_1^{1+\alpha}u^\rho\gamma(u)du\leq
2\sigma\gamma_1\alpha rV(r),
$$
уфх $\gamma_1=\max\{u^\rho\gamma(u):\:u\in[1,2]\}$. ╥хь ёрь√ь
єЄ\-тхЁ\-ц\-фх\-эшх  1) фю\-ърчр\-эю. ╙ЄтхЁц\-фх\-эш  2)-4) ыхуъю
ёых\-фє\-■Є шч єЄ\-тхЁ\-ц\-фх\-эш   1) ш ётющёЄт єЄюў\-э╕э\-эю\-ую
яю\-Ё ф\-ър. ─хщёЄтш\-Єхы№эю, яЁю\-тхЁшь, эряЁш\-ьхЁ,
єЄ\-тхЁ\-ц\-фх\-эшх  2). ╧єёЄ№ $r\geq 2r_0$. ╬яЁхфхышь ўшёыю  $n_0$
шч єёыю\-тш  $2^{-n_0}r\in [r_0,2r_0)$. ╬сючэрўшь $r_1=2^{-n_0}r$.
╚ьххь
$$
|F(r)|=\left|F(r_1)+\sum\limits_{n=1}^{n_0}
\left(F\left(2^nr_1\right)-F\left(2^{n-1}r_1\right)\right)\right|
$$
$$
\leq |F(r_1)| +M_1\sum\limits_{n=1}^{n_0} 2^{n-1}r_1
V\left(2^{n-1}r_1\right) =|F(r_1)|+M_1\sum\limits_{n=1}^{n_0}
\int\limits_{2^{n-1}r_1}^{2^nr_1} V\left(2^{n-1}r_1\right)dt
$$
$$
\leq |F(r_1)|+M_5\sum\limits_{n=1}^{n_0}
\int\limits_{2^{n-1}r_1}^{2^nr_1} V(t)dt = |F(r_1)|+ M_5
\int\limits_{r_1}^r V(t)dt.
$$
╚ч яюыєўхэ\-эю\-ую эх\-Ёртхэ\-ёЄ\-тр ёых\-фє\-хЄ єЄ\-тхЁ\-ц\-фх\-эшх
2). ═хЁртхэ\-ёЄтю $\rho(F)\leq \rho(f)+1$ ыхуъю ёых\-фє\-хЄ шч
єЄ\-тхЁ\-ц\-фх\-эш   4), хёыш $\rho$ $\neq -1$, ш
єЄ\-тхЁ\-ц\-фх\-эш  2), хёыш $\rho=-1$. ╦хььр фю\-ърчр\-эр.

\par ╦хььр  \ref{abel_lemma_5} ёыєцшЄ ьюЄштр\-Ўшхщ фы  ттхфх\-эш 
ёых\-фє\-■\-∙хую юяЁх\-фхых\-эш .

\par ╧єёЄ№ $f(r)$  -- ыю\-ъры№\-эю шэ\-Єху\-Ёш\-Ёєх\-ьр  Їєэъ\-Ўш 
ъю\-эхў\-эю\-ую яю\-Ё ф\-ър эр яюыє\-юёш $(0,\infty)$,
ёє\-∙хёЄ\-тєхЄ шэ\-Єху\-Ёры $\int\limits_0^1 f(t)dt$ ш $F(r)$  --
ърэю\-эш\-ўхё\-ър  яхЁтю\-юсЁрч\-эр  Їєэъ\-Ўшш  $f$ эр яюыє\-юёш
$(0,\infty)$.

\par ╙Єюў\-э╕э\-э√щ яю\-Ё фюъ $\rho(r)$ Їєэъ\-Ўшш $f$ эр\-ч√тр\-хЄ\-ё 
{\it єёЄющ\-ўш\-т√ь єЄюў\-э╕э\-э√ь яю\-Ё ф\-ъюь} Їєэъ\-Ўшш  $f$,
хёыш $\rho=\rho(\infty)$ $\neq -1$ ш
$$
\mathop{\overline{\lim}}\limits_{r\to\infty} \frac{|F(r)|}{rV(r)}>0.
$$
\par ╟рьхЄшь, ўЄю Їєэъ\-Ўшш $\cos r$, $\sin r$, $e^{ir}$ эх шьх■Є
єёЄющ\-ўш\-тюую єЄюў\-э╕э\-эю\-ую яю\-Ё ф\-ър.
%%%%%%%%%%%%%%%%%%%%%%%%%%%%%%%%%%%%%%%%%%%%%%%%%%%%%%%%%%%%%%%%%%%%%%%%%%%%%%%%%%%%%%%%%%%%%%%%%%%
\begin{theorem}\label{abel_th_9}\hskip-2mm{.}\:
%%%%%%%%%%%%%%%%%%%%%%%%%%%%%%%%%%%%%%%%%%%%%%%%%%%%%%%%%%%%%%%%%%%%%%%%%%%%%%%%%%%%%%%%%%%%%%%%%%%
╧єёЄ№ $f(r)$  -- ыю\-ъры№эю шэ\-Єху\-Ёш\-Ёєх\-ьр  Їєэъ\-Ўш 
ъю\-эхў\-эю\-ую яю\-Ё ф\-ър эр яюыє\-юёш $(0,\infty)$,
ёє\-∙хёЄ\-тєхЄ шэ\-ЄхуЁры $\int\limits_0^1 f(t)dt$ ш $F(r)$  --
ърэюэш\-ўхё\-ър  яхЁтю\-юсЁрч\-эр  Їєэъ\-Ўшш  $f$ эр яюыє\-юёш
$(0,\infty)$. ╧єёЄ№ $\rho(r)$, $\rho(\infty)$ $\neq -1,-2$,...  --
єёЄющ\-ўш\-т√щ єЄюў\-э╕э\-э√щ яю\-Ё фюъ Їєэъ\-Ўшш  $f$. ╥юуфр
єЄюў\-э╕э\-э√щ яю\-Ё фюъ $\rho(r)+1$ сєфхЄ єёЄющ\-ўш\-т√ь
єЄюў\-э╕э\-э√ь яю\-Ё ф\-ъюь Їєэъ\-Ўшш  $F$.
\end{theorem}
%%%%%%%%%%%%%%%%%%%%%%%%%%%%%%%%%%%%%%%%%%%%%%%%%%%%%%%%%%%%%%%%%%%%%%%%%%%%%%%%%%%%%%%%%%%%%%%%%%%%
\par {\sc  ─юърчрЄхы№ёЄтю.}\: ╚ч юяЁх\-фхых\-эш  єёЄющ\-ўш\-тюую
єЄюў\-э╕э\-эю\-ую яю\-Ё ф\-ър ш єЄ\-тхЁ\-ц\-фх\-эш   4) ыхьь√
\ref{abel_lemma_5}  ёых\-фє\-хЄ, ўЄю єЄюў\-э╕э\-э√щ яю\-Ё фюъ
$\rho(r)+1$ сєфхЄ єЄюў\-э╕э\-э√ь яю\-Ё ф\-ъюь Їєэъ\-Ўшш  $F$.
╬ёЄрыюё№ фю\-ърчрЄ№, ўЄю ¤Єю єёЄющ\-ўш\-т√щ єЄюў\-э╕э\-э√щ яю\-Ё фюъ
Їєэъ\-Ўшш  $F$. ┼ёыш ¤Єю эх Єръ, Єю, ъръ ёых\-фє\-хЄ шч ыхьь√
\ref{abel_lemma_5} ш юяЁх\-фхых\-эш  єёЄющ\-ўш\-тюую
єЄюў\-э╕э\-эю\-ую яю\-Ё ф\-ър, шьххЄ ьхёЄю Ёртхэ\-ёЄтю
% % % % % % % % % % % % % % % % % % % % % % % % % % % % % % % % % %
\begin{equation}\label{abel_for_th_9}
\lim\limits_{r\to\infty} \frac{F_1(r)}{r^2V(r)}=0,
\end{equation}
уфх $F_1(r)$  -- ърэюэш\-ўхё\-ър  яхЁтю\-юсЁрч\-эр  Їєэъ\-Ўшш $F(r)$
эр яюыє\-юёш $(0,\infty)$.

\par ╧Ёхф\-яюыю\-цшь тэрўр\-ых, ўЄю т√\-яюыэ \-хЄё 
эх\-Ёртхэ\-ёЄтю $\rho(\infty)>-2$. ╧юърцхь, ўЄю т  ¤Єюь ёыєўрх
шэ\-ЄхуЁры $\int\limits_r^\infty F(t)dt$ Ёрё\-їюфшЄ\-ё . ╥ръ ъръ
$\rho(r)+1$  -- єЄюў\-э╕э\-э√щ яю\-Ё фюъ Їєэъ\-Ўшш  $F$, Єю ьюц\-эю
ёўш\-ЄрЄ№, ўЄю ёє\-∙хёЄ\-тє■Є ўшёыю $m>0$ ш
яю\-ёыхфю\-тр\-Єхы№\-эюёЄ№ $r_n\to\infty$ Єръшх, ўЄю $\Re F(r_n)\geq
2mr_nV(r_n)$. ▌Єюую ьюц\-эю фю\-сшЄ№ё , чр\-ьхэ   т ёыєўрх
эхюс\-їюфш\-ьюёЄш Їєэъ\-Ўш■ $f$ эр $-f$ шыш эр $\pm if$. ╚ч
єЄ\-тхЁ\-ц\-фх\-эш   1) ыхьь√ \ref{abel_lemma_5} ёых\-фє\-хЄ, ўЄю
ёє\-∙хёЄ\-тєхЄ $\alpha_0>0$ Єръюх, ўЄю фы  тёхї фюёЄр\-Єюўэю
сюы№\-°шї $n$ ш фы  тёхї $r\in[r_n,$ $(1+\alpha_0)r_n]$ сєфхЄ
т√\-яюыэ Є№\-ё  эх\-Ёртхэ\-ёЄ\-тю $\Re F(r)\geq mr_nV(r_n)$. ╥юуфр
$$
\left|\int\limits_{r_n}^{(1+\alpha_0)r_n}F(t)dt\right|\geq
\left|\int\limits_{r_n}^{(1+\alpha_0)r_n}\Re F(t)dt\right|=
\int\limits_{r_n}^{(1+\alpha_0)r_n}\Re F(t)dt\geq m\alpha_0 r_n^2
V(r_n).
$$
╥ръ ъръ $\lim\limits_{r\to\infty} r^2 V(r)=\infty$, Єю яю ъЁшЄх\-Ёш■
╩ю°ш шэ\-ЄхуЁры $\int\limits_r^\infty F(t)dt$ Ёрё\-їюфшЄ\-ё .
╧ю¤Єю\-ьє, ёюуырёэю юяЁх\-фхых\-эш■ ърэюэш\-ўхёъющ
яхЁтю\-юсЁрч\-эющ, Ёртхэ\-ёЄ\-тю  (\ref{abel_for_th_9}) шьххЄ тшф
\begin{equation}\label{abel_lim_th_9}
\lim\limits_{r\to\infty} \frac{1}{r^2V(r)}\int\limits_0^r F(t)dt=0.
\end{equation}

\par ╘єэъЎш  $g(r)$ эр\-ч√тр\-хЄ\-ё  {\it ьхфыхэ\-эю ьхэ ■\-∙хщ\-ё }
юЄ\-эюёш\-Єхы№\-эю єЄюў\-э╕э\-эю\-ую яю\-Ё ф\-ър $\rho(r)$, хёыш
$$
\lim\limits_{{r\to\infty}\atop{\alpha\to 0}} \frac{g(r+\alpha
r)-g(r)}{V(r)}=0.
$$
╚ч єЄ\-тхЁ\-ц\-фх\-эш   1) ыхьь√  \ref{abel_lemma_5} ёых\-фє\-хЄ,
ўЄю Їєэъ\-Ўш  $F(r)$  ты \-хЄ\-ё  ьхфыхэ\-эю ьхэ ■\-∙хщ\-ё 
юЄ\-эюёш\-Єхы№\-эю єЄюў\-э╕э\-эю\-ую яю\-Ё ф\-ър $\rho(r)+1$.

\par ╤ыхфёЄ\-тш\-хь Єхю\-Ёх\-ь√  2  \cite{Grishin_2003}  ты \-хЄ\-ё 
ёых\-фє\-■\-∙хх єЄ\-тхЁ\-ц\-фх\-эшх.
%%%%%%%%%%%%%%%%%%%%%%%%%%%%%%%%%%%%%%%%%%%%%%%%%%%%%%%%%%%%%%%%%%%%%%%%%%%%%%%%%%%%%%%%%%%%%%%%%%%
\par ╧єёЄ№ Їєэъ\-Ўш  $g$  ты \-хЄё  ьхфыхэ\-эю ьхэ ■\-∙хщ\-ё 
юЄ\-эюёш\-Єхы№\-эю єЄюў\-э╕э\-эю\-ую яю\-Ё ф\-ър $\rho(r)$, яЁш\-ў╕ь
$\rho(\infty)>-1$. ╥юуфр, хёыш
$$
\lim\limits_{r\to\infty} \frac{1}{rV(r)} \int\limits_0^r g(t)dt=a,
$$
Єю
$$
\lim\limits_{r\to\infty} \frac{g(r)}{V(r)}=a(\rho+1).
$$
\par ╤ыхфёЄтшхь ёЇЁьєыш\-Ёю\-трэ\-эю\-ую єЄ\-тхЁ\-ц\-фх\-эш  ш
Ёртхэ\-ёЄ\-тр  (\ref{abel_for_th_9})  ты \-хЄ\-ё  Ёртхэ\-ёЄтю
\begin{equation}\label{abel_lim_F}
\lim\limits_{r\to\infty} \frac{F(r)}{rV(r)}=0.
\end{equation}
▌Єю яЁюЄштю\-ЁхўшЄ Єюьє, ўЄю $\rho(r)$  -- єёЄющ\-ўш\-т√щ
єЄюў\-э╕э\-э√щ яю\-Ё фюъ Їєэъ\-Ўшш  $f$.

\par ╥ръшь юсЁр\-чюь Ёртхэ\-ёЄ\-тю  (\ref{abel_for_th_9}) тхф╕Є ъ
яЁю\-Єштю\-Ёхўш■, ш Єхю\-Ёх\-ьр фю\-ърчр\-эр фы  ёыєўр , ъюуфр
$\rho(\infty)>-2$. ╧Ёхф\-яюыю\-цшь ЄхяхЁ№, ўЄю $\rho(\infty)<-2$. ┬
¤Єюь ёыєўрх Ёртхэ\-ёЄ\-тю  (\ref{abel_for_th_9}) чряшё√\-тр\-хЄё  т
тшфх
$$
\lim\limits_{r\to\infty} \frac{1}{r^2V(r)}\int\limits_r^\infty
F(t)dt=0.
$$
╥хяхЁ№ фюяюыэш\-Єхы№\-эюх єЄ\-тхЁ\-ц\-фх\-эшх фы  юсюёэютр\-эш 
Ёртхэ\-ёЄ\-тр  (\ref{abel_lim_F}) т√\-уы фшЄ Єръ.

\par ╧єёЄ№ Їєэъ\-Ўш  $g$  ты \-хЄё  ьхфыхэ\-эю ьхэ ■\-∙хщё  Їєэъ\-Ўшхщ
юЄ\-эюёш\-Єхы№\-эю єЄюў\-э╕э\-эю\-ую яю\-Ё ф\-ър $\rho(r)$, яЁш\-ў╕ь
$\rho(\infty)<-1$. ╥юуфр, хёыш
$$
\lim\limits_{r\to\infty} \frac{1}{rV(r)}\int\limits_r^\infty
g(t)dt=a,
$$
Єю
$$
\lim\limits_{r\to\infty} \frac{g(r)}{V(r)}=-a(\rho+1).
$$
▌Єю єЄ\-тхЁ\-ц\-фхэшх Єръ\-цх ёых\-фєхЄ шч Ёх\-чєы№\-ЄрЄют Ёрсю\-Є√
\cite{Grishin_2003}. ╥хь ёрь√ь ь√ фю\-ърчрыш, ўЄю ш т ёыєўрх
$\rho(\infty)<-2$ шч Ёртхэ\-ёЄ\-тр (\ref{abel_for_th_9}) ёых\-фєхЄ
Ёртхэ\-ёЄтю (\ref{abel_lim_F}), ъюЄюЁюх яЁюЄш\-тю\-ЁхўшЄ єёыю\-тш ь
Єхю\-Ёх\-ь√. ╥ръшь юсЁр\-чюь, яЁхф\-яюыюцх\-эшх ю Єюь, ўЄю
$\rho(r)+1$ эх  ты \-хЄё  єёЄющ\-ўш\-т√ь єЄюў\-э╕э\-э√ь яю\-Ё ф\-ъюь
Їєэъ\-Ўшш $F(r)$, яЁш\-тюфшЄ ъ яЁюЄш\-тю\-Ёхўш■. ╥хю\-Ёхьр
фю\-ърчр\-эр.

\par ─юърчрээр  Єхю\-Ёхьр фюёЄр\-Єюўэю cтюх\-юсЁрчэр. ┼╕ ьюц\-эю ЄЁръЄю\-трЄ№,
ъръ Єхю\-Ёх\-ьє ю ёЄрсш\-ыш\-чр\-Ўшш ётющ\-ёЄтр Їєэъ\-Ўшш шьхЄ№
єёЄющ\-ўш\-т√щ єЄюў\-э╕э\-э√щ яю\-Ё фюъ яЁш яю\-ёыхфю\-тр\-Єхы№\-эюь
шэ\-Єху\-Ёш\-Ёю\-тр\-эшш. ╧єёЄ№ $F_k(r)$  --
яю\-ёыхфю\-тр\-Єхы№\-э√х ърэюэш\-ўхё\-ъшх яхЁтю\-юсЁрч\-э√х эр
яюыє\-юёш $(0,\infty)$ Їєэъ\-Ўшш $f(r)$. ╚ч фю\-ърчрэ\-эющ
Єхю\-Ёх\-ь√ ёых\-фєхЄ, ўЄю хёыш Їєэъ\-Ўш  $F_k(r)$ шьххЄ
єёЄющ\-ўш\-т√щ єЄюў\-э╕э\-э√щ яю\-Ё фюъ $\rho_k(r)$ ш хёыш
$\rho_k(\infty)$ эх  ты \-хЄё  Ўхы√ь ёЄЁюую юЄЁшЎр\-Єхы№\-э√ь
ўшёыюь, Єю фы  ы■сюую $m\geq 1$ єЄюў\-э╕э\-э√щ яю\-Ё фюъ
$\rho_k(r)+m$ сєфхЄ єёЄющ\-ўш\-т√ь єЄюў\-э╕э\-э√ь яю\-Ё ф\-ъюь
Їєэъ\-Ўшш $F_{k+m}(r)$.
%%%%%%%%%%%%%%%%%%%%%%%%%%%%%%%%%%%%%%%%%%%%%%%%%%%%%%%%%%%%%%%%%%%%%%%%%%%%%%%%%%%%%%%%%%%%%%%%%%%
\begin{theorem}\label{abel_th_10}\hskip-2mm{.}\:
%%%%%%%%%%%%%%%%%%%%%%%%%%%%%%%%%%%%%%%%%%%%%%%%%%%%%%%%%%%%%%%%%%%%%%%%%%%%%%%%%%%%%%%%%%%%%%%%%%%
╧єёЄ№ $f$ -- ыю\-ъры№\-эю шэ\-ЄхЁуш\-Ёєх\-ьр  Їєэъ\-Ўш 
ъю\-эхў\-эю\-ую яю\-Ё ф\-ър эр яюыє\-юёш $(0,\infty)$ ш
ёє\-∙хёЄ\-тєхЄ шэ\-ЄхуЁры $\int\limits_0^1 f(t)dt$. ╧єёЄ№ Їєэъ\-Ўш 
$f$ эх шьххЄ єёЄющ\-ўш\-тю\-ую єЄюў\-э╕э\-эю\-ую яю\-Ё ф\-ър ш яєёЄ№
$\rho(r)$ -- єЄюў\-э╕э\-э√щ яю\-Ё фюъ Їєэъ\-Ўшш  $f$, яЁш\-ў╕ь
$\rho(\infty)$ $\neq -1$. ╧єёЄ№ $\lambda$  -- ьхЁр ё яыюЄ\-эюёЄ№■
$f(t)$. ╥юуфр $\lambda\in$ $\frak{M}_\infty(\rho(r)+1)$,
$Fr[\rho(r)+1,\lambda]=\{0\}$.
\end{theorem}
%%%%%%%%%%%%%%%%%%%%%%%%%%%%%%%%%%%%%%%%%%%%%%%%%%%%%%%%%%%%%%%%%%%%%%%%%%%%%%%%%%%%%%%%%%%%%%%%%%%%
\par {\sc  ─юърчрЄхы№ёЄтю.}\: ╚ч Єюую, ўЄю $\rho(r)$  ты \-хЄё 
єЄюў\-э╕э\-э√ь яю\-Ё ф\-ъюь Їєэъ\-Ўшш $f(r)$, ёых\-фє\-хЄ, ўЄю
$\lambda\in\frak{M}_\infty(\rho(r)+1)$. ╧єёЄ№  $F$  --
ърэюэш\-ўхё\-ър  яхЁтю\-юсЁрч\-эр  Їєэъ\-Ўшш  $f$ эр яюыє\-юёш
$(0,\infty)$. ╥ръ ъръ Їєэъ\-Ўш   $f$ эх шьххЄ єёЄющ\-ўш\-тюую
єЄюў\-э╕э\-эю\-ую яю\-Ё ф\-ър, Єю
\begin{equation}\label{abel_for_th_10}
\lim\limits_{r\to\infty} \frac{F(r)}{rV(r)}=0.
\end{equation}

\par ╧єёЄ№ $\varphi$  -- яЁю\-шч\-тюы№\-эр  ЇшэшЄ\-эр  схё\-ъю\-эхў\-эю
фшЇ\-ЇхЁхэ\-Ўш\-Ёєх\-ьр  Їєэъ\-Ўш  эр яюыє\-юёш $(0,\infty)$ ш яєёЄ№
$\supp\varphi\subset[a,b]\subset(0,\infty)$. ╚ьххь
$$
\frac{1}{rV(r)}\int\limits_0^\infty
\varphi\left(\frac{t}{r}\right)d\lambda(t)
=\frac{1}{rV(r)}\int\limits_0^\infty
\varphi\left(\frac{t}{r}\right)f(t)dt
$$
$$
= -\frac{1}{r^2V(r)}\int\limits_0^\infty
\varphi'\left(\frac{t}{r}\right)F(t)dt
=-\frac{1}{rV(r)}\int\limits_a^b \varphi'(u)F(ur)du.
$$
╚ч яюыєўхэ\-эю\-ую Ёртхэ\-ёЄ\-тр ш  (\ref{abel_for_th_10})
ёых\-фє\-хЄ, ўЄю
\begin{equation}\label{abel_lim_th_10}
\lim\limits_{r\to\infty} \frac{1}{rV(r)} \int\limits_0^\infty
\varphi\left(\frac{t}{r}\right)d\lambda(t)=0.
\end{equation}

\par ╧єёЄ№ ЄхяхЁ№ $\nu$  -- яЁю\-шч\-тюы№\-эр  ьхЁр шч
ьэю\-цх\-ёЄ\-тр  $Fr[\lambda]$. ╤є∙хёЄ\-тє\-хЄ
яю\-ёыхфю\-тр\-Єхы№\-эюёЄ№ $r_n\to\infty$ Єр\-ър , ўЄю
$\lambda_{r_n}\to\nu$. ┬ч т т Ёртхэ\-ёЄ\-тх  (\ref{abel_lim_th_10})
$r=r_n$, яюыєўшь, ўЄю $\int\limits_0^\infty \varphi(u)d\nu(u)$ $=0$.
╧ю¤Єю\-ьє яю Єхю\-Ёх\-ьх  \ref{measures_th_6} $\:\nu=0$, р чэрўшЄ
$Fr[\lambda]$ $=\{0\}$. ╥хю\-Ёх\-ьр фю\-ърчр\-эр.

%%%%%%%%%%%%%%%%%%%%%%%%%%%%%%%%%%%%%%%%%%%%%%%%%%%%%%%%%%%%%%%%%%%%%%%%%%%%%%%%%%%%%%%%%%%%%%%%%%%
\begin{theorem}\label{abel_th_11}\hskip-2mm{.}\:
%%%%%%%%%%%%%%%%%%%%%%%%%%%%%%%%%%%%%%%%%%%%%%%%%%%%%%%%%%%%%%%%%%%%%%%%%%%%%%%%%%%%%%%%%%%%%%%%%%%
╧єёЄ№ $f$ -- ыю\-ъры№\-эю шэ\-Єху\-Ёш\-Ёєх\-ьр  Їєэъ\-Ўш 
ъю\-эхў\-эю\-ую яю\-Ё ф\-ър эр яюыє\-юёш $(0,\infty)$,
ёє\-∙хёЄ\-тєхЄ шэ\-ЄхуЁры $\int\limits_0^1 f(t)dt$ ш $F$  --
ърэюэш\-ўхё\-ър  яхЁтю\-юсЁрч\-эр  Їєэъ\-Ўш   $f$ эр яюыє\-юёш
$(0,\infty)$. ╧єёЄ№ Їєэъ\-Ўш   $f$ эх шьххЄ єёЄющ\-ўш\-тю\-ую
єЄюў\-э╕э\-эю\-ую яю\-Ё ф\-ър, р $\rho(r)$  -- єЄюў\-э╕э\-э√щ
яю\-Ё фюъ Їєэъ\-Ўшш $f$, яЁш\-ў╕ь $\rho(\infty)\neq -1$. ╥юуфр, хёыш
яю\-Ё фюъ $F$ эх Ёртхэ $-\infty$, Єю ёє\-∙хёЄ\-тє\-хЄ єЄюў\-э╕э\-э√щ
яю\-Ё фюъ $\rho_1(r)$ ¤Єющ Їєэъ\-Ўшш Єр\-ъющ, ўЄю
$$
\lim\limits_{r\to\infty} \frac{V_1(r)}{rV(r)}=0.
$$
\end{theorem}
%%%%%%%%%%%%%%%%%%%%%%%%%%%%%%%%%%%%%%%%%%%%%%%%%%%%%%%%%%%%%%%%%%%%%%%%%%%%%%%%%%%%%%%%%%%%%%%%%%%%
\par {\sc  ─юърчрЄхы№ёЄтю.}\: ╙ЄтхЁц\-фх\-эшх юўхтшф\-эю,
хёыш яю\-Ё фюъ Їєэъ\-Ўшш  $F$ єфютыхЄ\-тю\-Ё \-хЄ эх\-Ёртхэ\-ёЄ\-тє
$\rho_1<\rho(\infty)+1$. ╧ю¤Єю\-ьє сєфхь ёўш\-ЄрЄ№, ўЄю
$\rho_1=\rho(\infty)+1$. ╥ръ ъръ Їєэъ\-Ўш   $f$ эх шьххЄ
єёЄющ\-ўш\-тюую єЄюў\-э╕э\-эю\-ую яю\-Ё ф\-ър, Єю т√\-яюыэ \-хЄ\-ё 
Ёртхэ\-ёЄтю (\ref{abel_for_th_10}). ╚ч ¤Єюую ш Єхю\-Ёх\-ь√ 5
\cite{Grishin_Malyutina} ёых\-фє\-хЄ, ўЄю ёє\-∙хёЄ\-тєхЄ
єЄюў\-э╕э\-э√щ яю\-Ё фюъ $\rho_2(r)$ Єр\-ъющ, ўЄю
т√\-яюы\-э \-■Є\-ё  Ёртхэ\-ёЄ\-тр
$$
\lim\limits_{r\to\infty} \frac{F(r)}{V_2(r)}=0,\qquad
\lim\limits_{r\to\infty} \frac{V_2(r)}{rV(r)}=0.
$$
\par ╘єэъ\-Ўш  $F_1(r)=\frac{F(r)}{V_2(r)}$  ты \-хЄ\-ё 
Їєэъ\-Ўш\-хщ эєых\-тюую яю\-Ё ф\-ър, шьх■\-∙хщ эєых\-тющ яЁх\-фхы эр
схё\-ъю\-эхў\-эюё\-Єш. ╧ю Єхю\-Ёх\-ьх  \ref{order_th_V:V} є эх╕ хёЄ№
єЄюў\-э╕э\-э√щ яю\-Ё фюъ  $\rho_3(r)$, ъюЄю\-Ё√щ тюч\-ЁрёЄрхЄ эр
яюыє\-юёш $(1,\infty)$, яЁш\-ў╕ь $\rho_3(\infty)=0$ (ёыєўрщ
$\rho_3(r)\equiv 0$ эх шёъы■ўр\-хЄё ). ╥юуфр єЄюў\-э╕э\-э√щ
яю\-Ё фюъ $\rho_1(r)=\rho_2(r)+\rho_3(r)$ сєфхЄ єЄюў\-э╕э\-э√ь
яю\-Ё ф\-ъюь Їєэъ\-Ўшш $F(r)$. ╚ьххь
$$
\lim\limits_{r\to\infty} \frac{V_2(r)V_3(r)}{rV(r)}=0.
$$
╥ръшь юсЁр\-чюь $\rho_1(r)$  -- ¤Єю шёъюь√щ єЄюў\-э╕э\-э√щ
яю\-Ё фюъ. ╥хю\-Ёхьр фю\-ърчрэр.

\par ╥хяхЁ№ тюч\-тЁрЄшь\-ё  ъ тюяЁюёє юс юяЁх\-фхых\-эшш яю\-Ё ф\-ър
ЁюёЄр Їєэъ\-Ўшш  $\Psi(r)$, юяЁх\-фхы х\-ьющ Ёртхэ\-ёЄ\-тюь
(\ref{abel_def_Psi}). ┼ёыш $\mu\in$ $\frak{M}_\infty$ $(\rho(r))$,
Єю яю Єхю\-Ёх\-ьх  \ref{abel_th_1} ьэю\-цх\-ёЄ\-тю $L(J,\infty)$
шьххЄ тшф
$$
L(J,\infty)=\left\{\int\limits_0^\infty K(u)d\nu(u):\; \nu\in
Fr[\mu]\right\}.
$$
┼ёыш $L(J,\infty)\neq \{0\}$, Єю Їєэъ\-Ўш  $\Psi$ Ёрё\-Є╕Є эр
схё\-ъю\-эхў\-эюё\-Єш ъръ Їєэъ\-Ўш  $V(r)$. ┼ёыш цх
$L(J,\infty)=\{0\}$, Єю ь√ шьххь ёююЄ\-эю°х\-эшх $\Psi(r)$
$=o(V(r))$, ъюЄю\-Ёюх эх юяЁх\-фхы \-хЄ яю\-Ё фюъ ЁюёЄр Їєэъ\-Ўшш
$\Psi(r)$ эр схё\-ъю\-эхў\-эюё\-Єш.

\par ╨ртхэ\-ёЄтю $L(J,\infty)=\{0\}$ ью\-цхЄ т√\-яюыэ Є№\-ё  яю Ёрч\-ышў\-э√ь
яЁш\-ўшэрь. ╨рё\-ёьюЄ\-Ёшь тэрўр\-ых ёыєўрщ, ъюуфр
$Fr[\mu]\neq\{0\}$, юфэръю  фЁю  $K$ Єръютю, ўЄю фы  ы■сющ ьхЁ√
$\nu\in Fr[\mu]$ т√\-яюыэ \-хЄё  Ёртхэ\-ёЄтю $\int\limits_0^\infty
K(u)d\nu(u)=0$. ╠√ шьххь ёыєўрщ, ъюуфр Єхю\-Ёх\-ьр  \ref{abel_th_1}
эх яюч\-тюы \-хЄ юяЁх\-фхы Є№ яю\-Ё фюъ ЁюёЄр Їєэъ\-Ўшш $\Psi(r)$ эр
схё\-ъю\-эхў\-эюё\-Єш. ═ряЁш\-ьхЁ, тюч\-ьюцхэ трЁш\-рэЄ
$Fr[\mu]=\{u^{\rho-1} du\}$, $\int\limits_0^\infty
K(u)u^{\rho-1}du=0$. ┬  ¤Єюь ёыєўрх эх Єюы№ъю эр° ьхЄюф эх фр╕Є
юЄтхЄр эр тюяЁюё ю яю\-Ё ф\-ъх ЁюёЄр Їєэъ\-Ўшш $\Psi(r)$ эр
схё\-ъю\-эхў\-эюё\-Єш, эю тююс∙х шьх■\-∙хщ\-ё  шэЇюЁ\-ьр\-Ўшш ю  фЁх
$K$ ш ьхЁх  $\mu$ (ь√ чэрхь, ўЄю т√\-яюыэ \-хЄё  ёююЄ\-эю°х\-эшх
$\int\limits_0^\infty K(u)d\nu(u)=0$ фы  $\nu\in Fr[\mu]$)
эх\-фюёЄр\-Єюўэю фы  юяЁх\-фхых\-эш  яю\-Ё ф\-ър ЁюёЄр Їєэъ\-Ўшш
$\Psi$ эр схё\-ъю\-эхў\-эюё\-Єш. ┬  ¤Єюь ёыєўрх фы  юяЁх\-фхых\-эш 
яю\-Ё ф\-ър ЁюёЄр Їєэъ\-Ўшш  $\Psi$ эєцэр фю\-яюыэш\-Єхы№\-эр 
шэЇюЁ\-ьр\-Ўш  ш ёыюц\-эюёЄ№ чрфрўш яю юяЁх\-фхых\-эш■ яю\-Ё ф\-ър
ЁюёЄр Їєэъ\-Ўшш  $\Psi$ чр\-тшёшЄ юЄ ¤Єющ шэЇюЁ\-ьр\-Ўшш. ┬
ўрёЄ\-эюь ёыєўрх $d\mu(t)=t^{\rho-1}dt$ шьххь $\Psi(r)\equiv 0$.

\par ╨рёёьюЄЁшь ЄхяхЁ№ ёыєўрщ $Fr[\mu]=\{0\}$. ┬  ¤Єюь ёыєўрх
Ёртхэ\-ёЄ\-тр  (\ref{abel_F_n}) ьюц\-эю шё\-яюы№чю\-трЄ№ фы 
єыєў°х\-эш  юЎхэ\-ъш $\Psi(r)=o(V(r))$.

\par ╧єёЄ№ Їєэъ\-Ўш  $F_m(t)$ шьххЄ єЄюў\-э╕э\-э√щ яю\-Ё фюъ $\rho_m(r)$.
╥юуфр ьхЁр  $\lambda_m$, $d\lambda_m(t)=F_m(t)dt$ яЁш\-эрф\-ыхцшЄ
ъырёёє  $\frak{M}_\infty(\rho_m(r)+1)$. ╬сю\-чэрўшь
$$
H_m=\left\{\int\limits_0^\infty K^{(m+1)}(u)d\nu(u): \; u\in
Fr[\rho_m(r)+1,\lambda_m]\right\}.
$$
╧ю Єхю\-Ёх\-ьх  \ref{abel_th_1} яЁх\-фхы№\-эюх ьэю\-цх\-ёЄ\-тю
Їєэъ\-Ўшш $(-1)^{m+1}r^{m+1}\Psi(r)/rV_m(r)$ яю эр\-яЁрт\-ых\-эш■
$r\to\infty$ ёют\-ярфр\-хЄ ё  $H_m$. ┼ёыш $H_m\neq \{0\}$, Єю
яю\-ыєўшь, ўЄю Їєэъ\-Ўш  $\Psi(r)$ Ёрё\-Є╕Є эр схё\-ъю\-эхў\-эюё\-Єш
ъръ Їєэъ\-Ўш  $\frac{V_m(r)}{r^m}$. ┼ёыш цх $H_m=\{0\}$, Єю
яю\-ыєўр\-хь ёююЄ\-эю°х\-эшх
$\Psi(r)=o\left(\frac{V_m(r)}{r^m}\right)$.

\par ┼ёыш $H_m=\{0\}$, $Fr[\lambda_m]\neq \{0\}$, Єю тэют№ яюыєўр\-хь
ёыєўрщ, ъюуфр Єхю\-Ёх\-ьр  \ref{abel_th_1}  эх фр╕Є юЄтхЄр эр тюяЁюё
ю яю\-Ё ф\-ъх ЁюёЄр Їєэъ\-Ўшш  $\Psi(r)$ эр схё\-ъю\-эхў\-эюё\-Єш.

\par ╧Ёхф\-яюыю\-цшь ЄхяхЁ№, ўЄю Їєэъ\-Ўш  $F_m(r)$ эх шьххЄ єёЄющ\-ўш\-тюую
єЄюў\-э╕э\-эю\-ую яю\-Ё ф\-ър. ┬  ¤Єюь ёыєўрх яЁш\-ьхэх\-эшх
Єхю\-Ёхь \ref{abel_th_1} ш  \ref{abel_th_10} ъ Їєэъ\-Ўшш
$F_{m+1}(r)$ ш єЄюў\-э╕э\-эюьє яю\-Ё ф\-ъє $\rho_m(r)+1$ фр╕Є, ўЄю
$$
\Psi(r)=o\left(\frac{V_m(r)}{r^m}\right).
$$
╚ч Єхю\-Ёх\-ь√  \ref{abel_th_11} ёых\-фє\-хЄ, ўЄю ёє\-∙хёЄ\-тєхЄ
єЄюў\-э╕э\-э√щ яю\-Ё фюъ $\rho_{m+1}(r)$ Їєэъ\-Ўшш $F_{m+1}(r)$
Єр\-ъющ, ўЄю $V_{m+1}(r)=o(rV_m(r))$. ╧Ёшьхэх\-эшх Єхю\-Ёх\-ь√
\ref{abel_th_1} ъ Їєэъ\-Ўшш $F_{m+1}(r)$ ш єЄюў\-э╕э\-эюьє яю\-Ё фъє
$\rho_{m+1}(r)$ фр╕Є ёююЄ\-эю°х\-эшх
$\Psi(r)=O(r^{-m-1}V_{m+1}(r))$, ъюЄю\-Ёюх ёшы№\-эхх ёююЄ\-эю°х\-эш 
$\Psi(r)=o(r^{-m}V_m(r))$.

\par ╬сючэрўшь ўхЁхч $\rho_n$ яю\-Ё фюъ Їєэъ\-Ўшш $F_n(r)$. CяЁртхф\-ыштр
ёых\-фє\-■\-∙р  Єхю\-Ёхьр.
%%%%%%%%%%%%%%%%%%%%%%%%%%%%%%%%%%%%%%%%%%%%%%%%%%%%%%%%%%%%%%%%%%%%%%%%%%%%%%%%%%%%%%%%%%%%%%%%%%%
\begin{theorem}\label{abel_th_12}\hskip-2mm{.}\:
%%%%%%%%%%%%%%%%%%%%%%%%%%%%%%%%%%%%%%%%%%%%%%%%%%%%%%%%%%%%%%%%%%%%%%%%%%%%%%%%%%%%%%%%%%%%%%%%%%%
╧єёЄ№ $K(t)$  -- ЇшэшЄ\-эюх схё\-ъю\-эхў\-эю
фшЇ\-ЇхЁхэ\-Ўш\-Ёєх\-ьюх  фЁю, $\mu$  -- Ёр\-фю\-эютр ьхЁр
ъю\-эхў\-эю\-ую яю\-Ё ф\-ър эр яюыє\-юёш $(0,\infty)$ эх
эр\-уЁєцр■\-∙р  яюыє\-шэ\-ЄхЁ\-тры $(0,1]$. ╧єёЄ№ Їєэъ\-Ўш 
$\Psi(r)$ юяЁх\-фхы \-хЄё  Ёртхэ\-ёЄ\-тюь  (\ref{abel_def_Psi}), р т
Ёртхэ\-ёЄ\-тх  (\ref{abel_F_n}) Їєэъ\-Ўшш $F_n(t)$ юфэю\-чэрў\-эю
юяЁх\-фхых\-э√ ё яюью∙№■ яЁхф\-ыюцхэ\-эю\-ую Ёрэхх рыую\-ЁшЄьр.
╧єёЄ№ $\rho_n$  -- яю\-Ё фюъ Їєэъ\-Ўшш  $F_n$. ┼ёыш
$$
\mathop{\underline{\lim}}\limits_{m\to\infty}(\rho_m-m)=-\infty,
$$
Єю Їєэъ\-Ўш  $\Psi(r)$ єс√трхЄ эр схё\-ъю\-эхў\-эюё\-Єш с√ёЄЁхх
ы■сющ ёЄхях\-эш  $r$.
\end{theorem}
%%%%%%%%%%%%%%%%%%%%%%%%%%%%%%%%%%%%%%%%%%%%%%%%%%%%%%%%%%%%%%%%%%%%%%%%%%%%%%%%%%%%%%%%%%%%%%%%%%%%
\par {\sc  ─юърчрЄхы№ёЄтю.}\: ╥хю\-Ёхьр  \ref{abel_th_1}, яЁш\-ьхэ╕э\-эр 
ъ Їєэъ\-Ўшш $F_m(r)$ ш єЄюў\-э╕э\-эюьє яю\-Ё ф\-ъє $\rho_m(r)$, фр╕Є
$$
\Psi(r)=O\left(\frac{V_m(r)}{r^m}\right).
$$
╚ч ¤Єюую ёых\-фєхЄ єЄ\-тхЁ\-ц\-фх\-эшх Єхю\-Ёх\-ь√. ╥хю\-Ёхьр
фюърчрэр.

\par ┬ючьюцхэ ёыєўрщ, ўЄю фы  эх\-ъюЄю\-Ёюую $n$ т√\-яюыэ \-хЄ\-ё 
Ёртхэ\-ёЄтю $\rho_n=-\infty$ (эр\-яЁш\-ьхЁ, хёыш
$d\mu(t)=te^{ie^t}dt$). ┬  ¤Єюь ёыєўрх шч Ёртхэ\-ёЄ\-тр
(\ref{abel_F_n}) ш Єхю\-Ёх\-ь√  \ref{abel_th_1} ёых\-фє\-хЄ, ўЄю
Їєэъ\-Ўш  $\Psi(r)$ єс√тр\-хЄ эр схё\-ъю\-эхў\-эюё\-Єш с√ёЄЁхх ы■сющ
ёЄхях\-эш  $r$.

\par ╨рёёьюЄЁшь ЄхяхЁ№ ёыєўрщ, ъюуфр т√\-яюыэ \-хЄ\-ё 
эх\-Ёртхэ\-ёЄ\-тю
\begin{equation}\label{abel_underline_lim}
\mathop{\underline{\lim}}\limits_{m\to\infty}(\rho_m-m)>-\infty.
\end{equation}
╚ч эх\-Ёртхэ\-ёЄ\-тр $\rho_{n+1}\leq \rho_n+1$, фю\-ърчрэ\-эю\-ую т
ыхььх  \ref{abel_lemma_5}, ёых\-фєхЄ, ўЄю $\rho_{n+1}=
\rho_n+1-\varepsilon_n$, уфх $\varepsilon_n\geq 0$. ╚ч
эх\-Ёртхэ\-ёЄ\-тр  (\ref{abel_underline_lim}) ёых\-фє\-хЄ, ўЄю Ё ф
$\sum\limits_{n=0}^\infty \varepsilon_n$ ёїю\-фшЄ\-ё . ╚ч ¤Єюую, т
ётю■ юўхЁхф№, ёых\-фє\-хЄ ёїюфш\-ьюёЄ№ яю\-ёыхфю\-тр\-Єхы№\-эюё\-Єш
$\rho_n-n$. ┬  ¤Єюь ёыєўрх шьххь $\rho_n=n+p+\delta_n$, уфх
$\delta_n\to 0$, $\rho_n>0$ яЁш $n>n_0$. ┼ёыш ёє\-∙хёЄ\-тє\-хЄ
$m>n_0$ Єръюх, ўЄю Їєэъ\-Ўш  $F_m(r)$ шьххЄ єёЄющ\-ўш\-т√щ
єЄюў\-э╕э\-э√щ яю\-Ё фюъ $\rho_m(r)$, Єю шч Єхю\-Ёх\-ь√
\ref{abel_th_9} ёых\-фє\-хЄ, ўЄю фы  ы■сюую $k\geq 0$ Їєэъ\-Ўш 
$F_{m+k}(r)$ сєфхЄ шьхЄ№ єёЄющ\-ўш\-т√щ єЄюў\-э╕э\-э√щ яю\-Ё фюъ
$\rho_{m+k}(r)=\rho_m(r)+k$. ┬  ¤Єюь ёыєўрх ь√ яю\-ыєўшь, ўЄю ышсю
Їєэъ\-Ўш  $\Psi(r)$ Ёрё\-Є╕Є эр схё\-ъю\-эхў\-эюё\-Єш ъръ Їєэъ\-Ўш 
$\frac{V_m(r)}{r^m}$, хёыш $E_m\neq \{0\}$, ышсю, хёыш $E_m=\{0\}$,
Єю ъЁюьх ёююЄ\-эю°х\-эш  $\Psi(r)=
o\left(\frac{V_{m}(r)}{r^m}\right)$ яЁхф\-ыюцхэ\-э√щ ьхЄюф эшўхую эх
фр╕Є.

\par ┼ёыш цх фы  ы■сюую $m>n_0$ Їєэъ\-Ўш  $F_m(r)$ эх шьххЄ
єёЄющ\-ўш\-тюую єЄюў\-э╕э\-эю\-ую яю\-Ё ф\-ър, Єю т  ¤Єюь ёыєўрх,
ъръ ¤Єю ёых\-фє\-хЄ шч Єхю\-Ёх\-ь√  \ref{abel_th_11}, єЄюў\-э╕э\-э√х
яю\-Ё ф\-ъш $\rho_m(r)$ Їєэъ\-Ўшщ $F_m(r)$ ьюц\-эю т√\-сЁрЄ№ Єръшь
юсЁр\-чюь, ўЄю ёююЄ\-эю°х\-эшх
$\Psi(r)=o\left(\frac{V_{m+1}(r)}{r^{m+1}}\right)$ сєфхЄ єЄюўэ Є№
яЁх\-ф√фє\-∙хх $\Psi(r)= o\left(\frac{V_m(r)}{r^m}\right)$. ═шўхую
сюы№\-°хую т  ¤Єюь ёыєўрх яЁхф\-ыюцхэ\-э√щ ьхЄюф эх фр╕Є.

\par ╠√ юс\-ёєфш\-ыш тюч\-ьюц\-эюё\-Єш ьхЄюфр фы  ёыєўр  ЇшэшЄ\-эю\-ую
схё\-ъю\-эхў\-эю фшЇ\-Їх\-Ёхэ\-Ўш\-Ёєх\-ьюую  фЁр  $K$. ╧Ёш
юяЁхфхы╕э\-э√ї єёыю\-тш ї ¤Єш Ёх\-чєы№\-ЄрЄ√ ьюц\-эю
Ёрё\-яЁюёЄ\-Ёр\-эшЄ№ эр эх\-ЇшэшЄ\-э√х  фЁр.

\par ┬ Єхю\-Ёх\-ьх  \ref{abel_th_5} ЄЁхсє\-хЄё , ўЄюс√ ЄЁющър
$(K,\rho(r),\mu)$ єфют\-ыхЄ\-тю\-Ё ыр єёыю\-тш\- ь
эхщ\-ЄЁр\-ыш\-чр\-Ўшш эєы  ш схё\-ъю\-эхў\-эюё\-Єш. ╤є∙хёЄ\-тє■Є
Єхю\-Ёх\-ь√ юс рёшья\-Єю\-Єш\-ўхё\-ъюь яю\-тхфх\-эшш Їєэъ\-Ўшш
$\Psi$ ш фы  ёыєўрхт, ъюуфр ¤Єю єёыю\-тшх эр\-Ёє°р\-хЄ\-ё .

%%%%%%%%%%%%%%%%%%%%%%%%%%%%%%%%%%%%%%%%%%%%%%%%%%%%%%%%%%%%%%%%%%%%%%%%%%%%%%%%%%%%%%%
\begin{theorem}\label{abel_th_12}\hskip-2mm{.}\:
%%%%%%%%%%%%%%%%%%%%%%%%%%%%%%%%%%%%%%%%%%%%%%%%%%%%%%%%%%%%%%%%%%%%%%%%%%%%%%%%%%%%%%%
╧єёЄ№ $\rho(r)$ -- эєыхтющ єЄюў\-э╕э\-э√щ яюЁ \-фюъ, $\mu\in$
$\frak{M}(\rho(r))$, яЁш\-ў╕ь ёє\-∙хёЄ\-тєхЄ яЁх\-фхы
$\lim\limits_{\varepsilon\to 0}\mu([\varepsilon,1])=\mu((0,1])$.
╧єёЄ№ $K(t)$  -- эхяЁх\-Ё√т\-эр  Їєэъ\-Ўш  эр яюыє\-юёш
$[0,\infty)$, $K_1(t)=K(t)-K(0) \chi_{[0,1]}(t)$ ш, ъЁюьх Єюую,
яєёЄ№ ЄЁющър $(K_1,\rho(r),\mu)$ єфют\-ыхЄ\-тю\-Ё \-хЄ єёыю\-тш\- ь
эхщ\-ЄЁр\-ыш\-чр\-Ўшш эєы  ш схё\-ъю\-эхў\-эюё\-Єш. ╥юуфр
т√\-яюыэ \-хЄ\-ё  ёююЄ\-эю°х\-эшх
$$
\int\limits_0^\infty K\left(\frac{t}{r}\right)d\mu(t)
=K(0)\mu((0,r])+\varphi(r),
$$
яЁш\-ў╕ь яЁх\-фхы№\-эюх ьэю\-цхёЄ\-тю Їєэъ\-Ўшш
$\frac{\varphi(r)}{V(r)}$ яю эр\-яЁрт\-ых\-эш■ $r\to\infty$ шьххЄ
тшф
$$
\left\{\int\limits_0^1 \left(K(t)-K(0)\right) d\nu_1(t)
+\int\limits_1^\infty K(t)d\nu_2(t): (\nu_1,\nu_2) \in
\widehat{Fr}(\mu)\right\}.
$$
\end{theorem}
%%%%%%%%%%%%%%%%%%%%%%%%%%%%%%%%%%%%%%%%%%%%%%%%%%%%%%%%%%%%%%%%%%%%%%%%%%%%%%%%%%%%%%%
\par {\sc  ─юърчрЄхы№ёЄтю.}\: ╚ьххь
$$
\int\limits_0^\infty K\left(\frac{t}{r}\right)d\mu(t)=
K(0)\mu([0,r])+\int\limits_0^\infty
K_1\left(\frac{t}{r}\right)d\mu(t).
$$
\par ─ы  ЇєэъЎшш $K_1$ ш ьхЁ√ $\mu$ т√\-яюыэ \-■Єё  тёх єёыю\-тш 
Єхю\-Ёх\-ь√  \ref{abel_th_6}, яЁш\-ьхэх\-эшх ъю\-Єю\-Ёющ
чр\-ърэўш\-тр\-хЄ фю\-ърчр\-Єхы№\-ёЄтю. ╥хю\-Ёхьр фюърчрэр.

\par ┬ю ьэюушї Єхю\-Ёхьрї т ърўхёЄ\-тх шёїюф\-эющ яюё√ыъш т√\-ёЄєяр\-хЄ
ёююЄ\-эю°х\-эшх
\begin{equation}\label{abel_Psi_V}
\Psi(r)=\int\limits_0^\infty K\left(\frac{t}{r}\right)d\mu(t)\sim
MV(r),\qquad r\to\infty.
\end{equation}
╧Ёш юяЁх\-фхы╕э\-э√ї юуЁрэш\-ўх\-эш ї эр Їєэъ\-Ўш■  $K$ ш ьхЁє $\mu$
Їєэъ\-Ўш   $\Psi(r)$ шьххЄ тшф $\Psi(r)=r^{\rho_1(r)}$, уфх
$\rho_1(r)$  -- єЄюў\-э╕э\-э√щ яю\-Ё фюъ.
%%%%%%%%%%%%%%%%%%%%%%%%%%%%%%%%%%%%%%%%%%%%%%%%%%%%%%%%%%%%%%%%%%%%%%%%%%%%%%%%%%%%%%%
\begin{theorem}\label{abel_th_13}\hskip-2mm{.}\:
%%%%%%%%%%%%%%%%%%%%%%%%%%%%%%%%%%%%%%%%%%%%%%%%%%%%%%%%%%%%%%%%%%%%%%%%%%%%%%%%%%%%%%%
╧єёЄ№ $\rho(r)$ -- эєыхтющ єЄюў\-э╕э\-э√щ яюЁ \-фюъ, $\mu$  --
яюыю\-цш\-Єхы№\-эр  ыю\-ъры№\-эю ъю\-эхў\-эр  ьхЁр эр яюыє\-юёш
$(0,\infty)$, $K(t)$  -- эх\-яЁх\-Ё√т\-эю фшЇ\-Їх\-Ёхэ\-Ўш\-Ёєх\-ьр 
ёЄЁюую єс√тр■\-∙р  Їєэъ\-Ўш  эр яюыє\-юёш $(0,\infty)$ ш Єр\-ър ,
ўЄю ЄЁющър $(tK'(t),\rho(t),\mu)$ єфют\-ыхЄ\-тю\-Ё хЄ єёыю\-тш\- ь
эхщ\-ЄЁр\-ыш\-чр\-Ўшш эєы  ш схё\-ъю\-эхў\-эюё\-Єш. ╧єёЄ№
т√\-яюыэ \-хЄ\-ё  ёююЄ\-эю°х\-эшх (\ref{abel_Psi_V}) ш яєёЄ№
яЁюшч\-тюф\-эр  $\Psi'(r)$ т√\-ўшёы \-хЄ\-ё  яю яЁртш\-ыє ╦хщс\-эшЎр
$$
\Psi'(r)=-\frac{1}{r^2} \int\limits_0^\infty
tK'\left(\frac{t}{r}\right)d\mu(t).
$$
╥юуфр, хёыш Їєэъ\-Ўш■ $\rho_1(r)$ юяЁх\-фхышЄ№ Ёртхэ\-ёЄ\-трьш
$r^{\rho_1(r)}=\Psi(r)$ яЁш $r\geq 1$,
$\rho_1\left(\frac{1}{r}\right)$ $=$ $-\rho_1(r)$, Єю $\rho_1(r)$ --
эєыхтющ єЄюў\-э╕э\-э√щ яю\-Ё фюъ.
\end{theorem}
%%%%%%%%%%%%%%%%%%%%%%%%%%%%%%%%%%%%%%%%%%%%%%%%%%%%%%%%%%%%%%%%%%%%%%%%%%%%%%%%%%%%%%%
\par {\sc  ─юърчрЄхы№ёЄтю.}\: ═ряюьэшь, ўЄю Ёртхэ\-ёЄтю
$\rho_1\left(\frac{1}{r}\right)$ $=$ $-\rho_1(r)$ тїюфшЄ т эр°х
юяЁх\-фхых\-эшх эєых\-тюую єЄюў\-э╕э\-эю\-ую яю\-Ё ф\-ър. ╚ьххь
$$
\Psi(2r)-\Psi(r)=\int\limits_0^\infty
\left(K\left(\frac{t}{2r}\right)-
K\left(\frac{t}{r}\right)\right)d\mu(t)
$$
\begin{equation}\label{abel_for_th_13}
\geq \int\limits_r^{2r} \left(K\left(\frac{t}{2r}\right)-
K\left(\frac{t}{r}\right)\right)d\mu(t)\geq m\mu([r,2r]),
\end{equation}
уфх $m=\min\left\{K\left(\frac{u}{2}\right)-K(u):\:
u\in[1,2]\right\}>0$. ╚ч ¤Єюую эх\-Ёртхэ\-ёЄ\-тр ш ёююЄ\-эю°х\-эш 
(\ref{abel_Psi_V}) ёых\-фєхЄ, ўЄю $Fr[\rho(r),\mu]=\{0\}$.

\par ╥хяхЁ№ шч єёыю\-тшщ Єхю\-Ёх\-ь√ ш Єхю\-Ёх\-ь√  \ref{abel_th_5} ёых\-фєхЄ,
ўЄю
$$
\lim\limits_{r\to\infty}\frac{r\Psi'(r)}{V(r)}=0,\qquad
\lim\limits_{r\to\infty}\frac{r\Psi'(r)}{\Psi(r)}=0.
$$
\par ┼ёыш ЄхяхЁ№ $\rho_1(r)$  -- Їєэъ\-Ўш  шч ЇюЁ\-ьє\-ыш\-Ёют\-ъш Єхю\-Ёх\-ь√,
Єю шч фю\-ърчрэ\-эю\-ую ёых\-фє\-хЄ, ўЄю $\rho_1(r)$  -- эєыхтющ
єЄюў\-э╕э\-э√щ яю\-Ё фюъ. ╥хю\-Ёхьр фюърчрэр.

%%%%%%%%%%%%%%%%%%%%%%%%%%%%%%%%%%%%%%%%%%%%%%%%%%%%%%%%%%%%%%%%%%%%%%%%%%
\begin{remark}\label{abel_rem_th_14}\hskip-2mm{.}\:
%%%%%%%%%%%%%%%%%%%%%%%%%%%%%%%%%%%%%%%%%%%%%%%%%%%%%%%%%%%%%%%%%%%%%%%%%%
╨рё\-ёьюЄЁшь Єхю\-Ёх\-ьє  108 шч ъэшуш ╒рЁфш  \cite{Hardy}. ┬ эхщ
єЄ\-тхЁ\-ц\-фр\-хЄ\-ё , ўЄю хёыш $\rho(r)$  -- эєыхтющ
єЄюў\-э╕э\-э√щ яю\-Ё фюъ, р $\mu$  -- яюыю\-цш\-Єхы№\-эр 
ыю\-ъры№\-эю ъю\-эхў\-эр  ьхЁр эр яюыє\-юёш $[0,\infty)$, Єю
ёююЄ\-эю°х\-эш 
$$
\Psi(r)=\int\limits_0^\infty e^{-\frac{t}{r}}d\mu(t)\sim V(r),\qquad
\mu([0,r])\sim V(r)\quad (r\to\infty)
$$
¤ътш\-тр\-ыхэЄ\-э√.
\end{remark}
%%%%%%%%%%%%%%%%%%%%%%%%%%%%%%%%%%%%%%%%%%%%%%%%%%%%%%%%%%%%%%%%%%%%%%%%%%
\par ╧юёьюЄЁшь, ўЄю фр╕Є т Ёрё\-ёьрЄЁш\-трх\-ьюь ёыєўрх
яЁш\-ьхэх\-эшх фю\-ърчрэ\-э√ї эрьш Єхю\-Ёхь.

\par ╚ч Єхю\-Ёх\-ь√  \ref{abel_th_13} ёых\-фє\-хЄ, ўЄю Їєэъ\-Ўш  $\Psi(r)$
яЁхф\-ёЄрт\-ы \-хЄ\-ё  т тшфх $\Psi(r)=r^{\rho_1(r)}=V_1(r)$, уфх
$\rho_1(r)$  -- эєых\-тющ єЄюў\-э╕э\-э√щ яю\-Ё фюъ. ╤є∙хёЄ\-тє\-хЄ
эєыхтющ єЄюў\-э╕э\-э√щ яю\-Ё фюъ $\rho_2(r)$ Єр\-ъющ, ўЄю
т√\-яюыэ \-■Є\-ё  ёююЄ\-эю°х\-эш 
$$
\mathop{\overline{\lim}}\limits_{r\to\infty}
\frac{V_1(2r)-V_1(r)}{V_2(r)} <\infty,\qquad
\lim\limits_{r\to\infty} \frac{V_2(r)}{V_1(r)}=0.
$$
┬ ърўхёЄтх $\rho_2(r)$ яЁю∙х тёхую тч Є№ Їєэъ\-Ўш■, юяЁх\-фхы х\-ьє■
Ёртхэ\-ёЄ\-тюь $r^{\rho_2(r)}$ $=$ $V_1(2r)$ $-V_1(r)$, хёыш Єр\-ър 
Їєэъ\-Ўш   ты \-хЄ\-ё  єЄюў\-э╕э\-э√ь яю\-Ё ф\-ъюь.

\par ╚ьххь
$$
V_1(2r)-V_1(r)=\int\limits_0^\infty e^{-\frac{t}{2r}}
\left(1-e^{-\frac{t}{2r}}\right)d\mu(t)
$$
$$
\geq\int\limits_r^{2r}
e^{-\frac{t}{2r}}\left(1-e^{-\frac{t}{2r}}\right)d\mu(t)\geq
e^{-1}(1-e^{-\frac{1}{2}})\mu((r,2r]).
$$
╚ч ¤Єюую эх\-Ёртхэ\-ёЄ\-тр ёых\-фє\-хЄ, ўЄю $\mu\in$
$\frak{M}_\infty(\rho_2(r))$. ╥хяхЁ№, яЁш\-ьхэ   Єхю\-Ёх\-ьє
\ref{abel_th_6} ъю тЄюЁюьє ёырурх\-ью\-ьє шч яЁртющ ўрёЄш
эшцх\-ёых\-фє\-■\-∙хую Ёртхэ\-ёЄ\-тр
$$
\Psi(r)=\mu([0,r])+\int\limits_0^\infty
\left(e^{-\frac{t}{r}}-\chi_{[0,r]}(t)\right)d\mu(t),
$$
яюыєўшь, ўЄю
$$
\mu([0,r])=\Psi(r)+O(V_2(r)).
$$
\par ▌Єю ёшы№эхх ўхь Єхю\-Ёхьр  108.

%%%%%%%%%%%%%%%%%%%%%%%%%%%%%%%%%%%%%%%%%%%%%%%%%%%%%%%%%%%%%%%%%%%%%%%%%%%%%%%%%%%%%%%%%
\section{╥рєсхЁют√ Єхю\-Ёх\-ь√ фы  шэ\-ЄхуЁрыют}
%%%%%%%%%%%%%%%%%%%%%%%%%%%%%%%%%%%%%%%%%%%%%%%%%%%%%%%%%%%%%%%%%%%%%%%%%%%%%%%%%%%%%%%%%
\qquad ═рўэ╕ь ёю ёых\-фє\-■∙хую юяЁх\-фхых\-эш .
\par  ╧єёЄ№ $\mu$  -- Ёр\-фю\-эю\-тр ьхЁр эр тх∙хёЄ\-тхэ\-эющ юёш,
єфют\-ыхЄ\-тю\-Ё ■\-∙р  єёыю\-тш■ $|\mu|([-t,t])$ $\leq $
$M(1+t^\alpha)$ ё эх\-ъюЄю\-Ё√ь $\alpha\geq 0$. {\it
╧Ёх\-юсЁрчю\-тр\-эшхь ╩рЁых\-ьрэр} ьхЁ√ $\mu$ эр\-ч√тр\-хЄ\-ё  ярЁр
Їєэъ\-Ўшщ $G(z)$ $=(G_+(z)$, $G_-(z))$, уфх
$$
G_+(z)=\int\limits_0^{\;\;\;\infty\;\;\;_\prime} e^{itz}d\mu(t),
\qquad \Im z>0,
$$
$$
G_-(z)=-\int\limits_{-\infty}^{\;\;\;0\;\;\;_\prime} e^{itz}d\mu(t),
\qquad \Im z<0.
$$
\par ╪ЄЁшїш эрф чэрърьш шэ\-ЄхуЁр\-ыют ючэрўр\-■Є, ўЄю шэ\-Єху\-Ёш\-Ёю\-тр\-эшх
яю яюыє\-юё ь $[0,\infty)$, $(-\infty,0]$ тхф╕Є\-ё  эх яю ьхЁх
$\mu$, р яю ьхЁх $\mu-\frac{1}{2}\mu(\{0\})\delta$, уфх $\delta$  --
ьхЁр ─шЁрър (хфшэшў\-эр  ьхЁр, ёю\-ёЁхфю\-Єю\-ўхэ\-эр  т эєых). ┬
ёыєўрх, хёыш $\mu(\{0\})=0$, °ЄЁшїш ьюц\-эю юяєё\-ЄшЄ№.

\par ╬ўхтшфэю, ўЄю $G(z)$  -- ыю\-ъры№эю уюыю\-ьюЁЇ\-эр  Їєэъ\-Ўш  эр
ьэю\-цх\-ёЄ\-тх $\mathbb{C}\setminus\mathbb{R}$ ($\mathbb{C}$  --
ъюья\-ыхъё\-эр  яыюё\-ъюёЄ№, $\mathbb{R}$  -- тх∙хёЄ\-тхэ\-эр  юё№).

\par ═р°х шчыюцх\-эшх Єрє\-сх\-Ёю\-т√ї Єхю\-Ёхь тю ьэю\-уюь юяшЁр\-хЄ\-ё  эр
ёых\-фє\-■\-∙є■ Єхю\-Ёхьє юс рэрыш\-Єш\-ўхё\-ъюь яЁю\-фюы\-цх\-эшш.

%%%%%%%%%%%%%%%%%%%%%%%%%%%%%%%%%%%%%%%%%%%%%%%%%%%%%%%%%%%%%%%%%%%%%%%%%%
\begin{theorem}\label{tauber_th_G_properties}\hskip-2mm{.}\:
%%%%%%%%%%%%%%%%%%%%%%%%%%%%%%%%%%%%%%%%%%%%%%%%%%%%%%%%%%%%%%%%%%%%%%%%%%
╧єёЄ№ $M>0$  -- Їшъёш\-Ёю\-трэ\-эюх ўшёыю, $\mu$  -- Ёр\-фю\-эю\-тр
ьхЁр эр тх∙хёЄ\-тхэ\-эющ юёш Єр\-ър , ўЄю $|\mu|([\alpha,\beta])\leq
M$, хёыш $\beta-\alpha\leq 1$. ╧єёЄ№ сюЁх\-ыхт\-ёър  Їєэъ\-Ўш 
$K(t)$ яЁш\-эрф\-ыхцшЄ ъырёёє $L_1(-\infty,\infty)$ ш
$\widehat{K}(\lambda)$ $=$ $\int\limits_{-\infty}^\infty
K(t)e^{-i\lambda t}dt$ $\neq 0$ яЁш $\lambda\in(a,b)$. ╧єёЄ№
\begin{equation}\label{tauber_K_d_mu}
\int\limits_{-\infty}^\infty K(t-u)d\mu(t)=0,\qquad
u\in(-\infty,\infty),
\end{equation}
$G(z)$ $=(G_+(z),G_-(z))$  -- яЁх\-юсЁрчю\-трэшх ╩рЁых\-ьрэр ьхЁ√
$\mu$.
╥юуфр т√\-яюыэ \-■Єё  ёых\-фє\-■\-∙шх єёыю\-тш :\\
1) \begin{equation}\label{tauber_G_leq} \left|G(z)\right|\leq M
\left(1+\frac{1}{|y|}\right), \qquad z=x+iy,
\end{equation}
 2) Їєэъ\-Ўш  $G_+(z)$ фю\-яєёърхЄ рэрыш\-Єш\-ўхё\-ъюх яЁю\-фюы\-цх\-эшх
ўхЁхч шэ\-ЄхЁ\-тры $(a,b)$, ш Ёх\-чєы№\-Єр\-Єюь рэрыш\-Єш\-ўхё\-ъюую
яЁю\-фюы\-цх\-эш   ты \-хЄё  Їєэъ\-Ўш  $G_-(z)$.
\end{theorem}
%%%%%%%%%%%%%%%%%%%%%%%%%%%%%%%%%%%%%%%%%%%%%%%%%%%%%%%%%%%%%%%%%%%%%%%%%%
\par {\sc  ─юърчрЄхы№ёЄтю.}\: ┬эрўрых ь√ фю\-ърцхь, ўЄю хёыш $d-c=1$,
Єю ёє\-∙хёЄ\-тєхЄ ўшёыю $M_1$ Єръюх, ўЄю фы  ы■сюую
$c\in(-\infty,\infty)$ сєфхЄ т√\-яюыэ Є№\-ё  эх\-Ёртхэ\-ёЄтю
\begin{equation}\label{tauber_I(c,d)}
I(c,d)=\int\limits_c^d \int\limits_{-\infty}^\infty
|K(t-u)|d|\mu|(t)du\leq M_1.
\end{equation}
╚ьххь
$$
I(c,d)= \int\limits_{-\infty}^\infty \int\limits_c^d
|K(t-u)|dud|\mu|(t) =\int\limits_{-\infty}^\infty
\int\limits_{t-d}^{t-c} |K(\tau)|d\tau d|\mu|(t)
$$
$$
= \sum\limits_{n=-\infty}^\infty \int\limits_n^{n+1}
\int\limits_{t-d}^{t-c} |K(\tau)|d\tau d|\mu|(t).
$$
\par ╧ю Єхю\-Ёхьх ю ёЁхф\-эхь чэрўх\-эшш фы  шэ\-Єху\-Ёр\-ыют эрщфєЄ\-ё 
Єюўъш $t_n\in[n,n+1]$ Єръшх, ўЄю сєфхЄ т√\-яюыэ Є№\-ё  Ёртхэ\-ёЄтю
$$
I(c,d)= \sum\limits_{n=-\infty}^\infty \int\limits_{t_n-d}^{t_n-c}
|K(\tau)|d\tau \int\limits_n^{n+1}d|\mu|(t).
$$
\par ┬ ёшыє юуЁрэш\-ўхэшщ эр ьхЁє $\mu$ юЄё■фр ёых\-фєхЄ, ўЄю
$$
I(c,d)\leq M\sum\limits_{n=-\infty}^\infty
\int\limits_{t_n-d}^{t_n-c} |K(\tau)|d\tau
=M\int\limits_{-\infty}^\infty a(\tau)|K(\tau)|d\tau,
$$
уфх $a(\tau)$  -- ¤Єю ўшёыю ёхуьхэ\-Єют $[t_n-d,t_n-c]$,
ёюфхЁцр\-∙шї Єюўъє  $\tau$. ╥хяхЁ№ шч эх\-Ёртхэ\-ёЄ\-тр $a(\tau)\leq
3$ ёых\-фє\-хЄ эх\-Ёртхэ\-ёЄтю  (\ref{tauber_I(c,d)}).

\par ─рыхх сєфхь юЎхэш\-трЄ№ Їєэъ\-Ўш■ $G(z)$. ╧єёЄ№ $\check{\mu}(t)$  --
Їєэъ\-Ўш  Ёрё\-яЁх\-фхых\-эш  ьхЁ√ $|\mu|$, эюЁьш\-Ёю\-трэ\-эр 
єёыю\-тшхь $\check{\mu}(0)=0$. ╚ч єёыю\-тшщ Єхю\-Ёх\-ь√ ёых\-фєхЄ
эх\-Ёртхэ\-ёЄтю $|\check{\mu}(t)|\leq M(1+|t|)$. ╧ю¤Єю\-ьє
$$
\left|G_+(z)\right|\leq \int\limits_0^\infty
e^{-ty}d\check{\mu}(t)=\left.e^{-ty}\check{\mu}(t)\right|_0^\infty+
y\int\limits_0^\infty \check{\mu}(t)e^{-ty}dt.
$$
╚ёяюы№чє  юЎхэъє фы  $\check{\mu}(t)$, эрїюфшь, ўЄю $|G_+(z)|\leq
M\left(1+\frac{1}{y}\right)$. └эрыюушў\-эю юЎхэш\-тр\-хЄё  $G_-(z)$.
╥хь ёрь√ь, єЄ\-тхЁ\-ц\-фх\-эшх  1) фю\-ърчрэю.

\par ╬сючэрўшь  $T_\varepsilon(t)= \left(sin\frac{\varepsilon}{2}t/
\frac{\varepsilon}{2}t\right)^2$ ш юяЁх\-фхышь ьхЁє
$\mu_\varepsilon$ Ёртхэ\-ёЄ\-тюь $d\mu_\varepsilon(t)$ $=$
$T_\varepsilon(t)d\mu(t)$. ╧єёЄ№ $G^\varepsilon(z)$ --
яЁх\-юсЁрчю\-тр\-эшх ╩рЁых\-ьрэр ьхЁ√ $\mu_\varepsilon$. ╬ЄьхЄшь,
ўЄю ьхЁр  $\mu_\varepsilon$ ъю\-эхўэр ш ўЄю тёыхф\-ёЄ\-тшх
эх\-Ёртхэ\-ёЄ\-тр $T_\varepsilon(t)\leq 1$ т эх\-Ёртхэ\-ёЄ\-тх
(\ref{tauber_G_leq}) Їєэъ\-Ўш■  $G$ ьюц\-эю чр\-ьхэшЄ№ эр
$G_\varepsilon$.

\par ┴єфхь ёўш\-ЄрЄ№, ўЄю ўшёыю $\varepsilon$ фю\-ёЄрЄюў\-эю ьрыюх,
ш яєёЄ№ $[a_1,b_1]$ Єр\-ъющ ёху\-ьхэЄ, ўЄю
$a+\varepsilon<a_1<b_1<b-\varepsilon$. ╧єёЄ№ $\xi$  --
яЁю\-шч\-тюы№\-эр  Єюўър ёхуьхэ\-Єр  $[a_1,b_1]$. ═рщф╕ь Їєэъ\-Ўш■
$K_1$ $\in$ $L_1(-\infty,\infty)$ Єръє■, ўЄюс√ т√\-яюыэ \-ыюё№
Ёртхэ\-ёЄтю
\begin{equation}\label{tauber_int_K_1}
\int\limits_{-\infty}^\infty K_1(u)K(t-u)du= T_\varepsilon(t)
e^{i\xi t}.
\end{equation}
─юяєёЄшь, ўЄю Єр\-ър  Їєэъ\-Ўш  ёє\-∙хёЄ\-тєхЄ. ╥юуфр
т√\-яюыэ \-хЄё  Ёртхэ\-ёЄтю
\begin{equation}\label{tauber_int_int_K_1}
\int\limits_{-\infty}^\infty \int\limits_{-\infty}^\infty K_1(u)
K(t-u)e^{-itx} dudt = \widehat{T}(\xi,x),
\end{equation}
уфх
$$
\widehat{T}(\xi,x)= \int\limits_{-\infty}^\infty T_\varepsilon(t)
e^{i\xi t}e^{-ixt}dt.
$$
\par ╚чтхёЄэю, ўЄю шэ\-ЄхуЁры т ыхтющ ўрёЄш
Ёртхэ\-ёЄ\-тр (\ref{tauber_int_int_K_1}) ёїю\-фшЄ\-ё  рсёюы■Є\-эю.
╧ю¤Єю\-ьє т  ¤Єюь шэ\-Єху\-Ёрых тюч\-ьюц\-эю шч\-ьхэх\-эшх
яю\-Ё ф\-ър шэ\-Єху\-Ёш\-Ёю\-тр\-эш . ▌Єю фр╕Є
$$
\widehat{T}(\xi,x)= \int\limits_{-\infty}^\infty
\int\limits_{-\infty}^\infty K(t-u) e^{-itx} dt
K_1(u)du=\widehat{K}(x)\widehat{K}_1(x),
$$
\begin{equation}\label{tauber_hat_K_1}
\widehat{K}_1(x)=\frac{\widehat{T}(\varepsilon,x)}{\widehat{K}(x)}.
\end{equation}
\par ╥хяхЁ№, чрс√тр  ю Єюь, ъръшь ёяюёюсюь яюыєўхэю
Ёртхэ\-ёЄтю  (\ref{tauber_hat_K_1}), Ёрё\-ёьюЄЁшь Їєэъ\-Ўш■
$\widehat{K}_1(x)$, юяЁх\-фхы \-хьє■ ¤Єшь Ёртхэ\-ёЄ\-тюь. ═юёшЄхы№
Їєэъ\-Ўшш $\widehat{T}(\varepsilon,x)$ Ёрё\-яюыю\-цхэ эр ёху\-ьхэЄх
$[\xi-\varepsilon,\xi+\varepsilon]$, ъюЄю\-Ё√щ  ты \-хЄё  ўрёЄ№■
шэ\-ЄхЁ\-трыр $(a,b)$. ═р шэ\-ЄхЁ\-трых $(a,b)$ Їєэъ\-Ўш 
$\widehat{K}(x)$ эх юсЁр∙р\-хЄё  т эюы№. ╥ръшь юсЁр\-чюь, Їєэъ\-Ўш 
$\widehat{K}_1(x)$ ъюЁЁхъЄ\-эю юяЁх\-фхыхэр эр ёхуьхэ\-Єх
$[\xi-\varepsilon,\xi+\varepsilon]$. ┴єфхь ёўш\-ЄрЄ№, ўЄю Їєэъ\-Ўш 
$\widehat{K}_1(x)$ Ёртэр эєы■ тэх ¤Єюую ёху\-ьхэЄр. ╥хь ёрь√ь
$\widehat{K}_1(x)$ -- ¤Єю эх\-яЁх\-Ё√т\-эр  ЇшэшЄ\-эр  Їєэъ\-Ўш .

\par ╚ч Єхю\-Ёх\-ь√ ┬шэхЁр ю фхых\-эшш ёых\-фє\-хЄ, ўЄю хёыш Їєэъ\-Ўш 
$\widehat{K}_1(x)$ юяЁх\-фхы \-хЄё  Ёртхэ\-ёЄ\-тюь
(\ref{tauber_hat_K_1}), Єю ёє\-∙хёЄ\-тє\-хЄ Їєэъ\-Ўш  $K_1$ $\in$
$L_1(-\infty,\infty)$ Єр\-ър , ўЄю х╕ яЁх\-юсЁрчю\-тр\-эшх ╘єЁ№х
ёют\-ярфрхЄ ё $\widehat{K}_1(x)$. ╥хяхЁ№ шч Ёртхэ\-ёЄ\-тр
$\widehat{K}_1(x)\widehat{K}(x)$ $=\widehat{T}(\varepsilon,x)$
яЁш\-ьхэх\-эшхь юсЁрЄ\-эю\-ую яЁх\-юсЁрчю\-трэш  ╘єЁ№х
яю\-ыєўр\-хЄё  Ёртхэ\-ёЄтю (\ref{tauber_int_K_1}).
╤є∙хёЄ\-тю\-тр\-эшх Їєэъ\-Ўшш $K_1$ фю\-ърчр\-эю.

\par ╙ьэюцшь юсх ўрёЄш Ёртхэ\-ёЄ\-тр  (\ref{tauber_K_d_mu}) эр  $K_1(u)$
ш яЁю\-шэ\-Єху\-Ёш\-Ёєхь юсх ўрёЄш яю\-ыєўхэ\-эю\-ую Ёртхэ\-ёЄ\-тр
яю юёш $(-\infty,\infty)$. ╠√ яюыєўшь
\begin{equation}\label{tauber_int_0}
\int\limits_{-\infty}^\infty \int\limits_{-\infty}^\infty K_1(u)
K(t-u)d\mu(t)du=0.
\end{equation}
\par ┴єфхь юЎхэш\-трЄ№ ёых\-фє\-■\-∙шщ шэ\-ЄхуЁры
$$
I=\int\limits_{-\infty}^\infty \int\limits_{-\infty}^\infty |K_1(u)|
|K(t-u)|dud|\mu|(t).
$$
\par ╚ьххь
$$
I=\sum\limits_{n=-\infty}^\infty \int\limits_{-\infty}^\infty
\int\limits_n^{n+1} |K_1(u)| |K(t-u)|dud|\mu|(t)
$$
$$
\leq \sum\limits_{n=-\infty}^\infty \max\{|K_1(u)|:\: u\in[n,n+1]\}
\int\limits_{-\infty}^\infty \int\limits_n^{n+1} |K(t-u)|du
d|\mu|(t).
$$
\par ╧Ёшьхэ   эх\-Ёртхэ\-ёЄтю  (\ref{tauber_I(c,d)}), яю\-ыєўшь
эх\-Ёртхэ\-ёЄтю
$$
I\leq M_1 \sum\limits_{n=-\infty}^\infty \max\{|K_1(u)|:\:
u\in[n,n+1]\}.
$$
╘єэъЎш  $K_1\in L_1(-\infty,\infty)$  ты \-хЄё 
яЁх\-юсЁрчю\-тр\-эшхь ╘єЁ№х эх\-яЁх\-Ё√т\-эющ ЇшэшЄ\-эющ Їєэъ\-Ўшш.
─ы  Єръшї Їєэъ\-Ўшщ, ъръ ёых\-фєхЄ шч ыхьь√  $6_7$ (\cite{Vinner},
уыртр  2, \S11), эр\-яшёрэ\-э√щ т√°х Ё ф ёїю\-фшЄ\-ё . ╥ръшь
юсЁр\-чюь, $I<\infty$. ╥хяхЁ№ шч Єхю\-Ёхь ╥юэхыыш ш ╘єсшэш
ёых\-фє\-хЄ, ўЄю т шэ\-Єху\-Ёрых, ёЄю \-∙хь т ыхтющ ўрёЄш
Ёртхэ\-ёЄ\-тр (\ref{tauber_int_0}), ьюц\-эю яю\-ьхэ Є№ яю\-Ё фюъ
шэ\-Єху\-Ёш\-Ёю\-тр\-эш . ▌Єю фр╕Є ё єў╕Єюь Ёртхэ\-ёЄ\-тр
(\ref{tauber_int_K_1}) Єръюх Ёртхэ\-ёЄ\-тю

\begin{equation}\label{tauber_int_T}
\int\limits_{-\infty}^\infty e^{i\xi t}T_\varepsilon(t) d\mu(t)=0,
\qquad \xi\in[a_1,b_1].
\end{equation}

\par ╚ч ъю\-эхў\-эюё\-Єш ьхЁ√ $\mu_\varepsilon$ ёых\-фєхЄ, ўЄю Їєэъ\-Ўш 
$G^\varepsilon_+(z)$  ты \-хЄ\-ё  уюыю\-ьюЁЇ\-эющ Їєэъ\-Ўшхщ т
тхЁї\-эхщ яюыє\-яыюё\-ъюёЄш эх\-яЁх\-Ё√т\-эющ тяыюЄ№ фю уЁрэш\-Ў√.
└эрыюушў\-эю, Їєэъ\-Ўш  $G^\varepsilon_-(z)$  ты \-хЄ\-ё  Їєэъ\-Ўшхщ
уюыю\-ьюЁЇ\-эющ т эшц\-эхщ яюыє\-яыюё\-ъюёЄш эх\-яЁх\-Ё√т\-эющ
тяыюЄ№ фю уЁрэш\-Ў√. ╥юуфр Ёртхэ\-ёЄтю  (\ref{tauber_int_T}) ьюц\-эю
яхЁх\-яшёрЄ№ т тшфх $G^\varepsilon_+(\xi)=G^\varepsilon_-(\xi)$,
$\xi\in[a_1,b_1]$. ╧ю Єхю\-Ёхьх ю ёЄшЁр\-эшш юёюсхэ\-эюё\-Єхщ
(\cite{Evgrafov}, Єхю\-Ёх\-ьр  2.2, уыртр  4) Їєэъ\-Ўшш
$G^\varepsilon_+(z)$ ш  $G^\varepsilon_-(z)$  ты \-■Єё 
рэрыш\-Єш\-ўхё\-ъшьш яЁю\-фюы\-цх\-эш \-ьш фЁєу фЁєур ўхЁхч
шэ\-ЄхЁ\-тры $(a_1,b_1)$.

\par ─рыхх Ёрё\-ёьюЄ\-Ёшь т ъюьяыхъё\-эющ яыюё\-ъюёЄш ътрф\-ЁрЄ $Q_1$, фы 
ъюЄю\-Ёюую ёху\-ьхэЄ $[a_1,b_1]$  ты \-хЄ\-ё  фшрую\-эры№■, ш
Ёрё\-ёьюЄ\-Ёшь ёхьхщ\-ёЄтю Їєэъ\-Ўшщ
$$
F_\eta(z)=(z-a_1)(z-b_1)G^\eta(z), \qquad \eta\in(0,\varepsilon],
$$
т  ¤Єюь ътрф\-ЁрЄх. ╬фэр шч ёЄюЁюэ ътрф\-ЁрЄр $Q_1$ шьххЄ
ярЁр\-ьхЄЁш\-чрЎш■ $z=a_1+te^\frac{i\pi}{4}$,
$t\in\left[0,\frac{b_1-a_1}{\sqrt{2}}\right]$. ═р ¤Єющ ёЄюЁюэх
ътрф\-ЁрЄр т√\-яюыэ \-хЄё  эх\-Ёртхэ\-ёЄтю
$$
|F_\eta(z(t))|\leq Mt\left|a_1-b_1+te^\frac{i\pi}{4}\right|
\left|1+\frac{\sqrt{2}}{t}\right|\leq M_2.
$$
\par ╬ўхтшфэю, ўЄю эр уЁрэш\-Ўх ътрф\-ЁрЄр  т√\-яюыэ \-хЄё  эх\-Ёртхэ\-ёЄ\-тю
$|F_\eta(z)|\leq M_2$.

\par ╧єёЄ№ $\delta\in\left(0,\frac{b_1-a_1}{4}\right)$ ш яєёЄ№ $Q_2$  --
ътрфЁрЄ, фы  ъюЄю\-Ёюую ёху\-ьхэЄ $[a_1+\delta,b_1-\delta]$
 ты \-хЄё  фшрую\-эр\-ы№■. ═р уЁрэшЎх ътрф\-ЁрЄр $Q_2$ сєфхЄ
т√\-яюыэ Є№\-ё  эх\-Ёртхэ\-ёЄтю $|G^\eta(z)|\leq M_3(\delta)$ ё
эх\-ъю\-Єю\-Ёющ тхышўш\-эющ $M_3$, чр\-тшё \-∙хщ юЄ  $\delta$.

\par ╧ю Єхю\-Ёхьх ╠юэЄхы  ёхьхщ\-ёЄтю Їєэъ\-Ўшщ $G^\eta(z)$ сєфхЄ
ъюь\-яръЄ\-э√ь тэєЄЁш $Q_2$.

\par ╧ю¤Єю\-ьє ёє\-∙хёЄ\-тєхЄ уюыю\-ьюЁЇ\-эр  тэєЄЁш $Q_2$ Їєэъ\-Ўш  $H(z)$ ш
яю\-ёыхфю\-тр\-Єхы№\-эюёЄ№ $\eta_n\to 0$ Єръшх, ўЄю
яю\-ёыхфю\-тр\-Єхы№\-эюёЄ№ $G^{\eta_n}(z)$ сєфхЄ Ёртэю\-ьхЁ\-эю эр
ъюьяръ\-Єрї, ыхцр\-∙шї тю тэєЄ\-Ёхэ\-эюё\-Єш $Q_2$, ёїю\-фшЄ\-ё  ъ
Їєэъ\-Ўшш  $H(z)$. ╧Ёш $\Im z>0$ яю\-ёыхфю\-тр\-Єхы№\-эюёЄ№
$G^{\eta_n}(z)$ ёїю\-фшЄ\-ё  ъ $G_+(z)$, р яЁш $\Im z<0$ юэр
ёїю\-фшЄ\-ё  ъ  $G_-(z)$. ╚ч ёърчрэ\-эю\-ую ёых\-фєхЄ, ўЄю Їєэъ\-Ўш 
$G_-(z)$ сєфхЄ рэрыш\-Єш\-ўхё\-ъшь яЁю\-фюы\-цх\-эшхь Їєэъ\-Ўшш
$G^+(z)$ ўхЁхч шэ\-ЄхЁ\-тры $(a_1+\delta,b_1-\delta)$, р т ёшыє
яЁю\-шч\-тюы№\-эюё\-Єш $\varepsilon$ ш  $\delta$ ўхЁхч шэ\-ЄхЁ\-тры
$(a,b)$. ╥хю\-Ёхьр фю\-ърчрэр.

%%%%%%%%%%%%%%%%%%%%%%%%%%%%%%%%%%%%%%%%%%%%%%%%%%%%%%%%%%%%%%%%%%%%%%%%%%%%%%%
\begin{remark}\label{tauber_rem_th_G}\hskip-2mm{.}\:
%%%%%%%%%%%%%%%%%%%%%%%%%%%%%%%%%%%%%%%%%%%%%%%%%%%%%%%%%%%%%%%%%%%%%%%%%%%%%%%
╥хю\-Ёхьр  \ref{tauber_th_G_properties}  ты \-хЄ\-ё  єёшых\-эш\-хь
ыхь\-ь√ ╩рЁых\-ьр\-эр юс рэр\-ыш\-Єш\-ўхё\-ъюь яЁю\-фюы\-цх\-эшш. ┬
юЁшуш\-эры№\-эющ ыхь\-ьх ╩рЁых\-ьр\-эр (фю\-ърчр\-Єхы№\-ёЄ\-тю
ыхь\-ь√ ╩рЁых\-ьр\-эр ьюцэю эрщЄш Єръ\-цх т \cite{Gurarii})
Ёрё\-ёьрЄ\-Ёш\-тр\-■Є\-ё  ьхЁ√ ё юуЁрэш\-ўхэ\-эющ яыюЄ\-эюё\-Є№■.
╧Ёш\-тхфхэ\-эюх фю\-ърчр\-Єхы№\-ёЄ\-тю Єхю\-Ёх\-ь√
\ref{tauber_th_G_properties} ўрё\-Єшў\-эю ёют\-ярфр\-хЄ ё
фю\-ърчр\-Єхы№\-ёЄ\-тюь ыхьь√ ╩рЁых\-ьрэр.
\end{remark}
%%%%%%%%%%%%%%%%%%%%%%%%%%%%%%%%%%%%%%%%%%%%%%%%%%%%%%%%%%%%%%%%%%%%%%%%%%%%%%%
% \par ═рь сєфхЄ эєцэю ёых\-фє\-■\-∙хх юяЁх\-фхых\-эшх.

\par ╧єёЄ№ $\mu$  -- Ёр\-фю\-эю\-тр ьхЁр ъю\-эхў\-эю\-ую яю\-Ё ф\-ър эр
тх∙хёЄ\-тхэ\-эющ юёш, $G(z)$ $=$ $(G_+(z),$ $G_-(z))$  --
яЁх\-юсЁрчю\-тр\-эшх ╩рЁых\-ьрэр ьхЁ√ $\mu$. ┬√сЁюёшь шч
тх∙хёЄ\-тхэ\-эющ юёш ьэю\-цх\-ёЄ\-тю тёхї Єръшї шэ\-ЄхЁ\-трыют
$(a,b)$, ўЄю Їєэъ\-Ўш   $G_+(z)$ рэрыш\-Єш\-ўхё\-ъш
яЁю\-фюыцр\-хЄ\-ё  ўхЁхч шэ\-ЄхЁ\-тры $(a,b)$ ш Ёх\-чєы№\-Єр\-Єюь
рэрыш\-Єш\-ўхё\-ъюую яЁю\-фюы\-цх\-эш  хёЄ№ Їєэъ\-Ўш  $G_-(z)$.
╬ёЄрт\-°р \-ё  ўрёЄ№ тх∙хёЄ\-тхэ\-эющ юёш эр\-ч√тр\-хЄ\-ё  {\it
ёяхъЄ\-Ёюь ╩рЁых\-ьрэр} ьхЁ√ $\mu$.

\par ┬ ёт чш ё яЁш\-тхфхэ\-э√ь юяЁх\-фхых\-эшхь ёЄюшЄ чр\-ьхЄшЄ№,
ўЄю Їєэъ\-Ўш  $G_+(z)$ ью\-цхЄ шьхЄ№ рэрыш\-Єш\-ўхё\-ъюх
яЁю\-фюы\-цх\-эшх ўхЁхч шэ\-ЄхЁ\-тры $(a,b)$, юфэръю, Ёх\-чєы№\-ЄрЄ
Єръюую яЁю\-фюы\-цх\-эш  ью\-цхЄ эх ёют\-ярфрЄ№ ё Їєэъ\-Ўшхщ
$G_-(z)$. ┬ ¤Єюь ёыєўрх шэ\-ЄхЁ\-тры $(a,b)$ яЁш\-эрф\-ыхцшЄ ёяхъЄЁє
╩рЁых\-ьрэр ьхЁ√ $\mu$.

\par ┬ ёт чш ё Єхю\-Ёх\-ьющ  \ref{tauber_th_G_properties} ш юяЁх\-фхых\-эшхь
ёяхъЄЁр ╩рЁых\-ьрэр шэ\-ЄхЁхё\-эю юЄ\-ьхЄшЄ№ ёых\-фє\-■\-∙хх
єЄ\-тхЁ\-ц\-фх\-эшх, юЄ\-эюё \-∙ххё  ъ яЁюс\-ыхьх урЁью\-эшўхё\-ъюую
ёшэЄхчр. ▌Єр яЁюс\-ых\-ьр юсёєц\-фр\-хЄё  т  \cite{Ross},
\cite{Gurarii}.

%%%%%%%%%%%%%%%%%%%%%%%%%%%%%%%%%%%%%%%%%%%%%%%%%%%%%%%%%%%%%%%%%%%%%%%%%%%%%%%
\begin{theorem}\label{tauber_th_karleman}\hskip-2mm{.}\:
%%%%%%%%%%%%%%%%%%%%%%%%%%%%%%%%%%%%%%%%%%%%%%%%%%%%%%%%%%%%%%%%%%%%%%%%%%%%%%%
╧єёЄ№ $\mu$  -- Ёр\-фю\-эю\-тр ьхЁр ъю\-эхў\-эю\-ую яю\-Ё ф\-ър эр
тх∙хёЄ\-тхэ\-эющ юёш ё юуЁрэш\-ўхэ\-э√ь ёяхъЄЁюь ╩рЁых\-ьрэр. ╥юуфр
ьхЁр  $\mu$ рсёюы■Є\-эю эх\-яЁх\-Ё√т\-эр ш х╕ яыюЄ\-эюёЄ№ $g(t)$
хёЄ№ юуЁрэш\-ўх\-эшх эр тх∙хёЄ\-тхэ\-эє■ юё№ Ўхыющ Їєэъ\-Ўшш
¤ъёяю\-эхэ\-Ўшры№\-эю\-ую Єшяр.
\end{theorem}
%%%%%%%%%%%%%%%%%%%%%%%%%%%%%%%%%%%%%%%%%%%%%%%%%%%%%%%%%%%%%%%%%%%%%%%%%%%%%%%
\par \quad {\sc ─юърчрЄхы№ёЄтю.}\: ╧єёЄ№ $G(z)$ $=$ $(G_+(z),G_-(z))$  --
яЁх\-юсЁрчю\-тр\-эшх ╩рЁых\-ьрэр ьхЁ√ $\mu$. ╧ю\-ёъюы№\-ъє ьхЁр
$\mu$ шьххЄ ъю\-эхў\-э√щ яю\-Ё фюъ, Єю ёє\-∙хёЄ\-тє\-■Є ўшёыр $M>0$
ш $\alpha>0$ Єръшх, ўЄю фы  $\check{\mu}(r)$ (Їєэъ\-Ўшш
Ёрё\-яЁх\-фхых\-эш  ьхЁ√ $|\mu|$) сєфхЄ т√\-яюыэ Є№\-ё 
эх\-Ёртхэ\-ёЄ\-тю $|\check{\mu}(r)|\leq M(1+r^\alpha)$. ╤ўшЄр , ўЄю
Їєэъ\-Ўш  $\check{\mu}(r)$ эюЁьш\-Ёю\-трэр єёыю\-тш\-хь
$\check{\mu}(+0)=0$, яю\-ыєўшь эх\-Ёртхэ\-ёЄтю
$$
\left|G_+(z)\right| \leq \frac{1}{2} \left|\mu(\{0\})\right|+
\int\limits_{(0,\infty)} e^{-ty}d\check{\mu}(t)= \frac{1}{2}
|\mu(\{0\})| + y\int\limits_0^\infty \check{\mu}(t)e^{-ty}dt
$$
$$
\leq \frac{1}{2}|\mu|(\{0\})+ My\int\limits_0^\infty (1+t^\alpha)
e^{-ty}dt \leq M_1\left(1+\frac{\Gamma(1+\alpha)}{y^\alpha}\right).
$$
\par └эрыюушўэю юЎхэш\-тр\-хЄё  Їєэъ\-Ўш  $G_-(z)$.
╧ю¤Єю\-ьє Їєэъ\-Ўш  $G(z)$ єфют\-ыхЄ\-тю\-Ё хЄ юЎхэ\-ъх

\begin{equation}\label{tauber_|G|}
|G(z)|\leq M_1\left(1+\frac{\Gamma(1+\alpha)}{|y|^\alpha}\right).
\end{equation}
╚ч юуЁрэш\-ўхэ\-эюё\-Єш ёяхъЄЁр ьхЁ√ $\mu$ ёых\-фє\-хЄ, ўЄю Єюўър
$\infty$  ты \-хЄ\-ё  шчюыш\-Ёю\-трэ\-эющ юёюсхэ\-эюёЄ№■ Їєэъ\-Ўшш
$G(z)$. ╨рёёєцфх\-эш  ё ътрфЁр\-Єрьш $Q_1$ ш  $Q_2$ т ЄхъёЄх
фю\-ърчр\-Єхы№\-ёЄтр Єхю\-Ёх\-ь√  \ref{tauber_th_G_properties},
яЁш\-ьхэ╕э\-э√х ъ Їєэъ\-Ўшш $(z-a_1)^n(z-b_1)^nG(z)$, фрфєЄ, ўЄю
ёє\-∙хёЄ\-тє■Є ўшёыр $\alpha>0$, $\beta>0$, $M_2>0$ Єръшх, ўЄю эр
ьэю\-цх\-ёЄ\-тх $\{z=x+iy:\: |x|\geq \alpha,\:|y|\leq \beta\}$ сєфхЄ
т√\-яюыэ Є№\-ё  эх\-Ёртхэ\-ёЄ\-тю $|G(z)|\leq M_2$. ╚ч ¤Єюую
эх\-Ёртхэ\-ёЄ\-тр ш эх\-Ёртхэ\-ёЄ\-тр  (\ref{tauber_|G|})
ёых\-фє\-хЄ, ўЄю Єюўър $\infty$  -- єёЄЁрэш\-ьр  юёюсхэ\-эюёЄ№
Їєэъ\-Ўшш $G(z)$.

\par ╦хуъю тшфхЄ№, ўЄю
$$
\lim\limits_{y\to +\infty} G_+(iy)=\lim\limits_{y\to +\infty}
\int\limits_0^{\;\;\;\infty\;\;\;_\prime}
e^{-ty}d\mu(t)=\frac{1}{2}\mu(\{0\}),
$$
$$
\lim\limits_{y\to -\infty} G_-(iy)=-\lim\limits_{y\to -\infty}
\int\limits_{-\infty}^{\;\;\;0\;\;\;_\prime}
e^{ty}d\mu(t)=-\frac{1}{2}\mu(\{0\}).
$$
╚ч ¤Єшї ЁртхэёЄт тшфэю, ўЄю ёююЄ\-эю°х\-эшх $\mu(\{0\})$ $\neq$ $0$
яЁю\-Єштю\-ЁхўшЄ Єюьє, ўЄю $\infty$  -- єёЄЁрэш\-ьр  юёюсхэ\-эюёЄ№
Їєэъ\-Ўшш $G(z)$. ╧ю¤Єю\-ьє $\mu(\{0\})=0$. ╥хь ёрь√ь ь√
фю\-ърчр\-ыш, ўЄю $G(\infty)=0$. ╧ю¤Єю\-ьє ёє\-∙хёЄ\-тє■Є ўшёыр
$M_3>0$ ш  $R_1>0$ Єръшх, ўЄю яЁш $|z|>R_1$ сєфхЄ т√\-яюы\-э Є№\-ё 
эх\-Ёртхэ\-ёЄтю

\begin{equation}\label{tauber_|G|_leq}
|G(z)|<\frac{M_3}{|z|}.
\end{equation}

\par ╬яЁхфхышь ЄхяхЁ№ Їєэъ\-Ўш■

\begin{equation}\label{tauber_def_g}
g(t)=-\frac{1}{2\pi}\int\limits_{\frak{L}} G(w)e^{-itw}dw,
\end{equation}
уфх $\frak{L}$  -- чрьъэє\-Єр  уырфър  яЁртшы№\-эю
юЁшхэЄш\-Ёю\-трэ\-эр  цюЁфрэю\-тр ъЁштр , юїтрЄ√\-тр■\-∙р  ёяхъЄЁ
ьхЁ√ $\mu$. ╥хь ёрь√ь Їєэъ\-Ўш  $G(w)$  ты \-хЄё  уюыю\-ьюЁЇ\-эющ т
чрь√ър\-эшш эх\-юуЁрэш\-ўхэ\-эющ юсырёЄш ё уЁрэш\-Ўхщ  $\frak{L}$.
╬ўхтшф\-эю, ўЄю $g(t)$  -- Ўхыр  Їєэъ\-Ўш  ¤ъёяю\-эхэ\-Ўшры№\-эю\-ую
Єшяр.

\par ─рыхх Ёрё\-ёьюЄЁшь яЁш $\Im z>0$ Їєэъ\-Ўш■
$$
G_1(z)=\int\limits_0^\infty g(t)e^{itz}dt=-\frac{1}{2\pi}
\int\limits_{\frak{L}} \int\limits_0^\infty e^{it(z-w)}dtG(w)dw.
$$
╧єёЄ№ $\Im z\geq 2h$, $h>0$ ш яЁюшч\-тюы№эю. ╧ю\-ёъюы№ъє ёяхъЄЁ ьхЁ√
$\mu$ ыхцшЄ эр тх∙хёЄ\-тхэ\-эющ юёш, ъюэЄєЁ  $\frak{L}$ ьюц\-эю
т√\-сЁрЄ№ Єръ, ўЄюс√ фы  ы■сюую $w\in\frak{L}$ т√\-яюыэ \-ыюё№
эх\-Ёртхэ\-ёЄтю $\Im w<h$. ╥юуфр
$$
G_1(z)=-\frac{1}{2\pi i} \int\limits_{\frak{L}} \frac{G(w)}{w-z}dw.
$$
╥хяхЁ№ шч эхЁртхэ\-ёЄ\-тр \ref{tauber_|G|_leq} ш Єхю\-Ёх\-ь√ ╩ю°ш
ёых\-фє\-хЄ, ўЄю т√\-яюыэ \-хЄ\-ё  Ёртхэ\-ёЄтю $G_1(z)$ $=$ $G(z)$
яЁш $\Im z\geq 2h$, р чэрўшЄ ш яЁш $\Im z>0$.

\par ╠√ фюърчрыш, ўЄю є ьхЁ√ $\mu_1$, $d\mu_1(t)=g(t)dt$,
ш ьхЁ√ $\mu$ ёют\-ярфр\-■Є яЁх\-юсЁрчю\-тр\-эш  ╩рЁых\-ьрэр т
яюыє\-яыюё\-ъюёЄш $\Im z>0$. └эрыюушў\-эю фю\-ърч√\-тр\-хЄё 
ёют\-ярфх\-эшх яЁх\-юсЁрчю\-тр\-эшщ ╩рЁых\-ьрэр ¤Єшї ьхЁ ш т
яюыє\-яыюё\-ъюёЄш $\Im z<0$. ╥ръшь юсЁр\-чюь, є ьхЁ $\mu$ ш  $\mu_1$
ёют\-ярфр\-■Є яЁх\-юсЁрчю\-тр\-эш  ╩рЁых\-ьрэр. ╤ыхфю\-тр\-Єхы№\-эю,
$\mu=\mu_1$. ╥хю\-Ёхьр фю\-ърчрэр.

%%%%%%%%%%%%%%%%%%%%%%%%%%%%%%%%%%%%%%%%%%%%%%%%%%%%%%%%%%%%%%%%%%%%%%%%%%%%%%%
\begin{theorem}\label{tauber_th_d_mu}\hskip-2mm{.}\:
%%%%%%%%%%%%%%%%%%%%%%%%%%%%%%%%%%%%%%%%%%%%%%%%%%%%%%%%%%%%%%%%%%%%%%%%%%%%%%%
╧єёЄ№ $\mu$  -- Ёр\-фю\-эю\-тр ьхЁр ъю\-эхў\-эю\-ую яю\-Ё ф\-ър
$\rho$ эр тх∙хёЄ\-тхэ\-эющ юёш. ╧єёЄ№ ёяхъЄЁ ╩рЁых\-ьрэр $\Lambda$
ьхЁ√ $\mu$ ъю\-эхўхэ. ╥юуфр ёє\-∙хёЄ\-тє■Є яюыш\-эюь√
$P_{\lambda}(t)$, $\lambda\in\Lambda$, $\deg{P_{\lambda}(t)}$ $\leq$
$\rho-1$, Єръшх, ўЄю
$$
d\mu(t)= \sum\limits_{\lambda\in\Lambda} P_{\lambda}(t) e^{-i\lambda
t}dt.
$$
\end{theorem}
%%%%%%%%%%%%%%%%%%%%%%%%%%%%%%%%%%%%%%%%%%%%%%%%%%%%%%%%%%%%%%%%%%%%%%%%%%%%%%%
\par {\sc  ─юърчрЄхы№ёЄтю.}\: ╧єёЄ№ $G(z)$  -- яЁх\-юсЁрчю\-тр\-эшх
╩рЁых\-ьрэр ьхЁ√ $\mu$. ╤яхъЄЁ ьхЁ√ $\mu$ ъю\-эхўхэ ш,
ёыхфю\-тр\-Єхы№\-эю, юуЁрэш\-ўхэ. ╩ръ ёых\-фєхЄ шч
фю\-ърчр\-Єхы№\-ёЄ\-тр Єхю\-Ёх\-ь√  \ref{tauber_th_karleman}, т ¤Єюь
ёыєўрх Їєэъ\-Ўш  $G(z)$ уюыю\-ьюЁЇ\-эр т эх\-ъю\-Єю\-Ёющ
юъ\-ЁхёЄ\-эюё\-Єш схё\-ъю\-эхў\-эюё\-Єш, яЁш\-ў╕ь $G(\infty)=0$.
┴√ыю Єръ\-цх фю\-ърчр\-эю, ўЄю т Ёрё\-ёьрЄЁш\-трх\-ьюь ёыєўрх
$d\mu(t)$ $=$ $g(t)dt$, уфх $g$  -- Ўхыр  Їєэъ\-Ўш 
¤ъёяю\-эхэ\-Ўшры№\-эю\-ую Єшяр. ╬сючэрўшь $g_1(t)$ $=$
$\int\limits_0^t |g(u)|du$. ╥юуфр яЁш\-ьхэх\-эшх ЇюЁ\-ьє\-ы√
шэ\-Єху\-Ёш\-Ёю\-тр\-эш  яю ўрёЄ ь фр╕Є
$$
|G_+(z)|=\left|\int\limits_0^\infty e^{itz}g(t)dt\right|\leq
\int\limits_0^\infty e^{ty}|g(t)|dt
=e^{-ty}\left.g_1(t)\right|_0^\infty+ y\int\limits_0^\infty
e^{-ty}g_1(t)dt.
$$
┬ ёыєўрх, хёыш $\rho<0$, т√\-яюыэ \-хЄ\-ё  эх\-Ёртхэ\-ёЄ\-тю
$|g_1(t)|\leq M$. ┬ ёыєўрх, хёыш $\rho\geq 0$, Єю фы  ы■сюую
$\varepsilon>0$ т√\-яюыэ \-хЄ\-ё  эх\-Ёртхэ\-ёЄ\-тю $|g_1(t)|\leq
M_\varepsilon|t|^{\rho+\varepsilon}$ яЁш $t\geq 1$. ╚ч ¤Єюую
ёых\-фє\-■Є юЎхэъш $|G(z)|\leq M$ яЁш $\rho<0$, $|G(z)|\leq
M_\varepsilon \left(\frac{1}{|y|}\right)^{\rho+\varepsilon}$ яЁш
$|y|\leq 1$ ш $\rho\geq 0$.

\par ╚ч Єюую, ўЄю ёяхъЄЁ ъю\-эхўхэ ёых\-фє\-хЄ, ўЄю Єюўъш ёяхъЄЁр  ты \-■Є\-ё 
шчюыш\-Ёю\-трэ\-э√ьш юфэю\-чэрў\-э√ьш юёюсхэ\-эюё\-Є \-ьш Їєэъ\-Ўшш
$G(z)$. ╚ч яю\-ыєўхэ\-э√ї юЎхэюъ $G(z)$ ёых\-фє\-хЄ, ўЄю яЁш
$\rho<1$ ¤Єш юёюсхэ\-эюё\-Єш  ты \-■Є\-ё  єёЄЁр\-эш\-ь√\-ьш. ┬  ¤Єюь
ёыєўрх $G(z)=0$, $\mu=0$. ┬ ёыєўрх $\rho\geq 1$ Єюўъш ёяхъЄЁр ьюуєЄ
с√Є№ яюы■ёр\-ьш Їєэъ\-Ўшш $G(z)$ яю\-Ё ф\-ър эх т√°х ўхь $\rho$.
╧ю¤Єю\-ьє фы  Їєэъ\-Ўшш $G(z)$ ёяЁр\-тхф\-ыштю яЁхф\-ёЄрт\-ых\-эшх
$$
G(z)=\sum\limits_{\lambda\in \Lambda} \sum\limits_{n=1}^\rho
\frac{a_{n,\lambda}}{(z-\lambda)^n}.
$$
╥хяхЁ№ шч ЇюЁ\-ьєы√
$$
\int\limits_{\frak{L}} \frac{1}{(w-\lambda)^n}e^{-itw}dw=2\pi
i(-it)^n e^{-i\lambda t}
$$
ш Ёртхэ\-ёЄ\-тр  (\ref{tauber_def_g}) ёых\-фєхЄ єЄ\-тхЁ\-цфх\-эшх
Єхю\-Ёх\-ь√.

%%%%%%%%%%%%%%%%%%%%%%%%%%%%%%%%%%%%%%%%%%%%%%%%%%%%%%%%%%%%%%%%%%%%%%%%%%%%%%%
\begin{theorem}\label{tauber_th_4}\hskip-2mm{.}\:
%%%%%%%%%%%%%%%%%%%%%%%%%%%%%%%%%%%%%%%%%%%%%%%%%%%%%%%%%%%%%%%%%%%%%%%%%%%%%%%
╧єёЄ№ $M>0$ -- Їшъёш\-Ёю\-трэ\-эюх ўшёыю, $\mu$  -- Ёр\-фю\-эю\-тр
ьхЁр эр тх∙хёЄ\-тхэ\-эющ юёш, єфютыхЄ\-тю\-Ё ■\-∙р 
эх\-Ёртхэ\-ёЄ\-тє $|\mu|([\alpha,\beta])\leq M$, хёыш
$\beta-\alpha\leq 1$. ╧єёЄ№ $K$ --  сюЁхыхт\-ёър  Їєэъ\-Ўш  шч
яЁюёЄ\-Ёрэ\-ёЄтр $L_1(-\infty,\infty)$ Єр\-ър , ўЄю ьэю\-цх\-ёЄтю
$\Lambda$ $=$ $\{\lambda\in(-\infty,\infty):\:$
$\widehat{K}(\lambda)=0\}$  ты \-хЄё  ъю\-эхў\-э√ь. ╧єёЄ№
$$
\int\limits_{-\infty}^\infty K(t-u)d\mu(t)=0, \qquad
u\in(-\infty,\infty).
$$
╥юуфр ёє\-∙хёЄтє■Є ўшёыр $c_{\lambda}$, $\lambda\in\Lambda$, Єръшх,
ўЄю $d\mu(t)= \sum\limits_{\lambda\in\Lambda} c_{\lambda}
e^{-i\lambda t}dt$.
\end{theorem}
%%%%%%%%%%%%%%%%%%%%%%%%%%%%%%%%%%%%%%%%%%%%%%%%%%%%%%%%%%%%%%%%%%%%%%%%%%%%%%%
\par {\sc  ─юърчрЄхы№ёЄтю.}\: ╚ч єёыю\-тшщ Єхю\-Ёх\-ь√ ш шч
Єхю\-Ёх\-ь√  \ref{tauber_th_G_properties} ёых\-фєхЄ, ўЄю ёяхъЄЁ
╩рЁых\-ьрэр ьхЁ√ $\mu$ ёюфхЁ\-цшЄ\-ё  т ъю\-эхў\-эюь ьэю\-цх\-ёЄ\-тх
$\Lambda$. ╥хяхЁ№ шч Єхю\-Ёх\-ь√  \ref{tauber_th_d_mu} ш
эх\-Ёртхэ\-ёЄ\-тр  $\rho\leq 1$ фы  яю\-Ё ф\-ър  $\rho$ ьхЁ√ $\mu$
ёых\-фєхЄ єЄ\-тхЁ\-ц\-фх\-эшх Єхю\-Ёх\-ь√. ╥хю\-Ёхьр фюърчрэр.

\par ╥хю\-Ёх\-ь√ рэрыюушў\-э√х Єхю\-Ёх\-ьрь  \ref{tauber_th_karleman}-\ref{tauber_th_4}
Ёрэхх фюърч√\-тры ╩юЁхэ\-сы■ь \cite{Korenblum}.

\par ╤ЇюЁ\-ьєышЁєхь х∙╕ ьєы№Єш\-яыш\-ър\-Єшт\-э√щ трЁш\-рэЄ яюёыхфэхщ Єхю\-Ёх\-ь√.
%%%%%%%%%%%%%%%%%%%%%%%%%%%%%%%%%%%%%%%%%%%%%%%%%%%%%%%%%%%%%%%%%%%%%%%%%%%%%%%
\begin{theorem}\label{tauber_th_5}\hskip-2mm{.}\:
%%%%%%%%%%%%%%%%%%%%%%%%%%%%%%%%%%%%%%%%%%%%%%%%%%%%%%%%%%%%%%%%%%%%%%%%%%%%%%%
╧єёЄ№ $M>0$ -- Їшъёш\-Ёю\-трэ\-эюх ўшёыю, $\mu$  -- Ёр\-фю\-эю\-тр
ьхЁр эр яюыє\-юёш $(0,\infty)$, єфютыхЄ\-тю\-Ё ■\-∙р 
эх\-Ёртхэ\-ёЄ\-тє $|\mu|([\alpha,\beta])\leq M$, хёыш
$\frac{\beta}{\alpha}\leq e$. ╧єёЄ№ $K$  --  сюЁхыхт\-ёър  Їєэъ\-Ўш 
эр яюыє\-юёш $(0,\infty)$ Єр\-ър , ўЄю $\frac{1}{t}K(t)$ $\in
L_1(0,\infty)$, ш яєёЄ№ $\Lambda$ $=$
$\{\lambda\in(-\infty,\infty):\:$ $\int\limits_0^\infty
\frac{1}{t}K(t)t^{-i\lambda }dt=0\}$ хёЄ№ ъю\-эхў\-эюх
ьэю\-цх\-ёЄтю. ╧єёЄ№
$$
\int\limits_0^\infty K\left(\frac{t}{r}\right)d\mu(t)=0, \qquad
r\in(0,\infty).
$$
╥юуфр ёє\-∙хёЄ\-тє■Є ўшёыр $c_{\lambda}$, $\lambda\in\Lambda$,
Єръшх, ўЄю $d\mu(t)= \frac{1}{t} \sum\limits_{\lambda\in\Lambda}
c_{\lambda} t^{-i\lambda}dt.$
\end{theorem}
%%%%%%%%%%%%%%%%%%%%%%%%%%%%%%%%%%%%%%%%%%%%%%%%%%%%%%%%%%%%%%%%%%%%%%%%%%%%%%%
\par {\sc  ─юърчрЄхы№ёЄтю.}\: ╬яЁх\-фхышь Ёр\-фюэютє ьхЁє $\nu$
эр тх∙хёЄ\-тхэ\-эющ юёш ёых\-фє\-■\-∙шь юсЁр\-чюь $\nu([a,b])$ $=$
$\mu\left(\left[e^a,e^b\right]\right)$. ┼ёыш $b\leq a+1$, Єю сєфхь
шьхЄ№
$$
|\nu|([a,b])=|\mu|\left(\left[e^a,e^b\right]\right)\leq
|\mu|\left(\left[e^a,e^{a+1}\right]\right)\leq M.
$$
╧єёЄ№ $K_1(x)=K(e^x)$. ╥юуфр $K_1$  -- сюЁхыхт\-ёър  Їєэъ\-Ўш  эр
тх∙хёЄ\-тхэ\-эющ юёш, ъюЄю\-Ёр  яЁш\-эрф\-ыхцшЄ яЁюёЄ\-Ёрэ\-ёЄтє
$L_1(-\infty,\infty)$. ┬√яюыэ \-■Єё  Ёртхэ\-ёЄ\-тр
$$
\int\limits_0^\infty \frac{1}{t} K(t)t^{-i\lambda}dt=
\int\limits_{-\infty}^\infty K_1(x)e^{-ix\lambda}dx,
$$
$$
\int\limits_0^\infty K\left(\frac{t}{r} \right) d\mu(t)=
\int\limits_{-\infty}^\infty K_1(x-u)d\nu(x)=0,\quad u=\ln r.
$$
\par ─ы   фЁр $K_1$ ш ьхЁ√ $\nu$ т√\-яюыэ \-■Є\-ё  тёх єёыю\-тш 
Єхю\-Ёх\-ь√  \ref{tauber_th_4}, яю¤Єю\-ьє
$$
d\nu(x)=\sum\limits_{\lambda\in\Lambda} c_\lambda e^{-i\lambda x}dx,
\qquad d\mu(t)=\frac{1}{t} \sum\limits_{\lambda\in\Lambda} c_\lambda
t^{-i\lambda}dt.
$$
╥хю\-Ёхьр фюърчрэр.

\par ┬ ёых\-фє\-■\-∙хщ Єхю\-Ёхьх юяшё√\-тр■Є\-ё  ьхЁ√ $\mu$, фы  ъюЄю\-Ё√ї
Їєэъ\-Ўш  $\Psi(r)$, юяЁх\-фхы х\-ьр  Ёртхэ\-ёЄ\-тюь
(\ref{introduction_Psi(r)}),  ты \-хЄё  яыюЄ\-эюё\-Є№■
Ёхує\-ы Ё\-эющ ьхЁ√. ═ряюь\-эшь, ўЄю Їєэъ\-Ўш  $\gamma(t)$
юяЁх\-фхы \-хЄё  Ёртхэ\-ёЄ\-тюь (\ref{def_gamma}).
%%%%%%%%%%%%%%%%%%%%%%%%%%%%%%%%%%%%%%%%%%%%%%%%%%%%%%%%%%%%%%%%%%%%%%%%%%%%%%%
\begin{theorem}\label{tauber_th_6}\hskip-2mm{.}\:
%%%%%%%%%%%%%%%%%%%%%%%%%%%%%%%%%%%%%%%%%%%%%%%%%%%%%%%%%%%%%%%%%%%%%%%%%%%%%%%
╧єёЄ№ $\rho(r)$  -- яЁюшч\-тюы№\-э√щ єЄюў\-э╕э\-э√щ яю\-Ё фюъ, $\mu$
-- Ёр\-фю\-эю\-тр ьхЁр эр яюыє\-юёш $(0,\infty)$ шч ъырёёр
$\frak{M}(\rho(r))$, $K$  -- сюЁхыхт\-ёър  Їєэъ\-Ўш  эр яюыє\-юёш
$(0,\infty)$ Єр\-ър , ўЄю $t^{\rho-1}\gamma(t)$ $K(t)$ $\in$
$L_1(0,\infty)$, $c_1$ $=$ $\int\limits_0^\infty K(t)t^{\rho-1}dt$
$\neq 0$. ╧єёЄ№ ьэю\-цх\-ёЄтю
$$
\Lambda=\left\{\lambda\in(-\infty,\infty):\:\int\limits_0^\infty
K(t)t^{\rho-1-i\lambda}dt=0\right\}
$$
 ты хЄё  ъю\-эхўэ√ь. ╧єёЄ№ Їєэъ\-Ўш  $\Psi$ юяЁхфхы \-хЄ\-ё 
Ёртхэ\-ёЄ\-тюь  (\ref{introduction_Psi(r)}). ╥юуфр хёыш ьхЁр $s$,
$ds(t)=\Psi(t)dt$,  ты \-хЄё  Ёхує\-ы Ё\-эющ ьхЁющ
юЄ\-эюёш\-Єхы№\-эю єЄюў\-э╕э\-эю\-ую яю\-Ё ф\-ър $\rho(r)+1$,
яЁш\-ў╕ь $Fr[s]=\{\sigma\}$, уфх $d\sigma(t)=ct^\rho dt$, Єю
яЁх\-фхы№\-эюх ьэю\-цх\-ёЄтю $Fr[\mu]$ ьхЁ√ $\mu$ ёюёЄю\-шЄ шч
рсёю\-ы■Є\-эю эх\-яЁх\-Ё√т\-э√ї ьхЁ $\nu$, ш яыюЄ\-эюёЄ№ $h(t)$
ърц\-фющ Єр\-ъющ ьхЁ√ шьххЄ тшф
$$
h(t)= \left(\frac{c}{c_1}+\sum\limits_{\lambda\in\Lambda}
c_{\lambda} t^{-i\lambda}\right)t^{\rho-1},
$$
уфх $c_{\lambda}$ -- эх\-ъюЄюЁ√х ъюьяыхъё\-э√х ўшёыр.
\end{theorem}
%%%%%%%%%%%%%%%%%%%%%%%%%%%%%%%%%%%%%%%%%%%%%%%%%%%%%%%%%%%%%%%%%%%%%%%%%%%%%%%
\par \quad {\sc ─юърчр\-Єхы№ёЄтю.}\: ╧ю Єхю\-Ёхьх  \ref{abel_th_8} фы  ы■сющ ьхЁ√
$\nu\in Fr[\mu]$ т√\-яюыэ \-хЄё  Ёртхэ\-ёЄ\-тю
\begin{equation}\label{tauber_int_1}
\int\limits_0^\infty K\left(\frac{t}{u}\right)d\nu(t)=cu^\rho.
\end{equation}
┼ёыш $\nu_2(t)=\nu(t)-\nu_1(t)$, уфх $d\nu_1(t)=$
$\frac{c}{c_1\rho}t^{\rho-1}dt$ яЁш $\rho$ $\neq$ $0$ ш $d\nu_1(t)=$
$\frac{c}{c_1}\frac{1}{t}dt$ яЁш  $\rho=0$, Єю
\begin{equation}\label{tauber_int_2}
\int\limits_0^\infty K\left(\frac{t}{u}\right)d\nu_2(t)=0.
\end{equation}

\par ╧єёЄ№ $K_1(t)=t^\rho K(t)$, $\nu_3$  -- Єр\-ър  ьхЁр,
ўЄю $d\nu_3(t)=$ $t^{-\rho}d\nu_2(t)$. ╥юуфр Ёртхэ\-ёЄ\-тю
(\ref{tauber_int_2}) ьюц\-эю яхЁх\-яшёрЄ№ т тшфх
$$
\int\limits_0^\infty K_1\left(\frac{t}{u}\right)d\nu_3(t)=0.
$$

\par ╚ч ёююЄ\-эю°х\-эш  $\mu\in$ $\frak{M}(\rho(r))$ ёых\-фєхЄ, ўЄю Їєэъ\-Ўш 
$$
N(\alpha)=\mathop{\overline{\lim}}\limits_{r\to\infty}
\frac{\check{\mu}(r+\alpha r)-\check{\mu}(r)}{V(r)},
$$
уфх $\check{\mu}(r)$  -- Їєэъ\-Ўш  Ёрё\-яЁх\-фхых\-эш  ьхЁ√ $|\mu|$,
 ты \-хЄ\-ё  юуЁрэш\-ўхэ\-эющ эр ёху\-ьхэЄх $[0,e]$. ╚ч ¤Єюую ш
Єхю\-Ёх\-ь√  \ref{azarin_th_6} ёых\-фєхЄ, ўЄю фы  ы■сющ ьхЁ√ $\nu\in
Fr[\mu]$ эр ьэю\-цх\-ёЄ\-тх $(0,\infty)\times[1,e]$ т√\-яюыэ \-хЄё 
юЎхэър $|\nu|([r,qr])\leq Mr^\rho$ ё эх\-ъю\-Єю\-Ёющ яю\-ёЄю э\-эющ
$M$. ╬ўхтшф\-эю, ўЄю Єр\-ър  юЎхэ\-ър, тюч\-ьюц\-эю, ё фЁєующ
яю\-ёЄю э\-эющ  $M$, т√\-яюыэ \-хЄ\-ё  ш фы  ьхЁ√ $\nu_2$.

\par ╤є∙хёЄ\-тєхЄ ўшёыю $\xi\in[\alpha,\beta]$ Єръюх, ўЄю т√\-яюыэ \-■Є\-ё 
Ёртхэ\-ёЄ\-тр
$$
\int\limits_\alpha^\beta d|\nu_3|(t)= \int\limits_\alpha^\beta
t^{-\rho}d|\nu_2|(t)=\xi^{-\rho}|\nu_2|([\alpha,\beta]).
$$
┼ёыш ЄхяхЁ№ $\beta\leq e\alpha$, Єю
$$
|\nu_3|([\alpha,\beta])=\xi^{-\rho}|\nu_2|([\alpha,\beta])\leq
M\left(\frac{\alpha}{\xi}\right)^\rho\leq Me^{|\rho|}.
$$
╤яЁр\-тхф\-ыштю Єръ\-цх ёююЄ\-эю°х\-эшх $\frac{1}{t}$ $K_1(t)$ $\in$
$L_1(0,\infty)$. ╥ръшь юсЁр\-чюь,  фЁю $K_1$ ш ьхЁр $\nu_3$
єфют\-ыхЄ\-тю\-Ё ■Є тёхь єёыю\-тш\- ь Єхю\-Ёх\-ь√ \ref{tauber_th_5}.
╧ю¤Єю\-ьє $d\nu_3(t)=\frac{1}{t}$ $\sum\limits_{\lambda\in\Lambda}
c_{\lambda} t^{-i\lambda}dt$, ўЄю ¤ътш\-тр\-ыхэЄ\-эю
єЄ\-тхЁ\-ц\-фх\-эш■ Єхю\-Ёх\-ь√. ╥хю\-Ёхьр фю\-ърчр\-эр.

\par ╤яЁр\-тхфышт ёых\-фє\-■∙шщ трЁш\-рэЄ Єхю\-Ёх\-ь√  \ref{tauber_th_6},
т ъюЄю\-Ёюь юуЁрэш\-ўх\-эшх $c_1$ $\neq 0$ юЄ\-ёєЄ\-ёЄ\-тєхЄ.
%%%%%%%%%%%%%%%%%%%%%%%%%%%%%%%%%%%%%%%%%%%%%%%%%%%%%%%%%%%%%%%%%%%%%%%%%%%%%%%
\begin{theorem}\label{tauber_th_7}\hskip-2mm{.}\:
%%%%%%%%%%%%%%%%%%%%%%%%%%%%%%%%%%%%%%%%%%%%%%%%%%%%%%%%%%%%%%%%%%%%%%%%%%%%%%%
╧єёЄ№ $\rho(r)$  -- яЁюшч\-тюы№\-э√щ єЄюў\-э╕э\-э√щ яю\-Ё фюъ, $\mu$
-- Ёр\-фю\-эю\-тр ьхЁр эр яюыє\-юёш $(0,\infty)$ шч ъырёёр
$\frak{M}(\rho(r))$, $K$  -- сюЁх\-ыхт\-ёър  Їєэъ\-Ўш  эр яюыє\-юёш
$(0,\infty)$ Єр\-ър , ўЄю $t^{\rho-1}\gamma(t)$ $K(t)$ $\in$
$L_1(0,\infty)$. ╧єёЄ№ ьэю\-цх\-ёЄтю
$$
\Lambda=\left\{\lambda\in(-\infty,\infty):\:\int\limits_0^\infty
K(t)t^{\rho-1-i\lambda}dt=0\right\}
$$
 ты \-хЄё  ъю\-эхў\-э√ь, яЁш\-ў╕ь $0\in\Lambda$. ╧єёЄ№ Їєэъ\-Ўш 
$\Psi$ юяЁх\-фхы \-хЄё  Ёртхэ\-ёЄ\-тюь (\ref{introduction_Psi(r)}).
╥юуфр, хёыш ьхЁр  $s$, $ds(t)=\Psi(t)dt$,  ты \-хЄё  Ёхує\-ы Ё\-эющ
ьхЁющ юЄ\-эюёш\-Єхы№\-эю єЄюў\-э╕э\-эю\-ую яю\-Ё ф\-ър $\rho(r)+1$,
Єю $Fr[s]=\{0\}$, яЁх\-фхы№\-эюх ьэю\-цх\-ёЄтю $Fr[\mu]$ ёюёЄюшЄ шч
рсёю\-ы■Є\-эю эх\-яЁх\-Ё√т\-э√ї ьхЁ $\nu$, ш яыюЄ\-эюёЄ№ $h(t)$
ърцфющ Єр\-ъющ ьхЁ√ шьххЄ тшф
$$
h(t)= \left(\sum\limits_{\lambda\in\Lambda} c_\lambda
t^{-i\lambda}\right)t^{\rho-1},
$$
уфх $c_{\lambda}$ -- эхъюЄю\-Ё√х ъюьяыхъё\-э√х ўшёыр.
\end{theorem}
%%%%%%%%%%%%%%%%%%%%%%%%%%%%%%%%%%%%%%%%%%%%%%%%%%%%%%%%%%%%%%%%%%%%%%%%%%%%%%%
\par {\sc  ─юърчрЄхы№ёЄтю.}\:  ╥ръ ъръ ьхЁр  $s$  ты \-хЄё  Ёхуєы Ё\-эющ,
Єю ъръ ш яЁш фю\-ърчр\-Єхы№\-ёЄтх яЁх\-ф√фє\-∙хщ Єхю\-Ёх\-ь√ яюыєўшь
Ёртхэ\-ёЄ\-тю  (\ref{tauber_int_1}). ─рыхх Ёрё\-ёєц\-фрхь Єръ. ╧єёЄ№
$\nu$ ш  $\sigma$  -- фтх ьхЁ√ шч $Fr[\mu]$, $\nu_3=\nu-\sigma$.
╚ьххь $\int\limits_0^\infty K\left(\frac{t}{u}\right)d\nu_3(t)$
$=0$. ─рыхх, яютЄю\-Ё   Ёрё\-ёєцфх\-эш  яЁхф√\-фє∙хщ Єхю\-Ёх\-ь√,
яю\-ыєўр\-хь $d\nu_3(t)=\left(\sum\limits_{\lambda\in\Lambda}
c_\lambda t^{-i\lambda}\right)t^{\rho-1}dt$. ╠√ фю\-ърчр\-ыш, ўЄю
фы  ы■сющ ьхЁ√ $\nu\in Fr[\mu]$ ёє\-∙хёЄ\-тєхЄ эрсюЁ ъюья\-ыхъё\-э√ї
яю\-ёЄю э\-э√ї $c_{\lambda}$ (тююс∙х уютюЁ , чр\-тшё \-∙шї юЄ ьхЁ√
$\nu$) Єр\-ъющ, ўЄю т√\-яюыэ \-■Є\-ё  Ёртхэ\-ёЄтю
$d\nu=d\sigma+d\tau$, уфх $\sigma$  -- яЁю\-шч\-тюы№\-эю
чр\-Їшъёш\-Ёю\-трэ\-эр  ьхЁр шч $Fr[\mu]$, р  $d\tau(t)$ $=$ $
\sum\limits_{\lambda\in\Lambda} \left(c_\lambda
t^{-i\lambda}\right)t^{\rho-1}dt$.

\par ─ы  ы■сюую $r>0$ т√\-яюыэ \-хЄ\-ё  Ёртхэ\-ёЄ\-тю
$d\nu_r=d\sigma_r+d\tau_r$. ╟рьхЄшь, ўЄю $d\tau_r(t)=$ \\ $\left(
\sum\limits_{\lambda\in\Lambda} c_{\lambda} r^{-i\lambda}
t^{-i\lambda}\right)$ $t^{\rho-1}dt$. ╤ фЁєующ ёЄюЁюэ√, яю\-ёъюы№ъє
ьхЁр $\nu_r\in Fr[\mu]$, Єю тхЁэю Ёртхэ\-ёЄ\-тю
$$
d\nu_r(t)=d\sigma(t)+ \left( \sum\limits_{\lambda\in\Lambda}
 c_\lambda(r) t^{-i\lambda}\right) t^{\rho-1}dt.
$$
┼ёыш  $\Lambda=\{\lambda_1,...,\lambda_n\}$, Єю шч ёърчрэ\-эю\-ую
ёых\-фєхЄ, ўЄю т√\-яюыэ \-хЄё  Ёртхэ\-ёЄтю
\begin{equation}\label{tauber_d_sigma}
d\sigma_r(t)-d\sigma(t)=\left(\sum\limits_{k=1}^n
\left(c_k(r)-c_kr^{-i\lambda_k}\right)
t^{-i\lambda_k}\right)t^{\rho-1}dt.
\end{equation}

\par ╚ч Ёртхэ\-ёЄ\-тр  (\ref{tauber_d_sigma}) ёых\-фє\-хЄ Ёртхэ\-ёЄтю
\begin{equation}\label{tauber_int_d_sigma}
\int\limits_0^\infty f(t)(d\sigma_r(t)-d\sigma(t)) =
\sum\limits_{k=1}^n \left(c_k(r)-c_kr^{-i\lambda_k}\right)
\int\limits_0^\infty f(t)t^{\rho-1-i\lambda_k}dt.
\end{equation}
╩юэъЁхЄэю Їєэъ\-Ўшш $f(t)$ сєфєЄ т√сшЁрЄ№ё  яючцх. ─рыхх шьххь
$$
\int\limits_0^\infty f(t)t^{\rho-1-i\lambda_k}dt=
\int\limits_{-\infty}^\infty e^{\rho x} f(e^x) e^{-i\lambda_k x}dx.
$$
┬√схЁхь ЄхяхЁ№ $f(x)=f_m(x)$ Єръ, ўЄюс√ фы  Їєэъ\-Ўшш
$\varphi_m(x)=e^{\rho x}f_m(e^x)$ т√\-яюыэ \-ыюё№ Ёртхэ\-ёЄтю
$$
\int\limits_{-\infty}^\infty \varphi_m(x)e^{-i\lambda x}dx=
\prod\limits_{k\neq m} (\lambda-\lambda_k)
\left(\frac{\sin\alpha\lambda}{\lambda}\right)^p,
$$
уфх  $p$  -- фю\-ёЄр\-Єюўэю сюы№°юх Ўхыюх ўшёыю, р тх∙хёЄ\-тхэ\-эюх
ўшёыю  $\alpha$ Єръютю, ўЄю $\sin\alpha\lambda_m$ $\neq 0$, хёыш
$\lambda_m$ $\neq 0$ ш  $\alpha=1$, хёыш $\lambda_m=0$. ╚ч
юяЁх\-фхых\-эш   $\varphi_m$ ш Єхю\-Ёх\-ь√ ┬шэхЁр-╧¤ыш ёых\-фєхЄ,
ўЄю $\varphi_m$  -- ЇшэшЄ\-эр  Їєэъ\-Ўш . ╚ч ЇюЁ\-ьєы√ юсЁр∙х\-эш 
фы  яЁх\-юсЁрчю\-тр\-эш  ╘єЁ№х ыхуъю ёых\-фєхЄ, ўЄю Їєэъ\-Ўш 
$\varphi_m$ шьххЄ $p-n-3$ эх\-яЁх\-Ё√т\-э√ї яЁю\-шч\-тюф\-э√ї. ╚ч
¤Єюую т ётю■ юўхЁхф№ ёых\-фєхЄ, ўЄю $\supp f_m$ хёЄ№ ъюь\-яръЄ,
яю\-ьх∙р■\-∙шщё  эр яюыє\-юёш $(0,\infty)$, ш ўЄю Їєэъ\-Ўш   $f_m$
шьххЄ эх\-яЁх\-Ё√т\-э√х яЁю\-шч\-тюф\-э√х фю яю\-Ё ф\-ър $p-n-3$
тъы■ўш\-Єхы№\-эю. ╧Ёш Єръюь т√\-сюЁх  $f$ Ёртхэ\-ёЄтю
(\ref{tauber_int_d_sigma}) яЁш\-эшьр\-хЄ тшф
$$
\frac{1}{r^\rho} \int\limits_0^\infty f_m
\left(\frac{t}{r}\right)d\sigma(t)-\int\limits_0^\infty
f_m(t)d\sigma(t)=A_m\left(c_m(r)-c_mr^{-i\lambda_m}\right),
$$
уфх $A_m\neq 0$. ╚ч эр\-яшёрэ\-эю\-ую Ёртхэ\-ёЄ\-тр ыхуъю
т√\-тхё\-Єш, ўЄю ърцфр  шч Їєэъ\-Ўшщ $c_k(r)$  ты \-хЄ\-ё 
схё\-ъю\-эхў\-эю фшЇ\-Їх\-Ёхэ\-Ўш\-Ёєх\-ьющ эр яюыє\-юёш
$(0,\infty)$, ш ўЄю т√\-яюы\-э \-хЄ\-ё  Ёртхэ\-ёЄ\-тю $c_m(1)=c_m$.
─рыхх, фы  юяЁх\-фхы╕э\-эюё\-Єш, сєфхь ёўш\-ЄрЄ№, ўЄю $\rho>0$.
┴єфхь Єръ\-цх ёўш\-ЄрЄ№, ўЄю Їєэъ\-Ўш  Ёрё\-яЁх\-фхых\-эш 
$\sigma(t)$ ьхЁ√ $\sigma$ эюЁьш\-Ёю\-трэр єёыю\-тшхь $\sigma(0)=0$.
╥юуфр $\sigma_r(t)$ $=$ $\frac{\sigma(rt)}{r^\rho}$. ─рыхх,
шэ\-Єху\-Ёш\-Ёє  юсх ўрёЄш Ёртхэ\-ёЄ\-тр  (\ref{tauber_d_sigma}) яю
ёху\-ьхэ\-Єє $[0,1]$, яЁш\-їюфшь ъ Ёртхэ\-ёЄтє
$$
\frac{\sigma(r)}{r^\rho}-\sigma(1)= \sum\limits_{k=1}^n
\left(c_k(r)-c_k r^{-i\lambda_k}\right) \frac{1}{\rho-i\lambda_k}.
$$
╚ч ¤Єюую Ёртхэ\-ёЄ\-тр ёых\-фє\-хЄ, ўЄю ьхЁр $\sigma$ шьххЄ
схё\-ъю\-эхўэю фшЇ\-Їх\-Ёхэ\-Ўш\-Ёєх\-ьє■ яыюЄ\-эюёЄ№. ╬сю\-чэрўшь
¤Єє яыюЄ\-эюёЄ№ ўхЁхч $h(t)$. ╥юуфр Ёртхэ\-ёЄтю
(\ref{tauber_d_sigma}) ьюц\-эю яЁх\-юсЁрчю\-трЄ№ ъ тшфє
$$
\frac{h(rt)}{(rt)^{\rho-1}}-\frac{h(t)}{t^{\rho-1}}=
\sum\limits_{k=1}^n \left(c_k(r)-c_kr^{-i\lambda_k}\right)
t^{-i\lambda_k}.
$$
┼ёыш юсючэрўшЄ№ ўхЁхч $H(r,t)$ яЁртє■ ўрёЄ№ яю\-ыєўхэ\-эю\-ую
Ёртхэ\-ёЄ\-тр, Єю шч ¤Єюую Ёртхэ\-ёЄ\-тр ёых\-фєхЄ, ўЄю Їєэъ\-Ўш 
$H(r,t)$ єфют\-ыхЄ\-тю\-Ё хЄ фшЇ\-ЇхЁхэ\-Ўшры№\-эюьє єЁрт\-эх\-эш■
$$
rH'_r(r,t)-tH'_t(r,t)=t\left(\frac{h(t)}{t^{\rho-1}}\right)'.
$$
▌Єю фр╕Є
$$
\sum\limits_{k=1}^n \left(rc'_k(r)+i\lambda_kc_k(r)\right)
t^{-i\lambda_k}=t\left(\frac{h(t)}{t^{\rho-1}}\right)'.
$$
╚ч ышэхщ\-эющ эх\-чр\-тшёш\-ьюёЄш Їєэъ\-Ўшщ $t^{-i\lambda_k}$
ёых\-фєхЄ, ўЄю ёє\-∙хёЄ\-тє\-■Є ўшёыр  $d_k$ Єръшх, ўЄю
т√\-яюыэ \-■Є\-ё  Ёртхэ\-ёЄ\-тр $rc'_k(r)+i\lambda_kc_k(r)=d_k$.
╥юуфр $\left(\frac{h(t)}{t^{\rho-1}}\right)'$ $=$
$\sum\limits_{k=1}^n d_kt^{-i\lambda_k-1}$. ╤ўшЄр , ўЄю
$\lambda_1=0$, яюыєўшь, ўЄю $h(t)=\left(d_1\ln t +
\sum\limits_{k=2}^n \frac{d_k}{-i\lambda_k}t^{-i\lambda_k}+d
\right)t^{\rho-1}$. ╥ръ ъръ $h(t)$  -- яыюЄ\-эюёЄ№ ьхЁ√ $\sigma\in
Fr[\mu]$, Єю ьхЁр $\sigma$ яю Єхю\-Ёх\-ьх \ref{azarin_th_7} фюыцэр
єфют\-ыхЄ\-тю\-Ё Є№ эх\-Ёртхэ\-ёЄ\-тє $|\sigma([r,er])|\leq \alpha
r^\rho$ ё эх\-ъюЄю\-Ё√ь  $\alpha$. ╧ю¤Єю\-ьє $d_1=0$. ╠√ фю\-ърчрыш,
ўЄю яыюЄ\-эюёЄ№ $h(t)$ яЁю\-шч\-тюы№\-эющ ьхЁ√ $\sigma\in Fr[\mu]$
шьххЄ тшф, єърчрэ\-э√щ т ЇюЁ\-ьє\-ыш\-Ёют\-ъх Єхю\-Ёх\-ь√. ─рыхх шч
Єхю\-Ёх\-ь√ \ref{abel_th_8} ёых\-фє\-хЄ, ўЄю $Fr[s]=\{0\}$.
╥хю\-Ёхьр фюърчрэр.

\par ┬ ёт чш ё Єхю\-Ёх\-ьющ  \ref{tauber_th_7} Ёрё\-ёьюЄ\-Ёшь ёых\-фє\-■\-∙шщ
яЁш\-ьхЁ. ╧єёЄ№ $K$  --  фЁю эр яюыє\-юёш $(0,\infty)$ Єръюх, ўЄю
$\frac{1}{t}K(t)$, $\frac{\ln t}{t}$ $K(t)$ $\in$ $L_1(0,\infty)$,
$\int\limits_0^\infty \frac{1}{t}K(t)dt=0$. ╚ьххь
$$
\Psi(r)=\int\limits_0^\infty K\left(\frac{t}{r}\right)\frac{\ln
t}{t}dt= \int\limits_0^\infty K(u)\frac{\ln u+\ln r}{u}du=
\int\limits_0^\infty K(u)\frac{\ln u}{u}du.
$$
╧єёЄ№ $\rho(r)=0$. ╠хЁр  $s$, $ds(t)=\Psi(t)dt$,  ты \-хЄё 
Ёхуєы Ё\-эющ ьхЁющ юЄ\-эюёш\-Єхы№\-эю єЄюў\-э╕э\-эю\-ую яю\-Ё ф\-ър
ЄюцфхёЄ\-тхээю Ёртэю\-ую хфшэшЎх. ╠хЁр  $\mu$, $d\mu(t)=$ $\frac{\ln
t}{t}dt$, эх яЁш\-эрф\-ыхцшЄ ъырёёє $\frak{M}_\infty(0)$, р Єхь
сюыхх ъырёёє $\frak{M}(0)$.

\par ╫рёЄэ√щ ёыєўрщ Єхю\-Ёх\-ь√  \ref{tauber_th_6}, ъюуфр $\Lambda=\emptyset$,
яЁш\-тюфшЄ ъ ёых\-фє\-■\-∙хщ Єрє\-сх\-Ёю\-тющ Єхю\-Ёх\-ьх.
%%%%%%%%%%%%%%%%%%%%%%%%%%%%%%%%%%%%%%%%%%%%%%%%%%%%%%%%%%%%%%%%%%%%%%%%%%%%%%%
\begin{theorem}\label{tauber_th_8}\hskip-2mm{.}\:
%%%%%%%%%%%%%%%%%%%%%%%%%%%%%%%%%%%%%%%%%%%%%%%%%%%%%%%%%%%%%%%%%%%%%%%%%%%%%%%
╧єёЄ№ $\rho(r)$  -- яЁюшч\-тюы№\-э√щ єЄюў\-э╕э\-э√щ яю\-Ё фюъ, $\mu$
-- Ёр\-фю\-эю\-тр ьхЁр эр яюыє\-юёш $(0,\infty)$ шч ъырёёр
$\frak{M}(\rho(r))$, $K$  -- сюЁхыхт\-ёър  Їєэъ\-Ўш  эр яюыє\-юёш
$(0,\infty)$ Єр\-ър , ўЄю $t^{\rho-1}\gamma(t)$ $K(t)$ $\in$
$L_1(0,\infty)$. ╧єёЄ№ Їєэъ\-Ўш  $\int\limits_0^\infty K(t)
t^{\rho-1-i\lambda}dt$ эх юсЁр∙р\-хЄ\-ё  т эюы№ эр тх∙хёЄ\-тхэ\-эющ
юёш, ш яєёЄ№ Їєэъ\-Ўш  $\Psi$ юяЁх\-фхы \-хЄё  Ёртхэ\-ёЄ\-тюь
(\ref{introduction_Psi(r)}). ╥юуфр, хёыш ьхЁр  $s$,
$ds(t)=\Psi(t)dt$,  ты \-хЄё  Ёхуєы Ё\-эющ ьхЁющ юЄ\-эюёш\-Єхы№\-эю
єЄюў\-э╕э\-эю\-ую яю\-Ё ф\-ър $\rho(r)+1$, Єю ьхЁр $\mu$  ты \-хЄё 
Ёхуєы Ё\-эющ ьхЁющ юЄ\-эюёш\-Єхы№\-эю єЄюў\-э╕э\-эю\-ую яю\-Ё ф\-ър
$\rho(r)$, яЁш\-ў╕ь хёыш $Fr[s]$ ёюёЄюшЄ шч ьхЁ√ ё яыюЄ\-эюёЄ№■
$ct^\rho$, Єю $Fr[\mu]$ ёюёЄюшЄ шч ьхЁ√ ё яыюЄ\-эюёЄ№■
$\frac{c}{c_1}t^{\rho-1}$, уфх $c_1$ $=$ $\int\limits_0^\infty
K(t)t^{\rho-1}dt$.
\end{theorem}
%%%%%%%%%%%%%%%%%%%%%%%%%%%%%%%%%%%%%%%%%%%%%%%%%%%%%%%%%%%%%%%%%%%%%%%%%%%%%%%
\par {\sc  ─юърчрЄхы№ёЄтю.}\: ╧ю Єхю\-Ёхьх  \ref{tauber_th_6} яЁх\-фхы№\-эюх
ьэю\-цх\-ёЄ\-тю $Fr[\mu]$ ёюёЄюшЄ шч хфшэ\-ёЄ\-тхэ\-эющ ьхЁ√.
╧ю¤Єю\-ьє ¤Єр ьхЁр Ёхуєы Ё\-эр. ┼╕ яЁх\-фхы№\-эюх ьэю\-цх\-ёЄ\-тю
юяшёрэю т Єхю\-Ёх\-ьх  \ref{tauber_th_6}. ╥хю\-Ёхьр фюърчрэр.

\par ╩ръ ёых\-фєхЄ шч Єхю\-Ёх\-ь√  \ref{azarin_th_10} фы  яюыю\-цш\-Єхы№\-эющ
ьхЁ√ $\mu$ х╕ Ёхуєы Ё\-эюёЄ№ юЄ\-эюёш\-Єхы№\-эю єЄюў\-э╕э\-эю\-ую
яю\-Ё ф\-ър $\rho(r)$ ¤ътш\-трыхэЄ\-эр т√\-яюыэх\-эш■ Ёртхэ\-ёЄ\-тр
$$
\lim\limits_{r\to\infty}
\frac{\mu([ar,br])}{V(r)}=c\frac{b^\rho-a^\rho}{\rho}, \qquad
0<a<b<\infty,
$$
ё эх\-ъю\-Єю\-Ёющ яюёЄю э\-эющ  $c$. ╧Ёш ¤Єюь  $\rho$ ью\-цхЄ с√Є№
яЁюшч\-тюы№\-э√ь тх∙хёЄ\-тхэ\-э√ь ўшёыюь (т ёыєўрх $\rho=0$ яЁртр 
ўрёЄ№ эр\-яшёрэ\-эю\-ую Ёртхэ\-ёЄ\-тр юяЁхфхы \-хЄё  яю
эх\-яЁх\-Ё√т\-эюё\-Єш ш Ёртэр $c\ln\frac{b}{a}$). ┬ ёыєўрх $\rho>0$
эр\-яшёрэ\-эюх Ёртхэ\-ёЄтю ¤ътш\-тр\-ыхэЄ\-эю Ёртхэ\-ёЄ\-тє
$\lim\limits_{r\to\infty} \frac{\mu((1,r])}{V(r)}$ $=$
$\frac{c}{\rho}$, р т ёыєўрх $\rho<0$ юэю ¤ътш\-трыхэЄ\-эю
Ёртхэ\-ёЄтє $\lim\limits_{r\to\infty} \frac{\mu((r,\infty))}{V(r)}$
$=-\frac{c}{\rho}$.

\par ╨хуєы Ё\-э√х чэръю\-яхЁх\-ьхэ\-э√х ьхЁ√ юяшёрэ√ т
Єхю\-Ёх\-ьх  \ref{azarin_th_11}.

\par ┬ ЄхюхЁх\-ьх  \ref{tauber_th_8} Ёрё\-ёьрЄЁш\-тр\-хЄё  яЁюшч\-тюы№\-э√щ
єЄюў\-э╕э\-э√щ яю\-Ё фюъ. ╩ръ єцх юЄ\-ьхўр\-ыюё№ тю тёЄєяых\-эшш,
єцх ўрёЄ\-э√щ ёыєўрщ Єхю\-Ёх\-ь√  \ref{tauber_th_8}, ъюуфр
$\rho(r)\equiv 1$,  ты \-хЄё  єёшых\-эшхь тЄюЁющ Єрє\-сх\-Ёю\-тющ
Єхю\-Ёх\-ь√ ┬шэхЁр яю эх\-ёъюы№\-ъшь эр\-яЁрт\-ыхэш\- ь.

%%%%%%%%%%%%%%%%%%%%%%%%%%%%%%%%%%%%%%%%%%%%%%%%%%%%%%%%%%%%%%%%%%%%%%%%%%%%%%%%%%%%%%%%%%%%%%%%%%%%%%%%%%%%%

%%%%%%%%%%%%%%%%%%%%%%%%%%%%%%%%%%%%%%%%%%%%%%%%%%%%%%%%%%%%%%%%%%%%%%%%%%%%%%%%%%%%%%%%%%%%%%%%%%%%%%%%%%%%%

 %\contentsname{╤яшёюъ ышЄхЁрЄєЁ√}{toc}
 %%%%%%%%%%%%%%%%%%%%%%%%%%%%%%%%%%%%%%%%%%%%%%%%%%%%%%%%%%%%%%%%%%%%%%%%%%%%
\vskip1cm {\bf └.\,╘.  ├Ёш°шэ (A.\,F.  Grishin)}
%\par { ╒рЁ№ъютёъшщ
%эрЎшюэры№э√щ єэштхЁёшЄхЄ шьхэш  ┬.\,═.  ╩рЁрчшэр}
\par {\it E-mail:}{\;grishin@univer.kharkov.ua} \vskip0.5cm
\par {\bf ╚.\,┬. ╧юхфшэЎхтр (I.\,V.  Poedintseva)}
%\par {╒рЁ№ъютёъшщ эрЎшюэры№э√щ
%руЁрЁэ√щ єэштхЁёшЄхЄ шьхэш  ┬.\,┬.  ─юъєўрхтр}
\par {\it E-mail}: {\,Irina.V.Poedintseva@univer.kharkov.ua}

\end{document}